%% file: Higher_Inf_Sym_Final.tex
\documentclass[twoside,pdftex, 11pt]{article}

\usepackage[latin1]{inputenc}
\usepackage[T1]{fontenc}
\usepackage[english]{babel}
\usepackage{color}
\usepackage{tocloft}
\usepackage{verbatim}
\usepackage{bbm}

\usepackage{tikz}
\usetikzlibrary{snakes}
\usetikzlibrary{arrows}
\usepackage{tikz-cd}

\usepackage{longtable}

\usepackage[
  hmarginratio={1:1},     
  vmarginratio={1:1},     
  textwidth=17cm,        
  textheight=24cm,
  heightrounded          
]{geometry}

\usepackage[noindentafter]{titlesec}
\titlespacing{\paragraph}{%
  0pt}{
0.25\baselineskip}{
1em}

\titlespacing{\section}{%
0pt}{
0.2cm}{
0em}

\titlespacing{\subsection}{%
0pt}{
0.2cm}{
0em}

\titlespacing{\subsubsection}{%
0pt}{
0cm}{
0em}

\usepackage{amsmath,amsthm,amssymb}
\usepackage{mathrsfs, mathtools, xpatch}
\usepackage{dsfont}

\allowdisplaybreaks

\usepackage[numbers,sort&compress]{natbib}
\makeatletter
\renewcommand{\@biblabel}[1]{[#1]\hfill}
\makeatother
\usepackage[mathscr]{euscript}
\usepackage{calrsfs}
\DeclareMathAlphabet{\pazocal}{OMS}{zplm}{m}{n}

\usepackage[shortlabels]{enumitem}

\usepackage[hidelinks]{hyperref}

\newtheoremstyle{thm}                                                           
{0.15cm}                                         
{0.15cm}                                         
{\itshape}      
{}                                      
{\bfseries}                             
{}                                      
{0.2cm}                                         
{\thmname{#1}~\thmnumber{#2}\thmnote{ (#3)}}%

\makeatletter
\xpatchcmd{\proof}{\topsep6\p@\@plus6\p@\relax}{}{}{}
\makeatother

\newtheoremstyle{rmk}                                                           
{0.15cm}                                         
{0.15cm}                                         
{}      
{}                                      
{\bfseries}                             
{}                                      
{0.2cm}                                         
{}                                      

\theoremstyle{thm}
\newtheorem{theorem}[equation]{Theorem}
\newtheorem{conjecture}[equation]{Conjecture}
\newtheorem{corollary}[equation]{Corollary}

\newtheorem{lemma}[equation]{Lemma}
\newtheorem{proposition}[equation]{Proposition}

\theoremstyle{rmk}
\newtheorem{example}[equation]{Example}
\newtheorem{remark}[equation]{Remark}
\newtheorem{definition}[equation]{Definition}

\numberwithin{equation}{section}

\DeclareMathAlphabet{\mathbbe}{U}{bbold}{m}{n}

\input{Higher_Inf_Sym_Defs.tex}

\mathtoolsset{showonlyrefs,showmanualtags}

\newlength{\Displayskip}
\begin{document}

\setlength{\abovedisplayskip}{\Displayskip}
\setlength{\belowdisplayskip}{\Displayskip}

\

\vspace{2cm}

\begin{center}
\LARGE{\textbf{Symmetries and Higher-Form Connections\\[0.1cm]in Derived Differential Geometry}}
\end{center}
\vspace{0.5cm}
\begin{center}
\large Severin Bunk \ , \ 
Lukas Müller \ , \ 
Joost Nuiten \ \ 
and \ \ Richard J.~Szabo
\end{center}

\vspace{0.75cm}

\begin{abstract}
\noindent
We introduce a general definition of higher-form connections on principal $\infty$-bundles in differential geometry.
These connections are built as infinitesimal parallel transports and thus satisfy the fake-curvature condition.
We achieve this by developing the formal differentiation and integration of maps from smooth manifolds to derived stacks with sufficient deformation theory.
That allows us to introduce the Atiyah $L_\infty$-algebroid of a principal $\infty$-bundle and establish its global sections as the $L_\infty$-algebra of the derived higher symmetry group of the bundle.
We  define the space of $p$-form connections on an $\infty$-bundle as the space of order~$p$ splittings of its Atiyah $L_\infty$-algebroid. This can be cast equivalently as lifting the classifying map of a bundle on a manifold to the order~$p$ truncation of the de Rham stack of the manifold.
We demonstrate that our new concept of derived geometric $p$-form connections recovers the known notion of connections on higher $\rmU(1)$-bundles defined via \v{C}ech--Deligne differential cocycles. 
We further relate the $L_\infty$-algebras of derived higher symmetries of higher $\rmU(1)$-bundles and higher Courant algebroids. Some applications in higher gauge theory and in supergravity are mentioned.
\end{abstract}

\vspace{0.75cm}
\hspace{0.72cm}%
\begin{tabular}{p{0.75cm} p{15cm}}
	(SB) & Department of Physics, Astronomy and Mathematics
	\newline
	University of Hertfordshire
	\newline
	College Lane, Hatfield, AL10 9AB, United Kingdom
	\newline
	\textit{Email:} \texttt{s.bunk@herts.ac.uk}
	\\[0.2cm]
	(LM) & Faculty of Physics,
	\newline
	Ludwig Maximilian University Munich \newline Theresienstraße 37, 80333 München-Maxvorstadt, Germany
	\newline
	\textit{Email:} \texttt{lukas.mueller2@lmu.de}
	\\[0.2cm]
	(JN) & Institut de Mathématiques de Toulouse
	\newline
	Université de Toulouse
	\newline
	118 route de Narbonne, F-31062 Toulouse Cedex 9, France
	\newline
	\textit{Email:} \texttt{joost.nuiten@math.univ-tolouse.fr}
	\\[0.2cm]
	(RJS) & Department of Mathematics
	\newline
	Heriot-Watt University
	\newline
	Colin Maclaurin Building, Riccarton, Edinburgh EH14 4AS, United Kingdom
	\newline
	and Maxwell Institute for Mathematical Sciences, Edinburgh, United Kingdom
	\newline
	\textit{Email:} \texttt{r.j.szabo@hw.ac.uk}
\end{tabular}

\newpage

{\baselineskip=18pt
\tableofcontents
}

\newpage


\section{Introduction}
\label{sec:Introduction}



\subsection{Towards higher-form connections}
\label{sec: non-technical overview}


Connections on principal bundles form a cornerstone of differential geometry and have been adopted in many other areas.
They describe, for instance, geodesics and curvature in Riemannian and Lorentzian geometry, gauge fields in  physics, the path signature in data science, and Chern--Weil characteristic classes in differential topology, among many further applications.
If $P  \longrightarrow M$ is a principal bundle, a connection on $P$ essentially describes how to compare the fibres $P_{|x}$ and $P_{|y}$ at the start and end point of any smooth path $\gamma$ in $M$.
More precisely, a connection describes this only for infinitesimally short paths---i.e.~on tangent vectors---and these data can then be integrated to the desired comparison along all paths.

Additional data is necessary when the fibres of $P$ have higher internal structure.
As a first example, suppose that $P$ is a bundle whose fibres are groupoids (which may be topological or smooth in some sense).
Then the comparison map $\PT_\gamma \colon P_{|x} \longrightarrow P_{|y}$ associated to $\gamma$ is a functor, and we are free to specify, for each homotopy $h \colon \gamma \longrightarrow \gamma'$ of paths with fixed endpoints, a natural isomorphism $\PT_h \colon \PT_\gamma \longrightarrow \PT_{\gamma'}$ as part of the connection data.
Infinitesimally, we now need to specify the action of the connection not only on tangent vectors (which yields $\PT_\gamma$), but also on pairs of tangent vectors (which yields $\PT_h$).
Geometrically, a classical $1$-form connection on a principal bundle $P$ consists of locally defined $1$-forms valued in the Lie algebra of the structure group of $P$.
In contrast, a 2-form connection on a categorified principal bundle $P$, whose fibres are groupoids, consists of locally defined 1-forms and 2-forms taking values in the Lie 2-algebra of the structure 2-group of $P$.
In order for these connection data to produce a parallel transport, it is necessary that the field strength of the 1-form part coincides with the internal differential of the 2-form part; this is called the \textit{fake-curvature condition} and goes back to pioneering works on higher connections~\cite[Theorem~23]{BS:HGT}.
Another direction of current research, on `adjusted connections', aims to find notions of connections which circumvent fake-curvature conditions (motivated, for instance, by heterotic supergravity and M-theory).
However, our deeper motivation for the present paper is to eventually lift differential geometric techniques based on parallel transport to higher-geometric contexts.
As a consequence, we focus on connections that necessarily satisfy higher versions of the fake-curvature condition.
We further comment on the relation to adjusted connections below in Remark~\ref{rmk: relation to adjusted connections}.

Connections on bundles whose fibres are smooth spaces with higher internal structure are by no means exotic, but appear in many instances in mathematics and its applications in physics, as well as increasingly other areas.
Examples include (differential) string structures in the study of topological modular forms, positive scalar curvature manifolds~\cite{Stolz:Pos_Ric, Berwick-Evans:SuSy_Localisation_and_Witten_genus} and supergravity~\cite{TD:Chern_corr_for_higher_PrBuns}, conjectural higher-dimensional versions of the path signature in data science~\cite{Kapranov:Free_LiAgds_and_space_of_paths, LO:Random_surfaces, BL:Surface_signature}, higher gauge theory~\cite{BS:HGT, BFJKNRSW:GHT}, and in gauging higher-form symmetries in quantum field theory~\cite{GKSW:Generalised_global_symmetries,Schafer-Nameki:2023jdn, BBFTGGPT:Notes_on_generalised_symmetries,Santilli:2024dyz}.

In their most general form, spaces which combine smoothness and higher internal structure form the objects of an $\infty$-category $\scX$.
Usually $\scX$ has particularly nice properties making it into an $\infty$-topos.
As a first example, one could take $\scX = \Sh_\infty(\Mfd)$ to be the $\infty$-category of sheaves of $\infty$-groupoids on the site of manifolds and open coverings.
This is sufficient for a range of questions whose answers can be found in higher, but non-derived, geometry.
In particular, in any $\infty$-topos $\scX$ we can find, for any group object $G$, a classifying object $\rmB G \in \scX$.
The $\infty$-groupoid of principal $\infty$-bundles on an object $X \in \scX$ with structure group $G$ is then given by the mapping space~\cite{NSS:oo-bundles_I, Bunk:Pr_oo-Bundles, Bunk:oo-bundles}
\begin{equation}
	\Bun(X; G) = \scX(X, \rmB G)\ .
\end{equation}
Given that a good theory of principal bundles exists at this level of generality, we arrive at the following very natural and mathematically pressing question:
\begin{itemize}
\item[]	\textsf{\large
		How can one define connections on principal $\infty$-bundles?
	}
\end{itemize}

Despite its broad range of applications---including but not limited to the examples mentioned above---no general theory of connections on principal $\infty$-bundles exists to date.
Various partial answers and definitions of such higher connections have been given in the literature, for instance in~\cite{BS:HGT,Kapranov:Free_LiAgds_and_space_of_paths,SW:PT_and_functors,SSS:L_infty-alg_connections, FSS:Cech_diff_char_classes_via_L_infty,AAS:A_infty-de_Rham_and_integrating_RUTHs,BS:Higher_Riemann-Hilbert,Kapranov:Membranes,Waldorf:PT_in_2-bundles,BO:Integrating_2-term_RUTHs,Chatterjee:2021sif}.
However, each of these setups make restrictive assumptions on the bundles or their structure groups, or only describe flat higher-form connections (i.e.~$\infty$-local systems).
One main upshot of the results in this paper is that they facilitate the definition of connections on maximally generic principal $\infty$-bundles in this setting, with the only condition on their classifying object being that it admits a sufficiently strong notion of tangent complex.
It does not need to come presented as a Lie $n$-group or crossed module, or be finite-dimensional.

Our approach to the definition of higher connections on a bundle is based on the relationship with the infinitesimal symmetries of the bundle.
To better illustrate our approach and results, we begin by briefly recalling this perspective in the case of classical principal bundles.

\paragraph{Symmetries and connections from Atiyah algebroids.}

Let $G$ be a connected Lie group and $P\longrightarrow M $ a principal $G$-bundle over a manifold $M$. A \emph{symmetry of the pair $(P,M)$} is a diffeomorphism $ \phi \colon M \longrightarrow M $ together with an isomorphism of principal bundles $\psi \colon P \longrightarrow \phi^* P$.
This defines an infinite-dimensional Lie group $\cSym(P)$ of symmetries of $P$, which is part of a short exact sequence of Lie groups
\begin{align}
	1 \longrightarrow \cAut(P) \longrightarrow \cSym(P) \longrightarrow \cDiff (M) \longrightarrow 1\ ,
\end{align} 
where $\cAut(P)$ is the group of automorphisms, or gauge transformations, of $P$ and $\cDiff(M)$ is the group of diffeomorphisms of $M$ (see for example~\cite[Section~3]{BMS:Sm2Grp} for more details). 

Differentiating at the identity leads to a short exact sequence of Lie algebras 
\begin{align}\label{eq: inf sym}
	0 \longrightarrow \Ad(P) \longrightarrow \At(P) \overset{\rho}{\longrightarrow}  TM \longrightarrow 0\ ,
\end{align} 
where $TM$ denotes the Lie algebra of vector fields on $M$, $\Ad(P)$ is the Lie algebra of sections of the adjoint bundle of $P$, and $\At(P)$ is the Lie algebra of $G$-invariant vector fields on $P$ (or equivalently sections of $TP/G$). The exact sequence~\eqref{eq: inf sym} controls infinitesimal symmetries of $P$. A crucial observation is that this infinitesimal sequence is \textit{local} on $M$, in the sense that it is even a short exact sequence of \textit{Lie algebroids} on $M$.
The Lie algebroid $\At(P)$ is called the \emph{Atiyah algebroid} of $P$~\cite{Atiyah}.

The set of connections on $P$ can then be expressed as~\cite{Kapranov:Free_LiAgds_and_space_of_paths}
\begin{equation}
	\Con(P) = \Lie\Agd_M \big( Q^{{(1)}} TM,\, \At(P) \big)\ ,
\end{equation}
the set of maps from the free Lie algebroid on the tangent bundle $TM$ to the Atiyah algebroid $\At(P)$.
This is the same as vector bundle morphisms $TM \longrightarrow  \fgt \At(P)$, where we forget the Lie algebroid structure on $\At(P)$.
A connection $A \colon Q^{{(1)}}TM \longrightarrow \At(P)$ is flat precisely if $A$ respects the brackets, i.e.~it even induces a map of Lie algebroids $TM \longrightarrow \At(P)$.
 
In the present paper we develop a theory of higher-form connections on higher bundles from first principles, taking as its foundations this Lie algebroid perspective on classical connections.
To every higher principal bundle we assign an \emph{Atiyah $L_\infty$-algebroid} encoding its infinitesimal symmetries and define connections in terms of certain partial splittings, as we explain in more detail below.
Defining higher connections using Lie algebroids has been attempted before, see for instance~\cite{Stevenson:Lie_2-algs_and_gerbes, Kapranov:Free_LiAgds_and_space_of_paths, Kapranov:Membranes}.


\subsection{Summary of main results}
\label{sec: summary of main results}


Heuristically, a $p$-form connection on a principal $\infty$-bundle on a manifold $M$ should define parallel transports along infinitesimal paths, their homotopies, homotopies of these  homotopies, and so on up to dimension $p$.
That is, such a $p$-form connection on a bundle $x \colon M \longrightarrow \rmB G$ can be understood as a map from an order~$p$ truncation of the infinitesimal path $\infty$-groupoid of $M$ to an infinitesimal version of the $\infty$-groupoid which describes equivalences between the fibres of $x$.
In this paper we make this picture rigorous and justify that it yields the correct notion of connections.
For $p = 1$ it matches the classical case described in Section~\ref{sec: non-technical overview} above:
the free Lie algebroid $Q^{{(1)}} TM$ can be seen as the $1$-truncation of the infinitesimal path groupoid of $M$, and	the Atiyah algebroid $\At(P)$ is exactly the object which describes infinitesimal equivalences between fibres of $P$. For $p>1$, the $p$-form connections we define are necessarily fake flat since our connections are really built as parallel transports along infinitesimal paths and their iterated infinitesimal homotopies (see also Example~\ref{eg: cons on trivial crossed module bundle},  Remark~\ref{rmk: fake flatness and adjusted connections} and the discussion of Section~\ref{sec: de rham description}).
The main results of this paper can be summarised as follows.

\paragraph{Derived differential geometry and Atiyah $\boldsymbol{L_\infty}$-algebroids.}

The first part of this paper develops a theory which allows us to form and compute the analogue of the Atiyah algebroid for a principal $\infty$-bundle $x \colon M \longrightarrow \rmB G$ in maximum generality.
We achieve this by passing to a derived differential geometric setting: using differential graded (dg) $C^\infty$-rings we set up a theory of $\infty$-sheaves on the $\infty$-category $\InfMfd$ of \textit{inf-manifolds}, a version of derived infinitesimally thickened manifolds we introduce, with the main aim of obtaining geometric objects which have a good deformation theory.

To that end, given a map $x \colon M \longrightarrow X$ from a smooth manifold to a suitable $\infty$-sheaf $X$, we first construct the tangent complex $T_x X$ on $M$.
We do this in a way that avoids using the cotangent complex and thus having to dualise (Definition~\ref{def:tangent}). 
For example, this allows to identify $T_e\cDiff(M)$ with the infinite-dimensional vector space of vector fields on $M$.
In the case where $X$ satisfies a certain Schlessinger condition and admits a well-behaved tangent complex for all maps from smooth manifolds (Definitions~\ref{def:schlessinger} and~\ref{def:qcoh def thy}), we say that $X$ \textit{has deformation theory}, or is a \textit{dt stack}.
Our first main result is Theorem~\ref{thm:formal stack vs lie}, which is the key to forming Atiyah algebroids of higher bundles. A slightly reduced version of Theorem~\ref{thm:formal stack vs lie} can be stated as

\begin{theorem}
	\emph{(Theorem~\ref{thm:formal stack vs lie})}
	\label{thm:formal stack vs lie -- intro}
	Let $M$ be a smooth manifold and let $x \colon M \longrightarrow X$ be a map to a stack with deformation theory.
	The relative tangent complex
	$$
		\label{eq: relative tangent complex -- intro}
		T(M/X) \coloneqq \hofib (TM \longrightarrow T_x X)\ ,
	$$
	which is a priori a dg $C^\infty(M)$-module, has a natural lift to a complete $L_\infty$-algebroid on $M$.
	This determines a right adjoint functor of $\infty$-categories
	$$
		T(M/-) \colon \DT_{M/} \longrightarrow L_\infty\Agd_M^{\cpl, \infty} \ ,
	$$
	from stacks with deformation theory under $M$ to complete $L_\infty$-algebroids on $M$.
\end{theorem}

In Theorem~\ref{thm:formal stack vs lie} we also show that the left adjoint of $T(M/-)$ is fully faithful and identify its essential image, as well as the global functions on any $\infty$-sheaf in this essential image.
We then introduce 

\begin{definition}
	Given a smooth manifold $M$ and a map $x \colon M \longrightarrow X$ to a stack with deformation theory, the $L_\infty$-algebroid $T(M/X)$ is the \textit{Atiyah $L_\infty$-algebroid} of $x$.
\end{definition}

Along the way we develop several important results about $\infty$-sheaves on $\InfMfd$, including formal versions of the Inverse Function and Submersion Theorems (Propositions~\ref{prop:inverse_function} and~\ref{prop:subm theorem}, respectively), as well as formally smooth maps (Definition~\ref{def:formally smooth}), formally smooth $\infty$-groupoids (Definition~\ref{def:form smooth oo-groupoid}) and their deformation theories.

\begin{remark}[\textbf{On the use of derived geometry}]
	It may seem unnecessary to pass to a fully derived geometric setting if one is interested only in $p$-form connections on principal bundles with a smooth, non-derived, higher structure group.
	Instead, for bundles with non-derived structure group, one might expect that the mere inclusion of infinitesimals in degree zero, i.e.~a `synthetic' setting, might suffice to construct Atiyah algebroids for higher principal bundles, and hence to define $p$-form connections.
	However, a key ingredient to Theorem~\ref{thm:formal stack vs lie -- intro}, and hence to our formalism, is a family version of the Lurie--Pridham Theorem in derived deformation theory~\cite{Lurie:DAG-X, Pridham:Unifying_DDTs, Nuiten:HoAlg_for_Lie_Algds}.
	We are not aware of a formalism for such a theorem which avoids derived techniques, even when restricted to higher but non-derived moduli problems (we are not saying that this would be impossible).
	In fact, introducing degree zero infinitesimals already leads to a need for a derived tangent complex: the inclusion of the first-order infinitesimal line into a higher-order one induces an isomorphism on classical tangent spaces, but not on derived tangent complexes.
	Additionally, passing to a derived geometric framework enlarges the scope of our theory: it also allows to define connections on bundles classified by derived stacks, like complexes of vector bundles.
	\qen
\end{remark}

\paragraph{Moduli of geometric structures on manifolds, globally.}

Theorem~\ref{thm:formal stack vs lie -- intro} contains, in particular, a classification of formal moduli problems over a manifold $M$ in terms of $L_\infty$-algebroids on $M$.
Indeed, if we view a dt stack $X$ as a moduli space of geometric data, then a map $x \colon M \longrightarrow X$ is an $M$-point in this moduli space.
By the internal hom adjunction, the same data is encoded by a map $* \longrightarrow \Hom(M,X)$, which is a single point in the moduli stack of $X$-structures on $M$.
However, we need to remember that some of the deformations of the map $x \colon M \longrightarrow X$ come purely from permuting points in $M$ by infinitesimal diffeomorphisms.
In this global picture that is manifested by the canonical action of the diffeomorphism group $\cDiff(M)$ on $\Hom(M,X)$.
We show that the moduli problem defined by the induced morphism $$x^\dashv \colon * \longrightarrow \Hom(M,X) \, \big/ \, \cDiff(M)$$ to the quotient is the correct global version of the moduli problem associated to $x \colon M \longrightarrow X$. Again a reduced version of our actual result can be stated as

\begin{theorem}
	\emph{(Theorem~\ref{thm:global})}
	\label{thm:global -- intro}
	Let $x \colon M \longrightarrow X$ be a map from a smooth manifold to a stack with deformation theory.
	\begin{myenumerate}
		\item The quotient $\Hom(M,X)/\cDiff(M)$ is a stack with deformation theory.
		
		\item The moduli problem associated to $x^\dashv \colon * \longrightarrow \Hom(M,X)/\cDiff(M)$ is controlled by an $L_\infty$-algebra with a canonical map to the Lie algebra of vector fields on $M$.
		It is equivalent over this Lie algebra to the global sections of the Atiyah $L_\infty$-algebroid $T(M/X)$ of $x$, which controls the local moduli problem (see Theorem~\ref{thm:formal stack vs lie -- intro}).
	\end{myenumerate}
\end{theorem}

Theorem~\ref{thm:global -- intro} has the following consequence, which provides further justification to view $T(M/X)$ as a derived generalisation of the classical  Atiyah Lie algebroid:
note that the loop stack of $x^\dashv \colon * \longrightarrow \Hom(M,X)/\cDiff(M)$ agrees with the group stack $\cSym(x)$ of derived higher symmetries of $x$.
Theorem~\ref{thm:global -- intro} then implies, in particular, that the $L_\infty$-algebra of $\cSym(x)$ consists precisely of the global sections of the Atiyah $L_\infty$-algebroid and comes with an appropriate morphism to the vector fields on $M$.
As in the classical case, we thus observe that the finite symmetries of $x$ are global objects on $M$, whereas upon passing to infinitesimal symmetries they enhance to local objects, amounting to the enhancement of structure from an $L_\infty$-algebra to an $L_\infty$-algebroid.
It is this enhancement of structure of infinitesimal symmetries which  allows us to define higher-form connections.

So far little is known about how to explicitly compute the $L_\infty$-algebra associated to a formal moduli problem.
Only passing in the opposite direction is computationally understood and explicit.
Here we make partial progress on this general problem by computing the induced brackets on the cohomology $\rmH^\bullet( T(*/Y))$ in all degrees via an intricate interaction between deformation theory and Whitehead brackets (Proposition~\ref{prop:bracket on homotopy}), where $y \colon * \longrightarrow Y$ is a point of a stack with deformation theory (in particular, we may take $Y = \Hom(M,X)/\cDiff(M)$).
This is likely to be of wider interest.

\paragraph{Higher-form connections and their curvatures.}

Now that we have a general notion of the Atiyah algebroid in place, the next step in making precise the heuristic perspective on $p$-form connections sketched at the beginning of this section is to find the order~$p$ truncation of $TM$:
at the level of cocommutative dg coalgebras (cdgcs), the $L_\infty$-algebroid $TM$ corresponds to its Chevalley--Eilenberg cdgc \smash{$\Sym_{\RN}(TM[1])$}.
This has a straightforward order~$p$ truncation \smash{$\Sym_{\RN}^{\leqslant p}(TM[1])$}: it is obtained by projecting away all tensor powers greater than $p$.
Dually, its Chevalley--Eilenberg cdga corresponds to the truncated de Rham algebra $\Omega^{\leqslant p}(M)$.
It turns out that there is another $L_\infty$-algebroid, denoted $Q^{{(p)}} TM$, which corresponds to \smash{$\Sym_{\RN}^{\leqslant p}(TM[1])$} (up to equivalence at the level of cdgcs).
We understand this $L_\infty$-algebroid as the order~$p$ truncation of $TM$. 
Dually, in terms of dg algebras, the Chevalley--Eilenberg algebra of $Q^{{(p)}} TM$ simply corresponds to the order~$p$ truncation $\Omega^{\leqslant p}(M)$ of the de Rham complex.

Given any morphism $x \colon M \longrightarrow X$ from a smooth manifold to a stack with deformation theory, we can now make the above heuristics fully rigorous through

\begin{definition}
	\label{def:l-cons -- intro}
	The \emph{$\infty$-groupoid of $p$-form connections on $x \colon M \longrightarrow X$} is the mapping space
	\begin{equation}
		\Con_p(x) \coloneqq L_\infty\Agd_M^\infty \big( Q^{(p)} TM,\, T(M/X) \big)
	\end{equation}
	in the $\infty$-category of $L_\infty$-algebroids on $M$.
\end{definition}

We then present several fundamental aspects of this notion of connection.
We define the field strength $(p{+}1)$-form of a $p$-form connection as the obstruction of the $p$-form connection to lift to a $(p{+}1)$-form connection (Definition~\ref{def:curvature of l-conn}) and check that this recovers the classical definition of field strength of a connection.
We also show that null homotopies of the field strength $(p{+}1)$-form correspond exactly to extensions of the $p$-form connection to a $(p{+}1)$-form connection (Proposition~\ref{st: extending l-con to (l+1)-con}), adding to the interpretation of the fake-curvature condition of higher connections, and show that when $x \colon M \longrightarrow \rmB G$ is a principal $\infty$-bundle its higher gauge group acts canonically on the space of $p$-form connections on $x$ (Definition~\ref{def: G-buns with l-con} and Remark~\ref{rmk: gauge action on l-cons}). We furthermore make precise the interpretation in terms of infinitesimal higher parallel transports by giving an equivalent description of $p$-form connections in terms of an order~$p$ truncation of the de Rham stack.
We expect that this notion of higher connections also generalises straightforwardly beyond the case where $M$ is a smooth manifold.

\paragraph{Computing spaces of connections on higher U(1)-bundles.}

In the second part of this paper we explore applications of our theory of $p$-form connections on higher bundles, and illustrate the relationship between derived symmetries of geometric structures and possible weakened versions of the notion of a connection.

First, we provide an independent verification of the derived geometric framework of connections in the case of higher $\rmU(1)$-bundles.
These are also called $(n{-}1)$-gerbes and are classified by maps $$x \colon M \longrightarrow \rmB^n\rmU(1)$$ to a smooth delooping of the circle group $\rmU(1)$.
For this special type of principal $\infty$-bundles there already exists a fully developed theory of higher-form connections which, crucially, is entirely independent of the derived geometric formalism we set up here.
Concretely, an $n$-form connection on $x \colon M \longrightarrow \rmB^n \rmU(1)$ is given by enhancing the map $x$ to a differential cocycle, which we model using \v{C}ech--Deligne cocycles (for background, see for instance~\cite{Brylinski:Loops_and_GeoQuan,Gajer:Geo_of_Deligne_Coho,  FSS:Cech_diff_char_classes_via_L_infty, Szabo:2012hc,Bunk:Gerbes_in_Geo_and_FT}).
In a similar vein, one can also endow $x$ with $p$-form connections for any $p \in \NN$.
It was observed in~\cite{BS:Higher_Syms_and_Deligne_Coho} that in this \v{C}ech--Deligne approach, for each $(n{-}1)$-gerbe on $M$ there is an entire $\infty$-groupoid of $p$-form connections, which we denote by \smash{$\Con_p^{\mathrm{\check{C}D}}(x)$}.
We compute the Atiyah $L_\infty$-algebroid $T(M/\rmB^n\rmU(1))$ in Corollary~\ref{st:At-L_oo-agd of n-gerbe from first principles}, strictify it to a dg Lie algebroid on $M$ in Theorem~\ref{st:strictification of At(CG)}, and then  by an explicit computation of mapping spaces of $L_\infty$-algebroids we prove

\begin{theorem}
	\emph{(Theorem~\ref{st:finite and infinitesimal l-cons on n-gerbes})}
	\label{st:finite and infinitesimal l-cons on n-gerbes -- intro}
	Let $x \colon M \longrightarrow \rmB^n\rmU(1)$ be an $(n{-}1)$-gerbe on a smooth manifold $M$.
	For any $p \in \NN$, there is a canonical equivalence
	\begin{equation}
		\Con_p^{\mathrm{\check{C}D}}(x) \simeq \Con_p(x)
	\end{equation}
	between the $\infty$-groupoids of $p$-form connections on $x$ in the \v{C}ech--Deligne sense and $p$-form connections on $x$ in the derived geometric sense introduced in the present paper.
\end{theorem}

In particular, we obtain an $L_\infty$-algebraic model for differential cohomology.

\paragraph{Higher Courant algebroids and weak connections.}

We conclude by exploring the relationship of the theory developed in this paper to higher Courant algebroids and symmetries of higher $\rmU(1)$-bundles with $p$-form connections. Higher Courant algebroids provide a geometric framework for unifying the diffeomorphism and higher-form gauge symmetries together with duality symmetries, as well as the different flavours of fields, in many supergravity theories. One can associate a higher Courant algebroid to any $(n{-}1)$-gerbe through a choice of $(n{-}1)$-form connection, extending Hitchin's generalised tangent bundle construction for $n=2$~\cite{Hitchin:2003cxu} to all $n\geqslant2$ (Remark~\ref{rem:higherCourantU1}); from this perspective, a choice of extension of the $(n{-}1)$-form connection to an $n$-form connection is a splitting of the short exact sequence of vector bundles in which the higher Courant algebroid sits. This association is crucial to the formulation of semi-classical flux quantization conditions in supergravity, and to higher geometric prequantization.

Of particular relevance to the present work is that the bracket of a higher Courant algebroid has a canonical extension to a Lie $n$-algebra that naturally projects to the Lie $n$-algebra underlying the Atiyah $L_\infty$-algebroid of a corresponding $(n{-}1)$-gerbe (Theorem~\ref{prop:Linftyproj}). Compared to the Atiyah $L_\infty$-algebra which controls the symmetries of an $(n{-}1)$-gerbe on its own, we argue that the Courant $L_\infty$-algebra controls the symmetries of the $(n{-}1)$-gerbe which also preserve a choice of $(n{-}1)$-form connection~(Conjectures~\ref{conj:Coursym} and~\ref{conj:CourAtsym}). For $n=2$, similar statements were  proven in~\cite{FRS}. Further relations of higher Courant algebroids with higher gerbes and dg geometry can be found in~\cite{LRU:TD_and_ExGeo_through_syms_of_dgMfds, Cueca:Geo_of_graded_cotangent_bundles}.

Because in this case the objects whose deformation theory we aim to compute contain differential form data, they do not have deformation theory and so only the global version of the moduli problem exists.
That is, we have to describe them globally as maps $y \colon * \longrightarrow \Hom_v(M,X)/\cDiff(M) = Y$, where $\Hom_v$ denotes a `concretified' version of the internal hom whereby families of differential forms on $M$ over a parameter space $U$ act only on tangent vectors along $M$ and not $U$.
We present a general way to compute the deformation theory of such maps in cases where $X$ has an $\bbE_\infty$-group structure (Proposition~\ref{st: computing sym(a) in A//Diff(M) case}).
This results in $L_\infty$-algebras $T(*/Y)$ over the Lie algebra of vector fields on $M$ which do not stem from $L_\infty$-algebroids.
Among other examples, we compute the deformation theory of a gerbe with $1$-form connection (i.e.~with a connective structure) and prove that the resulting Lie 2-algebra is equivalent to the Lie 2-algebra of global sections of the Courant algebroid associated to the gerbe (Example~\ref{eg: deforming gerbe with connective structure}), thus recovering the main theorem of~\cite{Collier:Inf_symmetries_of_gerbes}.

In the spirit of Definition~\ref{def:l-cons -- intro} we can, however, still consider order~$p$ splittings of the map of $L_\infty$-algebras $T(*/Y) \longrightarrow TM$ which are $C^\infty(M)$-linear and view these as a weakened notion of $p$-form connections.
We show that a $\rmU(1)$-bundle with connection admits a $1$-form connection in this weakened sense if and only if the original connection is flat, in which case the additional weak $1$-form connection agrees with the original connection (Example~\ref{eg: deforming U(1)-bundle with connection}).
We also show that a weak $1$-form connection on a gerbe with connective structure is the same as an enhanced curving in the sense of~\cite{TD:Chern_corr_for_higher_PrBuns} (Example~\ref{eg: deforming gerbe with connective structure}).
An enhanced curving canonically splits into a standard curving and a symmetric $2$-tensor on $M$, which in special cases will be a Riemannian or Lorentzian metric on $M$.
In particular, weak $1$-form connections on gerbes with connective structure may provide a unification of the $B$-field and metric in NSNS supergravity when the \v{S}evera class of the associated Courant algebroid is integral (see also~\cite{TD:Chern_corr_for_higher_PrBuns, BS:Higher_Syms_and_Deligne_Coho}).


\subsection{Outline of the paper}


For the convenience of the reader, as an aid in navigating this paper we provide a brief overview of the structure of its remaining contents.

In Section~\ref{sec: Rev of DefThy} we collect material on $L_\infty$-algebroids on manifolds and the geometry of derived stacks which serves as the foundation for our main results in later sections.
This section consists of both old and new material.
It includes a review of the homotopy theory of $L_\infty$-algebroids and dg-$C^\infty$-rings, focussing on square zero extensions.
In Section~\ref{sec:derived C^oo-geometry} we then introduce the $\infty$-category $\InfMfd$ of derived and infinitesimally thickened manifolds which serves as the foundation for our theory of derived stacks.
Section~\ref{sec:infinitesimal} covers tangent complexes and introduces stacks with (quasi-coherent) deformation theory, and Section~\ref{sec:Examples of stacks with defthy} contains various constructions and examples of stacks with (quasi-coherent) deformation theory.

In Section~\ref{sec:Lie diff} we recall key results about formal moduli problems over commutative differential graded algebras from~\cite{Nuiten:Thesis, Nuiten:Koszul_duality_for_Lie_algebroids} and build on these to establish the formal differentiation of maps $x \colon M \longrightarrow X$ from a smooth manifold to a stack with deformation theory.
Such differentiation results in an $L_\infty$-algebroid structure on the relative tangent complex $T(M/X)$.
We also provide a formal integration, thus identifying the essential image of the differentiation $\infty$-functor.
Moreover, we show that the global sections of the $L_\infty$-algebroid $T(M/X)$ agree with the $L_\infty$-algebra of the group of derived symmetries of $x \colon M \longrightarrow X$.
We include several computational tools and constructions in Section~\ref{sec:differentiation examples}, and show how the brackets on the cohomology of the tangent complex of a formal moduli problem can be computed explicitly via methods from rational homotopy theory.
Proofs of the main theorems in this section are deferred to Appendices~\ref{sec:lie diff proof} and~\ref{sec:global sections proof}.

In Section~\ref{sec: Higher connections} we recall the notion of an order $p$ morphism of $L_\infty$-algebroids from~\cite{Nuiten:HoAlg_for_Lie_Algds} and present our definition of a $p$-form connection on $L_\infty$-algebroids $\frg$ over a unital commutative dg algebra $A$ (Definition~\ref{def:k-cons general}).
We are then able to define the space of $p$-form connections on a morphism $x \colon M \longrightarrow X$ to a stack with deformation theory as the space of order $p$ connections on its associated Atiyah $L_\infty$-algebroid $T(M/X)$ (Definition~\ref{def: Connection}).
The section also contains model categorical tools for computing spaces of order $p$ connections and includes several examples.
With the definitions established, we extract the curvature of an order $p$ connection as its failure to already be an order $(p{+}1)$-connection (Definition~\ref{def:curvature of l-conn} and Proposition~\ref{st: extending l-con to (l+1)-con}).
Finally, we show that spaces of order $p$ connections are functorial with respect to automorphisms of the map $x \colon M \longrightarrow X$, and provide an equivalent description of these spaces in terms of stacks which makes precise the intuition that our higher-form connections really correspond to infinitesimal higher parallel transports.

Section~\ref{sec: cons on higher U(1)-bundles} contains two models of the Atiyah $L_\infty$-algebroid of higher $\rmU(1)$-bundles $x \colon M \longrightarrow \rmB^n\rmU(1)$, one as an $L_\infty$-algebroid concentrated in degrees $0$ and $-n +1$ with a non-trivial $(n{+}1)$-bracket, and one as a dg Lie algebroid concentrated in degrees $-n+1, \ldots, 0$ (Corollary~\ref{st:At-L_oo-agd of n-gerbe from first principles} and Theorem~\ref{st:strictification of At(CG)}).
The proof that these models are equivalent occupies Appendix~\ref{app:Pf of strictification}.
Subsequently, we establish in Theorem~\ref{st:finite and infinitesimal l-cons on n-gerbes} that the space of $p$-form connections on any $n$-bundle $x \colon M \longrightarrow \rmB^n \rmU(1)$ in our new derived geometric sense is canonically equivalent to the space of $p$-form connections on $x$ as obtained from the classical \v{C}ech--Deligne description of differential cocycles.
The proof of this theorem is an involved computation, which we carry out in Appendix~\ref{app:Pf of equivalence of connection spaces} and which relies on an explicit fibrant simplicial resolution of the target $L_\infty$-algebroid established in Appendix~\ref{app:simplicial resolution}.

Section~\ref{sec: Higher Courant algebroids} concerns the higher derived symmetries of higher $\rmU(1)$-bundles with fixed $p$-form connections and the question of which additional connection-like structures these symmetries induce.
Crucially, since we fix part of a connection as part of the data here, we are not able to differentiate the classifying map $x \colon M \longrightarrow \rmB^{n-p}\,\rmB_\nabla^p \rmU(1)$ of the geometric structure as before, but only the point $x^\dashv$ which $x$ defines in an (appropriately defined) moduli stack of $\rmU(1)$-principal $n$-bundles with $p$-form connection on $M$ (the key result being Proposition~\ref{st: computing sym(a) in A//Diff(M) case}).
Thus the morphism $x$ does not admit connections in the sense of Section~\ref{sec: Higher connections}, but we can still define a relaxed notion of connection using $x^\dashv$ instead.
In the case where $x$ classifies a gerbe with connective structure (i.e.~with a 1-form connection), these weakened 1-form connections on $x^\dashv$ consist of a curving on the gerbe and a symmetric 2-tensor on the base manifold $M$ (compare the  results in~\cite{TD:Chern_corr_for_higher_PrBuns}).
As first steps towards generalising these results, we survey higher Courant algebroids (Section~\ref{sub:highercourant}) as well as a higher version of Hitchin's generalised tangent bundle construction (Remark~\ref{rem:higherCourantU1}), and use the homotopy transfer theorem to generalise Zambon's $L_\infty$-algebras of higher Courant algebroids~\cite{Zambon:2010ka} to associate $L_\infty$-algebras to $\rmU(1)$-principal $n$-bundles with $(n{-}1)$-form connections (Section~\ref{sub:derivedbracket} and Proposition~\ref{prop:cour-intermediate}).
Finally, we show that there is a projection from Zambon's Courant $L_\infty$-algebra associated to such an $n$-bundle with $(n{-}1)$-form connection onto the $L_\infty$-algebra of global sections of the $n$-bundle's Atiyah $L_\infty$-algebroid.
This leaves a few further questions open, which we include as Conjectures~\ref{conj:Coursym}, \ref{conj:CourAtsym} and~\ref{conj:hetSuGra}.


\subsection{Acknowledgements}


We would like to thank David Carchedi, Owen Gwilliam, Alexander Schenkel and Pelle Steffens for helpful discussions. We are also grateful to Miquel Cueca and Urs Schreiber for comments on the manuscript.
SB is grateful to the Isaac Newton Institute for support through an INI Retreat.
LM gratefully acknowledges support of the Simons Collaboration on Global Categorical Symmetries. 
Research at Perimeter Institute is supported in part by the Government of Canada through the Department of Innovation, Science and Economic Development and by the Province of Ontario through the Ministry of Colleges and Universities. The Perimeter Institute is in the Haldimand Tract, land promised to the Six Nations. 
JN was supported by the CNRS through the program ``Premier soutien -- Jeunes chercheurs et jeunes chercheuses'', and by the ANR project ANR-24-CE40-5367-01 (``LieDG'').
RJS is a member of COST Actions CaLISTA CA21109 and THEORY-CHALLENGES CA22113 supported by COST (European Cooperation in Science and Technology).



\subsection{Notation and conventions}


This paper is engulfed in the language of derived geometry, $\infty$-categories and model categories. To aid the reader in following our exposition, below we summarise some of the underlying concepts as well as our notation and conventions.

\begin{myitemize}
	\item For $n\in\NN$, we write $\Sigma_n$ for the symmetric group of degree~$n$, and $\mathrm{Sh}(k,n-k)\subseteq\Sigma_n$ for the subset of $(k,n{-}k)$-shuffles with $0\leqslant k\leqslant n$, i.e. permutations $\sigma\in\Sigma_n$ such that $\sigma(1)<\cdots<\sigma(k)$ and $\sigma(k+1)<\cdots<\sigma(n)$. We put $\mathrm{Sh}(0,n) = \mathrm{Sh}(n,0) = \{\id\}$. We write $\epsilon(\sigma)$ for the sign of a permutation $\sigma$.

	\item We use cohomological conventions, so that a \emph{complex} will always mean a \emph{cochain complex}. Consequently, objects like the de Rham complex are in non-negative degrees while our main $L_\infty$-algebroids of interest will be in non-positive degrees.
	Explicitly, for a graded vector space $C$ we use the desuspension convention $(\sfOmega C)^i = (C[-1])^i = C^{i-1}$. We write $|c|\in\mathbbm{Z}$ for the degree of a homogeneous element~$c\in C$; that is, if $c\in C^i$ then $|c|=i$.
	
	\item We abbreviate `differential graded' to dg, `differential graded Lie algebra' to dgla, `commutative differential graded algebra' to cdga and `cocommutative differential graded coalgebra' to cdgc.
	
	\item We write $\bbDelta$ for the simplex category of finite ordinals $[p] = \{0,1,\dots,p\}$ and order-preserving maps. A simplicial object in a ($\infty$-)category $\scC$ is a functor $\bbDelta^\opp\longrightarrow\scC$, where the superscript ${}^\opp$ indicates the opposite category; a cosimplicial object is a functor $\bbDelta\longrightarrow\scC$. A simplicial set is a simplicial object in the category of sets $\mathscr{S}\mathrm{et}$, and the category of simplicial sets is denoted by $\sSet = \mathrm{Fun}(\bbDelta^\opp,\mathscr{S}\mathrm{et})$. For $p\in\NN_0$, the standard $p$-simplex is the simplicial set $\Delta^p=\bbDelta(-,[p])$ induced by the Yoneda embedding. We write $\Lambda_i^p\subset\partial\Delta^p\subset\Delta^p$ for the $i$-th horn. A Kan complex is a fibrant simplicial set, i.e. a simplicial set satisfying horn-filling conditions.

	\item An $\infty$-category is a simplicial set satisfying the inner horn-lifting conditions (see e.g.~\cite{Cisinski:HCats_HoAlg,Lurie:HTT}). If $\scC$ is an $\infty$-category and $c, d \in \scC$, we write $\scC(c, d)$ for the $\infty$-categorical mapping space in $\scC$.
	For an ordinary category $\scM$ with objects $x, y \in \scM$, we write $\scM(x, y)$ for the set of morphisms $x \longrightarrow y$ in $\scM$.
	
	\item A model category is a category with three specified classes of morphisms --- fibrations, cofibrations and weak equivalences --- which satisfy axioms reminiscent of Serre fibrations, relative cell complexes and weak homotopy equivalences of topological spaces. If $\scM$ is a model category, an object $x\in\scM$ is called fibrant if the unique map from $x$ to the terminal object is a fibration and cofibrant if the unique map from the initial object to $x$ is a cofibration. For any object $x\in\scM$, there is a cofibrant object $Q(x)$ and an acyclic fibration $Q(x)\longrightarrow x$ (a weak equivalence), called a cofibrant resolution of $x$. There is dual notion called fibrant resolution.
	
	\item If $\scM$ is a model category and $x,y\in\scM$, we write $\Map_\scM(x, y)$ for the model categorical mapping space. We sometimes write $\scM^\infty = \scM[W^{-1}] = \Ho\scM$ for its $\infty$-categorical localisation at the weak equivalences.
	Then there is an equivalence of spaces $\scM^\infty(x,y) \simeq \Map_\scM(x,y)$, for each $x, y \in \scM$.
	
	\item Throughout this paper we work over the field of real numbers $\RN$. Let $A$ be a cdga over $\RN$.
	We write $\Mod_A^\dg$ for the category of unbounded dg $A$-modules and \smash{$\Mod_A^\infty = \Mod_A^\dg[W^{-1}]$} for its $\infty$-categorical localisation at the quasi-isomorphisms.
	The morphism sets in \smash{$\Mod_A^\dg$} are denoted $\Mod_A^\dg(C,D)$, whereas we write the hom complex as $\Hom_A(C,D)$.
	The $\infty$-categorical mapping space is $\Mod_A^\infty(C,D)$.	
	
	\item We write $\Mfd$ for the category of smooth manifolds and smooth maps. For $M\in\Mfd$, we write $\mathrm{Open}(M)$ for the category of open subsets of $M$ and inclusions; the category $\mathrm{Open}(M)$ is a poset.
	
	\item We denote exponential objects inside a category or an $\infty$-category $\scC$ by $Y^X$.
	
	\item We use the $C^\infty$-version of the Serre--Swan theorem to identify vector bundles over a manifold $M$ with finitely generated projective modules over the ring of smooth functions $C^\infty(M)$.
	
	\item We sometimes use the algebraic notational convention for the tangent bundle of a manifold $M\in\Mfd$:
	we write $TM$ to denote the tangent complex of the $C^\infty$-ring $C^\infty(M)$; that is, $TM$ is the $C^\infty(M)$-module of smooth vector fields on $M$.
	
	\item When working with dg $C^\infty$-rings, we use $\otimes$ to denote the tensor product of the underlying cdgas and $\amalg$ to denote the coproduct, or $C^\infty$-algebraic tensor product, of the dg $C^\infty$-rings.
		
	\item In this paper we use three categories of $L_\infty$-algebroids over a cgda $A$:
	\begin{myitemize}
		\item[$\diamond$] the category $L_\infty\Agd_A^\dg$ of $L_\infty$-algebroids with morphisms the bracket-preserving maps of the underlying complexes,
		
		\item[$\diamond$] the category $L_\infty\Agd_A$ of $L_\infty$-algebroids with their $\infty$-morphisms, which are the same as morphisms of the associated cdgcs satisfying a certain compatibility with the anchor and an $A$-linearity condition (see Definition~\ref{def:oo-morphism of L_oo-algebroids} and Remark~\ref{rmk:constructing cdgc-maps into CE cdgc} for details), and
		
		\item[$\diamond$] the $\infty$-categorical localisation $L_\infty\Agd_A^\infty = L_\infty\Agd_A^\dg[W^{-1}] \simeq L_\infty\Agd_A[W'{}^{-1}]$, where the first localisation is at the quasi-isomorphisms and the second is at the $\infty$-equivalences.
	\end{myitemize}
	 
	\item We write $\DK \colon \mathrm{Ch}_{\leqslant0} \longrightarrow \Ab_\Delta$ for the Dold--Kan correspondence, where $\mathrm{Ch}_{\leqslant0}$ is the category of non-positively graded cochain complexes and $\Ab_\Delta$ is the category of simplicial abelian groups.
	For a more compact notation, if $\DK$ is applied to a complex which is not concentrated in non-positive degrees, we will implicitly let $\DK$ first truncate the complex at level zero and then apply the usual Dold--Kan construction.
		
	\item The $\infty$-category $\scS$ of spaces (or equivalently $\infty$-groupoids) is the homotopy coherent nerve of the simplicially enriched category of Kan complexes. It arises as the $\infty$-categorical localisation of the Kan--Quillen model structure on simplicial sets.
		
	\item If $\scC$ is an $\infty$-category and $x \in \scC$, we denote by $\scC_{/x}$ (resp. $\scC_{x/}$) the  $\infty$-category of objects in $\scC$ over (resp. under) $x$, i.e. with a map to (resp. from) $x$.

	\item In this paper, by a stack we will always mean an $\infty$-stack. We write $\Sh_\infty$ (resp. $\PSh_\infty$)  for the $\infty$-category of stacks (resp. prestacks) on a site.

	\item We denote $\infty$-morphisms between $L_\infty$-algebroids by a squigly arrow $\longleadsto$. We write $\simeq$ for weak equivalences of objects, and $\cong$ for isomorphisms.
	
	\item We write $\int\mathcal{F}\longrightarrow \scC$ for the left fibration classified by an $\infty$-functor $\mathcal{F}:\scC\longrightarrow\scS$.
\end{myitemize}



\section{Deformation theory}
\label{sec: Rev of DefThy}


To introduce our definition of connections on $\infty$-bundles, we need a formalism of derived $C^\infty$-geometry which admits  a classification theorem of formal moduli problems over manifolds.
In this section we begin by developing some foundational material on $C^\infty$-rings and $L_\infty$-algebroids as well as their homotopy theory.
We introduce a notion of infinitesimally thickened manifolds suitable for our purposes, $\infty$-sheaves and formally smooth $\infty$-groupoids thereon, as well as their tangent complexes. This follows ideas from derived algebraic geometry \cite{GR:Study_Derived_Geom,TV:Homotopical_Alg_Geom_II}, but some additional care is needed in the differential-geometric setting.
In Section~\ref{sec:Lie diff} below we build on these foundations to prove classification results for deformation problems over manifolds.


\subsection{Homotopy theory of $L_\infty$-algebroids}
\label{sec: HoThy of L_oo-agds}


The purpose of this section is to give a brief review of the homotopy theory of $L_\infty$-algebroids over a smooth manifold $M$. To be able to use model-categorical tools, it will be useful to first consider $L_\infty$-algebroids over $M$ from a purely algebraic perspective in terms of their $C^\infty(M)$-modules of global sections. Since vector bundles over $M$ can be identified with finitely generated projective $C^\infty(M)$-modules by the $C^\infty$-version of the Serre--Swan theorem~\cite[Chapter 11]{Nestruev:Smooth_mfds_and_observables}, this purely algebraic approach specialises to more geometric definitions of $L_\infty$-algebroids over $M$ in terms of vector bundles, or more generally sheaves; we will come back to this point in Definition~\ref{def:complete}.

Let us fix a cdga $A$ over $\RN$ and let $T_A$ denote its complex of derivations.
We will write $\pounds_X\colon A\longrightarrow A$ for the action of a derivation $X\in T_A$ on $A$. The commutator bracket and pointwise multiplication by an element in $A$ endow the complex $T_A$ with the structure of a dg Lie algebra and a dg $A$-module.
When $A = C^\infty(M)$, this reproduces the usual Lie algebra of vector fields on $M$ which we also denote by $TM = T_{C^\infty(M)}$; then $\pounds_X$ is the Lie derivative along a vector field $X\in TM$.

\begin{definition}\label{def:L-infinity algebroid}
	Let $A$ be a cdga over $\RN$.
	An \emph{$L_\infty$-algebroid} over $A$ is a dg $A$-module $\frg$ together with
	\begin{myenumerate}
		\item a map of dg $A$-modules $\rho\colon \frg\longrightarrow T_A$ called the \emph{anchor map}, and
		
		\item an $\RN$-linear $L_\infty$-algebra structure on $\frg$ such that $\rho$ is a \emph{strict} map of $L_\infty$-algebras (in particular, all $n$-ary brackets with $n\geqslant 3$ are sent to zero) and such that for all $x_i\in \frg$ and $a\in A$:
		\begin{align*}
			[x_1, a\, x_2] &= (-1)^{|a|\,|x_1|}\, a\, [x_1, x_2]+ \pounds_{\rho(x_1)}(a)\, x_2 \ , \\[4pt]
			[x_1, \dots, x_{n-1}, a\, x_n] &=(-1)^{|a|\,(|x_1|+\dots+|x_{n-1}|)}\,a\, [x_1, \dots, x_n] \quad \text{for }n\geqslant 3 \ .
		\end{align*}
	\end{myenumerate}
	An $L_\infty$-algebroid is \emph{fibrant} if its anchor map $\rho:\frg\longrightarrow T_A$ is a fibration, i.e. it is  surjective.
	
	A (\emph{strict}) \emph{map of $L_\infty$-algebroids} is  a map of dg $A$-modules $\frg\longrightarrow\frh$ that commutes with the anchor maps and preserves the $L_\infty$-structure.
	We denote the category of $L_\infty$-algebroids over $A$ and strict maps by $L_\infty\Agd_A^\dg$.
	A morphism of $L_\infty$-algebroids is  a \emph{quasi-isomorphism} if its underlying map of complexes is a quasi-isomorphism.
\end{definition}

Recall that an $L_\infty$-structure is a collection of $n$-ary brackets $$[-,\dots,-]:\midwedge_\RN^n\, \frg\longrightarrow\frg$$ of degree $2{-}n$ satisfying higher homotopy Jacobi identities for all $n$.
We refer to~\cite{KS:Intro_to_L_oo-algebras} for background on $L_\infty$-algebras.

\begin{example}
	An $L_\infty$-algebroid whose $n$-ary brackets vanish for all $n\geqslant 3$ is a \emph{dg Lie algebroid}. One can show that every $L_\infty$-algebroid is connected by a zig-zag of quasi-isomorphisms of $L_\infty$-algebroids to a dg Lie algebroid~\cite[Corollary 3.10]{Nuiten:HoAlg_for_Lie_Algds}.
	\qen
\end{example}

In Definition~\ref{def:L-infinity algebroid} we allow for $L_\infty$-algebroids that are unbounded both in the positive and in the negative direction. However, the positive and non-positive parts play a very different role: the positive part encodes a \emph{derived} structure while the non-positive part encodes a (formally smooth) \emph{stacky} structure. An $L_\infty$-algebroid whose underlying complex is concentrated in degrees $[-(n {-} 1), \dots , 0]$ is sometimes called a \emph{Lie $n$-algebroid}.

\begin{proposition}[{\cite[Theorem 3.1]{Nuiten:HoAlg_for_Lie_Algds}}]\label{prop:semi-model}
	The category of $L_\infty$-algebroids over $A$ carries a cofibrantly generated left semi-model structure in which a map is a weak equivalence if it is a quasi-isomorphism and a fibration if it is degreewise surjective.
\end{proposition}

\begin{definition}\label{def:Loo algd oo-cat}
	We denote by
	\begin{equation}
		L_\infty\Agd_A^\infty
		= L_\infty\Agd_A^\dg[W^{-1}]
	\end{equation}
	the $\infty$-category obtained as the ($\infty$-categorical) localisation of the category of $L_\infty$-algebroids at the weak equivalences.
\end{definition}

\begin{remark}
	\label{rem:nonpositively graded case}
	When $A=C^\infty(M)$ for a manifold $M$, the category of \emph{non-positively graded} $L_\infty$-algebroids carries a (full) model structure in which the weak equivalences are the quasi-isomorphisms, and the fibrations are maps that are surjective in negative degrees. Indeed, this follows from the transfer argument from~\cite{Vezzosi:Model_structure_on_Lie_algds}: every non-positively graded $L_\infty$-algebroid $\frg$ is fibrant and admits a functorial path object, obtained from the $L_\infty$-algebroid $\frg\boxtimes \Omega^*[\Delta^1]$ from~\cite[Lemma 2.5]{Vezzosi:Model_structure_on_Lie_algds} by taking its canonical truncation ($0$-cycles and elements in negative degree). 
	
	The evident inclusion into all $L_\infty$-algebroids is then a left Quillen functor,\footnote{A \emph{left Quillen functor} is a left adjoint functor that preserves (acyclic) fibrations and weak equivalences.} whose derived unit map is a weak equivalence. Consequently, the model category of non-positively graded $L_\infty$-algebroids over $A=C^\infty(M)$ (with strict maps) describes the full $\infty$-subcategory of $L_\infty\Agd_A^\infty$ on the $L_\infty$-algebroids with vanishing cohomology in positive degrees.
	\qen
\end{remark}

Proposition~\ref{prop:semi-model} implies that the $\infty$-category $L_\infty\Agd_A^\infty$ can be studied using (most of) the usual computational tools from model category theory. For example, the mapping space $L_\infty\Agd_A^\infty\big(\frg, \frh)$ can be computed using a cofibrant replacement of $\frg$ and a fibrant simplicial resolution of $\frh$. In fact, there is a rather explicit model for this mapping space in terms of \emph{$\infty$-morphisms} of $L_\infty$-algebroids.

To describe this, let us start by recalling that any $\RN$-linear $L_\infty$-algebra has an associated \emph{Chevalley--Eilenberg coalgebra}, given by the reduced%
\footnote{The \textit{reduced} symmetric coalgebra arises from the full symmetric coalgebra by projecting out the degree zero part.}
cofree non-unital cocommutative differential graded coalgebra
\begin{equation}
	\ChEil_*(\frg) = \Big(\ol{\Sym}_{\RN}(\frg[1]), \delta_{\ChEil}\Big)\ ,
\end{equation}
with the shuffle coproduct (see for instance~\cite[Equations~(2.28) and~(2.31)]{KS:Intro_to_L_oo-algebras}), and differential $\delta_{\ChEil}$ given onto cogenerators by the differential of $\frg$ as well as all higher brackets~\cite[Proposition~10.1.20]{LodayVallette:Algebraic_Operads}.
Given two $\RN$-linear $L_\infty$-algebras, an $\infty$-morphism $\frg \longleadsto \frh$ is then given by a map of (non-unital) cocommutative dg coalgebras 
\begin{equation}
	\phi\colon \ChEil_*(\frg)\longrightarrow \ChEil_*(\frh)\ .
\end{equation}
Note that without differential $\phi$ is uniquely determined by an $\RN$-linear map $\ol{\Sym}_{\RN}(\frg[1])\longrightarrow \frh$ of degree~$+1$.

Compatibility with the differential then translates into the following condition:
the complex of $\RN$-linear maps $\Hom_{\RN}(\ChEil_*( \frg), \frh\big)$ has the structure of an $L_\infty$-algebra, with $n$-ary bracket $[\phi_1, \dots, \phi_n]_{\Hom,n}$ computed as the composition
\begin{equation}
	\label{eq:hom brackets}
	[\phi_1,\ldots,\phi_n]_{\Hom, n}
	= [-, \dots, -]_{\frh, n} \circ (\phi_1\otimes \dots \otimes \phi_n)\circ \sfDelta^{n-1}_\frg
\end{equation}
with the $n$-ary bracket of $\frh$ and the $(n{-}1)$-fold coproduct of $\ChEil_*(\frg)$.
Then a degree one $\RN$-linear map $\phi\colon \ChEil_*(\frg)\longrightarrow \frh$ is an $\infty$-morphism $\frg \longleadsto \frh$ if and only if it satisfies the \textit{Maurer--Cartan equation} (see for example~\cite[Lemma 2.9]{DHR:Homotopy_Algebras})
\begin{equation}
	\label{eq:MC eqn for cdgc mps into CE cdgc}
	\dd_\Hom \phi + \sum_{n \geqslant 2}\, \frac{[\phi, \dots, \phi]_{\Hom, n}}{n!} = 0\ .
\end{equation}
Here $\dd_\Hom$ is the differential induced on the complex $\Hom_\RN(\frg, \frh)$ from the differentials on $\frg$ and $\frh$.

Finally, we extend the notion of an $\infty$-morphism from $L_\infty$-algebras to $L_\infty$-algebroids in this vein.

\begin{definition}
	\label{def:oo-morphism of L_oo-algebroids}
	Let $\frg$ and $\frh$ be two $L_\infty$-algebroids over a cdga $A$.
	An \textit{$\infty$-morphism} $\phi \colon \frg \longleadsto \frh$ consists of a map of $\RN$-linear cocommutative dg coalgebras
	\begin{equation}
		\begin{tikzcd}[column sep = 0.cm, row sep = 1cm]
			\ChEil_*(\frg)
			\arrow[rd, "\ChEil_*(\rho_{\frg})"{swap}] \arrow[rr, "\phi"]
			& & \ChEil_*(\frh)\arrow[ld, "\ChEil_*(\rho_{\frh})"]
			\\
			& \ChEil_*(T_A) &
		\end{tikzcd}
	\end{equation}
	such that the underlying map of non-unital graded coalgebras $\ol{\Sym}_{\RN}(\frg[1])\longrightarrow \ol{\Sym}_{\RN}(\frh[1])$ descends to an $A$-linear map of graded coalgebras $\ol{\Sym}_{A}(\frg[1])\longrightarrow \ol{\Sym}_{A}(\frh[1])$. The category of $L_\infty$-algebroids over $A$ and $\infty$-morphisms is denoted $L_\infty\Agd_A$.
\end{definition}

\begin{remark}
\label{rmk:constructing cdgc-maps into CE cdgc}
	Equivalently, an $\infty$-morphism of $L_\infty$-algebroids is given by the data of a degree one map of graded $A$-modules $\phi\colon \ol{\Sym}_{A}(\frg[1])\longrightarrow \frh$ such that
	\begin{myenumerate}
		\item the $\RN$-linear map underlying $\phi$ defines a Maurer--Cartan element in $\Hom_{\RN}\big(\ol{\Sym}_\RN(\frg[1]), \frh[1]\big)$ (see~\eqref{eq:MC eqn for cdgc mps into CE cdgc}), and
		
		\item the composition \smash{$\ol{\Sym}_{A}(\frg[1]) \longrightarrow \frh\xrightarrow{ \ \rho_{\frh} \ } T_A$} is given by $\rho_{\frg}$ on the linear copy $\frg[1]$ and vanishes on all higher tensor powers of $\frg[1]$.
	\qen
	\end{myenumerate}
\end{remark}

\begin{lemma}[{\cite[Lemmas 5.8 and 5.10]{Nuiten:HoAlg_for_Lie_Algds}}]\label{lem:lie algebroid resolution}
	For any $L_\infty$-algebroid $\frg$, there exists another $L_\infty$-algebroid $Q\frg$ with the following universal property: there is a natural bijection between strict maps of $L_\infty$-algebroids $Q\frg\longrightarrow \frh$ and $\infty$-morphisms $\frg\longleadsto \frh$; that is
	\begin{equation}
		L_\infty\Agd_A^\dg(Q\frg, \frh) \cong L_\infty\Agd_A(\frg, \frh)\ .
	\end{equation}
	Furthermore, $Q\frg\longrightarrow \frg$ is a cofibrant resolution if $\frg$ is cofibrant as a dg $A$-module.
\end{lemma}

We will give a more explicit description of $Q\frg$ in Section~\ref{sec:l-conn}.

\begin{remark}
	Suppose that $A$ is concentrated in degree zero, and recall that any bounded-above complex of projective $A$-modules is cofibrant.
	In the particular case where $A=C^\infty(M)$, this means that any bounded-above, and in particular every finite, complex of finite-rank vector bundles on $M$ determines such a cofibrant module.
\qen
\end{remark}

\begin{corollary}
	Let $\frg$ and $\frh$ be  $L_\infty$-algebroids and suppose that $\frg$ is cofibrant as a dg $A$-module.
	Then the derived mapping space can be computed by the simplicial set of $\infty$-morphisms
	\begin{equation}
		L_\infty\Agd_A^\infty\big(\frg, \frh) \simeq L_\infty\Agd^\dg_A \big( Q\frg,  \widehat{\frh}\, \big)
		\simeq L_\infty\Agd_A \big( \frg,  \widehat{\frh}\, \big)\ ,
	\end{equation}
	where $\widehat{\frh}$ is any Reedy fibrant simplicial resolution of $\frh$.
\end{corollary}

\begin{remark}\label{Rem: Reed fibrant resolutions <1}
Depending on the situation, there are various ways to produce a Reedy fibrant simplicial resolution $\widehat{\frh}$:
\begin{myenumerate}
\item If $\frh$ is fibrant, i.e. $\frh\longrightarrow T_A$ is surjective, then one can construct an explicit fibrant simplicial resolution using polynomial differential forms~\cite[Lemma 5.24]{Nuiten:HoAlg_for_Lie_Algds}.

\item If $A=C^\infty(M)$, and $\frg$ and $\frh$ are both concentrated in non-positive degrees, then it suffices to construct a Reedy fibrant resolution of $\frh$ in the model structure on non-positively graded $L_\infty$-algebroids from Remark~\ref{rem:nonpositively graded case}.

\item If $\frh\longrightarrow T_A$ is a \emph{dg abelian extension} of $T_A$ in the sense of Definition~\ref{def:ab ext}, then one can use the simplicial resolution from Proposition~\ref{st:hat(frg) is dg Lie agd}. \qen
\end{myenumerate}
\end{remark}

The usual de Rham algebra has a natural analogue for $L_\infty$-algebroids.

\begin{definition}
	If $\frg$ is an $L_\infty$-algebroid  over $A$, its \emph{Chevalley--Eilenberg algebra} $\ChEil^*(\frg)$ (or \emph{de Rham algebra}) is the cdga given as follows:
	\begin{myitemize}
		\item the underlying graded algebra is given by the $A$-linear dual of the $A$-linear cofree cdgc on $\frg[1]$:
		\begin{equation}
			\ChEil^*(\frg)=\Hom_{A}\big(\Sym_{A} (\frg[1]), A\big)\ .
		\end{equation}
		\item the differential is given by
		\begin{align*}
			 (\dd_\ChEil\alpha)(x_1, \dots, x_n) & = \dd_A\big(\alpha(x_1, \dots, x_n)\big) +  \sum_{k=1}^n \hspace{2pt} \pm\, \pounds_{\rho(x_k)}\alpha(x_{1}, \dots, \widehat{x_k}, \dots , x_n)\\
			&\quad\, - \sum_{k\geqslant 1} \ \sum_{\sigma\in \mathrm{Sh}(k, n-k)} \pm \hspace{2pt} \alpha \Big([x_{\sigma(1)}, \dots , x_{\sigma(k)}], x_{\sigma(k+1)}, \dots , x_{\sigma(n)}\Big)\ .
		\end{align*}
		Here the permutation $\sigma$ ranges over all $(k,n{-}k)$-shuffles, where $\pm$ denotes the Koszul sign associated to the permutation of the elements $x_i\in \frg[1]$, and $[x_i] = \dd_{\frg[1]}(x_i)$ denotes the differential of $\frg[1]$. The hat indicates omission of the corresponding entry.
	\end{myitemize}
\end{definition}

There is a natural bijection between $n$-cocycles $\alpha\in Z^n\ChEil^*(\frg)$ and $\infty$-morphisms
$$
	\begin{tikzcd}
		(\rho, \alpha)\colon \frg\arrow[r, rightsquigarrow] & T_A\oplus A[n-1]
	\end{tikzcd}
$$
where $T_A\oplus A[n-1]$ is the dg Lie algebroid with zero differential and bracket $$\big[(X, a)\,,\, (Y, b)\big] = \big([X, Y]\,,\, \pounds_X(b) - \pounds_Y(a)\big) \ . $$
Indeed, the $T_A$-component of this $\infty$-morphism is fixed to be the anchor map $\rho$, and the second component is given by an $A$-linear map
$$
	\begin{tikzcd}
		\alpha\colon \Sym_{A}(\frg[1]) \arrow[r] & A[n]\ .
	\end{tikzcd}
$$
Unravelling the definitions as in~\cite[Example 6.18]{Nuiten:HoAlg_for_Lie_Algds}, one sees that $\alpha$ is an $n$-cocycle in the Chevalley--Eilenberg complex if and only if $(\rho, \alpha)$ is an $\infty$-morphism. One can refine this to an equivalence between mapping spaces.

\begin{lemma}
	Let $\frg$ be an $L_\infty$-algebroid which is cofibrant as a dg $A$-module. Then there is a natural equivalence of mapping spaces
	$$
		L_\infty\Agd_A^\infty\big(\frg, T_A\oplus A[n-1]\big)\simeq \Mod^\infty_\RN\big(\RN[-n], \ChEil^*(\frg)\big)\ .
	$$
\end{lemma}

The right mapping space can be presented explicitly as the image of a complex under the Dold--Kan correspondence as
$$
	\Mod^\infty_\RN \big(\RN[-n], \ChEil^*(\frg)\big)
	\simeq \DK\big( Z^n\ChEil^*(\frg) \longleftarrow \ChEil^{n-1}(\frg) \longleftarrow \cdots \big)\ .
$$
Indeed, this simplicial set is isomorphic to the simplicial set of maps $\Mod_{\RN}^{\dg}\big(C_*(\Delta^\bullet)[-n], \ChEil^*(\frg)\big)$ out of a cofibrant cosimplicial resolution of the complex $\RN[-n]$, where $C_*(\Delta^p)$ denotes the (normalised Moore) complex of singular chains on the $p$-simplex $\Delta^p \in \sSet$, for $p \in \NN_0$.

\begin{proof}
	It follows from~\cite[Corollary 6.16]{Nuiten:HoAlg_for_Lie_Algds} that $\Mod^\infty_{\RN}\big(\RN[-n], \ChEil^*(\frg)\big)$ is naturally equivalent to the space of sections of the projection 
	$$
		\pi\colon \frg\oplus A[n-1] = \frg\times_{T_A} \big(T_A\oplus A[n-1]\big) \longrightarrow \frg\ .
	$$
	Since $\pi$ is the homotopy pullback of $T_A\oplus A[n-1]\longrightarrow T_A$ along $\rho\colon \frg\longrightarrow T_A$ and $T_A$ is the terminal $L_\infty$-algebroid, this space of sections is equivalent to the space of maps from $\frg$ to $T_A\oplus A[n-1]$.
\end{proof}

Finally, let us consider the special case where $A=C^\infty(M)$ is the ring of smooth functions on a smooth manifold $M$. In this case, the definition of $L_\infty$-algebroid is a bit too algebraic: a general $L_\infty$-algebroid over $C^\infty(M)$ cannot be studied locally on $M$ using descent methods. For this reason, it will often be convenient to add the condition below.
In the case where $A = C^\infty(M)$ we also write
\begin{equation}
	L_\infty\Agd_M \coloneqq L_\infty\Agd_{C^\infty(M)}\ ,
	\quad
	L_\infty\Agd_M^\dg \coloneqq L_\infty\Agd_{C^\infty(M)}^\dg
	\qandq
	L_\infty\Agd_M^\infty \coloneqq L_\infty\Agd_{C^\infty(M)}^\infty\ .
\end{equation}

\begin{definition}
	\label{def:complete}
	Let $M$ be a smooth manifold and $E$  a dg $C^\infty(M)$-module. Let $\mathcal{E}$ be the induced presheaf of complexes on $M$ given by
	$$
		\mathcal{E}(U) = C^\infty(U) \otimes_{C^\infty(M)} E\ ,
	$$
	for any open subset $U \subset M$, and let $E^\cpl$ be the associated $\infty$-sheaf of complexes. The dg $C^\infty(M)$-module $\Gamma(M, E^\cpl)$ of global sections of $E^\cpl$ is the \emph{completion} of $E$, and  $E$ is \emph{complete} if the map $E\longrightarrow \Gamma(M, E^\cpl)$ is a quasi-isomorphism.
	We denote the full subcategory of \smash{$\Mod_{C^\infty(M)}^\dg$} on the complete modules by \smash{$\Mod_{C^\infty(M)}^{\dg, \cpl}$}.
	
	An $L_\infty$-algebroid is \textit{complete} if its underlying complex of $C^\infty(M)$-modules is complete, and we denote the full subcategory of complete $L_\infty$-algebroids by \smash{$L_\infty\Agd_M^{\dg,\cpl}$}.
\end{definition}

One can think of complete $L_\infty$-algebroids as (the global sections of) \emph{sheaves} of $L_\infty$-algebroids over $M$.

\begin{lemma}[{\cite[Lemma 5.1.33]{Nuiten:Thesis}}]
	\label{lem:complete modules vs sheaves}
	For any smooth manifold $M$, the functors sending $E\longmapsto E^\cpl$ and taking global sections give rise to an adjunction
	$$
		\begin{tikzcd}[column sep = 1cm]
			(-)^\cpl : \Mod_{C^\infty(M)}^\infty \arrow[r, yshift = 1ex, "\perp"' yshift = 0ex] \arrow[r, yshift = -1ex, hookleftarrow]
			& \Mod_{C^\infty_M}^\infty : \Gamma(M, -)
		\end{tikzcd}
	$$
	between the $\infty$-categories of dg $C^\infty(M)$-modules and complexes of sheaves on $M$ that are dg modules over the sheaf $C^\infty_M$ of $C^\infty$-functions on $M$.
	The right adjoint $\Gamma(M, -)$ is fully faithful with essential image consisting precisely of the complete $C^\infty(M)$-modules.
\end{lemma}

In particular, this shows that completion is an idempotent operation.

\begin{remark}\label{rem:complete}
	Since the sheaf of $C^\infty$-functions on $M$ has partitions of unity, its module sheaves have vanishing higher sheaf cohomology. Consequently, for any $-\infty\leqslant i\leqslant j\leqslant \infty$, both functors in the adjunction of Lemma~\ref{lem:complete modules vs sheaves} preserve complexes with cohomology in the range $[i, j]$. This implies the following pleasant properties of completion:
	\begin{myenumerate}
	\item If $E^\cpl$ is the completion of $E$, then $\rmH^i(E^\cpl)$ is the completion of $\rmH^i(E)$; in particular, $\rmH^i(E^\cpl)$ is complete. Consequently, $E$ is complete if and only if each cohomology module $\rmH^i(E)$ is complete.
	
	\item If $E$ is a complex of $C^\infty(M)$-modules given in each degree by a complete module in the sense of Joyce~\cite[Definition 5.25]{Joyce:Alg_Geom_Coo}, then $E$ is complete in the sense of Definition~\ref{def:complete}. Most importantly, this implies that any complex of finite rank projective modules (i.e.~vector bundles) is complete.
	
	\item Let $f\colon M\longrightarrow N$ be a map between smooth manifolds. Then there is a natural adjoint pair between their categories of complete modules
	$$
		\begin{tikzcd}
			f^*: \Mod_{C^\infty(N)}^{\infty, \cpl} \arrow[r, yshift = 1ex, "\perp"']
			& \Mod_{C^\infty(M)}^{\infty, \cpl}: f_* \arrow[l, yshift=-1ex]
		\end{tikzcd} \ .
	$$
	Here $f_*$ is given by restriction of scalars along $C^\infty(N)\longrightarrow C^\infty(M)$ and $f^*(E)$ is the completion of the module $C^\infty(M)\otimes^h_{C^\infty(N)} E$. Under the identification from Lemma~\ref{lem:complete modules vs sheaves}, this adjunction is identified with the direct image and inverse image of module sheaves; the functor $f_*$ does not need to be derived because it is exact (by partitions of unity). \qen
	\end{myenumerate}
\end{remark}


\subsection{Homotopy theory of dg $C^\infty$-rings}
\label{sec:HoTheoryrings}


Perhaps the simplest non-trivial example of a Lie algebroid on a smooth manifold $M$ is the following: let $\pi\colon M\longrightarrow N$ be a submersion and consider the subbundle
$$
T(M/N)=\ker\big(\dd \pi\colon TM\longrightarrow \pi^*TN\big)
$$
of vector fields on $M$ that are tangent to the fibres of $\pi$. Then $T(M/N)$ has the natural structure of a Lie algebroid, since vector fields tangent to the fibres of $\pi$ are closed under the commutator bracket. 
The purpose of this and the following sections is to describe a general method for associating an $L_\infty$-algebroid (well-defined up to quasi-isomorphism) to \emph{any} map from a manifold $M$ to a higher differentiable stack which admits sufficient deformation theory.

To do this in a functorial way, we will need a sufficiently formal way to define `infinitesimal objects' such as the tangent bundle of a smooth manifold.
One way to achieve this is via synthetic differential geometry, by enlarging the site of smooth manifolds to also include `infinitesimally thickened manifolds'.
In fact, to ensure sufficient functoriality of our constructions, we will need to go one step further and also include `\emph{derived} infinitesimally thickened manifolds' to our site.
Since each smooth manifold $M$ is uniquely determined by its $C^\infty$-ring of smooth functions, we will use the algebraic formalism of dg $C^\infty$-rings from~\cite{CR:HoAlg_Superalg} to describe such objects.

\begin{definition}
	A \textit{dg $C^\infty$-ring} $A$ is an $\RN$-linear commutative dg algebra concentrated in non-positive degrees, equipped with a compatible $C^\infty$-ring structure on its degree~zero part $A^0$.
\end{definition}

Recall that a $C^\infty$-ring is a set $A$ equipped with an operation $f_*\colon A^{\times n}\longrightarrow A$ for every $C^\infty$-function $f\colon \RN^n\longrightarrow \RN$ such that composition is respected (see e.g.~\cite{MR:Models_Smooth_Inf_Analysis} or~\cite{Joyce:Alg_Geom_Coo} for detailed accounts). In particular, the operations associated to polynomial maps $\RN^n\longrightarrow \RN$ endow any dg $C^\infty$-ring $A$ with the structure of a cdga over $\RN$.

\begin{example}
	Let us mention the following examples of dg $C^\infty$-rings:
	\begin{myenumerate}
		\item For every smooth manifold $M$, its smooth functions form a dg $C^\infty$-ring concentrated in degree zero, which we again denote by $C^\infty(M)$.
		
		\item If $A$ is a dg $C^\infty$-ring and $I \subseteq A$ is any dg ideal in the usual algebraic sense, then $A/I$ inherits the structure of a dg $C^\infty$-ring by~\cite[Proposition 1.2]{MR:Models_Smooth_Inf_Analysis}.
		
		\item Consider a dg Artin (or Weil) algebra over $\RN$, i.e.\ a finite-dimensional cdga of the form $A = \RN[x_1, \dots, x_n]/I$ where $I \subseteq (x_1, \dots, x_n)$ is an ideal such that for all $i$, $x_i^{k_i} \in I$ for some $k_i$. Here the $x_i$ can have any non-positive degree, and there can be a non-zero differential.
		Then $A$ has the canonical structure of a dg $C^\infty$-ring (the $C^\infty$-ring structure on $A^0$ arises from~\cite[Theorem~3.17]{MR:Models_Smooth_Inf_Analysis}).
		
		\item The tensor product $C^\infty(M) \otimes \RN[x_1, \dots, x_n]/I$ of a dg Artin algebra and the function algebra of a smooth manifold is again a dg $C^\infty$-ring.
		
		\item Suppose that $\frg$ is an $L_\infty$-algebroid over a smooth manifold $M$, and let
		$$
			\tau^{\leqslant 0}\ChEil^*(\frg) = \Big[Z^0\big(\ChEil^*(\frg)\big)\longleftarrow \ChEil^{-1}(\frg)\longleftarrow \cdots\Big]
		$$
		be the truncation of its Chevalley--Eilenberg algebra to non-positive degrees. Then $\tau^{\leqslant 0}\ChEil^*(\frg)$ has the structure of a dg $C^\infty$-ring~\citep[Lemma 4.1.20]{Nuiten:Thesis}.
		
		\item If $A$ is a dg $C^\infty$-ring and $M$ is a non-positively graded dg $A$-module, then the square-zero extension $A\oplus M$ has a natural $C^\infty$-ring structure by~\cite[Proposition 2.42]{CR:Theories_Superalg}. Explicitly, the $C^\infty$-structure on $A^0\oplus M^0$ is given by
		$$
			f_*\big(a_1+m_1, \dots, a_n+m_n\big)
			= f_*\big(a_1, \dots, a_n\big)+\sum_{i=1}^n\, \Big(\frac{\partial f}{\partial x_i}\Big)_*(a_1, \dots, a_n)\cdot m_i\ .
		$$
		
		\item For any dg $C^\infty$-ring $A$, its zeroth cohomology $\rmH^0(A)$ has the structure of an ordinary $C^\infty$-ring. We define the \emph{reduction} $A_\red$ of $A$ to be the quotient of $\rmH^0(A)$ by its ideal of nilpotent elements (we will not need any more subtle notion of reduction):
		$$
			A_{\red}\coloneqq \rmH^0(A)\,\big/\,\mathrm{nilradical}\ .
		$$
		There is a natural map of dg $C^\infty$-rings $A\longrightarrow A_\red$.
\qen
	\end{myenumerate}
\end{example}

The category of dg $C^\infty$-rings supports a model structure~\cite[Theorem 6.10]{CR:HoAlg_Superalg} whose weak equivalences are the quasi-isomorphisms and fibrations are maps that are surjective in negative degree. The cofibrant objects are retracts of dg $C^\infty$-rings of the form
$$
	\Big(C^\infty(\RN^S)[x_\alpha], \dd \Big)
$$
with arbitrary differential $\dd$.
Here the $x_\alpha$ are generators in negative degree, and for any set $S$ we denote by $C^\infty(\RN^S)$ the ring of functions $f \colon \prod_{s \in S}\, \RN \longrightarrow \RN$ that factor as $$\prod_{s \in S}\, \RN \longtwoheadrightarrow \RN^{\times k}\longrightarrow \RN \ . $$
The first map projects onto finitely many factors and the second map is $C^\infty$.
Any dg $C^\infty$-ring is a quotient of such a quasi-free dg $C^\infty$-ring.

\begin{definition}
We write $\CRing^\dg$ for the category of dg $C^\infty$-rings endowed with the above model structure, and $\CRing^\infty = \CRing^\dg[W^{-1}]$ for its associated $\infty$-category.
\end{definition}

\begin{remark}\label{rem:ho pushout of coo-rings}
	For any diagram of dg $C^\infty$-rings $$\begin{tikzcd} & \arrow[dl] \arrow[dr] A & \\ A' & & B \end{tikzcd}$$ in which $A \longrightarrow A'$ is surjective on $\rmH^0$, the map $A'\otimes^h_A B\longrightarrow A'\amalg^h_A B$ from the derived tensor product of the underlying cdgas to the homotopy pushout of dg $C^\infty$-rings is an equivalence. Indeed, in this case one can model $A \longrightarrow A'$ by a cofibration in which one only adds generators in negative degree, so that $A$ and $A'$ are isomorphic $C^\infty$-rings in degree $0$; in this case, the pushout of dg-$C^\infty$-rings is simply given by the tensor product.
	\qen
\end{remark}

For our purposes it will suffice to consider dg $C^\infty$-rings that arise as nilpotent extensions of the function ring of a smooth manifold.
To make this more precise, let us recall how to describe square-zero extensions of dg $C^\infty$-rings (see for instance~\cite[Definition 2.2.20]{Nuiten:Thesis} or~\cite[Definition 1.2.1.6]{TV:Homotopical_Alg_Geom_II}).

\begin{definition}\label{def:square-zero}
	Let $A$ and $B$ be dg $C^\infty$-rings, and let $I$ be a dg $A$-module in non-positive degrees. A homotopy pullback diagram in the $\infty$-category $\CRing^\infty$ of the form
	\begin{equation}\label{eq:square-zero}
			\begin{tikzcd}[row sep = 1cm,column sep=1cm]
			B_\eta\arrow[d]\arrow[r] & A\arrow[d, "{(\id, 0)}"]\\
			B\arrow[r, "\eta"] & A\oplus I[1]
		\end{tikzcd}
	\end{equation}
	exhibits $B_\eta$ as a \emph{square-zero extension} of $B$ by the dg $A$-module $I$.
\end{definition}

\begin{remark}
	\label{rem:square_zero}
	The homotopy pullback in Definition~\ref{def:square-zero} can be identified more explicitly as follows:
	suppose that the map $\eta$ in $\CRing^\infty$ arises from a (strict) map of dg $C^\infty$-rings (one can always assume this by replacing $B$ with a cofibrant resolution) $$\eta = \eta_A \oplus \eta_I \colon B\longrightarrow A\oplus I[1] \ . $$
	Such a map of dg $C^\infty$-rings corresponds to a map of dg $C^\infty$-rings $\eta_A:B\longrightarrow A$, so that we can view $I$ as a dg $B$-module, together with a \emph{derivation} $\eta_I \colon B\longrightarrow I[1]$, i.e.\ a map such that
	$$
		\eta_I(b_1\, b_2) = b_1\cdot \eta_I(b_2) + (-1)^{|b_1|\, |b_2|}\, b_2 \cdot \eta_I(b_1)\ ,
		\quad
		\eta_I \big( f_*(b_1, \dots, b_n) \big) = \sum_{i=1}^n\, \Big( \frac{\partial f}{\partial x_i} \Big)_* (b_1, \dots, b_n) \cdot \eta_I(b_i)\ ,
	$$
	where the last equation only applies to elements $b_i \in B$ in degree zero.
	One can then compute the homotopy pullback by resolving the right vertical map by
	$$
		\begin{tikzcd}
			A\oplus \text{Cone}(\id_I)\arrow[r] & A\oplus I[1]\ ,
		\end{tikzcd}
	$$
	whose kernel is $I$; here $\Cone$ denotes the mapping cone of a map of complexes.
	The homotopy pullback is then given by $B \oplus I$ with differential $$\dd (b+m) = \dd_B(b) + \eta_I(b) + \dd_I(m)$$ for all $b \in B$ and $m \in I$.
	
	Conversely, suppose that $I$ is a square-zero ideal in a dg $C^\infty$-ring $B'$, and let $B\longrightarrow B'/I$ be a quasi-free resolution. Then the base change $B\times_{B'/I} B' \longtwoheadrightarrow B$
is a (strict) square-zero extension of dg $C^\infty$-rings which is quasi-isomorphic to $B'\longrightarrow B'/I$.
	That $B$ is quasi-free allows us to lift generators, and so one can identify \smash{$B\times_{B'/I} B'\cong B\oplus I$} with differential determined by a certain derivation $\eta_I\colon B\longrightarrow I[1]$. In other words, up to quasi-isomorphism every square-zero extension arises from a homotopy pullback as in Definition~\ref{def:square-zero}.
	\qen
\end{remark}

\begin{lemma}
	\label{lem:coprod of square-zero}
	Fix a homotopy pullback square \eqref{eq:square-zero} and a dg $C^\infty$-ring $C$, and consider the induced square of homotopy coproducts
	$$
		\begin{tikzcd}[row sep = 1cm,column sep=1cm]
			C \amalg^h B_\eta \arrow[r] \arrow[d]
			& C \amalg^h A\arrow[d, "{\id\amalg (\id, 0)}"]
			\\
			C \amalg^h B \arrow[r, "\id \amalg \eta"]
			& C \amalg^h (A\oplus I[1])
		\end{tikzcd}
	$$
	If the map $\rmH^0 \eta \colon \rmH^0(B) \longrightarrow \rmH^0(A \oplus I[1]) \cong \rmH^0(A)$ is surjective, then this square is a homotopy pullback that exhibits $C\amalg^h B_\eta$ as the square-zero extension of $C \amalg^h B$ by the dg $C^\infty$-ring  $(C \amalg^h A) \otimes^h_A I$.
\end{lemma}

\begin{proof}
	We first show that the square is a homotopy pullback, for which it suffices to work at the level of cdgas.
	Since $\eta$ is surjective on $\rmH^0$, all maps in the square~\eqref{eq:square-zero} are surjective on $\rmH^0$. By Remark~\ref{rem:ho pushout of coo-rings}, the square in the statement of the lemma is then obtained from~\eqref{eq:square-zero} by applying the functor
	\begin{equation}
		(C \amalg^h B_\eta) \otimes^h_{B_\eta}(-)
		\simeq (C\amalg^h B_\eta) \amalg^h_{B_\eta}(-)\ .
	\end{equation}
	It thus suffices to check that the derived tensor product preserves homotopy pullbacks of maps that induce a surjection on $\rmH^0$. This can be verified at the level of non-positively graded $B_\eta$-modules, where a square is a pullback if it is a pushout, and the converse holds whenever two parallel maps in the square are surjective on $\rmH^0$ (see \cite[Corollary C.1.2.6]{Lurie:SAG}); it then follows since the derived tensor product preserves homotopy pushouts.
	
	A similar argument shows that at the level of dg algebras, the map
	\begin{equation}
		\id \amalg (\id, 0) \colon C \amalg^h A \longrightarrow C\amalg^h (A\oplus I[1])
	\end{equation}
	is equivalent to the canonical map
	\begin{equation}
		(\id, 0) \colon C \amalg^h A \longrightarrow (C \amalg^h A) \oplus (C\amalg^h A) \otimes^h_A I[1]
	\end{equation}
	to the trivial square-zero extension.
	It follows from~\cite[Theorem 9.3]{Carchedi:D_modules} that this is also an equivalence as maps of dg $C^\infty$-rings. This shows that the square of coproducts indeed classifies a square-zero extension by $(C\amalg^h A)\otimes^h_A I$.
\end{proof}

\begin{lemma}
	\label{lem:nil ext}
	Let $B\longrightarrow A$ be a morphism of dg $C^\infty$-rings with the following two properties:
	\begin{myenumerate}
	\item[(a)] the induced morphism $\rmH^0(B) \longrightarrow \rmH^0(A)$ is surjective with nilpotent kernel $J$;
	
	\item[(b)] there exists an $N\in \NN$ such that $\rmH^{-l}(B) \longrightarrow \rmH^{-l}(A)$ is an isomorphism for all $l>N$.
	\end{myenumerate}
	Then:
	\begin{myenumerate}
		\item The map $B\longrightarrow A$ decomposes (up to equivalence) as a finite sequence of square-zero extensions by shifted $A_{\rm red}$-modules.
		
		\item If $C\longrightarrow C'$ is a map of dg $C^\infty$-rings under $B$ such that the map of pushouts $A\amalg_B C\longrightarrow A\amalg_B C '$ is an equivalence, then $C\longrightarrow C'$ is an equivalence.
	\end{myenumerate}
\end{lemma}

\begin{proof}
	For assertion~(1), we may assume that $B\longrightarrow A$ is surjective, with kernel $I$. For each $n\geqslant 0$, let us consider the following two sub-complexes of $I$:
	\begin{equation}
		\begin{tikzcd}[row sep = 0.1cm, column sep = 0.75cm]
			I(n) = \big[ 0
			& \cdots \ar[l]
			& 0 \ar[l]
			& Z^{-n} I \ar[l]
			& I^{-n-1} \ar[l]
			& \cdots \big]\ , \ar[l]
			\\[4pt]
			I'(n) = \big[ 0
			& \cdots \ar[l]
			& \dd(I^{n-1}) \ar[l]
			& I^{-n} \ar[l]
			& I^{-n-1} \ar[l]
			& \cdots \big]\ . \ar[l]
		\end{tikzcd}
	\end{equation}
	Both define dg ideals of $B$, so that we obtain a tower of dg $C^\infty$-rings
	$$
		\begin{tikzcd}
			B \arrow[r]
			& \cdots \arrow[r]
			& B/I'(n) \arrow[r]
			& B/I(n)\arrow[r, "\sim"]
			& B/I'(n-1)\arrow[r]
			& \cdots \arrow[r]
			& B/I(0) = A\ .
		\end{tikzcd}
	$$
	Each map $B/I(n)\longrightarrow B/I'(n-1)$ is a quasi-isomorphism because it is surjective with kernel given by the acyclic complex $I^{-n}/Z^{-n}I\longrightarrow \dd(I^{-n})$. Furthermore, the map $B\longrightarrow B/I(n)$ is a quasi-isomorphism for all sufficiently large $n$, as $\rmH^{-l}(I) = 0$ for all $l \in \NN$ sufficiently large (since $\rmH^{-l}(f)$ is an isomorphism for all $l\in\NN$ sufficiently large). 
	
	It therefore suffices to verify that each map $B/I'(n)\longrightarrow B/I(n)$ decomposes into a finite sequence of square-zero extensions by $\rmH^0(A)$-modules. Let us write $B'\longrightarrow B/I(n)$ for this map and $K=\rmH^{-n}(I)[n]$ for its kernel, which is an $\rmH^0(B)$-module. Since the kernel $J$ of the map $\rmH^0(B)\longrightarrow \rmH^0(A)$ satisfies $J^p=0$ for some $p\in\NN$, we then obtain a finite sequence of dg $C^\infty$-rings
	$$
		\begin{tikzcd}
			B/I'(n) = B'/J^p\cdot K\arrow[r]
			& B'/J^{p-1} \cdot K \arrow[r]
			& \cdots\arrow[r] &
			B'/J\cdot K \arrow[r]
			& B'/K=B/I(n) \ .
		\end{tikzcd}
	$$
	Each map in this sequence is a (strict) square-zero extension by a shift of an $\rmH^0(A)$-module of the form $J^k\cdot K/J^{k+1}\cdot K$. By Remark~\ref{rem:square_zero}, these are also square-zero extensions in the sense of Definition~\ref{def:square-zero}.
	
	For assertion~(2), note that $A\amalg^h_B C\simeq A\otimes^h_B C'$ since $B \longrightarrow A$ is surjective on $\rmH^0$.
	It therefore suffices to show that $A\otimes^h_B (-)$ detects equivalences between non-positively graded dg $B$-modules.
	This follows because $B\longrightarrow A$ is the limit of a finite tower of square-zero extensions by $A$-modules $I_k$ in increasingly negative degrees:
	tensoring up, this induces a natural tower for each non-positively graded $B$-module $E$ such that $E=\lim_{k\in\NN_0} E_k$ and equipped with natural equivalences 
	$$\fib(E_k \longrightarrow E_{k-1}) \simeq I_k \otimes_{A} (A\otimes^h_B E).
	$$
	So if $f \colon E \longrightarrow E'$ is a morphism of non-positively graded dg $B$-modules such that $(A\otimes^h_B f)$ is an equivalence, then $f$ induces an equivalence \smash{$\fib(E_k \longrightarrow E_{k-1}) \eq \fib(E'_k \longrightarrow E'_{k-1})$} on the homotopy fibre, for each $k \in \NN_0$.
	Thus $f$ is itself an equivalence.
\end{proof}

\begin{remark}
		\label{rem:extensions linear}
	The same argument shows that, in the situation of Lemma~\ref{lem:nil ext}, any non-positively graded dg $B$-module $E$ can be obtained as an iterated extension of dg $A$-modules.
	\qen
\end{remark}

\begin{corollary}
	\label{cor:inf manifolds retract}
	Let $A$ be a dg $C^\infty$-ring such that $A_{\red}\cong C^\infty(M)$ for a smooth manifold $M$. The natural map $\epsilon\colon A\longrightarrow A_{\red}=C^\infty(M)$ admits a (non-canonical) section up to homotopy.
\end{corollary}

\begin{proof}
	The map $\epsilon$ decomposes into a tower of maps
	$$\begin{tikzcd}
	\epsilon\colon A=\underset{n\in\NN_0}{\lim}\, \tau^{\geqslant -n}(A)\arrow[r] & \cdots \arrow[r] & \tau^{\geqslant -n}(A)\arrow[r] & \cdots\arrow[r] & \tau^{\geqslant 0}(A)\arrow[r] & A_{\red} \ .
	\end{tikzcd}$$
	It will suffice to construct a compatible family of sections $s_n\colon A_{\red}\longrightarrow \tau^{\geqslant -n}(A)$ for each $n\geqslant 0$;
	the desired section is then given by taking the limit.
	We will produce this family of sections by induction, using Lemma \ref{lem:nil ext}, which implies that the maps $\tau^{\geqslant 0}(A)\longrightarrow A_{\red}$ and $\tau^{\geqslant -n}(A)\longrightarrow \tau^{\leqslant 1-n}(A)$ all decompose into finite sequences of square-zero extensions by $A_{\red}$-modules.
	
	 For the inductive step, it therefore suffices to check that we can lift these sections against square-zero extensions. For this, let us consider a diagram
	 $$\begin{tikzcd} [column sep = 1cm , row sep = 1cm]
& B_\eta\arrow[r]\arrow[d] & A_{\red}\arrow[d, "{(\id, 0)}"]\\
A_{\red}\arrow[r, "s"]\arrow[ru, dashed, "\tilde{s}"] & B\arrow[r, "\eta"] & A_{\red}\oplus I[1] 	 
	 \end{tikzcd}$$
where the square is a pullback exhibiting $B_\eta$ as a square-zero extension of $B$. To find the dashed lift, we have to show that the map $\eta\circ s\colon A_{\red}\longrightarrow A_{\red}\oplus I[1]$ factors over $(\id, 0)\colon A_{\red}\longrightarrow A_{\red}\oplus I[1]$.

	Now recall from Remark~\ref{rem:square_zero} that the map $\eta\circ s$ is classified by a derivation $A_{\red}=C^\infty(M)\longrightarrow I[1]$. 
	In turn, this is determined uniquely by a map $\phi\colon \Omega^1(M)\longrightarrow I[1]$ in the $\infty$-category of $C^\infty(M)$-modules (see Example~\ref{ex:corep tangent} below, as well as~\cite[Example~2.2.18]{Nuiten:Thesis}). 
	Under these identifications, the space of factorisations of $\eta\circ s$ over $(\id, 0)$ is equivalent to the space of null-homotopies of $\phi$. 
	Since $\Omega^1(M)$ is a projective $C^\infty(M)$-module (hence cofibrant) and $I[1]$ is concentrated in negative degrees, the map $\phi$ is necessarily null-homotopic. We can therefore factor $\eta\circ s$ over $(\id, 0)$ and find a lift $\tilde{s}$.	
\end{proof}


\subsection{Background from derived $C^\infty$-geometry}
\label{sec:derived C^oo-geometry}


Using the elements from the homotopy theory of dg $C^\infty$-rings from Section~\ref{sec:HoTheoryrings}, we will now set up a version of derived $C^\infty$-algebraic geometry for infinitesimal thickenings of smooth manifolds.
Ideally, one would like to say that for every dg $C^\infty$-ring $A$ such that $A_\red\cong C^\infty(M)$ for a smooth manifold $M$, there is an infinitesimal thickening $\Spec(A)$ of $M$. However, such dg $C^\infty$-rings need not behave well with respect to descent, and we need to impose a completeness condition similar to Definition~\ref{def:complete}.

\begin{definition}
	\label{def:complete2}
	Let $A$ be a dg $C^\infty$-ring such that $\rmH^0(A)$ is a nilpotent extension of $C^\infty(M)$ for some smooth manifold $M$. Consider the ($\infty$-)presheaf of dg $C^\infty$-rings
	\begin{equation}
			\begin{tikzcd}
			\calO_{A, \mathrm{pre}} \colon \mathrm{Open}(M)^{\mathrm{op}}\arrow[r]
			& \CRing^\infty\ ,
		\end{tikzcd}
	\end{equation}
	sending each open subset $U \subset M$ to the dg $C^\infty$-ring obtained from $A$ by inverting all elements in $\rmH^0(A)$ whose image under $$\rmH^0(A)\longrightarrow C^\infty(M)\longrightarrow C^\infty(U)$$ is invertible~\cite[Lemma 5.1.13]{Nuiten:Thesis}. Let $\calO_A$ be the associated sheaf. Then $A$ is \emph{complete} if the natural map 
	\begin{equation}
		A = \Gamma(M, \calO_{A, \mathrm{pre}}) \longrightarrow \Gamma(M, \calO_A)
	\end{equation}
	of global sections is an equivalence of dg $C^\infty$-rings.
	Similarly, a dg $A$-module $E$ is \emph{complete} if the natural map $E\longrightarrow \Gamma(M, \mathcal{E})$ is an equivalence of dg $A$-modules, where $\mathcal{E}$ is the associated sheaf of~$\calO_A\otimes_A E$.
\end{definition}

This coincides with the notion of completeness for dg $C^\infty$-rings and their modules from~\cite{Joyce:Alg_Geom_Coo, Nuiten:Thesis, Steffens:Der_Coo}. As shown there, the functors $A\longmapsto \Gamma(M, \calO_A)$ and $E\longmapsto \Gamma(M, \mathcal{E})$ define left adjoints to the inclusions of the full subcategories of complete objects; we will refer to them as \emph{completion}.

\begin{remark}
	The presheaf $\calO_{A, \mathrm{pre}}$ can be made more explicit via a choice of section $C^\infty(M)\longrightarrow A$ as in Corollary~\ref{cor:inf manifolds retract}.
	Indeed, such a section gives rise to natural maps
	\begin{equation}
		\mu \colon C^\infty(U) \amalg_{C^\infty(M)} A \longrightarrow \calO_{A, \mathrm{pre}}(U)
		\qqandqq
		\mu' \colon C^\infty(U) \otimes_{C^\infty(M)} A \longrightarrow \calO_{A, \mathrm{pre}}(U)
	\end{equation}
	of $C^\infty$-rings and cdgas, respectively (recall that when we take the algebraic tensor product $\otimes$ of $C^\infty$-rings, rather than their pushouts, we only remember their underlying cdga structures). The maps $\mu$ and $\mu'$ are both equivalences, since they become equivalences after taking the pushout (tensor product) along the map $A\longrightarrow C^\infty(M)$ (see Lemma~\ref{lem:nil ext}). In particular, it follows that $A$ is complete if and only if it is complete as a dg $C^\infty(M)$-module (for any choice of section), and similarly for a dg $A$-module. Remark~\ref{rem:complete} then shows that completion has good homological properties; for example, the cohomology groups of the completion of $A$ can be identified with the completions of the $C^\infty(M)$-modules $\rmH^i(A)$.
	\qen
\end{remark}

\begin{definition}\label{def:inf-manifold}
	The $\infty$-category of \emph{inf-manifolds},
	$$
		\InfMfd \subseteq \CRing^{\infty, \opp}\ ,
	$$
	is the opposite of the full subcategory of $\CRing^\infty$ on those dg $C^\infty$-rings $A$ such that 
	\begin{myenumerate}
		\item $A_{\red}=C^\infty(M_{\red})$ for a smooth manifold $M_{\red}$,
		
		\item $\rmH^0(A)\longrightarrow A_{\red}$ is a nilpotent extension,
		
		\item $A$ is complete, and
		
		\item $\rmH^*(A)$ is concentrated in finitely many degrees.
	\end{myenumerate}
	For any such dg $C^\infty$-ring $A$, we  write $\Spec(A)$ for the corresponding inf-manifold and for any inf-manifold $M$, we  write $\calO(M)$ for the corresponding dg $C^\infty$-ring. We also write
	$$
		M_{\red}=\Spec(A_\red)
	$$
	for the smooth manifold underlying an inf-manifold $M=\Spec(A)$.
\end{definition}

\begin{example}
	Let $M$ be a smooth manifold and $$A=\RN[x_1, \dots, x_n]\,\big/\,I$$ a dg Artin algebra. Then the dg $C^\infty$-ring $C^\infty(M)\otimes_\RN A$ is a nilpotent extension of $C^\infty(M)$ and complete (since it is a finite complex of free $C^\infty(M)$-modules). Geometrically, 
	$$
		\Spec \big( C^\infty(M)\otimes_\RN A \big)
		\simeq M \times \Spec(A)
	$$
	looks like the product of a smooth manifold and an infinitesimally thickened point.
	\qen
\end{example}

\begin{example}
	\label{ex:p-th neighbourhood}
	Let $i\colon M\longhookrightarrow N$ be a closed embedding of smooth manifolds, and let $I$ denote the ideal of the surjection of $C^\infty$-rings $i^*\colon C^\infty(N)\longtwoheadrightarrow C^\infty(M)$.
	The \emph{$p$-th order infinitesimal neighbourhood} of $M$ in $N$ is the inf-manifold
	\begin{equation}\label{eq:p-th neighbourhood}
	N^{(p)}_M = \Spec \big( C^\infty(N)\,\big/\,I^{p+1} \big) \ .
	\end{equation}
	In particular, the underlying smooth manifold of $N^{(p)}_M$ is $M$. 
	In the special case where $i\colon M\longhookrightarrow N$ is given by $\RN^m\times \{0\}\longhookrightarrow \RN^{m+k}$, the $p$-th order infinitesimal neighbourhood is given by
	\begin{equation}
	\Spec \big( C^\infty(\RN^{m+k})\,\big/\,(y_1, \dots, y_k)^{p+1}\big)\simeq 
	\Spec \big( C^\infty(\RN^m)\otimes_{\RN} \RN[y_1, \dots y_k]\,\big/\,(y_1, \dots, y_k)^{p+1}\big) \ ,
	\end{equation}
	where $y_1, \dots, y_k$ are the last $k$ coordinates of $\RN^{m+k}$. For a general closed embedding, the inf-manifold $N^{(p)}_M$ is locally of this form.
	\qen
\end{example}

\begin{remark}
	\label{rem:almost inf-manifolds}
	Removing the boundedness condition~(4) from Definition~\ref{def:inf-manifold} leads to another notion of nilpotent thickenings of smooth manifolds; let us refer to these as \emph{almost inf-manifolds}. The $\infty$-category of almost inf-manifolds $\InfMfd^\mathrm{alm}$ has the advantage of being closed under the following operations:
	\begin{myenumerate}
		\item Let us say that a sequence of almost inf-manifolds $$M_0 \longrightarrow M_1\longrightarrow \cdots$$ is \emph{almost eventually constant} if for every $n\geqslant 0$, the tower of dg $C^\infty$-rings $$\tau^{\geqslant -n}\, \calO(M_0)\longleftarrow \tau^{\geqslant -n}\, \calO(M_1) \longleftarrow \cdots$$ is eventually constant.
		The basic example is (the opposite of) the Postnikov tower. Then $\colim_{k\in\NN_0} M_k$ is an almost inf-manifold, with $$\calO\Big(\underset{k\in\NN_0}\colim\, M_k\Big)=\underset{k\in\NN_0}\lim\, \calO(M_k) \ . $$
		
		\item Let $N\longrightarrow M$ be a submersion between smooth manifolds and $M'\longrightarrow M$ a map from an almost inf-manifold. Then there exists a fibre product $M' \times_M N$ in $\InfMfd^\mathrm{alm}$, which coincides with the usual fibre product if $M'$ is a smooth manifold:
		indeed, $\calO(M' \times_M N)$ is simply the completion of the dg $C^\infty$-ring $A = \calO(M') \amalg^h_{\calO(M)} \calO(N)$, which has $A_{\red} = C^\infty(M'_{\red} \times_M N)$.
	\end{myenumerate}
	The restriction functor $\PSh_\infty(\InfMfd^\mathrm{alm})\longrightarrow \PSh_\infty(\InfMfd)$ has a fully faithful right adjoint given by right Kan extension, which can also be presented as sending $X\colon \InfMfd\longrightarrow \scS$ to the presheaf  $\Spec(A)\longmapsto \lim_{n\in\NN_0} X(\Spec(\tau^{\geqslant -n}A))$. The essential image of this right adjoint consists of  presheaves $Y$ such that $Y(\colim_{k\in\NN_0} M_k)\simeq \lim_{k\in\NN_0} Y(M_k)$ for every almost eventually constant sequence of almost inf-manifolds. Essentially all presheaves of geometric interest have this property.
	Furthermore, this embedding restricts to the level of $\infty$-sheaves with respect to the Grothendieck topology below.
	\qen
\end{remark}

\begin{definition}
	A map of inf-manifolds $j \colon U = \Spec(B) \longrightarrow \Spec(A ) =M$ is an \emph{open embedding} if $U_{\red}\longhookrightarrow  M_{\red}$ is an open embedding of smooth manifolds and the map $A \longrightarrow B$ exhibits $B$ as the initial complete dg $C^\infty$-ring under $A$ with the following property: each element in $\rmH^0(A)$ whose image under $\rmH^0(A)\longrightarrow C^\infty(M_{\red})\longrightarrow C^\infty(U_\red)$ is invertible is sent to an invertible element of $\rmH^0(B)$.
\end{definition}

For every inf-manifold $M$, the full subcategory of $\InfMfd_{/M}$ on the open embeddings is equivalent to the poset of open subsets of the smooth manifold $M_{\red}$: for every open subspace $U_{\red}\subseteq M_{\red}$, the corresponding open subspace of $M$ is $\Spec(\calO_A(U_{\red}))$ where $\calO_A$ is the sheaf from Definition~\ref{def:complete2}.

We endow $\InfMfd$ with a Grothendieck topology by declaring that a family of open embeddings $\{U_a \longrightarrow M\}_{a \in \varLambda}$ is a \textit{cover} if and only if $\{U_{a, \red} \longrightarrow M_{\red}\}_{a \in \varLambda}$ is an open cover of the underlying smooth manifold in the usual sense. There is then an adjoint pair
\begin{equation}
	\begin{tikzcd}
		\iota : \Mfd\arrow[r, hookrightarrow, yshift=1ex, "\perp"']
		& \InfMfd\arrow[l, yshift=-1ex] : r\ ,
	\end{tikzcd}
\end{equation}
where $\iota$ is the natural inclusion (which is, in particular, fully faithful) and $r$ sends each inf-manifold $M$ to its underlying smooth manifold $M_{\red}$.
Both functors preserve open covers and hence induce adjunctions between $\infty$-categories of $\infty$-sheaves:
\begin{equation}
	\label{diag:infinitesimal cohesion}
	\begin{tikzcd}[column sep = 2.5cm]
		\Sh_\infty(\Mfd)
		\ar[r, hookrightarrow, bend left = 30, "\iota_!" {description, pos = 0.45}, "\perp"' {yshift = -0.1cm, pos = 0.405}]
		\ar[r, shift right = 0.25cm, "\iota_* = r^*" description]
		& \Sh_\infty(\InfMfd) \ .
		\ar[l, shift right = 0.25cm, "\iota^* = r_!" description, "\perp" yshift = -0.05cm]
		\ar[l, hookrightarrow, bend left = 30, "r_*" {description, pos = 0.52}, "\perp"' {yshift = 0.1cm, pos = 0.595}]
	\end{tikzcd}
\end{equation}

Because $\Mfd$ has finite products and $\iota$ preserves these (as a consequence of Remark~\ref{rem:almost inf-manifolds}), the functor $\iota_!$ preserves finite products as well.
Furthermore, $r^*$ has a further right adjoint $r_*$, since $r$ has the covering lifting property (see for example~\cite[Appendix~A]{Pstragowski:Synthetic_spectra} or~\cite[Section~A.2]{ADH:Differential_Cohomology} for background).
In the terminology of~\cite[Definition 4.2.1]{Schreiber:DCCT}, this equips the $\infty$-topos $\Sh_\infty(\Mfd)$ with \emph{differential cohesion}.
In practise, we will often leave the functor $\iota$ implicit and treat $\Mfd$ as an $\infty$-subcategory of $\InfMfd$.

We will always consider $\infty$-sheaves on smooth manifolds as $\infty$-sheaves on inf-manifolds via the functor $\iota_!$. For example, $\rmB^n \rmU(1)$ can be viewed as a higher stack on inf-manifolds, given by the $n$-fold delooping of $\rmU(1)$ viewed as an inf-manifold that happens to be an ordinary smooth manifold.

The adjunctions \eqref{diag:infinitesimal cohesion} give rise to two idempotent endofunctors on $\mathrm{Sh}_\infty(\InfMfd)$ whose value on $X$ we will denote by
\begin{equation}
	\label{eq:X_red --> X --> X_dR}
	\begin{tikzcd}
		X_\red = \iota_!\, \iota^*X \arrow[r]
		& X \arrow[r]
		& r^*\, \iota^*X = X_{\dR} \ ,
	\end{tikzcd}
\end{equation}
and refer to respectively as the \emph{reduced stack} and the \emph{de Rham stack} of $X$. The de Rham stack is given explicitly by
\begin{equation}
	X_{\dR}(M)=X(M_\red)\ ,
\end{equation}
for each inf-manifold $M$.
As a pleasant property, forming the de Rham stack $X \longmapsto X_{\dR}$ preserves both limits and colimits.
Moreover, if $M$ is a smooth manifold with a map $x \colon M \longrightarrow X$, then there is a canonical map
\begin{equation}
	X(M)
	\simeq \Sh_\infty(\InfMfd)(\iota_! M, X)
	 \longrightarrow \Sh_\infty(\InfMfd)(\iota_! M, \iota_!\, \iota^* X)
	\simeq X_\red(M)\ ,
\end{equation}
where we used that, for any functor of $\infty$-categories $F \colon \scC \longrightarrow \scD$, the Yoneda embeddings $h$ satisfy  $h_{F(c)} \simeq F_!\, h_c$ for each $c \in \scC$, and that $\InfMfd$ is subcanonical (Proposition~\ref{prop:subcanonical} below).
In other words, any point of $X$ from a smooth manifold factors through $X_\red$.

\begin{proposition}
	\label{prop:subcanonical}
	The Grothendieck topology on $\InfMfd$ is subcanonical, i.e.\ every representable presheaf is a sheaf.
\end{proposition}

\begin{proof}
	Let $M$ and $N$ be inf-manifolds and let $\CU$ be a covering sieve of $M$.
	We have to verify that
	\begin{equation}
		\begin{tikzcd}
			\InfMfd(M, N)\arrow[r] & \underset{U \in \CU}{\lim} \ \InfMfd(U, N)
		\end{tikzcd}
	\end{equation}
	is an equivalence. In terms of dg $C^\infty$-rings, this corresponds to the assertion that
	\begin{equation}
		\begin{tikzcd}
			\CRing^\infty \big( \calO(N), \calO(M) \big) \arrow[r]
			& \underset{U\in \CU}{\lim} \ \CRing^\infty \big(\calO(N), \calO(U)\big)
		\end{tikzcd}
	\end{equation}
	is an equivalence.
	It therefore suffices to verify  $\calO(M)\simeq \lim_{U\in \CU}\,\calO(U)$, which holds exactly by completeness of $\calO(M)$ (Definition~\ref{def:complete2}).
\end{proof}

\begin{example}
	\label{ex:function algebras}
	Every $\infty$-sheaf $X\in \Sh_\infty(\InfMfd)$ has an associated \emph{function algebra} $\Gamma(X, \calO_X)$, given by an---a priori \emph{unbounded}---commutative dg algebra over $\RN$ (well-defined up to quasi-isomorphism).
	In abstract terms, this is defined as follows.
	Let $\cdga^\infty_{\RN}$ denote the $\infty$-category associated to the model category \smash{$\cdga_\RN^\dg$} of unbounded cdgas (with its model structure inherited from that on \smash{$\Mod_\RN^\dg$} along the forgetful functor).
	Then there is a functor
	\begin{equation}
		\calO\colon \InfMfd^{\opp} \longrightarrow \cdga_\RN^\infty\ ,
		\qquad
		\Spec(A) \longmapsto A\ ,
	\end{equation}
	which sends each inf-manifold to the corresponding dg $C^\infty$-ring and then applies the forgetful functor from dg $C^\infty$-rings to cdgas (which we suppress in the notation).
	
	Recall that any open cover $\{U_a \longrightarrow M\}_{a \in \varLambda}$ of an inf-manifold is determined by an open cover of the underlying manifold.
	Thus, for each cover $\{U_a\longrightarrow M\}_{a \in \varLambda}$, the map
	\begin{equation}
		\calO(M) \longrightarrow \underset{[n]\in \bbDelta}{\holim} \ \prod_{a_0, \dots, a_n \in \varLambda}\calO(U_{a_0 \cdots a_n})
	\end{equation}
	is a quasi-isomorphism by completeness of $\calO(M)$.
	Here we have used the notation $U_{a_0 \cdots a_n} = \bigcap_{i = 0}^n\, U_{a_i}$.
	This implies that $\calO$ extends to an adjoint pair
	\begin{equation}
		\label{eq: coSpec-functions adjunction}
		\begin{tikzcd}
			\coSpec : \cdga_\RN^\infty \arrow[r, yshift=1ex, "\perp"']
			&\Sh_\infty(\InfMfd)^\opp : \Gamma(-, \calO)\ . \arrow[l, yshift=-1ex]
		\end{tikzcd}
	\end{equation}
	
	Explicitly, if $X$ can be written as a homotopy colimit $$X = \underset{\alpha}\hocolim\ \Spec(A_\alpha)$$ of inf-manifolds, then $$\Gamma(X, \calO_X)\simeq \underset{\alpha}\holim\ A_\alpha \ . $$ The left adjoint $\coSpec$ can be thought of as sending a cdga to the corresponding affine stack in the sense of To\"en~\cite{Toen:Champs_Affines} (called coaffine stack by Lurie \cite{Lurie:DAG-VIII}).
	\qen
\end{example}

\begin{example}
	\label{ex:function algebras and maps to B^n RN}
	The affine stack $\coSpec(\RN[x_{n}])$ associated to a polynomial algebra on a single degree~$n$ generator is given by the stack $\rmB^n\RN$. It then follows from the adjunction~\eqref{eq: coSpec-functions adjunction} that for each $n\geqslant 0$,
	\begin{align}
		\Sh_\infty(\InfMfd) \big( X, \rmB^n \RN \big)
		&\simeq \Map_{\cdga_{\RN}} \big( \RN[x_n], \calO(X) \big)
		\\[4pt]
		&\simeq \Map_{\Mod^\dg_\RN} \big( \RN[-n], \calO(X) \big)
		\\[4pt]
		&\simeq \DK \Big( \tau^{\leqslant 0} \big( \Gamma(X, \calO_X)[n] \big) \Big)
	\end{align}
	is the space corresponding to the non-positively graded complex
	\begin{equation}
		0 \longleftarrow Z^n \Gamma(X, \calO_X)
		\longleftarrow \Gamma(X, \calO_X)^{n-1}
		\longleftarrow \cdots
	\end{equation}
	under the Dold--Kan correspondence, where $Z^n \Gamma(X, \calO_X)$ sits in degree zero.
	\qen
\end{example}


\subsection{Infinitesimal structure}\label{sec:infinitesimal}


The infinitesimal thickenings of smooth manifolds included in $\InfMfd$ allow us to study the infinitesimal structure of a stack on $\InfMfd$. For this to be well-behaved, we will impose an analogue of the Schlessinger condition, as studied for example in~\cite{GR:Study_Derived_Geom} in the algebro-geometric setting.

\begin{definition}
	\label{def:schlessinger}
	Let $S$ be a class of morphisms in $\InfMfd$, and let $S_0$ be the class of morphisms $g \colon \Spec(A) \longrightarrow \Spec(A')$ in $\InfMfd$ such that the induced map $\rmH^0(A') \longrightarrow \rmH^0(A)$ is surjective with nilpotent kernel.
	A prestack $X \colon \InfMfd^{\opp} \longrightarrow \scS$ satisfies the \emph{Schlessinger condition along morphisms in $S$} if $X$ sends any pushout square
	\begin{equation}
		\label{eq:pushout_sqz}
		\begin{tikzcd}[column sep = 1.25cm, row sep = 1cm]
			\Spec(A) \arrow[d, "g "'] \arrow[r, "f"]
			& \Spec(B)\arrow[d]
			\\
			\Spec(A')\arrow[r]
			& \Spec(B')
		\end{tikzcd}
	\end{equation}
 in $\InfMfd$, with $g \in S_0$ and $f \in S$, to a pullback square of spaces.
\end{definition}

\begin{remark}
	\label{rmk: inf-cohesive prestack}
	There are various choices for the class of maps along which we impose the Schlessinger condition; see Definition~\ref{def:qcoh def thy} for two choices that will be of interest to us.
	For now, we let $S = S_0$ be the class of morphisms; prestacks satisfying the Schlessinger condition along all maps in $S_0$ are often called \emph{infinitesimally cohesive}.
	
	In this situation, since $\ker (\rmH^0 (f))$ and $\ker (\rmH^0 (g))$ are nilpotent, all corners of~\eqref{eq:pushout_sqz} have the same underlying smooth manifold.
	One can therefore think of the Schlessinger condition as an infinitesimal gluing condition. 
	\qen
\end{remark}

One of the Schlessinger condition's main features is that it gives rise to a well-behaved notion of tangent space.
To see this, recall that the tangent space $T_xM$ of a smooth manifold can be identified with the set of first order infinitesimal paths in $M$ around $x$.
The $\infty$-category $\InfMfd$ allows us to describe such first order infinitesimal paths as maps out of the inf-manifold $\Spec(\RN[\epsilon]/\epsilon^2)$. More generally, by giving $\epsilon$ a negative cohomological degree, one obtains spaces of derived first order infinitesimal paths.
If $X$ is infinitesimally cohesive, the resulting spaces of first order paths can be organised into a \emph{linear} object---just as for ordinary smooth manifolds---using the $\infty$-categorical construction we now describe.

\begin{definition}
	\label{def: reduced excisive}
	Let $F \colon \scC \longrightarrow \scD$ be an $\infty$-functor.
	\begin{myenumerate}
		\item $F$ is  \emph{excisive} if it sends all pushout squares in $\scC$ to pullback squares in $\scD$.
		We write $\mathrm{Exc}(\scC, \scD)$ for the full $\infty$-subcategory of the $\infty$-functor category $\Fun(\scC, \scD)$ on the excisive functors.
		
		\item If $\scC$ also has a final object, then $F$ is  \emph{reduced excisive} if it additionally preserves the final object.
		We write $\mathrm{Exc}_*(\scC, \scD)$ for the full $\infty$-subcategory of $\Fun(\scC, \scD)$ on the reduced excisive functors.
	\end{myenumerate}
\end{definition}

\begin{remark}
	\label{rem:reduced excisive}
	Let $A$ be a non-positively graded cdga and let $\Perf^{\leqslant 0}_A$ denote its $\infty$-category of non-positively graded \emph{perfect} modules, i.e.\ the smallest subcategory of $\Mod_A^\infty$ containing $A$ which is closed under retracts and finite homotopy colimits.
	It follows from Remark~\ref{rem:complete}(2) that if $A = C^\infty(M)$, for a smooth manifold $M$, then each perfect $A$-module is also complete.
	
	Every $A$-module $E$ determines a functor
	$$
		\Phi_E\colon \Perf_A^{\leqslant 0} \longrightarrow \scS\ ,
		\qquad
		\Phi_E(I) = \DK \big( E \otimes_A I \big)\ .
	$$
	The discussion in~\cite[Construction 17.5.1.1]{Lurie:SAG} shows that this determines an equivalence 
	\begin{equation}
		\label{eq: modules as red exc fctrs}
		\Mod_A^\infty\ \eq\ \mathrm{Exc}_\ast \big( \Perf^{\leqslant 0}_A, \scS \big)
		\subseteq \Fun\big(\Perf^{\leqslant 0}_A, \scS\big)\ ,
		\qquad E \longmapsto \Phi_E
	\end{equation}
	between the $\infty$-category of $A$-modules and the $\infty$-category of reduced excisive functors $\Perf^{\leqslant 0}_A\longrightarrow \scS$.~\qen
\end{remark}

Now let $x \colon M \longrightarrow X$ be a map from an inf-manifold $M$ to a prestack $X$.
Let $\Mod_A^{\leqslant 0, \infty}$ denote the $\infty$-category of non-positively graded dg $A$-modules.
For any non-positively graded, bounded and complete $\calO(M)$-module $I$, we write
\begin{equation}
	\label{eq:derivations}
	\Der_X(x; I) = X\big( \Spec(\calO(M) \oplus I) \big) \times_{X(M)} \{x\}
\end{equation}
for the space of first order paths around $x$ parametrised by $\Spec(\calO(M)\oplus I)$.

When $X$ is infinitesimally cohesive, the functor \smash{$\Der_X(x; -)\colon \Mod_{\calO(M)}^{\leqslant 0, \infty} \longrightarrow \scS$} is reduced excisive:
first, note that $\Der_X(x;-)$ is reduced.
Then by~\cite[Proposition~1.4.2.13]{Lurie:HA} it suffices to check that it sends suspension pushouts to pullbacks.
We can see this as follows:
in $\Mod^{\leqslant 0, \infty}_{\calO(M)}$, the canonical square
\begin{equation}
	\begin{tikzcd}[row sep = 1cm,column sep=1cm]
		I \ar[r] \ar[d]
		& 0 \ar[d]
		\\
		0 \ar[r]
		& I[1]
	\end{tikzcd}
\end{equation}
is both cartesian and cocartesian.
Taking square-zero extensions by the corners of this square produces a commutative square of dg $C^\infty$-rings which is a pullback  whose  legs are surjective on $\rmH^0$ with nilpotent kernel, and so its image under $\Der_X(x;-)$ is a pullback square of $\infty$-groupoids, as desired.

Thus there exists a unique (in the $\infty$-categorical sense) $\calO(M)$-module $T_x X$ which presents the restriction of $\Der_X(x; -)$ to $\Perf_{\calO(M)}^{\leqslant 0}$ under the equivalence~\eqref{eq: modules as red exc fctrs}.
We can therefore make 

\begin{definition}
	\label{def:tangent}
	Let $X$ be an infinitesimally cohesive prestack, and let $ x\colon M\longrightarrow X$ be a map from an inf-manifold.
	The \emph{tangent complex} of $X$ at the point $x$ is the unique object $T_x X \in \Mod_{\calO(M)}^\infty$ equipped with a natural equivalence
	$$
		\DK \big(T_xX\otimes_{\calO(M)} I \big)
		\simeq \Der_X(x; I)\ ,
	$$
	for $I \in \Perf_{\calO(M)}^{\leqslant 0}$.
\end{definition}

Taking $I=\calO(M)[n]$, one sees that, up to equivalence, the complex
$$\begin{tikzcd}
Z^nT_xX & \arrow[l] (T_xX)^{n-1} & \arrow[l] \cdots
\end{tikzcd}$$
corresponds under the Dold--Kan correspondence to the space $\Der_X(x; \calO(M)[n])$ of shifted tangent vectors at $x$.

\begin{example}
	\label{ex:corep tangent}
	If the functor \smash{$\Der_X(x; -) \colon \Perf_{\calO(M)}^{\leqslant 0} \longrightarrow \scS$} is \emph{corepresented} by a dg $\calO(M)$-module $L$, then $T_xX$ is the $\calO(M)$-linear dual $L^\vee$.
	For example, for any map $x\colon M\longrightarrow N$ to a smooth manifold, the functor $\Der_N(x; -)$ is corepresented by the module $\Omega^1(N)\otimes_{C^\infty(N)} C^\infty(M)\simeq \Gamma(M, x^*T^*N)$ of $1$-forms. It follows that $T_xN$ is the module of sections of the pullback of the tangent bundle $TN$ along~$x$.
	\qen
\end{example}

To study the dependence of the tangent complex $T_xX$ on the basepoint $x$, it will be convenient to use a more global version of the above definitions. Restricting attention to basepoints $x\colon M\longrightarrow X$ where $M$ is a smooth manifold (which suffices for our purposes), let us consider the  diagram of $\infty$-categories
\begin{equation}
	\label{eq: GrConstr over Mfd for T_x X}
	\begin{tikzcd}
		\Perf^{\leqslant 0} \arrow[r, hookrightarrow, "\kappa"]
		& \Mod^{\leqslant 0, \infty} \arrow[r]
		& \Mfd^{\opp}\ .
	\end{tikzcd}
\end{equation}
Here $\Mod^{\leqslant 0, \infty}$ is the Grothendieck construction of the $\infty$-functor $\Mfd^\opp \longrightarrow \Cat_\infty$ which sends a manifold $M$ to its $\infty$-category of bounded, complete, and non-positively graded dg $C^\infty(M)$-modules \smash{$\Mod^{\leqslant 0, \infty}_{C^\infty(M)}$}; that is, it is the $\infty$-category of pairs $(M, I)$ of a smooth manifold $M$ and a bounded, complete, and non-positively graded dg $C^\infty(M)$-module $I$, with morphisms $(M, I)\longrightarrow (N, J)$ given by a smooth map $f \colon N \longrightarrow M$ and a $C^\infty(M)$-linear map \smash{$I\longrightarrow C^\infty(M) \otimes_{C^\infty(N)} J$}.
Likewise, $\Perf^{\leqslant 0}$ is the full $\infty$-subcategory of $\Mod^{\leqslant 0, \infty}$ on those objects $(M, I)$ where $I$ is perfect.

Every prestack $X$ gives rise to a functor
$$
	\Der_X\colon \Mod^{\leqslant 0,\infty} \longrightarrow \scS\ ,
	\qquad
	(M, I) \longmapsto X \big( \Spec(C^\infty(M) \oplus I) \big)\ .
$$
If $X$ is an $\infty$-sheaf, then $\Der_X$ is an $\infty$-sheaf for the following Grothendieck topology on $(\Mod^{\leqslant 0, \infty})^\opp$:
the open covers of $(M, I)$ are of the form $(U_a, j_a^*I)$ with $j_a\colon U_a\longrightarrow M$ an open cover in the usual sense and $j_a^*I$ the pullback of the complete module $I$ (see Remark~\ref{rem:complete}(3)).

The fibre of $\Der_X(M, I) \longrightarrow \Der_X(M, 0)\simeq X(M)$ over a point $x$ is precisely the space $\Der_X(x; I)$ from \eqref{eq:derivations}. In particular, the restriction of $\Der_X$ to the full subcategory $\Perf^{\leqslant 0}$ encodes the data of all the tangent complexes of $X$ at basepoints parametrised by smooth manifolds. In addition, one can use the tangent complex to approximate the values of $\Der_X$ on non-perfect modules:
recall the inclusion $\infty$-functor $\kappa$ introduced in~\eqref{eq: GrConstr over Mfd for T_x X}.

\begin{lemma}
	\label{lem:qcdt}
	Let $X$ be an infinitesimally cohesive prestack, $M$ a smooth manifold and $I$ a bounded, complete, and non-positively graded $C^\infty(M)$-module.
	Consider the diagram
	$$
		\begin{tikzcd}[column sep = {2.5cm,between origins}, row sep = 1cm]
			\kappa_!\,\kappa^*\Der_X(M, I) \arrow[rd] \arrow[rr]
			& & \Der_X(M, I)\arrow[ld]
			\\
			& \Der_X(M, 0)=X(M)
		\end{tikzcd}
	$$
	Then for any point $x \in X(M)$, the induced map between the fibres over $x$ of the diagonal arrows can be identified with a natural map
	\begin{equation}
		\label{eq:qcdt comparison map}
		\DK \big(T_x X\otimes_{C^\infty(M)} I\big) \longrightarrow \Der_X(x; I)\ .
	\end{equation}
\end{lemma}

\begin{proof}
	It follows by the pointwise formula for the left Kan extension that
	$$
		\kappa_!\,\kappa^*\Der_X(M, I) \simeq \underset{(N, J)\in \Perf^{\leqslant 0}_{/(M, I)}}{\colim} \, \Der_X(N, J)
		\simeq \underset{J\in (\Perf^{\leqslant 0}_{C^\infty(M)})_{/I}}{\colim} \, \Der_X(M, J)\ ,
	$$
	where we have used that the inclusion
	\begin{equation}
		\begin{tikzcd}
			\big( \Perf^{\leqslant 0}_{C^\infty(M)} \big)_{/I} \ar[r, hookrightarrow]
			& \Perf^{\leqslant 0}_{/(M, I)}
		\end{tikzcd}
	\end{equation}
	is cofinal (see the proof of~\cite[Lemma 4.3.3.9]{Lurie:HTT}).
	The first overcategory is taken with respect to the inclusion \smash{$\Perf^{\leqslant 0}_{C^\infty(M)} \longhookrightarrow  \Mod^{\leqslant 0, \infty}_{C^\infty(M)}$}.
	Since \smash{$\big(\Perf^{\leqslant 0}_{C^\infty(M)}\big)_{/I}$} has finite colimits, it is a filtered $\infty$-category.
	It follows that the  colimit is a filtered colimit and commutes with taking fibres. Consequently, the fibre of $\kappa_!\,\kappa^*\Der_X(M, I)$ over a point $x$ can be identified with
	$$
		\underset{J \in (\Perf^{\leqslant 0}_{C^\infty(M)})_{/I}}{\colim} \, \Der_X(x; J)
		\simeq \underset{J \in (\Perf^{\leqslant 0}_{C^\infty(M)})_{/I}}{\colim} \, \DK \big( T_x X \otimes_{C^\infty(M)} J \big)
		\simeq \DK \big(T_x X \otimes_{C^\infty(M)} I \big)\ ,
	$$
	where the last equivalence uses that the Dold--Kan correspondence and the tensor product preserve filtered colimits.
\end{proof}

Using this, we now specify two classes of (pre)stacks whose infinitesimal structure is controlled by the tangent complex.

\begin{definition}
	\label{def:qcoh def thy}
	A prestack $X$ has \emph{deformation theory} (resp.\ \emph{quasi-coherent deformation theory}) if it satisfies:
	\begin{myenumerate}
		\item (\emph{Schlessinger condition}) It satisfies the Schlessinger condition from Definition~\ref{def:schlessinger} along the class of maps $f \colon \Spec(A) \longrightarrow \Spec(B)$ in $\InfMfd$ such that $\rmH^0(B) \longrightarrow \rmH^0(A)$ is surjective (resp.\ along all maps $f \colon \Spec(A) \longrightarrow \Spec(B)$).
		
		\item (\emph{Tangent representability}) The map $\kappa_!\, \kappa^* \Der_X \longrightarrow \Der_X$ from Lemma~\ref{lem:qcdt} induces an equivalence on associated $\infty$-sheaves on $\Mod^{\leqslant 0, \infty}$.
	\end{myenumerate}
	We will abbreviate the term `prestacks with (quasi-coherent) deformation theory' to (\emph{qc})\emph{dt} prestacks.
	
	The prestack $X$ has \textit{strictly} (\emph{quasi-coherent}) \emph{deformation theory}, or is a \emph{strict} (\emph{qc})\emph{dt prestack}, if it satisfies condition~(1) and the map in~(2) is an equivalence even before sheafifying.
\end{definition}

\begin{remark}\label{rem:def thy on sqz}
	By Lemma~\ref{lem:nil ext}, every map in the class $S$ from Definition~\ref{def:schlessinger} decomposes into a finite sequence of square-zero extensions.
	In turn, every square-zero extension $\Spec(A) \longrightarrow \Spec(A_\eta)$ of inf-manifolds arises as the pushout of the canonical projection $\Spec(A \oplus I[1]) \longrightarrow \Spec(A)$.
	To check that a prestack has deformation theory, it therefore suffices to verify that it sends a pushout \eqref{eq:pushout_sqz} to a pullback in the case where $A'\longrightarrow A=A'\oplus I[1]$ is the zero map and $\rmH^0(B)\longrightarrow \rmH^0(A)$ is surjective. The same argument applies to quasi-coherent deformation theory.
	\qen
\end{remark}

	We will be mainly interested in the case where $X$ is a stack. In this case, the tangent representability condition (2) in Definition~\ref{def:qcoh def thy} implies that the tangent complex controls the infinitesimal structure of $X$ through
	
\begin{proposition}
	\label{prop:qcdt stack 2}
	Let $X$ be a dt stack and let $x\colon M\longrightarrow X$ be a map from a smooth manifold.
	Then for every bounded, complete, and non-positively graded $C^\infty(M)$-module $I$, there is a natural equivalence
	$$
		\DK \big((T_xX\otimes_{C^\infty(M)} I)^\cpl\big)
		\simeq \Der_X(x; I)
	$$
	where the domain is the Dold--Kan image of the completion of $T_x X \otimes_{C^\infty(M)} I$. 
\end{proposition}

That is, while the comparison morphism~\eqref{eq:qcdt comparison map} generally fails to be an equivalence if $I$ is not perfect, if $X$ is a dt stack the morphism still becomes an equivalence upon completing the module $T_x X \otimes_{C^\infty(M)} I$.

\begin{proof}
The values of the comparison map~\eqref{eq:qcdt comparison map} at the open subspaces $j\colon U\longhookrightarrow  M$ determine a natural transformation of presheaves on $M$ of the form
$$
	\begin{tikzcd}
		\DK \big(T_{j^*x}X\otimes_{C^ \infty(U)} j^*I\big) \arrow[r]
		& \Der_X(j^*x; j^*I) \ .
	\end{tikzcd}
$$
Since $X$ is a stack, the codomain of this map is an $\infty$-sheaf on $M$. By tangent representability, the codomain is the associated sheaf of the domain. The associated sheaf of the domain can be identified explicitly: since $\DK$ preserves associated sheaves (as a consequence of Remark~\ref{rem:complete}), it corresponds to the associated sheaf of the presheaf of modules $U \longmapsto T_{j^*x} \otimes_{C^ \infty(U)} j^*I$, i.e.\ the tensor product of the two sheaves sending $U$ to $T_{j^*x}$ and $j^*I$. By Lemma~\ref{lem:complete modules vs sheaves}, the global sections of this tensor product of sheaves is the \emph{completion} of $T_xX\otimes_{C^{\infty}(M)} I$.
\end{proof}

\begin{remark}
The proof of Proposition~\ref{prop:qcdt stack 2} explains why we do not require the map $\kappa_!\,\kappa^*\Der_X \longrightarrow \Der_X$ to be an equivalence in Definition~\ref{def:qcoh def thy}: the comparison maps \eqref{eq:qcdt comparison map} are seldom equivalences when $X$ is a stack. Indeed, their codomain satisfies descent, but the domain need not since the tensor product does not commute with infinite homotopy limits.
\qen
\end{remark}

The Schlessinger condition~(1) in Definition~\ref{def:qcoh def thy} allows us to relate tangent complexes at different points of a given dt stack $X$ through

\begin{proposition}
	\label{prop:qcdt stack 1}
	Let $X$ be a dt stack and let $x \colon M \longrightarrow X$ be a map from a smooth manifold. Then for every map of smooth manifolds $f \colon N\longrightarrow M$, there is a natural map of complete $C^\infty(N)$-modules
	$$
		\begin{tikzcd}
			f^*T_x X = \big( T_x X\otimes_{C^\infty(M)} C^{\infty}(N)  \big)^\cpl \arrow[r]
			& T_{f^*x}X\ ,
		\end{tikzcd}
	$$
	which is an equivalence if $f$ is an open or closed embedding.
	If $X$ is even a qcdt stack, then this map is an equivalence for all $f \colon N\longrightarrow M$.
\end{proposition}

	In other words, when $X$ is a qcdt stack, the tangent complexes $T_x X$, where $x$ ranges over all morphisms from a smooth manifold to $X$, define a quasi-coherent sheaf on $X_\red$ (the restriction of $X$ along the inclusion $\Mfd \subset \InfMfd$).

\begin{proof}
	Let $I$ be a perfect non-positively graded $C^\infty(N)$-module and let $f_*I$ be its underlying $C^\infty(M)$-module. Using Proposition~\ref{prop:qcdt stack 2}, we obtain a natural map 
	$$
		\begin{tikzcd}
			\DK \big( (T_x X \otimes_{C^\infty(M)} f_*I)^\cpl \big)
			\simeq \Der_X(x; f_*I)
			\arrow[r]
			& \Der_X \big( f^*x ; I \big)
			\simeq \DK (T_{f^*x} X \otimes_{C^\infty(N)} I)\ ,
		\end{tikzcd}
	$$
	which arises by restriction along $\Spec(C^\infty(N)\oplus I)\longrightarrow \Spec(C^\infty(M)\oplus f_*I)$. This natural transformation of reduced excisive functors on \smash{$\Perf^{\leqslant 0}_{C^\infty(N)}$} is equivalent to a map of $C^\infty(N)$-modules of the form $f^*T_xX\longrightarrow T_{f^*x}X$ (see Remark~\ref{rem:reduced excisive}).
	
	The Schlessinger condition implies that $\Der_X(x; f_*I) \longrightarrow \Der_X\big(f^*x ; I\big)$ is an equivalence whenever $X$ is a qcdt stack or whenever $X$ is a dt stack and $f$ is a closed embedding (so that $f^*:C^\infty(M)\longrightarrow C^\infty(N)$ is surjective). In both cases it follows that $f^*T_xX\simeq T_{f^*x}X$. Proposition~\ref{prop:qcdt stack 2} and its proof imply that if $X$ is a dt stack, then also $f^*T_xX\simeq T_{f^*x}X$ whenever $f$ is an open embedding.
\end{proof}

In particular, it follows that dt stacks admit a formal version of the Inverse Function Theorem given by

\begin{proposition}
	\label{prop:inverse_function}
	Let $f\colon Y\longrightarrow X$ be a map of dt stacks. Then $f$ is an equivalence if and only if
	\begin{myenumerate}
		\item the map $Y_{\red}\longrightarrow X_{\red}$ is an equivalence, and
		
		\item for every map $y\colon M\longrightarrow Y$ from a smooth manifold, the map $T_yY\longrightarrow T_{f(y)}X$ is an equivalence.
	\end{myenumerate}
\end{proposition}

\begin{remark}
	\label{rmk: detecting equivalences for qcdt}
	By Proposition~\ref{prop:qcdt stack 1}, as soon as condition (2) holds for a point $y\colon M\longrightarrow Y$, it automatically also holds for the restrictions of $y$ to all open and all closed submanifolds of $M$.
	If $X$ is a qcdt stack, then it automatically holds for all points of $X$ from manifolds which factor through $y$.
	\qen
\end{remark}

\begin{proof}
	Let $\mathcal{K}$ be the class of inf-manifolds $U$ such that $Y(U)\longrightarrow X(U)$ is an equivalence. Since $Y_{\red}\simeq X_{\red}$ by assumption, this class contains all smooth manifolds.
	Together with the fact that $f$ induces equivalences on tangent complexes, this implies that the induced map $\kappa^*\Der_Y\longrightarrow \kappa^*\Der_X$ is a natural equivalence (because the map is an equivalence on all fibres by Lemma~\ref{lem:qcdt}).
	Since $Y$ and $X$ are dt stacks, $\Der_Y$ and $\Der_X$ are the associated sheaves of $\kappa_!\,\kappa^*\Der_Y$ and $\kappa_!\,\kappa^*\Der_X$ by tangent representability. Consequently $\Der_Y\longrightarrow \Der_X$ is an equivalence as well, which means that $\mathcal{K}$ contains all square-zero extensions $\Spec(C^\infty(M)\oplus I)$ of a smooth manifold $M$ by a bounded, complete, and non-positively graded $C^\infty(M)$-module $I$. 
	 
	Because $X$ and $Y$ have deformation theory, the class $\mathcal{K}$ is then stable under square-zero extensions by bounded, complete, and non-positively graded $C^\infty(M)$-modules, for any smooth manifold $M$. Using Lemma~\ref{lem:nil ext} (and taking completions of all dg $C^\infty$-rings involved), one finds that any inf-manifold can be obtained by such a sequence of square-zero extensions, so that $\mathcal{K}$ contains all inf-manifolds.
\end{proof}


\subsection{Examples and constructions of stacks with deformation theory}
\label{sec:Examples of stacks with defthy}


Let us now turn to examples of stacks with (quasi-coherent) deformation theory, as well as several ways of constructing such stacks that will be important later in this paper.

\begin{example}\label{ex:mfd are qcdt}
	Every smooth manifold $N$ has strictly quasi-coherent deformation theory (see Definition~\ref{def:qcoh def thy}).
	Indeed, any representable stack satisfies the Schlessinger condition along all maps $\Spec(A) \longrightarrow \Spec(B)$, since $\InfMfd(-, N)$ sends all pushouts to pullbacks.
	
	For condition~(2) of Definition~\ref{def:qcoh def thy}, we note that for any map of smooth manifolds $x \colon M \longrightarrow N$, the functor $\Der_N(x; -)$ is corepresented by the finite-type projective module $L = \Omega^1(N) \otimes_{C^\infty(N)} C^\infty(M)$.
	By Example~\ref{ex:corep tangent}, the map \eqref{eq:qcdt comparison map} is then the natural map $L^\vee\otimes_{C^\infty(M)} I \longrightarrow \Hom_{C^\infty(N)}(L, I)$.
	It is an equivalence because finite-type projective modules are dualisable.
	\qen
\end{example}

\begin{example}\label{ex:de rham}
	Recall the notion of  de Rham stack from~\eqref{eq:X_red --> X --> X_dR}.
	For each stack $X \in \Sh_\infty(\InfMfd)$, the de Rham stack $X_{\dR}$ has strictly quasi-coherent deformation theory. The Schlessinger condition is satisfied because $\Spec(A)_\red = \Spec(A')_\red$, so that $X_{\dR}(\Spec(A))\simeq X_{\dR}(\Spec(A'))$, and likewise for $B$ and $B'$ (in the notation of Definition~\ref{def:schlessinger}).
	From~\eqref{eq:derivations} we  see  $\Der_{X_\dR}(x; I) \simeq *$, for any point $x$ and  any non-positively graded, bounded, and perfect module $I$, since $\Spec(\calO(M) \oplus I)_\red = \Spec(\calO(M))$.
	Therefore the domain and codomain of \eqref{eq:qcdt comparison map} are both contractible, implying condition~(2) of Definition~\ref{def:qcoh def thy}.
	\qen
\end{example}

To produce more examples of (pre)stacks with quasi-coherent deformation theory and help with computations later on, we record the following properties of (qc)dt prestacks and stacks:
we write $\Hom$ for the internal hom, or mapping (pre)stack, in both $\PSh_\infty(\InfMfd)$ and $\Sh_\infty(\InfMfd)$.

\begin{proposition}
	\label{prop:def thy}
	\begin{myenumerate}
		\item If $X$ is a dt (resp.~qcdt) prestack, the associated stack $\widetilde{X}$ is also dt (resp.~qcdt).
		
		\item The classes of prestacks and stacks with (quasi-coherent) deformation theory are each closed under finite limits and filtered colimits.
		
		\item If $X$ is a (pre)stack with deformation theory and $M$ is a smooth manifold, then the mapping stack $\Hom(M, X)$ has deformation theory as well.
		For every morphism $x^\dashv \colon N \longrightarrow \Hom(M, X)$ from a smooth manifold, adjoint to a morphism $x \colon M \times N \longrightarrow X$, there is an equivalence of complete $C^\infty(N)$-modules
		$$
			T_{x^\dashv}\Hom(M, X)
			\simeq \pi_{N*}\big(T_xX\big)\ ,
		$$
	where $\pi_N\colon M \times N\longrightarrow N$ is the canonical projection.
	\end{myenumerate}
\end{proposition}

Assertion~(3) does not hold for qcdt prestacks, because Lemma~\ref{lem:coprod of square-zero} does not hold without the condition that $\rmH^0(B) \longrightarrow \rmH^0(A)$ is surjective.

\begin{proof}
	For (1), let us first consider the Schlessinger condition. Given a pushout square of inf-manifolds as in~\eqref{eq:pushout_sqz}, together with an open subspace $U \longrightarrow \Spec(B')$, restricting to $U$ gives another pushout square of inf-manifolds. Furthermore, if the map $\rmH^0(B)\longrightarrow \rmH^0(A)$ was surjective in the original square, then the resulting map for the square of restrictions is surjective as well. 
	Varying $U$ and applying $X$ to each of the four corners of the square gives a commuting square of presheaves on the underlying manifold $\Spec(B')_{\red}$. Since $X$ has (quasi-coherent) deformation theory, this square of presheaves is a pullback square. Applying the same construction to $\widetilde{X}$ provides a commuting square of $\infty$-sheaves on $\Spec(B')_{\red}$, which is obtained from the square for $X$ by taking associated sheaves. The result now follows from the fact that $\infty$-sheafification is left exact. (See~\cite[Lemma 6.1.6]{Nuiten:Thesis} for a similar proof.)
	
	For tangent representability, note that the map $\Der_X\longrightarrow \Der_{\widetilde{X}}$ induces an equivalence on associated sheaves. Since $\kappa$ preserves covers and has the covering lifting property, it follows that $\kappa_!\,\kappa^*\Der_X\longrightarrow \kappa_!\,\kappa^*\Der_{\widetilde{X}}$ induces an equivalence on associated sheaves as well. The result then follows from the two-out-of-three property.
	
	For assertion (2), note that the Schlessinger condition is a pullback condition. Consequently, any limit or filtered colimit of prestacks satisfying the Schlessinger condition will do so as well. Tangent representability follows from the fact that $\kappa_!\,\kappa^*\Der_X$ and $\Der_X$ (as well as their associated sheaves) both preserve fibre products and pointwise colimits in the variable $X$.
	
	For (3), we will use the fully faithful inclusion $\PSh_\infty(\InfMfd)\longhookrightarrow  \PSh_\infty(\InfMfd^\mathrm{alm})$ from Remark~\ref{rem:almost inf-manifolds}. Since this inclusion preserves mapping stacks, we may view $X$ and $M$ as (representable) presheaves on almost inf-manifolds and work with $\Hom(M, X)$ in $\PSh_\infty(\InfMfd^\mathrm{alm})$. The Schlessinger condition for $\Hom(M, X)$ then follows by adjunction: for any square of the form \eqref{eq:pushout_sqz}, taking the product of all its objects with the smooth manifold $M$ produces another square of the form \eqref{eq:pushout_sqz}, by Lemma~\ref{lem:coprod of square-zero}. The result then follows from the fact that $X$ sends such squares to pullback squares.
	
	For the tangent representability condition, let us fix a map $x^\dashv \colon N \longrightarrow \Hom(M, X)$, adjoint to a map $x \colon M {\times} N \longrightarrow X$. Using Lemma~\ref{lem:coprod of square-zero}, for each $C^\infty(N)$-module $I$ there is a natural equivalence
	$$
		\Der_{\Hom(M, X)}(x^{\dashv}; I)
		\simeq \Der_{X} \big( x; C^\infty (M \times N) \otimes^h_{C^\infty(N)} I \big)\ .
	$$
	By Remark~\ref{rem:reduced excisive} and Definition~\ref{def:tangent} this implies that $T_{x^\dashv}\Hom(M, X)\simeq \pi_{N*}\big(T_xX\big)$, and that there are equivalences
	$$
		\begin{tikzcd}[row sep = 1.5cm,column sep=1cm]
			\DK \big(T_{x^{\dashv}} \Hom(M, X) \otimes^h_{C^\infty(N)} I \big) \arrow[r] \arrow[d, "\sim"{swap}]
			& \Der_{\Hom(M, X)}(x^{\dashv} ; I) \arrow[d, "\sim"]
			\\
			\DK \big( (\pi_{N*}T_{x} X) \otimes^h_{C^{\infty}(N)} I \big) \arrow[r]
			& \Der_X \big( x; C^\infty(M \times N) \otimes^h_{C^\infty(N)} I \big)
		\end{tikzcd}
	$$
	The bottom map induces equivalences on $\infty$-sheaves on $N$ because $X$ is a dt prestack, so that the same holds for the top map.
	This completes the proof by Lemma~\ref{lem:qcdt}.
\end{proof}

\begin{example}
	\label{ex:vertical forms def thy}
	Let $M$ be a smooth manifold, and consider the functor $$\Omega^p_v(M) \colon \InfMfd^{\opp} \longrightarrow \scS$$ sending each inf-manifold $U$ to the space of vertical $p$-forms on $U\times M$, i.e.\ the $p$-forms with legs along $M$ only. Algebraically, this space can be identified under the Dold--Kan correspondence as
	$$
		\Omega^p_v(M)(U) \simeq \DK \big( \calO (U\times M) \otimes^h_{C^\infty(M)} \Omega^p(M) \big)\ .
	$$
	Then $\Omega^p_v(M)$ is a stack with deformation theory (but not with quasi-coherent deformation theory). Indeed, to verify the Schlessinger condition along maps $\Spec(A)\longrightarrow \Spec(B)$ such that $\rmH^0(B)\longrightarrow \rmH^0(A)$ is surjective, we use Remark~\ref{rem:def thy on sqz} to reduce to the case where $A=A'\oplus I[1]$, so that $B'=B_\eta$ is a square-zero extension of $B$. We then have to verify that the square of complexes
	is a homotopy pullback, which follows from Lemma~\ref{lem:coprod of square-zero} and the fact that $\Omega^p(M)$ is a projective $C^\infty(M)$-module.
	
	For tangent representability, since $\Omega^p_v(M)$ has the natural structure of an abelian group (in fact an $\RN$-module) and unravelling the definitions, one sees that for any point $\alpha \colon N \longrightarrow \Omega^p_v(M)$ from a smooth manifold $N$, there are canonical equivalences
	$$
	T_\alpha\Omega^p_v(M)\simeq T_0\Omega^p_v(M)\simeq \Omega^p_v(N\times M)= C^\infty(N\times M)\otimes_{C^\infty(M)} \Omega^p(M)\ .
	$$
	Similarly
	\begin{equation}
		\Der_{\Omega^p_v(M)}(\alpha; I)
		\simeq \DK \big(I\otimes_{C^\infty(N)} C^\infty(N\times M)\otimes_{C^\infty(M)} \Omega^p(M)\big) \ ,
	\end{equation}
	which shows that the maps \eqref{eq:qcdt comparison map} are equivalences.
\qen
\end{example}

\begin{example}
	\label{eg: complex of vertical forms as dt stack}
	For every smooth manifold $M$, consider the cochain complex of stacks in cochain degrees $[0, \dim(M)-p]$ given by
	$$
		\begin{tikzcd}
			\Omega^p_v(M) \ar[r, "\dd_v"]
			& \Omega^{p+1}_v(M) \ar[r, "\dd_v"]
			& \cdots \ar[r, "\dd_v"]
			& \Omega^{\dim(M)}_v(M)\ ,
		\end{tikzcd}
	$$
	where $\dd_v$ is the de Rham differential along $M$.
	This associates a cochain complex of simplicial $\RN$-modules to each inf-manifold $U$. Under the cosimplicial Dold--Kan correspondence, this corresponds to a (coskeletal) cosimplicial diagram in $\Sh_\infty(\InfMfd)$, and we define $\Omega^{\geqslant p}_v(M)\colon \InfMfd^{\opp}\longrightarrow \scS$ to be its limit. Explicitly, $\Omega^{\geqslant p}_v(M)$ acts as
	$$
		U \longmapsto \DK \Big( \tau^{\leqslant 0}\, \mathrm{Tot} \big( \calO(U\times M) \otimes^h_{C^\infty(M)} \Omega^{\geqslant p}(M), \dd_v \big) \Big)\ ,
	$$
	where $\Omega^{\geqslant p}(M)$ denotes the complex of differential forms on $M$ in form degrees at least $p$.
	For example, its value on a \emph{smooth} manifold $U$ is the set of vertically closed vertical $p$-forms on $U\times M \longrightarrow U$. 
	
	The stack $\Omega^{\geqslant p}_v(M)$ has deformation theory, since it is a finite limit of dt stacks. Similarly to Example~\ref{ex:vertical forms def thy}, the tangent complex at a map $\alpha\colon N\longrightarrow \Omega^{\geqslant p}_v(M)$ is equivalent to the complex $\Omega^{\geqslant p}_v(N\times M)$ with the vertical de Rham differential.
	\qen
\end{example}

\begin{example}
	\label{eg: forms on M have no DefThy}
	In contrast to Examples~\ref{ex:vertical forms def thy} and~\ref{eg: complex of vertical forms as dt stack}, now consider the functor $\Omega^1 \colon \InfMfd \longrightarrow \scS$ which assigns to each inf-manifold its space of $1$-forms. This functor does \emph{not} have deformation theory. Indeed, consider the pullback square of square-zero extensions of dg $C^\infty$-rings
	\begin{equation}
		\label{eq:square-zero pullback}
		\begin{tikzcd}[row sep = 1cm,column sep=1cm]
			\RN[\epsilon_0]\,\big/\,\epsilon_0^2 \ar[r] \ar[d]
			& \RN \ar[d]
			\\
			\RN \ar[r] 
			& \RN[\epsilon_{-1}]\,\big/\,\epsilon_{-1}^2
		\end{tikzcd}
	\end{equation}
	The bottom right corner is the free dg $C^\infty$-ring on a generator of degree $-1$ (which squares to zero by the Koszul sign rule). Its space of $1$-forms corresponds under the Dold--Kan correspondence to the complex
	\begin{equation}
		\begin{tikzcd}
			\cdots \arrow[r]
			& 0\arrow[r]
			& \RN \, \epsilon_{-1}\, \dd \epsilon_{-1} \arrow[r, "0"]
			& \RN \, \dd \epsilon_{-1} \arrow[r]
			& 0\ ,
		\end{tikzcd}
	\end{equation}
	concentrated in degrees $-1$ and $-2$.
	Since the space of $1$-forms on $\ast=\Spec(\RN)$ is zero, the homotopy pullback of spaces of $1$-forms is given by the two-term complex $\RN \, \epsilon_{-1} \, \dd \epsilon_{-1} \longrightarrow \RN \, \dd\epsilon_{-1}$ in degrees $0$ and $-1$. This is not weakly equivalent to the space of $1$-forms on $\Spec(\RN[\epsilon_0] / \epsilon_0^2)$, whose set of path components is the two-dimensional vector space $(\RN[\epsilon_0] / \epsilon_0^2) \, \dd \epsilon_0$.
	
	An analogous argument for the square
	\begin{equation}
		\begin{tikzcd}[row sep = 1.5cm,column sep=1cm]
			\RN[\epsilon_{1,0}, \ldots, \epsilon_{p,0}]\,\big/\,\{\epsilon_{i,0}\, \epsilon_{j,0}\, | \, i,j = 1, \ldots, p\} \ar[r] \ar[d]
			& \RN \ar[d]
			\\
			\RN \ar[r] 
			& \RN[\epsilon_{1,-1}, \ldots, \epsilon_{p,-1}]\,\big/\,\{\epsilon_{i,-1}\, \epsilon_{j,-1}\, | \, i,j = 1, \ldots, p\}
		\end{tikzcd}
	\end{equation}
	shows that $$\Omega^p \colon \InfMfd \longrightarrow \scS$$ does not have deformation theory for any $p\in\NN$.
	\qen
\end{example}

Finally, the class of qcdt prestacks also includes quotients by formally smooth higher groupoids.

\begin{definition}
	\label{def:formally smooth}
	A map between prestacks $f \colon Y\longrightarrow X$ is \emph{formally smooth} if in any square
		\begin{equation}
			\label{eq:formally_smooth}
			\begin{tikzcd}[column sep = 1.5cm, row sep = 1cm]
				\Spec(A) \arrow[d] \arrow[r]
				& Y \arrow[d,"f"]
				\\
				\Spec(A') \arrow[r] \arrow[ru, dashed]
				& X
			\end{tikzcd}
		\end{equation}
		where $\rmH^0(A')\longrightarrow \rmH^0(A)$ is surjective with nilpotent kernel, there exists a dashed lift as indicated.
\end{definition}

\begin{remark}
	\label{rem:formally smooth}
	By Lemma~\ref{lem:nil ext}, a map between prestacks with deformation theory is formally smooth if and only if it has the right lifting property against all maps $\Spec(A\oplus I[1])\longrightarrow \Spec(A)$. In fact, it suffices to check this when $A=C^\infty(M)$ for a smooth manifold: Remark~\ref{rem:extensions linear} shows that any morphism $\Spec(A\oplus I[1])\longrightarrow \Spec(A)$ can be further decomposed into a sequence of square-zero extensions by modules over $A_{\red}=C^\infty(M)$, i.e.\ into pushouts of maps of the form $\Spec(C^\infty(M)\oplus J[1]) \longrightarrow \Spec(C^\infty(M))$ for a $C^\infty(M)$-module $J$.
	\qen
\end{remark}

We obtain a formal version of the Submersion Theorem given by

\begin{proposition}
	\label{prop:subm theorem}
	Let $f\colon Y\longrightarrow X$ be a map between  qcdt stacks. Then $f$ is formally smooth if and only if for any map from a smooth manifold $y\colon M\longrightarrow Y$, the fibre of the map of tangent complexes $T_yY\longrightarrow T_{f(y)}X$ has vanishing cohomology in positive degrees.
\end{proposition}

\begin{proof}
	By Remark~\ref{rem:formally smooth}, $f$ is formally smooth if and only if there exists a lift
	$$
		\begin{tikzcd}[row sep = 1cm,column sep=1cm]
			\Spec(C^\infty(M) \oplus I[1]) \arrow[r] \arrow[d]
			& Y \arrow[d]
			\\
			M \arrow[r] \arrow[ru, dashed]
			& X
		\end{tikzcd}
	$$
	for every smooth manifold $M$ and complete non-positively graded $C^\infty(M)$-module $I$. In this case, the lifting condition is equivalent to the map 
	$$
		\Der_Y(M, 0) \longrightarrow \Der_Y(M, I[1]) \times_{\Der_X(M,I[1])} \Der_X(M, 0)
	$$
	having non-empty fibres. Unravelling the definitions using Proposition~\ref{prop:qcdt stack 2}, this means that the map
	$$
		\big( T_y Y \otimes_{C^\infty(M)} I[1] \big)^\cpl
		 \longrightarrow \big( T_{f(y)} X \otimes_{C^\infty(M)} I[1] \big)^\cpl
	$$
	has fibre with cohomology concentrated in strictly negative degrees, for every point $y$ and every choice of $I$. Because completion can be computed at the level of cohomology groups (Remark~\ref{rem:complete}), this is equivalent to the fibre of $T_yY\longrightarrow T_{f(y)}X$ having cohomology concentrated in non-positive degrees.
\end{proof}

\begin{example}
	\label{ex:submersions are formally smooth}
	Any submersion between smooth manifolds $\pi\colon M\longrightarrow N$ is formally smooth: this follows directly from Proposition~\ref{prop:subm theorem} and the fact that for any map of smooth manifolds $f\colon V\longrightarrow M$, the map of tangent bundles $f^*TM\longrightarrow f^*\,\pi^*TN$ is surjective.
	\qen
\end{example}

\begin{definition}
	\label{def:form smooth oo-groupoid}
	Let $X_\bullet$ be a simplicial diagram of prestacks with deformation theory. Then $X_\bullet$ is a \emph{formally smooth $\infty$-groupoid} if for each horn inclusion the map 
	\begin{equation}
			X_n \longrightarrow X(\Lambda^n_i) = \underset{\Delta^k\longrightarrow \Lambda^n_i}{\holim}\ X_k
	\end{equation} 
	is formally smooth.
\end{definition}

\begin{proposition}
	\label{prop:form smooth oo-groupoid qcdt}
	Let $X_\bullet$ be a formally smooth $\infty$-groupoid. If each $X_n$ has (quasi-coherent) deformation theory, then the colimit $|X_\bullet|$ has (quasi-coherent) deformation theory as well and
	$$
		T_x\big|X_\bullet\big|\simeq \big|T_xX_\bullet\big|
	$$
for all $x\colon M\longrightarrow X_0\longrightarrow |X_\bullet|$.
\end{proposition}

\begin{proof}
	First, recall that a morphism $K \longrightarrow L$ of simplicial spaces is called a \textit{realisation fibration} if, for every map of simplicial spaces $S \longrightarrow L$, the diagram
	\begin{equation}
		\begin{tikzcd}[row sep = 1cm,column sep=1cm]
			{|S \times_L K|} \ar[r] \ar[d]
			& {|K|} \ar[d]
			\\
			{|S|} \ar[r]
			& {|L|}
		\end{tikzcd}
	\end{equation}
	is homotopy cartesian in $\scS$~\cite{Rezk:hocolims_and_hopullbacks}.
	In other words, the canonical morphism $|S \times_L K| \longrightarrow |S| \times_{|L|} |K|$ is an equivalence of spaces.
	
	For the Schlessinger condition, consider a pushout square \eqref{eq:pushout_sqz} whose image under each $X_n$ is a pullback square. It remains to verify that this pullback is preserved under taking geometric realisations. This follows from the formal smoothness conditions: these imply that the map of simplicial spaces $X_\bullet(\Spec(A')) \longrightarrow X_\bullet (\Spec(A))$ is a Kan fibration of simplicial spaces~\cite[Definition~3.2]{BBP:Classifying_Spaces}, and thus also a realisation fibration~\cite[Theorem~3.17]{BBP:Classifying_Spaces}.
	Tangent representability and the description of the tangent complex follows from $\Der_{|X_\bullet|}\simeq |\Der_{X_\bullet}|$ and $\kappa_!\,\kappa^*\Der_{|X_\bullet|}\simeq |\kappa_!\,\kappa^*\Der_{X_\bullet}|$.
\end{proof}

\begin{example}
	\label{ex:tangent of diff stacks}
	Every (non-derived, higher) differentiable stack $X$ has deformation theory: indeed, it arises as the geometric realisation (computed in stacks) of a Lie $\infty$-groupoid $X_\bullet$, each of whose horn-filling maps $X_n\longrightarrow X(\Lambda^n_i)$ is a submersion and hence formally smooth.
	
	One can use Proposition~\ref{prop:form smooth oo-groupoid qcdt} and Example~\ref{ex:mfd are qcdt} to compute (locally) the tangent complex of higher differentiable stacks. For example, let$\begin{tikzcd}[column sep=0.6cm] \mathcal{G}\arrow[r, yshift=0.5ex] \arrow[r,yshift=-0.5ex]
			& M \end{tikzcd}$be a Lie groupoid and $x\colon M\longrightarrow [M/\mathcal{G}]$ the canonical map to the associated quotient stack. Then $T_x[M/\mathcal{G}]$ is the complex of vector bundles over $M$ corresponding under the Dold--Kan correspondence to the simplicial vector bundle
	$$	
		\begin{tikzcd}
			\cdots \ \arrow[r, yshift=1.5ex] \arrow[r,yshift=0.5ex] \arrow[r,yshift=-0.5ex] \arrow[r, yshift=-1.5ex] & e^*T \big( \mathcal{G} \times_M \mathcal{G}) \arrow[r, yshift=1ex] \arrow[r,yshift=-1ex] \arrow[r]
			& e^*T \mathcal{G} \arrow[r, yshift=0.5ex] \arrow[r,yshift=-0.5ex]
			& TM\ ,
		\end{tikzcd}
	$$
	where $e \colon M\longrightarrow \mathcal{G}\times_M\dots\times_M\mathcal{G}$ denotes the iterated degeneracy maps. Under the Dold--Kan correspondence, this corresponds precisely to the two-term complex $$\frg\longrightarrow TM$$ given in degree $-1$ by the Lie algebroid $\frg$ of $\mathcal{G}$ and with differential given by the anchor map.
	\qen
\end{example}


\section{Lie differentiation and infinitesimal symmetries}
\label{sec:Lie diff}


In this section we discuss how one can use derived $C^\infty$-geometry to associate an $L_\infty$-algebroid to a map $x \colon M \longrightarrow X$ from a smooth manifold to a dt stack.
One can think of the morphism $x \colon M \longrightarrow X$ as classifying a certain type of geometric structure on $M$.
The $L_\infty$-algebroid we produce can then be understood as the $L_\infty$-version of its Atiyah algebroid.
This perspective is the starting point for the theory of connections on higher bundles we develop in the later sections.
Furthermore, the $L_\infty$-algebra of global sections of this $L_\infty$-algebroid can be thought of as the $L_\infty$-algebra of infinitesimal symmetries of the map $x$. We will start by stating the main (formal) results and illustrating how they can be used in examples. The technical details behind the main theorems are discussed in Appendix~\ref{sec:lie diff proof} (Theorem~\ref{thm:formal stack vs lie}) and in Appendix~\ref{sec:global sections proof} (Theorem~\ref{thm:global}).



\subsection{Formal differentiation and integration}
\label{sec: infinitesimal symmetries}


The relevance of the derived geometric formalism we developed Section~\ref{sec: Rev of DefThy} to the problem of differentiating higher stacks rests in the following result.

\begin{definition}
	Let $M$ be a smooth manifold.
	We  write $$\QCDT_{M/} \subseteq \DT_{M/} \subseteq \Sh_\infty(\InfMfd)_{M/}$$ for the full $\infty$-subcategories on those morphisms $M \longrightarrow X$ where $X$ is a dt (resp.~qcdt) stack.
\end{definition}

Given a dg $C^\infty$-ring $A$, let $\Mod^{\cpl, \infty}_A \subset \Mod^\infty_A$ denote the full $\infty$-subcategory on the complete dg $A$-modules, and let \smash{$L_\infty\Agd_A^{\cpl, \infty} \subset L_\infty\Agd_A^\infty$} denote the full $\infty$-subcategory on those $L_\infty$-algebroids whose underlying dg $A$-module is complete.

\begin{theorem}
	\label{thm:formal stack vs lie}
	Let $M$ be a smooth manifold and let $x \colon M \longrightarrow X$ be a map to a dt stack.
	Then the \emph{relative tangent complex}
	$$
	\label{eq: relative tangent complex}
		T(M/X) \coloneqq \hofib (TM \longrightarrow T_x X)
		\
		\in \ \big( \Mod^{\cpl, \infty}_{C^\infty(M)} \big)_{/TM}
	$$
	has a natural lift to an object in the $\infty$-category $L_\infty\Agd_M^{\cpl, \infty}$ of complete $L_\infty$-algebroids over $M$. This determines a right adjoint functor of $\infty$-categories
	$$
		T(M/-) \colon \DT_{M/} \longrightarrow L_\infty\Agd_M^{\cpl, \infty}\ ,
	$$
	whose left adjoint is fully faithful and has essential image given by the full subcategory of $\DT_{M/}$ on those morphisms $M \longrightarrow X$ where $X$ is a qcdt stack such that $M\longrightarrow X_{\red}$ is an equivalence.
	
	Furthermore, for any morphism $M \longrightarrow X$ in this essential image, there is an equivalence
	$$
		\Gamma(X, \mathcal{O}_X) \simeq \ChEil^* \big( T(M/X) \big)
	$$
	between the function cdga of $X$ (Example~\ref{ex:function algebras}) and the Chevalley--Eilenberg cdga of $T(M/X)$.
\end{theorem}

An algebro-geometric version of part of this theorem appears as~\cite[Proposition~V.1.4.2]{GR:Study_Derived_Geom}; combined with~\cite[Theorem~8.2]{Nuiten:Koszul_duality_for_Lie_algebroids}, this identifies the notion of $L_\infty$-algebroid used there with the one from Definition~\ref{def:Loo algd oo-cat}.
When $M=\ast$ is a point, this recovers the principle that for any point $x\colon \ast\longrightarrow X$ of a suitable stack, the shifted tangent complex $T_xX[-1]$ has the structure of an $L_\infty$-algebra (see for example~\cite{Hinich:FormalStacks, Lurie:DAG-X}).
Here the shift arises from the equivalences
\begin{equation}
	T(*/X) \coloneqq \hofib (0 \longrightarrow T_x X) \simeq T_x X[-1]\ .
\end{equation}

When $G$ is a (finite-dimensional) Lie $\infty$-group, an explicit description of the Lie differentiation of the map $\ast\longrightarrow \rmB G$ is given in \cite{Rogers:Differentiation}. When $M$ is a manifold, its seems more difficult to describe the Lie differentiation of a map $M\longrightarrow X$ explicitly, because such maps are typically presented as maps out of a hypercover of $M$. In the special case where $M\longrightarrow X$ arises from a map $M\longrightarrow X_\bullet$ of Lie $\infty$-groupoids, without taking a hypercover of $M$, then the Lie differentiation given by Theorem \ref{thm:formal stack vs lie} probably coincides with the one proposed by \v{S}evera \cite{Severa:Differentiation,LRWZ:Differentiation}.

\begin{example}
	\label{ex:tangent}
	The de Rham stack $M_{\dR}$ is (by definition) the terminal stack equipped with a map $M\longrightarrow M_{\dR}$ inducing an equivalence on reductions. In addition, $M_{\dR}$ has strictly quasi-coherent deformation theory by Example~\ref{ex:de rham}. It follows that $M_{\dR}$ corresponds under the equivalence of Theorem~\ref{thm:formal stack vs lie} to the terminal Lie algebroid $TM$; that is, there is a canonical equivalence $$T(M/M_\dR) \simeq TM \ . $$
	
	Similarly, the initial Lie algebroid $0$ corresponds to $M$ itself, i.e.~it is the relative tangent complex of the identity map on $M$.
	\qen
\end{example}

\begin{corollary}
	\label{cor:loo algebra}
	For any smooth manifold $M$ and  retract diagram 
\begin{equation}\label{eq:retractLinfty}	
	\begin{tikzcd}[column sep=1cm] M\arrow[r, yshift=0.8ex] 
			& \arrow[l,yshift=-0.8ex] X \end{tikzcd}
			\end{equation}
 where $X$ is a dt stack, the relative tangent bundle $T(M/X)$ has the natural structure of a complete $C^\infty(M)$-linear $L_\infty$-algebra. In fact, there is an equivalence between complete $C^\infty(M)$-linear $L_\infty$-algebras and retract diagrams \eqref{eq:retractLinfty} where $X$ is a qcdt stack such that $M \longrightarrow X_{\red}$ is an equivalence.
\end{corollary}

\begin{proof}
Theorem~\ref{thm:formal stack vs lie} and Example~\ref{ex:tangent} show that $T(M/X)$ defines a complete $L_\infty$-algebroid over the zero Lie algebroid, that is, a complete $C^\infty(M)$-linear $L_\infty$-algebra.
\end{proof}

The second part of Theorem~\ref{thm:formal stack vs lie} quantifies how much of the structure of $X$ can be recovered from the $L_\infty$-algebroid $T(M/X)$. It essentially asserts that $T(M/X)$  encodes the infinitesimal structure of $X$ supported at $M$, or more precisely its formal completion.

\begin{definition}
	\label{def: de Rham stack}
	Let $M\longrightarrow X$ be any map in $\mathrm{Sh}_\infty(\InfMfd)$. The \emph{formal completion} of $X$ along $M$ is the fibre product
	\begin{equation}
		X^\wedge_M \coloneqq X\times_{X_\dR} M_\dR\ .
	\end{equation}
\end{definition}

Evaluating $X^\wedge_M$ at an inf-manifold $U$, one sees that
$$
	X^\wedge_M(U) = X(U)\times_{X(U_{\red})} M(U_{\red})
$$
describes the infinitesimal structure of $X$ around $M$. In particular $X^\wedge_M\simeq X$ as soon as $M_{\red}\simeq X_{\red}$.

\begin{corollary}\label{cor:lie diff and completion}
	Let $x\colon M\longrightarrow X$ be a map from a smooth manifold to a qcdt stack. Then the counit map of the adjunction from Theorem \ref{thm:formal stack vs lie} is equivalent to the map $X^\wedge_M\longrightarrow X$. 
\end{corollary}

In particular, the data of the $L_\infty$-algebroid $T(M/X)$ completely recovers precisely the formal completion $X^\wedge_M$. When $X$ is not a qcdt stack, the counit map is much less explicit.

\begin{proof}
	Since $X$ has quasi-coherent deformation theory, $X^\wedge_M$ is a stack with quasi-coherent deformation theory whose reduction is equivalent to $M$. Since $T(M/X^\wedge_M)\longrightarrow T(M/X)$ is an equivalence as a direct consequence of~\eqref{eq:derivations} and Definition~\ref{def: de Rham stack}, the result follows from Theorem~\ref{thm:formal stack vs lie}.
\end{proof}

\begin{remark}
	As an application of Corollary~\ref{cor:lie diff and completion}, suppose that $\frg$ is a complete $L_\infty$-algebroid over $M$ and that $f\colon N\longrightarrow M$ is a map of smooth manifolds. Then $f^*\frg\times^h_{f^*TM} TN$ can be endowed with a natural $L_\infty$-algebroid structure: if $M\longrightarrow X$ is the qcdt stack associated to $\frg$ by Theorem~\ref{thm:formal stack vs lie}, then $f^*\frg\times_{f^*TM} TN$ is the underlying complex of the $L_\infty$-algebroid of $X^\wedge_N$.
	\qen
\end{remark}

Every $L_\infty$-algebroid $\frg$ over a smooth manifold $M$ gives rise to an $\RN$-linear $L_\infty$-algebra $\Gamma(M, \frg)$ by taking global sections; in terms of the Serre--Swan perspective we are employing, this amounts to forgetting the $C^\infty(M)$-linear structure on $\frg$ while remembering its $\RN$-linear structure and brackets. For example, when $\frg = TM$ is the tangent bundle, we recover the Lie algebra of vector fields on $M$.
At the same time, this $\RN$-linear Lie algebra should also arise in a different way: if we consider the diffeomorphism group $\cDiff(M)$ as an object in $\Sh_\infty(\InfMfd)$, then the Lie algebra of vector fields should be its Lie differentiation. One can make this more precise through

\begin{definition}\label{def:diff}
	Let $$\cDiff(M)\colon \InfMfd^{\opp}\longrightarrow \scS$$ be the stack sending each inf-manifold $U$ to the space of equivalences relative to $U$, i.e.~the space of commutative triangles in $\InfMfd$ of the form
	$$
		\begin{tikzcd}[row sep = 1cm,column sep=1cm]
			U \times M \arrow[rr, "\sim"', "{(\pi_U, \phi)}"] \arrow[rd, "\pi_U"{swap}]
			& & U \times M \arrow[ld, "\pi_U"]
			\\
			& U
		\end{tikzcd}
	$$
\end{definition}

\begin{lemma}
	\label{lem:diff def thy}
	The stack $\cDiff(M)$ has deformation theory and is formally smooth, i.e.\ $\cDiff(M)\longrightarrow \ast$ is formally smooth in the sense of Definition~\ref{def:formally smooth}.
\end{lemma}

\begin{proof}
By Lemma~\ref{lem:nil ext}(2), a map $(\pi_U, \phi)\colon U\times M\longrightarrow U\times M$ is an equivalence (over $U$) if and only if its pullback $U_{\red}\times M\longrightarrow U_{\red}\times M$ is a diffeomorphism of smooth manifolds. This implies that the square
	\begin{equation} \label{eq:pullbackDiff}
		\begin{tikzcd}[row sep = 1cm,column sep=1cm]
			\cDiff(M) \ar[r] \ar[d]
			& \Hom(M, M) \ar[d]
			\\
			\cDiff(M)_\dR \ar[r]
			& \Hom(M, M)_\dR
		\end{tikzcd}
	\end{equation}
	is a pullback. Examples~\ref{ex:mfd are qcdt} and \ref{ex:de rham}  together with Proposition \ref{prop:def thy} then imply that $\cDiff(M)$ is a pullback of qcdt stacks, and hence a qcdt stack itself. 
	
	It follows from Proposition~\ref{prop:subm theorem} that $\cDiff(M)$ is formally smooth if for each map $\phi^{\dashv}\colon U\longrightarrow \cDiff(M)$ from a smooth manifold $U$, the tangent complex $T_{\phi^{\dashv}}\cDiff(M)$ has vanishing positive cohomology. Using the pullback square \eqref{eq:pullbackDiff} and Proposition~\ref{prop:def thy}(3), one then deduces
	\[
	T_{\phi^{\dashv}}\cDiff(M)\simeq T_{\phi^{\dashv}}\Hom(M, M)\simeq \pi_{U*}(T_{\phi}M) \ .
	\]
	 In other words, this is the complete $C^\infty(U)$-module underlying the $C^\infty(U\times M)$-module $\phi^*TM$, where $\phi\colon U\times M\longrightarrow M$ is adjoint to $\phi^{\dashv}$. This module is concentrated in degree $0$, so the result follows.
\end{proof}

Composition of equivalences endows $\cDiff(M)$ with the structure of an $\infty$-group, and we will write $\rmB\cDiff(M)$ for the corresponding classifying stack. Similarly, each mapping stack $\Hom(M, X)$ carries a natural right action $\varrho$ of $\cDiff(M)$ by pre-composition, encoded by a simplicial diagram
$$
	\begin{tikzcd}
		\Hom(M, X)
		&\Hom(M, X) \times \cDiff(M) \arrow[l, yshift=-0.5ex, "\varrho"] \arrow[l, yshift=0.5ex, "\pi_1"{swap}]
		& \cdots  \arrow[l]\arrow[l, yshift=1ex] \arrow[l, yshift=-1ex]
	\end{tikzcd}
$$
where $\pi_1$ is the projection to the first factor.
Lemma~\ref{lem:diff def thy} and Proposition~\ref{prop:def thy} imply that for any dt stack $X$ and any manifold $M$, this is a formally smooth $\infty$-groupoid. Proposition~\ref{prop:form smooth oo-groupoid qcdt} then implies that the corresponding quotient stack $\Hom(M, X)/\cDiff(M)$~\cite[Definition~3.1]{NSS:oo-bundles_I} is a dt stack as well. In particular, taking $X=\ast$ shows that $\rmB\cDiff(M)$ is a dt stack itself. We give a more precise description of the stacks $\rmB\cDiff(M)$ and $\Hom(M, X)/\cDiff(M)$ in Appendix~\ref{sec:global sections proof}.

Every morphism $x \colon M \longrightarrow X$ induces a basepoint $$x^\dashv \colon \ast \longrightarrow \Hom(M, X)\,\big/\,\cDiff(M) \ . $$ We may then give a stacky interpretation of the functor taking the $L_\infty$-algebra of global sections of an $L_\infty$-algebroid.
In the following, if $\scC$ is an $\infty$-category and $x, z \in \scC$, we write $\scC_{x \dslash z}$ for the $\infty$-category of objects under $x$ and over $z$, i.e.~diagrams $x \longrightarrow y \longrightarrow z$ in $\scC$. We also write $L_\infty\Alg^\infty_\RN$ for the $\infty$-category of $\RN$-linear $L_\infty$-algebras.

\begin{theorem}
	\label{thm:global}
	There is a commutative diagram of right adjoint functors
	\begin{equation}
		\label{eq: global symmetries}
		\begin{tikzcd}[column sep = 2cm, row sep = 1.25cm]
			\DT_{M/} \arrow[r, "p_*\,t^*"] \arrow[d, "T(M/-)"{swap}]
			& \DT_{\ast \dslash \rmB \cDiff(M)} \arrow[d, "T(\ast/-)"]
			\\
			L_\infty\Agd_M^{\cpl, \infty} \arrow[r, "\fgt"']
			& \big(L_\infty\Alg^\infty_\RN \big)_{/\Gamma(M,TM)} 
		\end{tikzcd}
	\end{equation}
	where the vertical maps arise from Theorem~\ref{thm:formal stack vs lie}, the bottom functor forgets the $C^\infty(M)$-linear structure, and the top functor $p_*\,t^*$ sends $x \colon M \longrightarrow X$ to the formal completion of $\Hom(M, X)/\cDiff(M)$ at the basepoint $x^\dashv$.
\end{theorem}

\begin{example}
	The Lie differentiation of $\ast\longrightarrow \rmB\cDiff(M)$ in the sense of Theorem~\ref{thm:formal stack vs lie} recovers the Lie algebra of vector fields on $M$, as one would expect: this follows by comparing the two images of the terminal object of $\DT_{M/}$ in \eqref{eq: global symmetries} using Example~\ref{ex:tangent}.
	\qen
\end{example}

\begin{remark}
	Theorem~\ref{thm:global} provides the following unified derived perspective on symmetries of geometric structures on a manifold $M$:
	first, we may consider the global symmetries of $x \colon M \longrightarrow X$.
	These are not local on $M$, as $U \longmapsto \cDiff(U)$ does not define a sheaf (nor even a presheaf) on $M$, and so the global symmetries are encoded in a group object in $\Sh_\infty(\InfMfd)$, or equivalently its classifying object.
	Differentiating, we obtain the $L_\infty$-algebra of the symmetry group.
	Theorem~\ref{thm:global} states, in particular, that the \textit{infinitesimal} symmetries \textit{are} local, i.e.~that the $L_\infty$-algebra of the higher derived symmetry group has an enhancement to an $L_\infty$-algebroid on $M$, and that the functor $T(M/-)$ computes precisely this local enhancement.
	
	Geometrically, the left-hand vertical functor in~\eqref{eq: global symmetries} treats the deformation theory of the map $x \colon M \longrightarrow X$ as a family of deformation problems in $X$, parametrised by the points of $M$.
	It keeps track of deformations which merely arise from permuting the indexing set by derived diffeomorphisms of $M$.
	In contrast, the right-hand vertical functor treats the deformation theory of the map $x \colon M \longrightarrow X$ as the deformations of a single point $x^\dashv \colon * \longrightarrow \Hom(M,X)/\cDiff(M)$ in the moduli stack of geometric structures on $M$ classified by $X$, modulo their equivalences and diffeomorphisms of $M$.
	Theorem~\ref{thm:global} implies 
	\begin{equation}
		T \big( *\,\big/\,(\Hom(M,X)\,/\,\cDiff(M)) \big)
		\simeq \Gamma \big( M, T(M/X) \big)
	\end{equation}
	and makes precise the relationship between these two deformation problems.
	\qen
\end{remark}

The first part of Theorem~\ref{thm:formal stack vs lie} allows us, in particular, to associate  an $L_\infty$-algebroid on $M$ which controls the deformation theory of any morphism $x \colon M \longrightarrow X$ from a manifold to a dt stack.
Viewing $X$ as the classifying stack of a certain type of geometric structure on $M$, Theorem~\ref{thm:global} suggests that the $L_\infty$-algebroid $T(M/X)$ should then be understood as the generalisation of the Atiyah algebroid of a classical principal bundle.
This motivates

\begin{definition}
	\label{def:Atiyah L_oo-algebroid from defthy}
	Given a morphism $x \colon M \longrightarrow X$ from a smooth manifold to a dt stack $X$, the relative tangent complex $T(M/X)$ is its \textit{Atiyah $L_\infty$-algebroid}.
\end{definition}

A key insight and justification of Definition~\ref{def:Atiyah L_oo-algebroid from defthy} is, as we will show in the later sections of this paper, that it provides access to what we believe is the most general definition of higher-form connections on higher geometric structures, such as principal $\infty$-bundles.


\subsection{Properties and examples}
\label{sec:differentiation examples}


Let us now record some properties and examples of the Lie differentiation procedure provided by Theorem~\ref{thm:formal stack vs lie}, starting with

\begin{proposition}
	The Lie differentiation from Theorem~\ref{thm:formal stack vs lie} satisfies the following properties:
	\begin{myenumerate}
		\item Lie differentiation preserves homotopy limits: for any diagram $M \longrightarrow X_\alpha$ of dt stacks under $M$,
		\begin{equation}
			T\Big( M\,\big/\,\underset{\alpha}\holim\ X_\alpha\Big) \simeq \underset{\alpha}\holim\ T( M/X_\alpha) \ .
		\end{equation} 
		
		\item Lie diffentiation preserves filtered homotopy colimits.
		
		\item Given a formally smooth $\infty$-groupoid $X_\bullet \colon \bbDelta^{\opp} \longrightarrow \Sh_\infty(\InfMfd)$ and a map $M\longrightarrow X_0$, there is an equivalence of $L_\infty$-algebroids
		\begin{equation}
			T \big( M/|X_\bullet| \big) \simeq \big| T(M/X_\bullet) \big|\ .
		\end{equation}
	\end{myenumerate}
\end{proposition}

\begin{proof}
	Since the forgetful functor from $L_\infty$-algebroids to complexes respects homotopy limits and sifted\footnote{A colimit is \emph{sifted} if it commutes with finite coproducts.} homotopy colimits, the properties readily follow from the fact that taking tangent complexes (Definition~\ref{def:tangent}) preserves limits and filtered colimits (by~\eqref{eq:derivations} and Theorem~\ref{eq: relative tangent complex}) as well as geometric realisations of formally smooth $\infty$-groupoids (see Proposition~\ref{prop:form smooth oo-groupoid qcdt}).
\end{proof}

\begin{example}
	\label{ex:loop lie algebra is trivial}
	Let $M$ be a smooth manifold and consider a retract diagram
	$$\begin{tikzcd}
	M\arrow[r, "x"] & X\arrow[r, "i", yshift=1ex] & Y\arrow[l, "r", yshift=-1ex]
	\end{tikzcd}$$
 of dt stacks under $M$. This gives rise to a retract diagram of $L_\infty$-algebroids $$\begin{tikzcd}[column sep=1cm]\frg\arrow[r,yshift=0.8ex] &  \mathfrak{h}=\frg\oplus \mathfrak{n}\arrow[l,yshift=-0.8ex] \end{tikzcd}$$ on $M$, where $\mathfrak{n}$ denotes the kernel of $\mathfrak{h}\longrightarrow\frg$. Then \smash{$T(M/X\times^h_Y X)$} can be identified as the complex $$\frg \oplus \mathfrak{n}[-1]\longrightarrow \frg \longrightarrow TM \ , $$ with the $L_\infty$-structure
	$$
		\big[ (x_1, y_1)\,,\, \dots\,,\, (x_n, y_n) \big] = \Big( [x_1, \dots, x_n]\,,\, \sum\limits_{k=1}^n\, [x_1, \dots, y_k, \dots, x_n] \Big)\ .
	$$
	
	Indeed, we can compute $T(M/X\times^h_Y X) \simeq \frg \times^h_{\mathfrak{h}} \frg$ as a homotopy pullback of $L_\infty$-algebroids. To do this, we  factor the morphism $\frg \longrightarrow \frh$ as the composition $\frg \longrightarrow \tilde{\frg}\longrightarrow \mathfrak{h}$ of a weak equivalence followed by a fibration, as follows. Let $$\tilde{\frg} = \frh \oplus \mathfrak{n}[-1]$$ with differential $$\dd(z, y) = (\dd z, \pi(z) + \dd y) \ , $$ where $\pi \colon \mathfrak{h} = \frg \oplus \mathfrak{n} \longrightarrow \mathfrak{n}$ is the projection. This comes equipped with an $L_\infty$-algebroid structure where 
	$$
		\big[ (z_1, y_1)\,,\, \dots\,,\, (z_n, y_n) \big] = \Big([z_1, \dots, z_n]\,,\, \sum\limits_{k=1}^n\, [z_1, \dots, y_k, \dots, z_n] \Big)\ .
	$$
	One readily checks that the inclusion of $\frg$ and the projection to $\mathfrak{h}$ preserve the brackets. The pullback $\frg\times_{\mathfrak{h}} \tilde{\frg}$ then coincides with the $L_\infty$-algebroid given above.
	
	If $X = M$ and $x = \id_M$, we can think of $$\begin{tikzcd}[column sep=1cm]M \arrow[r,yshift=0.8ex]& Y \arrow[l,yshift=-0.8ex] \end{tikzcd}$$ as a bundle over $M$ with a specified global section.
	In that case, $M \times^h_{Y} M$ produces the stack of fibrewise loops in $Y$ based at the global section.
	Then $\frg = 0$ is the trivial $L_\infty$-algebroid (by Example~\ref{ex:tangent}), and so $T(M/M \times^h_Y M) \simeq \frg \times^h_\frh \frg$ is equivalent to $\mathfrak{n}[-1]$ with trivial brackets.
	\qen
\end{example}

\begin{example}
	\label{ex:lie diff group}
	Suppose that $G$ is a dt stack that admits the structure of an $\infty$-group. Let $e\colon \ast\longrightarrow G$ be the unit. Then the Atiyah $L_\infty$-algebroid $$T(\ast/G) \simeq T_e G[-1]$$ is a \emph{trivial} $L_\infty$-algebra~\cite[Chapter~6, Proposition~1.7.2]{GR:Study_Derived_Geom} (where $T(*/G)$ is the relative tangent complex of the canonical morphism $e \colon * \longrightarrow G$). Indeed, we can identify $T(\ast/G)\simeq T\big(\ast/(\ast \times_{\rmB G} \ast)\big)$ with the homotopy pullback $0\times^h_{T(\ast/\rmB G)} 0$ of $L_\infty$-algebras. Example~\ref{ex:loop lie algebra is trivial} shows that this homotopy pullback is trivial.
\qen
\end{example}

\begin{example}
	\label{ex: Lie diff group under M}
	The same analysis as in Example~\ref{ex:lie diff group} applies when $M$ is a manifold and $e\colon M\longrightarrow G$ can be identified with the unit map for some $\infty$-group in the $\infty$-category $\Sh_\infty(\InfMfd)_{/M}$ of stacks over $M$. In this case, the $C^\infty(M)$-linear $L_\infty$-algebra $T(M/G)$ (Corollary~\ref{cor:loo algebra}) is the trivial $L_\infty$-algebra.
	
	For instance, this shows that for any closed embedding of smooth manifolds $M\longhookrightarrow  N$, the shifted normal bundle $$T(M/N)\simeq (T_MN/TM)[-1]$$ is a trivial $L_\infty$-algebroid. Indeed, $T(M/N)$ only depends on an open neighbourhood of $M$ inside $N$, and so using a tubular neighbourhood we can reduce to the case where $M\longhookrightarrow  N$ is the zero section of a vector bundle on $M$ (which has a fibrewise group structure over $M$).
	\qen
\end{example}

\begin{remark}
	\label{rem:ab group}
	Example~\ref{ex:lie diff group} shows that for a dt stack $G$ that admits some group structure, the \emph{shifted} tangent complex $T_eG[-1]$ is a trivial $L_\infty$-algebra. The usual tangent complex $T_eG$ normally has a \emph{non-trivial} $L_\infty$-structure that depends on the group structure on $G$: this $L_\infty$-structure arises by differentiating the map $\ast \longrightarrow \rmB G$ to the classifying stack of $G$. It reproduces, for example, the usual Lie bracket on the tangent space $T_eG$ of a Lie group $G$ (by Example~\ref{ex:classical lie} below). 
	
	Conversely, this gives the following informal interpretation of the Lie differentiation of $x \colon \ast \longrightarrow X$: it can be identified with the $L_\infty$-algebra $$T(\ast/X) \simeq T_x(\sfOmega_xX)$$ of the based loop space $\sfOmega_x X$ endowed with the group structure given by composition of loops. If $X$ has itself the structure of a group, then $\sfOmega_x X$ has an $\mathbbm{E}_2$-structure, which trivialises all brackets.
	\qen
\end{remark}

\begin{remark}
	\label{rem:loop space relative diffeo}
	In light of Remark~\ref{rem:ab group} we revisit the interpretation of Theorem~\ref{thm:global} in terms of infinitesimal symmetries. Suppose that $x\colon M\longrightarrow X$ is a map to a dt stack. The $L_\infty$-algebra underlying $T(M/X)$ can be identified with the Lie differentiation of the basepoint inclusion $x^\dashv:\ast \longrightarrow \Hom(M, X)/\cDiff(M)$. The loop space of $\Hom(M, X)/\cDiff(M)$ is precisely the group stack
	\begin{equation}
		\cSym(x) = \cDiff_{/X}(M, M)
	\end{equation}
	of derived diffeomorphisms of $M$ relative to $X$.
	(In the present setting this is a derived enhancement of the smooth higher symmetry groups investigated in~\cite{BMS:Sm2Grp, Bunk:Pr_oo-Bundles, BS:Higher_Syms_and_Deligne_Coho}.)
	By the above arguments one can think of the $L_\infty$-algebra underlying $T(M/X)$ as the tangent Lie algebra of the group stack $\cDiff_{/X}(M, M)$. 
	
	In particular, if $X = \rmB G$ for a smooth group stack $G$, so that $x \colon M \longrightarrow X$ classifies a $G$-bundle $P \longrightarrow M$, then a point in $\cDiff_{/X}(M, M)$ corresponds to a diffeomorphism $\phi \colon M \longrightarrow M$ together with an equivalence of $G$-bundles $P\simeq \phi^*P$.
	\qen
\end{remark}

Let us now describe a way to identify the Lie bracket at the level of homotopy groups. To this end, suppose that $x \colon \ast \longrightarrow X$ is a point of a dt stack and let us write $\sfOmega_x^m X$ for the $m$-fold loop stack at the basepoint $x$. For any $a\leqslant 0$, let us write $$D^1_a= \Spec\big(\RN[\epsilon_a]\,\big/\,\epsilon_a^2\big)$$ for the infinitesimal line, where $\epsilon_a$ has degree $a$. There is a canonical basepoint $\ast\longrightarrow D^1_a$, and the tangent complex classifies maps out of such infinitesimal lines: we can identify
$$
	\rmH^{-a-m}\big(T(\ast/X)\big) \cong \big[ (D^1_a, \ast)\,,\, (\sfOmega^{m+1}_xX, \ast) \big]
$$
with the set of homotopy classes of pointed maps $D^1_a \longrightarrow \sfOmega_x^{m+1}X$ (since this set of homotopy classes agrees with $\pi_0 \Der_{\sfOmega_x^{m+1}X}(x; \RN[a])$ and Definition~\ref{def:tangent}).

Given two infinitesimal lines we will write
$$
	D^1_a \vee D^1_b = \Spec\big(\RN[\epsilon_a, \epsilon_b]\,\big/\,\epsilon_a^2, \epsilon_b^2, \epsilon_a\,\epsilon_b\big)
	\qquad \text{and} \qquad
	D^1_a \times D^1_b = \Spec\big(\RN[\epsilon_a, \epsilon_b]\,\big/\,\epsilon_a^2, \epsilon_b^2\big)
$$
for the wedge sum and the product, respectively. 
Given two pointed maps $\alpha \colon D^1_a \longrightarrow \sfOmega^{m+1}_xX$ and $\beta \colon D^1_b \longrightarrow \sfOmega^{n+1}_xX$, let us now consider the  diagram
\begin{equation} \label{eq:Whproduct}
	\begin{tikzcd}[column sep = 1.25cm, row sep = 1cm]
		D^1_a \vee D^1_b \arrow[r] \arrow[d]
		& D^1_a \times D^1_b \arrow[r, "\alpha \times \beta"] \arrow[d]
		& \sfOmega^{m+1}_x X \times \sfOmega_x^{n+1} X \arrow[d, "\mathrm{Wh}"]
		\\ 
		\ast \arrow[r]
		& D^1_{a+b} \arrow[r, dashed, "{[\alpha, \beta]_{\mathrm{Wh}}}"]
		& \sfOmega_x^{m+n+1}X
	\end{tikzcd}
\end{equation}
Here $\mathrm{Wh}$ denotes the Whitehead product, induced by the attaching map $S^{m+n+1}\longrightarrow S^{m+1}\vee S^{n+1}$ for the $(m{+}n{+}2)$-cell of $S^{m+1} \times S^{n+1}$. The composition $D^1_a \vee D^1_b \longrightarrow \sfOmega_x^{m+n+1}X$ is equivalently described by an $(m{+}n{+}1)$-fold loop in the space of pointed maps $D^1_a \vee D^1_b\longrightarrow X$ given as the Whitehead product of an $(m{+}1)$-fold loop and an $(n{+}1)$-fold loop.

However, this mapping space is an infinite loop space, so all Whitehead products vanish. By the Schlessinger condition, we therefore obtain a map $[\alpha, \beta]_{\mathrm{Wh}}\colon D^1_{a+b}\longrightarrow \sfOmega_x^{m+n+1}X$ out of the pushout.

\begin{proposition}
	\label{prop:bracket on homotopy}
	For each $a, b \leqslant 0$, the map $[\alpha, \beta]_{\mathrm{Wh}}\colon D^1_{a+b}\longrightarrow \sfOmega_x^{m+n+1}X$ represents the image of $(\alpha, \beta)$ under the Lie bracket on cohomology:
	\begin{equation}
		[-, -] \colon \rmH^{-a-m} \big( T(\ast/X) \big)\, \times\, \rmH^{-b-n} \big( T(\ast/X) \big)
		 \longrightarrow \rmH^{-a-b-m-n} \big( T(\ast/X) \big)\ .
	\end{equation}
	That is, the formula described by \eqref{eq:Whproduct} computes the graded Lie algebra structure on all cohomology groups of $T(*/X)$.
\end{proposition}

\begin{proof}
	Let us write $\frg$ for a dg Lie algebra modelling $T(\ast/X)$. It follows from (the proof of) Theorem~\ref{thm:formal stack vs lie} and~\cite[Example~7.20]{CCN:Moduli_Operads} that for any non-positively graded dg Artin algebra $A$ with maximal ideal $\mathfrak{m}_A$, there is a natural equivalence of spaces
	$$
		\DT\big(\Spec(A), X\big)\times_{\DT(\ast, X)} \{x\}\simeq \mathrm{MC}_\bullet\big(\frg\otimes \mathfrak{m}_A)\ .
	$$
	Here the right hand side denotes the Maurer--Cartan space from~\cite{Hinich:FormalStacks}: it is obtained by tensoring with the simplicial cdga $\Omega^*[\Delta^\bullet]$ of polynomial differential forms on simplices and taking Maurer--Cartan elements.
	
Using this, the construction of \eqref{eq:Whproduct} boils down to the following: we can identify $\alpha$ and $\beta$ with maps
	$$
		\begin{tikzcd} [row sep = 1.2cm]
			& \mathrm{MC}_\bullet \big( \frg\otimes (\epsilon_a, \epsilon_b, \epsilon_a\,\epsilon_b) \big) \arrow[dr, hookleftarrow] 
			& \\
			\mathrm{MC}_\bullet\big(\frg \otimes (\epsilon_a)\big) \arrow[ur, hookrightarrow] & & \mathrm{MC}_\bullet\big(\frg\otimes (\epsilon_b)\big)  \\
			S^{m+1 }\arrow[u, "\alpha = \alpha' \cdot \epsilon_a"] & &   [1.4pc] S^{n+1} \arrow[u, "\beta = \beta'\cdot\epsilon_b",swap]
		\end{tikzcd}
	$$
	where $(\epsilon_a)$ is the maximal ideal of $\calO(D^1_a)$, while $\alpha'$ and $\beta'$ correspond to cocycles in $\frg$ via the Dupont contraction~\cite[Section~11.2]{BFMT:Lie_Models}. Next we take the Whitehead product of the resulting pair of elements in the homotopy groups of the space $\mathrm{MC}_\bullet\big(\frg\otimes (\epsilon_a, \epsilon_b, \epsilon_a\,\epsilon_b)\big)$, or more precisely of its connected component containing the basepoint. This connected component is a nilpotent rational space, whose Lie model is the nilpotent dg Lie algebra $\tau^{\leqslant 0}\big(\frg\otimes(\epsilon_a, \epsilon_b, \epsilon_a\,\epsilon_b)\big)$. 
	
It now follows from~\cite[Theorems~12.41 and~12.42]{BFMT:Lie_Models} that the Whitehead product of $\alpha$ and $\beta$ is given by the Lie bracket in $\frg\otimes (\epsilon_a, \epsilon_b, \epsilon_a\,\epsilon_b)$. In other words, the Whitehead product is given up to homotopy by
	$$
		[\alpha', \beta']_{\frg} \cdot \epsilon_a\, \epsilon_b \colon S^{m+n+1}
		 \longrightarrow \mathrm{MC}_\bullet \big( \frg \otimes (\epsilon_a, \epsilon_b, \epsilon_a\, \epsilon_b) \big)\ ,
	$$
	where $[-, -]_\frg$ is the Lie bracket of $\frg$. This class is indeed null modulo $\epsilon_a\,\epsilon_b$, and we find that the class $[\alpha, \beta]_{\mathrm{Wh}}$ can be represented up to homotopy by
	$$
		[\alpha, \beta]_{\mathrm{Wh}} = [\alpha', \beta']_{\frg} \cdot \epsilon_a\, \epsilon_b \colon S^{m+n+1}
		 \longrightarrow \mathrm{MC}_\bullet \big( \frg\otimes (\epsilon_a\,\epsilon_b) \big)\ .
	$$
	It follows that $[\alpha, \beta]_{\mathrm{Wh}}$ indeed presents the Lie bracket on $\rmH^*(\frg)$.
\end{proof}

Proposition~\ref{prop:bracket on homotopy} is particularly useful for determining the Lie bracket on the non-positive cohomology groups of $T(\ast/X)$: in this case we can take $\epsilon_a$ and $\epsilon_b$ to be of degree zero and work in a (non-derived) setting of synthetic differential geometry.

\begin{example}
	\label{ex:synthetic lie diff}
	Suppose that $G\in \Sh_\infty(\InfMfd)$ is a formally smooth dt stack equipped the structure of an $\infty$-group. By Proposition~\ref{prop:form smooth oo-groupoid qcdt} the classifying stack $\rmB G$ has deformation theory, and the tangent complex $T(\ast/\rmB G) \simeq T_e(\sfOmega \rmB G) = T_eG$ is concentrated in non-positive degrees. The Lie bracket $\rmH^m(T_eG) \times \rmH^n(T_eG)\longrightarrow \rmH^{m+n}(T_eG)$ on cohomology can be obtained entirely using the Whitehead product of maps $\alpha\colon D^1\longrightarrow \sfOmega_e^mG$ and $\beta\colon D^1\longrightarrow \sfOmega_e^nG$.
	
	This is well-known for $m=n=0$: in this case there are two (non-derived) infinitesimal paths $\alpha, \beta\colon D^1 \longrightarrow G$ at the identity and their Lie bracket is obtained as the unique factorisation
	$$
		\begin{tikzcd}[column sep = 1cm, row sep = 1cm]
			D^1\times D^1 \arrow[r, "\alpha\times \beta"] \arrow[d, "{(\epsilon_1, \epsilon_2) \longmapsto \epsilon_1\,\epsilon_2}"{swap}]
			& G\times G \arrow[d, "{(g, h)\longmapsto g\,h\,g^{-1}\,h^{-1}}"]\arrow[d]
			\\
			D^1\arrow[r, dashed, "{[\alpha, \beta]}"]
			& G
		\end{tikzcd}
	$$
 so that the Whitehead product is the commutator~\cite[Theorem~V.1.6]{MR:Models_Smooth_Inf_Analysis}.
	\qen
\end{example}

\begin{example}
	\label{ex:classical lie}
	Let$\begin{tikzcd}[column sep=0.6cm] \mathcal{G}\arrow[r, yshift=0.5ex] \arrow[r,yshift=-0.5ex]
			& M \end{tikzcd}$be a Lie groupoid and let $x\colon M\longrightarrow X=[M/\mathcal{G}]$ be the canonical map to the associated stack. Example~\ref{ex:tangent of diff stacks} identifies the morphism $TM\longrightarrow T_xX$ with the evident inclusion $TM \longrightarrow (\frg\longrightarrow TM)$, where $\frg$ is the Lie algebroid of $\mathcal{G}$ and the differential is its anchor map. Consequently the map $T(M/X) \longrightarrow TM$ coincides with the anchor map $\rho \colon \frg\longrightarrow TM$. 
	
	For degree reasons, an $L_\infty$-algebroid structure on $\frg$ is just the datum of a Lie bracket. To show that this recovers the usual Lie algebroid structure on $\frg$, we  use Theorem~\ref{thm:global} and consider the Lie differentiation of $x^\dashv\colon \ast\longrightarrow\Hom(M, X)/\cDiff(M)$. Since we already know that the underlying module is $\frg$ (concentrated in degree $0$), we can use Proposition~\ref{prop:bracket on homotopy} to identify the Lie bracket exactly as in Example~\ref{ex:synthetic lie diff}.
	
	By Remark~\ref{rem:loop space relative diffeo}, the loop stack $\cDiff_{/X}(M, M) = \sfOmega_{x^\dashv} \big(\Hom(M, X)/\cDiff(M)\big)$ can be identified with the stack sending each inf-manifold $U$ to the space of dashed maps in the square
	$$
		\begin{tikzcd}[column sep = 1.25cm, row sep = 1cm]
			U
			& [3pc] U\times M \arrow[l, "\pi_1"{swap}] \arrow[d, "x\circ \pi_2"] \arrow[Leftarrow, ld, shorten=1.5ex, end anchor={[xshift=8ex, yshift=0ex]}, start anchor={[xshift=3ex, yshift=-1.5ex]}, "g"]
			\\
			U\times M\arrow[u, "\pi_1"]\arrow[ru, dashed, "\phi"{description}]\arrow[r, "x\circ \pi_2"{swap}]
			& {[M/\mathcal{G}]}
		\end{tikzcd}
	$$
  which are equivalences and make both triangles commute up to homotopy.	When $U$ is non-derived, the space of such lifts is a set: for the top left triangle, commutativity is simply a condition expressing $\phi$ as a map of the form $(\pi_1, \phi')$ for $\phi'\colon U\times M\longrightarrow M$. Since the quotient stack $[M/\mathcal{G}]$ sends each non-derived $U$ to a groupoid, there is a set of homotopies making the bottom right triangle commute: they correspond to the set of \emph{bisections} $$g\colon U\times M\longrightarrow \mathcal{G}$$ of the Lie groupoid$\begin{tikzcd}[column sep=0.6cm] s,t:\mathcal{G}\arrow[r, yshift=0.5ex] \arrow[r,yshift=-0.5ex]
			& M \end{tikzcd}$over \smash{$\phi'$}, i.e.\ maps such that $s \circ g = \pi_2$ and \smash{$t \circ g  =\phi'$}.
	
	Proceeding as in Example~\ref{ex:synthetic lie diff}, it now suffices to take two first order infinitesimal families of bisections $g, h\colon D^1\times M\longrightarrow \mathcal{G}$ (sending the basepoint to the identity) and consider the map 
	$$
		D^1\times D^1\times M \longrightarrow \mathcal{G}\ ,
		\qquad
		(\epsilon_1, \epsilon_2, m) \longmapsto g(\epsilon_1, m)\, h(\epsilon_2, m)\, g(\epsilon_1, m)^{-1}\, h(\epsilon_2, m)^{-1}\ .
	$$
	The $\epsilon_1\, \epsilon_2$ term precisely gives the usual Lie bracket of the Lie algebroid $\frg$ (see~\cite[Section~15.4]{CW:Geometric_Models_for_Noncommutative_Algebras}). For an alternative proof of the fact that Theorem~\ref{thm:formal stack vs lie} indeed reproduces the usual Lie algebroid associated to a Lie groupoid, see~\cite[Corollary 6.3.33]{Nuiten:Thesis}.
	\qen
\end{example}


\section{Connections and curvature on principal $\infty$-bundles}
\label{sec: Higher connections}


In this section we present a general derived differential geometric approach to higher-form connections and curvatures on geometric structures on manifolds.
More precisely, we describe geometric structures on a manifold $M$ in terms of morphisms $x \colon M \longrightarrow X$ in $\Sh_\infty(\InfMfd)$ from $M$ to a dt stack $X$ (see Definition~\ref{def:qcoh def thy}), where we  understand $X$ as classifying a certain type of geometric structure on $M$.
In particular, for  the classifying stack of a Lie group $X = \rmB G$ we recover classical connections on principal $G$-bundles, while in Section~\ref{sec: cons on higher U(1)-bundles} we compute explicitly how our new definition of higher-form connections recovers all spaces of $p$-form connections on higher principal bundles classified by $X=\rmB^n \rmU(1)$ for any $n \in \NN$ (i.e.~$(n{-}1)$-gerbes).


\subsection{$p$-form connections on geometric structures}
\label{sec:l-conn}


Let $M$ be a smooth manifold and $X$ a dt stack.
We will define the space of $p$-form connections on $x \colon M \longrightarrow X$ in terms of the associated Atiyah $L_\infty$-algebroid $T(M/X)$ (see Theorem~\ref{thm:formal stack vs lie} and Definition~\ref{def:Atiyah L_oo-algebroid from defthy}).
It is worth emphasising that once this $L_\infty$-algebroid has been associated to $x$, our definition is purely algebraic.
We thus present it first at the level of $L_\infty$-algebroids.

To that end, let  $A$ be a cdga over $\RN$.
We recall from~\cite[Section~5]{Nuiten:HoAlg_for_Lie_Algds} (see in particular the proof of Proposition~5.17 therein) the chain of functors and natural transformations
\begin{equation}
\label{eq:sequence of Q^(k)s}
	Q^{(1)} \longrightarrow Q^{(2)} \longrightarrow \cdots \longrightarrow Q^{(\infty)} = Q\ ,
\end{equation}
where each $Q^{(p)}$ is a functor (of ordinary categories)
\begin{equation}
	Q^{(p)} \colon L_\infty\Agd_A^\dg \longrightarrow L_\infty\Agd_A^\dg\ .
\end{equation}
The functor $Q$ is the left adjoint in an adjunction~\cite[Lemma~5.8]{Nuiten:HoAlg_for_Lie_Algds}
\begin{equation}
	\begin{tikzcd}
		Q : L_\infty\Agd_A \ar[r, shift left=0.15cm, "\perp"' yshift=0.05cm]
		& L_\infty\Agd_A^\dg : \iota\ . \ar[l, shift left=0.15cm]
	\end{tikzcd}
\end{equation}

The functors $Q^{(p)}$ can be described as follows~\cite[Section~5]{Nuiten:HoAlg_for_Lie_Algds}:
let
\begin{equation}
	F \colon \big(\Mod_A^\dg\big)_{/T_A} \longrightarrow L_\infty\Agd_A^\dg
\end{equation}
denote the functor which takes the free $L_\infty$-algebroid on a dg $A$-module with its anchor map to $T_A$ (see e.g.~\cite{Kapranov:Free_LiAgds_and_space_of_paths}).
Then the underlying dg $A$-module and brackets of $Q^{(p)} \frg$ are the same as those of
\begin{equation}
	F \Big( \Sym^{1 \leqslant \bullet \leqslant p}_A \big( \frg[1] \big) [-1] \Big)\ ,
\end{equation}
but the differential is not the differential produced by $F$.
Instead, the differential on $Q^{(p)} \frg$ is obtained by observing that the differential on $Q \frg$ restricts to $Q^{(p)} \frg \subset Q \frg$, for each $p \in \NN$.
On $Q \frg$ the differential reads as~\cite[Equations~(5.11) and (5.12)]{Nuiten:HoAlg_for_Lie_Algds}
\begin{equation}
	\dd_{Q\frg} (x_1 \otimes \cdots \otimes x_n)
	= \dd_\frg (x_1 \otimes \cdots \otimes x_n)
	+ \kappa(x_1 \otimes \cdots \otimes x_n)\ ,
\end{equation}
where
\begin{align}\label{eq:kappa}
\begin{split}
	\dd_\frg (x_1 \otimes \cdots \otimes x_n)
	&= \sum_{k = 1}^n\,\pm \, x_1 \otimes \cdots \otimes \dd_\frg (x_k) \otimes \cdots \otimes x_n\ ,
	\\[4pt]
	\kappa(x_1 \otimes \cdots \otimes x_n)
	&= - \sum_{k \geqslant 2} \ \sum_{\sigma \in {\rm Sh}(k, n-k)}\,\pm \, [x_{\sigma(1)}, \ldots, x_{\sigma(k)}]_{\frg,k} \otimes x_{\sigma(k+1)} \otimes \cdots \otimes x_{\sigma(n)}
	\\
	&\quad\, + \sum_{k \geqslant 2}\, \frac{1}{k!} \, \big[ \sfDelta^k (x_1 \otimes \cdots \otimes x_n) \big]_{F,k+1}\ .
	\end{split}
\end{align}
Here $[-]_{\frg,k}$ is the $k$-ary bracket on $\frg$, $\sfDelta^k$ is the $k$-fold comultiplication in $\overline{\Sym}_A(\frg[1])$ and $[-]_{F,k}$ is the $k$-ary formal bracket on the free $L_\infty$-algebroid.
Note that $\dd_\frg$ and the double sum  in~\eqref{eq:kappa} combine into the Chevalley--Eilenberg differential on the cdgc $\overline{\Sym}_{\RN}(\frg[1])$. 

\begin{example}
	\label{eg:Q^(1) frg for non-dg Lie algebroid}
	If $A$ and $\frg$ are both concentrated in degree zero, then
	\begin{equation}
		Q^{(1)} \frg = \big( F(\frg), \dd_{\overline{\Sym}_A(\frg)| \leqslant 1} \big)
		= F(\frg)
	\end{equation}
	agrees with the free Lie algebroid over $A$ on $\frg$.
	This will allow us to recover Kapranov's approach to connections~\cite{Kapranov:Free_LiAgds_and_space_of_paths}.
	\qen
\end{example}

Let $\frh$ be an $L_\infty$-algebroid over $A$.
We let $\ChEil_{\leqslant p}(\frh)\subset \ChEil_*(\frh)$ denote the sub-dg coalgebra on polynomials of degree at most $p$.
Since morphisms out of $L_\infty$-algebroids of the form $Q^{(p)}\frh$ will be central to the computations in the remainder of this paper, we briefly describe them explicitly through

\begin{lemma}
	\label{st:describing order~l-mps of L_oo-Agds}
	Let $\frg$ and $\frh$ be $L_\infty$-algebroids over $A$.
	For each $p \in \NN$, there are canonical natural bijections between
	\begin{myenumerate}
		\item the set of morphisms $L_\infty\Agd_A^\dg (Q^{(p)} \frh, \frg)$,
		
		\item the set of morphisms of $\RN$-linear cdgcs
		\begin{equation}
			\begin{tikzcd}[column sep = 0.cm, row sep = 1cm]
				\ChEil_{\leqslant p}(\frh)
				\arrow[rd, "\ChEil_{\leqslant p}(\rho_{\frh})"{swap}] \arrow[rr, "\phi"]
				& & \ChEil_*(\frg)\arrow[ld, "\ChEil_*(\rho_{\frg})"]
				\\
				& \ChEil_*(T_A) &
			\end{tikzcd}
		\end{equation}
		such that the underlying map of non-unital graded coalgebras $\ol{\Sym}_\RN(\frh[1]) \longrightarrow \ol{\Sym}_\RN (\frg[1])$ descends to an $A$-linear map of graded coalgebras $\ol{\Sym}_A(\frh[1])\longrightarrow \ol{\Sym}_A(\frg[1])$, and
		
		\item the set of degree one maps of graded $A$-modules $\phi \colon \ol{\Sym}^{\,\leqslant p}_A(\frh[1])\longrightarrow \frg$ such that 
		\begin{mysubenumerate}
			\item the  underlying $\RN$-linear map defines a Maurer--Cartan element in $\Hom_\RN\big(\ol{\Sym}_\RN(\frh[1]), \frg[1])\big)$ (see~\eqref{eq:MC eqn for cdgc mps into CE cdgc}), and
			
			\item the composition \smash{$\ol{\Sym}_A(\frh[1]) \longrightarrow \frg \xrightarrow{ \ \rho_{\frg} \ } T_A$} is given by $\rho_{\frh}$ on the linear copy $\frh[1]$ and vanishes on all higher tensor powers of $\frh[1]$.
		\end{mysubenumerate}
	\end{myenumerate}
\end{lemma}

\begin{proof}
	This follows readily by the description of the differential on $Q^{(p)} \frh$ and the fact that $\ChEil_*(\frg)$ is a cofree cocommutative dg coalgebra, together with Definition~\ref{def:oo-morphism of L_oo-algebroids} and Remark~\ref{rmk:constructing cdgc-maps into CE cdgc}.
\end{proof}

\begin{definition}
	\label{def:order~l-mp of Loo_Agds}
	For $p \in \NN \cup \{\infty\}$, an \textit{order~$p$ morphism of $L_\infty$-algebroids} $\frh \longrightarrow \frg$ is an element of
	\begin{equation}
		L_\infty\Agd_A^{(p)} (\frh, \frg)
		\coloneqq L_\infty\Agd_A^\dg \big( Q^{(p)} \frh, \frg \big)\ .
	\end{equation}
\end{definition}

We can now introduce an algebraic definition of $p$-form connections through

\begin{definition}
	\label{def:k-cons general}
	Let $A$ be a unital commutative dg algebra (cdga) over $\RN$.
	Let $\frg$ be an $L_\infty$-algebroid over $A$.
	The \textit{space of $p$-form connections on $\frg$} is the mapping space
	\begin{equation}
		\Con_p(\frg) \coloneqq \Map_{L_\infty\Agd_A^\dg} \big( Q^{(p)} T_A, \frg \big)
		\simeq L_\infty\Agd_A^\infty \big( Q^{(p)} T_A, \frg \big)\ .
	\end{equation}
\end{definition}

This definition becomes particularly accessible for $A$-cofibrant $L_\infty$-algebroids~\cite[Definition~5.16]{Nuiten:HoAlg_for_Lie_Algds}.

\begin{definition}
	An $L_\infty$-algebroid $\frh$ over $A$ is  \textit{$A$-cofibrant} if its underlying dg $A$-module is cofibrant.
\end{definition}

\begin{remark}
	\label{rmk:computing Con_l(frg)}
	In the situation of Definition~\ref{def:k-cons general}, if $T_A$ is $A$-cofibrant then we can compute the space of $p$-form connections on $\frg$ as
	\begin{equation}
		\label{eq:computation of Con_l(frg)}
		\Con_p(\frg) \cong L_\infty\Agd_A^\dg \big( Q^{(p)} T_A, \widehat{\frg} \big)
		\cong L_\infty\Agd_A \big( Q^{(p)} T_A, \widehat{\frg} \big)\ ,
	\end{equation}
	where the isomorphisms are in the homotopy category $\Ho\sSet$ and $\widehat{\frg}$ is any fibrant simplicial resolution of $\frg$ in \smash{$L_\infty\Agd_A^\dg$}.
	Indeed, given an $A$-cofibrant $L_\infty$-algebroid \smash{$\frh \in L_\infty\Agd_A^\dg$}, the $L_\infty$-algebroid $Q^{(p)} \frh$ is cofibrant for each $p \in \NN \cup \{\infty\}$ (see the proof of~\cite[Proposition~5.17]{Nuiten:HoAlg_for_Lie_Algds}).
	This implies the first isomorphism in~\eqref{eq:computation of Con_l(frg)}, since we can now compute the mapping space by choosing a fibrant simplicial resolution $\widehat{\frg}$ of $\frg$.
	To pass from strict to $\infty$-morphisms of $L_\infty$-algebroids, we can further apply the functor $Q$ to the source (this produces a weekly equivalent cofibrant $L_\infty$-algebroid by~\cite[Proposition~5.22(c)]{Nuiten:HoAlg_for_Lie_Algds}) to obtain the second isomorphism in $\Ho \sSet$, since
	\begin{equation}
		L_\infty\Agd^\dg_A \big( Q Q^{(p)} T_A, \widehat{\frg} \big)
		\cong L_\infty\Agd_A \big( Q^{(p)} T_A, \widehat{\frg} \big)
	\end{equation}
	by the defining property of $Q$.
	\qen
\end{remark}

\begin{example}
	\label{eg:T_(C^oo(M)) is C^oo(M)-cofibrant}
	Let $M$ be a connected manifold, and set $A = C^\infty(M)$ (regarded as a cdga concentrated in degree zero).
	For each vector bundle $E \longrightarrow M$, the $C^\infty(M)$-module of global sections $\Gamma(M,E)$ is finitely generated and projective~\cite[Theorem~12.32]{Nestruev:Smooth_mfds_and_observables}.
	Since the projective model structure on $\Mod^\dg_A$ is inherited from that on $\Mod^\dg_\RN$ along the forgetful functor, each projective $A$-module is cofibrant.
	In particular, $TM = T_{C^\infty(M)}$ is a $C^\infty(M)$-cofibrant $L_\infty$-algebroid.
	\qen
\end{example}

Let again $M$ be a smooth manifold and $x \colon M \longrightarrow X$ a morphism to a dt stack (see Definition~\ref{def:qcoh def thy}).
Combining the algebraic Definition~\ref{def:k-cons general} with the derived  geometric differentiation of maps $x \colon M \longrightarrow X$ from Theorem~\ref{thm:formal stack vs lie}, we introduce a derived geometric definition of higher-form connections on geometric structures through

\begin{definition}
	\label{def: Connection}
	Let $x \colon M \longrightarrow X$ be an object in $\DT_{M/}$.
	The \textit{space of $p$-form connections on $x$} is
	\begin{equation}
		\Con_p(x) \coloneqq \Con_p \big( T(M/X) \big)\ .
	\end{equation}
\end{definition} 

\begin{example}
	\label{eg:BG and qcoh defthy}
	An important class of target stacks consists of the classifying stacks for higher principal bundles:
	let $G_\bullet$ be a group object in $\Sh_\infty(\InfMfd)$ whose underlying object $G = G_1 \in \Sh_\infty(\InfMfd)$ has (quasi-coherent) deformation theory and is formally smooth.
	Recall that a group object in an $\infty$-topos is in particular a groupoid object in the sense of~\cite[Definition~6.1.2.7]{Lurie:HTT}.
	By~\cite[Example~A.5.2.5]{Lurie:SAG}, for any groupoid object $X_\bullet$  the morphisms
	\begin{equation}
		X_n \longrightarrow X(\Lambda^n_i) = \underset{\Delta^k\longrightarrow \Lambda^n_i}{\holim}\ X_k
	\end{equation}
	are equivalences.
	Since any equivalence in $\Sh_\infty(\InfMfd)$ is formally smooth (Definition~\ref{def:formally smooth}), $G_\bullet$ is a formally smooth $\infty$-groupoid (Definition~\ref{def:form smooth oo-groupoid}).
	Moreover, since $G_n \simeq (G_1)^{\times n}$ for each $n \in \NN_0$, by Proposition~\ref{prop:def thy}(2) each $G_n$ for $n \in \NN_0$ as well as the homotopy colimit
	\begin{equation}
		\rmB G = \underset{[n] \in \bbDelta^\opp}{\hocolim}\ G_n
	\end{equation}
	are stacks with (quasi-coherent) deformation theory as well.
	\qen 
\end{example}

The stack $\rmB G$ classifies $G$-principal $\infty$-bundles in $\Sh_\infty(\InfMfd)$~\cite{NSS:oo-bundles_I, NSS_oo-Bundles_II} (see also~\cite{Bunk:oo-bundles} for an overview and~\cite{Sati:2021eyj} for a more in-depth account).
In Example~\ref{Example: Lie group} below we consider the cases where $G$ is an ordinary Lie group and show how we recover the classical theory of connections.
In Section~\ref{sec: cons on higher U(1)-bundles} we will compute the spaces of $p$-form connections on morphisms $x \colon M \longrightarrow \rmB^n \rmU(1)$ and show that these recover the higher groupoids of $p$-form connections on higher $\rmU(1)$-bundles.

\begin{example}
	\label{Example: Lie group}
	Let $G$ be a Lie group and $$x_P \colon M \longrightarrow \rmB G$$ the classifying map for a principal $G$-bundle $P\longrightarrow M$.
	This is, a priori, a morphism in $\Sh_\infty(\Mfd)$.
	It induces a morphism $\iota_!\, x_P \colon M \longrightarrow \iota_! \rmB G$ in $\Sh_\infty(\InfMfd)$ under the functor $\iota_!$ from~\eqref{diag:infinitesimal cohesion}.%
	\footnote{
		Since $\iota_!$ preserves finite products and is a left adjoint, it preserves group objects and there is a canonical equivalence $\iota_! \rmB G \simeq \rmB(\iota_! G)$ for each group object $G$ in $\Sh_\infty(\Mfd)$.
		In particular, it commutes with forming classifying objects.
	}
	In computations involving objects in $\Sh_\infty(\Mfd)$ we will often omit the functor $\iota_!$ for ease of notation.
	
	The Atiyah $L_\infty$-algebroid of $x_P:M\longrightarrow \rmB G$ is the traditional Atiyah algebroid of $P$~\cite[Corollary~6.4.33]{Nuiten:Thesis}
	\begin{align}
		T(M/\rmB G) = \mathfrak{X}(P)^G \longrightarrow TM\ ,
	\end{align}
	where $\mathfrak{X}(P)^G$ denotes the $G$-invariant vector fields on $P$.
	We also write
	\begin{equation}
		\At(P) = T(M/\rmB G)\ .
	\end{equation}
	The space of 1-form connections on $P$ in Definition~\ref{def: Connection} is the mapping space
	\begin{equation}
		\Con_1(P):=\Con_1(x_P) = \Map_{L_\infty\Agd_{M}^\dg} \big( Q^{(1)} TM, \At(P) \big)\ .
	\end{equation}
	
	The $L_\infty$-algebroid $Q^{(1)} TM$ was computed in Example~\ref{eg:Q^(1) frg for non-dg Lie algebroid}.
	Both $Q^{(1)} TM$ and $T(M/\rmB G)$ are concentrated in degree zero.
	Hence we can compute the derived mapping space in the model structure on $L_\infty$-algebroids concentrated in non-positive degrees from Remark~\ref{rem:nonpositively graded case}.
	To do this, we have to pick a Reedy fibrant simplicial resolution of $\At(P)$.
	As a consequence of Remark~\ref{Rem: Reed fibrant resolutions <1} we can choose the constant resolution in this case.
	By Remark~\ref{rmk:computing Con_l(frg)} and Example~\ref{eg:T_(C^oo(M)) is C^oo(M)-cofibrant} it is not necessary to choose a cofibrant resolution of the source $L_\infty$-algebroid $Q^{(1)} TM$.
	Altogether we find 
	\begin{align}
		\Con_1(P) \cong L_\infty\Agd_M^\dg \big( Q^{(1)} TM, \At(P) \big)
		\cong \Mod^\dg_{C^\infty(M)}\big(TM,\At(P)\big) \ ,
	\end{align}  
	where for the last isomorphism  we used that $Q^{(1)} TM$ is the free Lie algebroid on $TM$, see Example~\ref{eg:Q^(1) frg for non-dg Lie algebroid}.
	Hence our notion of a 1-form connection on a principal $G$-bundle $P\longrightarrow M$ from Definition~\ref{def: Connection} reduces to the usual definition of a connection on $P$ in the formulation of Kapranov~\cite{Kapranov:Free_LiAgds_and_space_of_paths}, i.e.~using the free Lie algebroid on the tangent bundle of $M$.
		
	By the same arguments, for $p \geqslant 2$ the space of $p$-form connections on $P$ is
	\begin{equation}
		\Con_p(P) = L_\infty\Agd_M^\dg \big( Q^{(p)} TM, \At(P) \big)\ .
	\end{equation}
	An element of this set is the same as a morphism $TM \longrightarrow \At(P)$ of Lie algebroids on $M$, which are in one-to-one correspondence with \textit{flat} connections on $P$.
	\qen
\end{example}

\begin{example}
	\label{eg: cons on trivial crossed module bundle}
	Another important class of examples consists of principal 2-bundles whose structure 2-group is a crossed module $$G = \big(G_{-1}\,,\, G_0\,,\, \delta\colon G_{-1} \longrightarrow G_0\,,\, \alpha \colon G_{-1} \longrightarrow \cAut(G_0)\big)$$ of Lie groups.
	Using that all manifolds are qcdt stacks by Example~\ref{ex:mfd are qcdt} and that all submersions of manifolds are formally smooth (Example~\ref{ex:submersions are formally smooth}), together with Proposition~\ref{prop:form smooth oo-groupoid qcdt}, we see that any such crossed module gives rise first to a (2-coskeletal) simplicial Lie group and then to a qcdt stack $\rmB G$.
	Thus for each classifying map $$x \colon M \longrightarrow \rmB G$$ of a crossed module 2-bundle, there is an associated Atiyah $L_\infty$-algebroid $$\At(x)=T(M/\rmB G) \ . $$
	
	Here we consider the case where $x$ classifies the trivial $G$-bundle and check that our new notion of $2$-form connection on $x$ recovers the known one in the literature (see for instance~\cite{Waldorf:Global_perspective}).
	Since $x$ is trivial, the Atiyah $L_\infty$-algebroid is
	\begin{equation}
		\At(x) = \big( \frg_{-1} \otimes C^\infty(M) \longrightarrow \frg_0 \otimes C^\infty(M) \oplus TM\big)\ ,
	\end{equation}
	which lies in degrees $-1$ and $0$,  where the brackets are induced from those of the corresponding Lie 2-algebra $\frg$ and the canonical actions of $TM$ on itself as well as on $C^\infty(M)$.
	Using Lemma~\ref{st:describing order~l-mps of L_oo-Agds}, it follows that the vertices of $\Con_2(x)$ are given by  pairs $$\big(A^{(1)} \in \Omega^1(M; \frg_0)\,,\, A^{(2)} \in \Omega^2(M;\frg_{-1})\big)$$ satisfying the fake-curvature condition
	\begin{equation}
		\delta A^{(2)}  = \dd A^{(1)} + [A^{(1)}, A^{(1)}] \ ,
	\end{equation}
	as expected.
	Flat 2-form connections have also been described as morphisms of higher Lie algebroids in~\cite[Section~3.6.5]{AC:Homological_vfs_on_stacks}.
	\qen
\end{example}

\begin{remark}
	\label{rmk: fake flatness and adjusted connections}
	In Example~\ref{eg: cons on trivial crossed module bundle} one also sees that our connections are fake flat.
	This is consistent with the perspective that Definition~\ref{def:k-cons general} really defines a $p$-form connection as a morphism out of the infinitesimal path $p$-groupoid of $M$, i.e.~as an infinitesimal order~$p$ version of a parallel transport; we make this statement precise in Section~\ref{sec: de rham description} below.
	Such higher parallel transports necessarily satisfy the fake-curvature condition~\cite{BS:HGT} so our connections will also satisfy this condition.
	Current discussions surrounding a modified notion of connection which does not impose the fake-curvature condition, which are crucial for applications in supergravity and string theory contexts, can be found in e.g.~\cite{SSS:L_infty-alg_connections, SS:Towards_M5-brane_model_II, KS:Adjusted_PT, FFKS:Adjusted_connections_I, BFJKNRSW:GHT}.
	\qen
\end{remark}

\begin{remark} \label{rem:String2group}
	In principle, we can compute the Atiyah algebroid and space of connections for a higher principal bundle $P \colon M\longrightarrow \rmB G$ with general higher structure group $G$ from a Whitehead-like decomposition of $\rmB G$ and the compatibility of all our constructions with pullbacks.
	For example, if $K$ is a compact, simple and simply-connected Lie group, then the string 2-group $\String(K)$ fits into the pullback
	\begin{equation}
		\begin{tikzcd}[row sep = 1cm,column sep=1cm]
			\rmB\String(K) \ar[r] \ar[d]
			& * \ar[d]
			\\
			\rmB K \ar[r]
			& \rmB^3 \rmU(1)
		\end{tikzcd}
	\end{equation}
	where the map $\rmB K \longrightarrow \rmB^3 \rmU(1)$ classifies the Chern--Simons 2-gerbe on $\rmB K$.
	From this description we see that a $\String(K)$-bundle $P_{\String(K)}\colon M \longrightarrow \rmB \String(K)$ on a manifold $M$ is the same as a principal $K$-bundle $P_K$ on $M$ together with a trivialisation of the pullback of the Chern--Simons 2-gerbe to its total space $P_K$.    
	By Theorem~\ref{thm:formal stack vs lie} the Atiyah $L_\infty$-algebroid $T(M / \rmB\String(K))$ is equivalent to the pullback $T(M / \rmB K)\times_{T(M / \rmB^3 \rmU(1))}*$. 
	
	It follows that a 1-form connection on $P_{\String(K)}$ is specified by a connection on $P_K$ together with a compatible trivialisation of the induced 1-form connection on the Chern--Simons 2-gerbe on $M$. If we consider instead $p$-form connections for $p>1$, then the induced connection on $P_K$ must be flat, which presents an obstruction to their existence.
	In particular, we recover the well-known fact that geometric string structures on $P_K$ are not the same as 2-form connections on the underlying string structure~$P_{\String(K)}$.
	\qen
\end{remark}

Definition~\ref{def: Connection} also allows us to define $p$-form connections only along group actions or foliations:
let $x \colon M \longrightarrow X$ be a morphism from a smooth manifold to a dt stack.
Let $\frh$ be an $L_\infty$-algebroid on $M$ with a morphism $\frh \longrightarrow TM$ in $L_\infty\Agd_M$.
We can then form the homotopy pullback
\begin{equation}
	\begin{tikzcd}[row sep = 1cm,column sep=1cm]
		\frh \times^h_{TM} T(M/X) \ar[r] \ar[d]
		& T(M/X) \ar[d]
		\\
		\frh \ar[r]
		& TM
	\end{tikzcd}
\end{equation}

\begin{definition}
	\label{def: Connection along h}
	The space of \textit{$p$-form connections on $x$ along $\frh$} is
	\begin{align}
		\Con_p(x; \frh) &= (L_\infty\Agd_M^\infty)_{/\frh} \big( Q^{(p)} \frh,\, \frh \times^h_{TM} T(M/X) \big)
		\\[4pt]
		&\simeq L_\infty\Agd_M^\infty \big( Q^{(p)} \frh,\, T(M/X) \big)\ .
	\end{align}
\end{definition}

\begin{example}\label{ex:p-connG}
	If $G$ is a higher smooth group acting on $M$ with associated $L_\infty$-algebra $\frg$, there is an induced infinitesimal action $\frg \otimes \calO(M) \longrightarrow TM$.
	A $p$-form connection on $x$ along $\frg \otimes \calO(M)$ can be viewed as an order~$p$ infinitesimal lift of the $G$-action on $M$ to the geometric structure on $M$ classified by $x$.
	\qen
\end{example}

\begin{example}
	Another source of $L_\infty$-algebroids on $M$ consists of foliations of $M$, or equivalently integrable distributions in $TM$.
	Definition~\ref{def: Connection along h} allows us to define $p$-form connections on $x$ along any integrable distribution in $TM$.
	
	A particularly interesting case is the following:
	suppose that $M$ comes endowed with a submersion $s \colon M \longrightarrow B$.
	Then the fibres of $s$ are submanifolds of $M$ and $\ker(s_*) \subset TM$ is an integrable distribution.
	If we understand $s \colon M \longrightarrow B$ as encoding a $B$-parametrised family of smooth manifolds $(M|b)_{b \in B}$ and the morphism $x \colon M \longrightarrow X$ as a $B$-parametrised family of geometric structures $x_{M|b} \colon M|b \longrightarrow X$, then a $p$-form connection along $\ker(s_*)$ can be viewed as a $B$-family of $p$-form connections on $x_{M|b}$.
	In particular, this provides a general and perhaps surprisingly simple description of concretification (see~\cite{Schreiber:DCCT, BSS:YM_fields_on_Lorentz_mfds}), i.e.~a description of smooth families of geometric structures with connections on manifolds where the connection data is defined not on the total space of a family of manifolds but only along its fibres.
	Finally, if the family $M \longrightarrow B$ is a bundle with fibre $F$, then this example is a sub-case of Example~\ref{ex:p-connG} where $G = \cDiff(F)$ acts canonically on $M$ along the fibres.
	\qen
\end{example}

\begin{remark}
	On every smooth manifold $M$ there is a universal $p$-form connection in the sense of Definition~\ref{def:k-cons general}, for each $p \in \NN$.
	It consists of the identity morphism
	\begin{equation}
		\CA_{p, \mathrm{univ}} = \id \colon Q^{(p)} TM \longrightarrow Q^{(p)} TM\ .
	\end{equation}
	We  expect that these universal connections (or translation-invariant versions thereof) furnish the classical path signature for $p=1$ and higher-dimensional versions thereof for higher values of $p$ (see e.g.~\cite{LO:Random_surfaces, BL:Surface_signature}).
	\qen
\end{remark}


\subsection{Curvature forms}
\label{sec: curvature of higher connections}


We now derive the curvature of the $p$-form connections we introduced in Definitions~\ref{def:k-cons general} and~\ref{def: Connection}.
To that end, fix a morphism $x \colon M \longrightarrow X$ in $\Sh_\infty(\InfMfd)$ from a smooth manifold to a dt stack, and consider a $p$-form connection $\tau$ on $x$.
Its curvature appears as the obstruction to lifting $\tau$ to a $(p{+}1)$-form connection.
To see this, we study the data required to extend an order~$p$ morphism out of an $L_\infty$-algebroid $\frh$ to an order~$(p{+}1)$ morphism, or equivalently a map out of $Q^{(p)}\frh$ to a map out of $Q^{(p+1)}\frh$. 
This is accomplished by means of~\cite[Proposition 5.17]{Nuiten:HoAlg_for_Lie_Algds} (or rather its proof), which describes $Q^{(p+1)}\frh$ as a pushout of $Q^{(p)}\frh$.
Recall the free-forgetful adjunction
\begin{equation}
	\begin{tikzcd}
		F : (\Mod_A^\dg)_{/T_A} \ar[r, shift left = 0.15cm, "\perp"' yshift = 0.05cm]
		& \mathrm{Lie}\Agd_A^\dg : U \ar[l, shift left = 0.15cm]\ ,
	\end{tikzcd}
\end{equation}
where  the right-hand side is the category of dg Lie algebroids over $A$ (the full subcategory of $L_\infty\Agd_A^\dg$ on those $L_\infty$-algebroids whose $m$-ary brackets vanish for all $m > 2$).

\begin{lemma}
	[{\cite[Proof~of Proposition~5.17]{Nuiten:HoAlg_for_Lie_Algds}}]
	\label{st:pushout square for Q^(l)}
	Let $A$ be a cdga over $\RN$, and $\frh$ an $L_\infty$-algebroid over $A$.
	For each $p\in \NN_0$ there is a pushout of $L_\infty$-algebroids
	\begin{equation}
		\label{eq:pushout square for Q^(l)}
		\begin{tikzcd}[column sep = 1.5cm, row sep = 1cm]
			F \big( \Sym_A^{p+1}(\frh[1]) [-2] \big) \ar[r, "\kappa"] \ar[d]
			& Q^{(p)} \frh \ar[d]
			\\
			F \big( \Sym_A^{p+1}(\frh[1]) [-2,-1] \big) \ar[r]
			& Q^{(p+1)} \frh
		\end{tikzcd}
	\end{equation}
	where the bottom-left vertex is the free $L_\infty$-algebroid on a shifted cone of the dg $A$-module $\Sym_A^{p+1}(\frh[1])$.
	Both free dg Lie algebroids on the left are formed using the zero anchor map, and the map $\kappa$ was defined in~\eqref{eq:kappa}.
	If $\frh$ is $A$-cofibrant, then this diagram is even a homotopy pushout.
\end{lemma}

To describe the curvature of a $p$-form connection, we will now present a dual version of this statement as follows.
Let $\tau \colon \frh \longrightarrow \frg$ be an order~$p$ morphism of $L_\infty$-algebroids over $A$. 
To extend $\tau$ to an order~$(p{+}1)$ morphism $\tau' \colon \frh \longrightarrow \frg$, by Lemma~\ref{st:describing order~l-mps of L_oo-Agds} it suffices to specify a map
\begin{equation}
	\tau'_{p+1}\colon \Sym_{A}^{p+1} (\frh[1]) \longrightarrow \ker (\rho_{\frg})[1]
\end{equation}
such that $\tau + \tau'_{p+1}$ is a Maurer--Cartan element in $\Mod_\RN^\dg(\ChEil_{\leqslant p+1}(\frh),\frg)$. 
Observe from~\eqref{eq:hom brackets} and~\eqref{eq:MC eqn for cdgc mps into CE cdgc} that in the Maurer--Cartan equation in question all the brackets involving $\tau'_{p+1}$ vanish, since $\ChEil_{\leqslant p+1}(\frh)$ vanishes in degrees less than $-(p{+}1)$.
Hence the Maurer--Cartan equation becomes
\begin{align}
	\dd_\frg \circ \tau'_{p+1} - \tau'_{p+1} \circ \dd_\frh
	= - \left( \dd_\Hom \tau + \sum_{m = 2}^p\, \frac{1}{m!}\, [\tau, \dots, \tau]_{\Hom, m} \right)\ ,
\end{align}
where the differential $\dd_\frg$ on the left-hand side is that on the complex $\frg$, while $\dd_\frh$ is---by a slight abuse of notation---induced from the extension of the differential of $\frh$ to \smash{$\Sym^{\leqslant p+1}_A(\frh[1])$}.
In other words, the left-hand side can be written as
\begin{equation}
	\dd_\frg \circ \tau'_{p+1} - \tau'_{p+1} \circ \dd_\frh = \dd_\Hom \tau'_{p+1}\ .
\end{equation}

Since by assumption $\tau$ is a Maurer--Cartan element in $\Mod_\RN(\ChEil_{\leqslant p}(\frh),\frg)$, the only non-trivial part of the Maurer--Cartan equation involves the component $\Sym^{p+1}_A( \frh[1]) \longrightarrow \frg$.
We write $F^p(\tau)$ for this component of the right-hand side of the equation.
Explicitly
\begin{equation}
	\label{eq:curvature of l-conn algebraically}
	F^p(\tau) = \tau \circ \dd_\frh + \sum_{m = 2}^p\, \frac{1}{m!}\, [\tau, \dots, \tau]_{\Hom, m| \Sym^{\leqslant p+1}_A(\frh[1])}\ .
\end{equation}

\begin{definition}
	\label{def:curvature of l-conn algebraically}
	Given an order~$p$ morphism $\tau \colon \frh \longrightarrow \frg$ of $L_\infty$-algebroids over $A$, the element
	\begin{equation}
		F^p(\tau) \ \in \ \Mod_A^\dg \left( \Sym^{p+1}_A \big( \frh[1] \big),\, \ker(\rho_\frg)[2] \right)
	\end{equation}
	defined in~\eqref{eq:curvature of l-conn algebraically} is the \emph{curvature} of $\tau$.
\end{definition}

The construction of the curvature $F^p$ gives rise to a map
\begin{equation}
	\label{eq:curvature map}
	F^p \colon L_\infty\Agd_A^\dg (Q^{(p)} \frh, \frg)
	 \longrightarrow \Mod_A^\dg \left( \Sym_{A}^{p+1} \big( \frh[1] \big)[-2], \ker(\rho_\frg) \right)\ .
\end{equation}
Furthermore, from every order~$(p{+}1)$ morphism $\tau' \in L_\infty\Agd_A^\dg (Q^{(p+1)} \frh, \frg)$ we can extract a map of dg $A$-modules
\begin{align}
	(\tau'_{p+1}, \dd_\Hom \tau'_{p+1}) \colon  \big( \Sym_{A}^{p+1}( \frh[1]) \big) [-2,-1]
	 \longrightarrow \ker (\rho_\frg)
\end{align} 
which extracts the component \smash{$\tau'_{p+1} \colon \Sym_{A}^{p+1}( \frh[1]) \longrightarrow \ker (\rho_\frg)[1]$} of $\tau'$ together with its differential (taken in the hom complex).
This defines a map
\begin{align}
\begin{split}
	\label{eq:tau'_(l+1), d_Hom tau'_(l+1)}
	\psi^p \colon L_\infty\Agd_A^\dg (Q^{(p+1)}\frh, \frg)
	& \longrightarrow \Mod_A^\dg \left( \Sym_{A}^{p+1} \big( \frh[1] \big) [-2,-1],  \ker (\rho_\frg) \right)\ ,
	\\
	\tau' &\longmapsto (\tau'_{p+1}, \dd_\Hom \tau'_{p+1})\ .
	\end{split}
\end{align}

The maps~\eqref{eq:curvature map} and~\eqref{eq:tau'_(l+1), d_Hom tau'_(l+1)} fit into a commutative diagram of sets
\begin{equation}
	\label{eq: pullback for curvature}
	\begin{tikzcd}[column sep = 1.5cm, row sep = 1cm]
		L_\infty\Agd_A^\dg (Q^{(p+1)} \frh, \frg) \ar[r, "\psi^p"] \ar[d]
		&  \Mod_A^\dg \left(\Sym_{A}^{p+1} \big( \frh[1] \big) [-2,-1], \ker(\rho_\frg) \right) \ar[d]
		\\ 
		L_\infty\Agd_A^\dg (Q^{(p)} \frh, \frg) \ar[r, swap, "F^p"]
		& \Mod_A^\dg \left( \Sym_{A}^{p+1} \big( \frh[1] \big) [-2], \ker (\rho_\frg) \right)
	\end{tikzcd}
\end{equation}
where the left-hand vertical morphism forgets the order~$(p{+}1)$ component of a morphism.
The commutative square~\eqref{eq: pullback for curvature} is obtained from the square~\eqref{eq:pushout square for Q^(l)} by taking morphisms to $\frg$ (this uses that the $L_\infty$-algebroids on the left-hand side of~\eqref{eq:pushout square for Q^(l)} are free, together with Lemma~\ref{st:describing order~l-mps of L_oo-Agds}(3)).
In particular, by Lemma~\ref{st:pushout square for Q^(l)} it follows that the square~\eqref{eq: pullback for curvature} is cartesian.

With these preparations we can now give

\begin{definition}
	\label{def:curvature of l-conn}
	Let $x \colon M \longrightarrow X$ be a morphism from a smooth manifold to a dt stack.
	Suppose that this morphism is equipped with a $p$-form connection $\tau$.
	The \emph{curvature} or \emph{field strength} of $\tau$ is the morphism of dg $C^\infty(M)$-modules
	\begin{equation}
		F^p(\tau) \colon \left( \Sym_{C^\infty(M)}^{p+1} \big( TM[1] \big) \right) [-2]
		 \longrightarrow \ker(\rho_{T(M/X)})\ .
	\end{equation}   
\end{definition}

\begin{remark}
	Motivated by the theory of classical principal bundles, we can think of $\ker(\rho_{T(M/X)})$ as the module of global sections of the adjoint bundle of $x$.
	In this sense the curvature $F^p(\tau)$ can be interpreted as a (derived) $(p{+}1)$-form with values in the adjoint bundle, generalising the classical picture.
	\qen
\end{remark}

We record the relation between curvature and lifts to $(p{+}1)$-form connections in

\begin{proposition}
	\label{st: extending l-con to (l+1)-con}
	Let $x \colon M \longrightarrow X$ be a morphism from a smooth manifold to a dt stack, and suppose that $\tau$ is a $p$-form connection on $x$.
	There is a one-to-one correspondence between the set of lifts of $\tau$ to a $(p{+}1)$-form connection on $x$ and the set of maps 
	\begin{equation}
		\tau'_{p+1} \colon \left( \Sym_{C^\infty(M)}^{p+1} \big( TM[1] \big) \right)[-1]
		 \longrightarrow \ker(\rho_{T(M/X)})
	\end{equation}
	such that $\dd_\Hom \tau'_{p+1} = F^p(\tau)$.
\end{proposition} 

\begin{proof}
	This is exactly the statement that the diagram~\eqref{eq: pullback for curvature} is a pullback. 
\end{proof}

\begin{example}
	\label{ex:conn-principal}
	For ordinary principal bundles, Definition~\ref{def:curvature of l-conn} exactly recovers the usual definition of the curvature 2-form with values in the adjoint bundle.
	Let $G$ be a Lie group with Lie algebra $\frg$, and let $x \colon M \longrightarrow \rmB G$ be a map in $\Sh_\infty(\Mfd)$.
	We know from Example~\ref{Example: Lie group} that the Atiyah $L_\infty$-algebroid associated to the morphism $x \colon M \longrightarrow \rmB G$ is the classical Atiyah Lie algebroid of the principal $G$-bundle $P\longrightarrow M$ classified by $x$.
	The morphism $x \colon M \longrightarrow \rmB G$ further has a presentation by a morphism $$g \colon \cC\, \CU\, \longrightarrow \rmB G$$ of simplicial presheaves on $\Mfd$, where $\cC\, \CU$ is the \v{C}ech nerve of a good open covering $\CU = \{U_a\}_{a \in \varLambda}$ of the manifold $M$.
	In other words, $g$ is a \v{C}ech 1-cocycle with values in~$G$.
	
	The Atiyah algebroid of $x$ can then be expressed as
	\begin{align}
		T(M/\rmB G) = \At(P)
		&= \big\{ (X, \xi) \ | \ X \in TM \ ,\ \xi = (\xi_a \colon U_a \longrightarrow \frg)_{a \in \varLambda}\ ,
		\\
		&\hspace{2cm} \xi_b = \Ad_{g_{ab}^{-1}}(\xi_a) + g_{ab}^*\,\mu_G(X) \ ,\ \forall\, a,b \in \varLambda \big\}\ ,
	\end{align}
	where $\Ad \colon G \longrightarrow \cAut(\frg)$ is the adjoint representation of $G$ and $\mu_G$ is its Maurer--Cartan 1-form.
	The anchor map $\rho$ is the projection $$\rho(X, \xi) = X \ , $$ and the bracket reads as
	\begin{equation}
		\big[ (X, \xi), (Y, \eta) \big]
		= \big( [X,Y],\, \pounds_X \eta - \pounds_Y \xi + [\xi, \eta] \big)\ .
	\end{equation}
	Note that $\ker(\rho)$ consists precisely of descent data for sections of the adjoint bundle of $P$.
	
	A 1-form connection on $x$ is then a morphism of $C^\infty(M)$-modules
	\begin{equation}
		\tau = (\id_{TM}, A) \colon TM \longrightarrow \At(P)\ ,
		\quad
		X \longmapsto \big( X, A(X) \big)\ .
	\end{equation}
	We compute the curvature of $\tau$:
	we find
	\begin{align}
		\dd_\Hom\, \tau (X_1 \wedge X_2)
		&= \tau \circ \dd_\ChEil (X_1 \wedge X_2)
		= - \tau \big( [X_1, X_2] \big)
		= - \big( [X_1, X_2],\, A([X_1, X_2]) \big)\ ,
\end{align}
and
\begin{align}
		\tfrac{1}{2}\, [\tau, \tau]_{\Hom,2}
		&= \tfrac{1}{2}\, [-]_{\At(P), 2} \circ (\tau \otimes \tau) \circ \sfDelta_{TM} (X_1 \wedge X_2)
		\\[4pt]
		&= \Big( [X_1, X_2],\, \pounds_{X_1} \big( A(X_2) \big) - \pounds_{X_2} \big( A(X_1) \big) - \big[ A(X_1), A(X_2) \big] \Big)\ .
	\end{align}
	Thus since $$\dd A(X_1, X_2) = \pounds_{X_1} \big(A(X_2)\big) - \pounds_{X_2} \big(A(X_1)\big) - A\big([X_1, X_2]\big) \ , $$ we indeed see that
	\begin{equation}
		F^1(A) (X_1, X_2)
		= \dd A(X_1, X_2) - \big[ A(X_1), A(X_2) \big]
	\end{equation}
	recovers the classical expression for the curvature 2-form of the 1-form connection defined by $A$.
	\qen
\end{example}

\begin{example}
	We revisit Example~\ref{eg: cons on trivial crossed module bundle}, where we recovered the notion of 2-form connections on trivial bundles for a crossed module 2-group.
	Here we obtain for the field strength the single 3-form
	\begin{equation}
		F^2(A^{(1)}, A^{(2)}) = \dd A^{(2)} + [A^{(1)}, A^{(2)}] \ ,
	\end{equation}
	which is the well-known expression for the curvature of a crossed module 2-form connection (compare for instance~\cite[Equation~(5.1.7)]{Waldorf:Global_perspective}).
	\qen
\end{example}

\begin{remark}
	\label{rmk: relation to adjusted connections}
	We comment further on the relationship between the new approach to fake-flat connections presented here and the parallel research on adjusted connections.
	By results of Tellez-Dominguez~\cite{TD:Chern_corr_for_higher_PrBuns} and Waldorf~\cite{Waldorf:Adjusted_connections}, 2-bundles with adjusted 2-connections for a large class of crossed module 2-groups can be recast in terms of trivialisations of 2-gerbes and ordinary (i.e.~fake-flat) connections.
	We will briefly sketch the main arguments here, referring the reader to~\cite{TD:Chern_corr_for_higher_PrBuns, Waldorf:Adjusted_connections} for the background on adjusted connections we are using.
	
	Let $G$ be a smoothly separable, central crossed module with a splitting and adapted adjustment.
	Let $\pi_0 G$ be the smooth group of isomorphism classes of objects in $G$.
	Then, for any fixed principal $\pi_0 G$-bundle $P \longrightarrow M$ with connection, by~\cite[Theorem~6.2.3]{Waldorf:Adjusted_connections} there is a canonical equivalence between enhancements of $P$ to a $G$-principal 2-bundle with adjusted connection and trivialisations of the Chern--Simons 2-gerbe $\mathbb{CS}_P$ with connection induced by the connection on $P$ and the adjustment data on $G$.
	
	In this latter picture, all connections can be treated by means of the formalism in this paper.
	To see this, let $H$ be any smooth higher group such that $\rmB H$ has deformation theory, and observe that a morphism of principal $H$-bundles on $M$ is classified by a map
	\begin{equation}
		E \colon M \longrightarrow (\rmB H)^{\Delta^1} \ .
	\end{equation}
	The case relevant in the present setting is $H = \rmB \rmU(1)$.
	Suppose that we are given such a morphism $E$ between principal $H$-bundles $E_0, E_1 \colon M \longrightarrow \rmB H$.
	Further, suppose that we are given $p$-form connections
	\begin{equation}
		A_i \colon Q^{(p)} TM \longrightarrow \At(E_i)\ ,
		\quad   i = 0,1\ .
	\end{equation}
	A connection on $E$ compatible with $A_0$ and $A_1$ is a morphism of complexes
	\begin{equation}
		A_{01} \colon Q^{(p-1)} TM \longrightarrow E
	\end{equation}
	such that
	\begin{equation}
		\dd_\Hom A_{01} + \sum_{n \geqslant 2}\, \frac{[A_{01}, \dots, A_{01}]_{\Hom, n}}{n!}
		= A_{1|Q^{(p-1)} TM} - A_{0|Q^{(p-1)} TM}\ ,
	\end{equation}
	i.e.~in degrees $0, \ldots, -(p-1)$, and
	\begin{equation}
		F^p(A_{01}) = A_1^{(p)} - A_0^{(p)}
	\end{equation}
	as maps 
	\begin{equation}
		\left( \Sym_{C^\infty(M)}^{p+1} \big( TM[1] \big) \right) [-2]
		\longrightarrow \ker(\rho_{T(M/X)})\ ,
	\end{equation}
	i.e.~in degree $-p$.
	Setting $H = \rmB \rmU(1)$ and $p = 2$ here recovers the notion of a connection-preserving morphism of 2-gebres with connection (using also the explicit comparison from Theorem~\ref{st:finite and infinitesimal l-cons on n-gerbes} below).
	A trivialisation of a 2-gerbe $\CG$ with 2-form connection is then such a connection-preserving morphism $\CI_0 \longrightarrow \CG$, where $\CI_0$ denotes the trivial 2-gerbe with the trivial 2-form connection.
	Thus, in the setting described at the beginning of this remark, Waldorf's result~\cite[Theorem~6.2.3]{Waldorf:Adjusted_connections} allows us to recover adjusted 2-connections through the derived geometric formalism developed here.
	
	It should also be possible to further extend the procedure in this remark to define $m$-morphisms of bundles with $p$-form connections for any $m \in \NN$.
	\qen
\end{remark}


\subsection{Functoriality of spaces of connections}
\label{sec: functoriality of conenctions}


We now discuss the functoriality of spaces of $p$-form connections on maps $x \colon M \longrightarrow X$ from a smooth manifold to a dt stack.
This will imply, for instance, a canonical action of automorphisms of $x$ on the space of $p$-form connections.
We first need some further background on $L_\infty$-algebroids and to see that $Q^{(p)}$ induces an $\infty$-functor.

As before, let $A$ be a cdga over $\RN$ (for instance arising as the underlying cdga of a dg $C^\infty$-ring).
Recall that we denote by
\begin{equation}
	L_\infty\Agd_A^\infty
	\coloneqq L_\infty\Agd_A^\dg [W^{-1}]
\end{equation}
the $\infty$-category obtained as the ($\infty$-categorical) localisation of $L_\infty\Agd_A^\dg$ at the weak equivalences of $L_\infty$-algebroids over $A$.
This is canonically equivalent to the localisation of $L_\infty\Agd_A$ at the $\infty$-morphisms which are weak equivalences.

\begin{lemma}
	\label{st:detecting weqs on ass grds}
	Let $\frg_\bullet, \frh_\bullet \colon \NN \longrightarrow L_\infty\Agd_A^\dg$ be functors which send each morphism in the poset $\NN$ to a cofibration between cofibrant objects (equivalently, $\frg_\bullet$ and $\frh_\bullet$ are cofibrant in the Reedy model structure on \smash{$\Fun(\NN, L_\infty\Agd_A^\dg)$)}.
	Suppose that a morphism $f_\bullet \colon \frg_\bullet \longrightarrow \frh_\bullet$ of diagrams induces weak equivalences of associated graded objects on the level of dg $\RN$-modules.
	Then $f_\bullet$ induces a weak equivalence $\colim_\NN f_\bullet \colon \colim_\NN \frg_\bullet \longrightarrow \colim_\NN \frh_\bullet$ in $L_\infty\Agd_A^\dg$.
\end{lemma}

\begin{proof}
	A morphism $f \colon \frg \longrightarrow \frh$ in \smash{$L_\infty\Agd_A^\dg$} is a weak equivalence if and only if it is a quasi-isomorphism of the underlying dg $A$-modules, since by~\cite[Theorem~3.1]{Nuiten:HoAlg_for_Lie_Algds} the semi-model structure on $L_\infty\Agd_A^\dg$ is transferred along the free-forgetful adjunction
	\begin{equation}
	\begin{tikzcd}
		F : (\Mod_A^\dg)_{/T_A} \ar[r, shift left=  0.15cm, "\perp"' yshift=0.05cm]
		& L_\infty\Agd_A^\dg : U\ . \ar[l, shift left=0.15cm]
	\end{tikzcd}
	\end{equation}
	As the model structure on $\Mod_A^\dg$ is itself transferred from that on $\Mod_\RN^\dg$, it follows  that $f$ is a weak equivalence if and only if it is a quasi-isomorphism of the underlying dg $\RN$-modules.
	
	The poset $\NN$ is filtered, and hence homotopy sifted (see for instance~\cite[Proposition~5.3.1.20]{Lurie:HTT}).
	Thus by~\cite[Theorem~3.3]{Nuiten:HoAlg_for_Lie_Algds} the canonical morphism
	\begin{equation}
		\underset{\NN}{\hocolim}\, (U \circ \frg_\bullet)
		 \longrightarrow U\, \underset{\NN}{\hocolim}\, \frg_\bullet
	\end{equation}
	is a weak equivalence of dg $A$-modules.
	We obtain a commutative square
	\begin{equation}
	\begin{tikzcd}[row sep = 1cm,column sep=1cm]
		\underset{\NN}{\hocolim}\, (U \circ \frg_\bullet) \ar[r, "\sim"] \ar[d, "\sim"']
		& U\, \underset{\NN}{\hocolim}\, \frg_\bullet \ar[d]
		\\
		\underset{\NN}{\colim}\, (U \circ \frg_\bullet) \ar[r, "\sim"']
		& U\, \underset{\NN}{\colim}\, \frg_\bullet
	\end{tikzcd}
	\end{equation}
	showing that the canonical morphism $\hocolim_\NN \frg_\bullet \longrightarrow \colim_\NN \frg_\bullet$ is a weak equivalence in $L_\infty\Agd_A^\dg$, for each Reedy cofibrant functor \smash{$\frg_\bullet \colon \NN \longrightarrow L_\infty\Agd_A^\dg$}.
	That is, the functor $U$ commutes with colimits of diagrams of this type up to weak equivalence.
	Consequently, to check whether $\colim_\NN f_\bullet$ is a weak equivalence in \smash{$L_\infty\Agd_A^\dg$}, it suffices to check whether $\colim_\NN U\, f_\bullet$ is a weak equivalence in \smash{$\Mod_A^\dg$}, or equivalently in \smash{$\Mod_\RN^\dg$}.
	This however is equivalent to $U\, f_\bullet$ inducing quasi-isomorphisms on the associated graded objects.
\end{proof}

\begin{remark}
\label{rmk:detecting weqs on ass grds}
Suppose that $f_\bullet \colon \frg_\bullet \longrightarrow \frh_\bullet$ is as in Lemma~\ref{st:detecting weqs on ass grds}, but suppose that there is a positive integer $N \in \NN$ such that the restriction of $f_\bullet \colon \frg_\bullet \longrightarrow \frh_\bullet$ to $\NN_{\geqslant N}$ is a constant natural transformation between constant diagrams.
Then  \smash{$\colim_\NN \frg_\bullet \cong \frg_N$} and $\colim_\NN \frh_\bullet \cong \frh_N$.
It follows that $\colim_\NN f_\bullet = f_N$, and $f_N \colon \frg_N \longrightarrow \frh_N$ is a weak equivalence in \smash{$L_\infty\Agd_A^\dg$} whenever $f_\bullet$ induces weak equivalences on associated graded objects in degrees $0,1, \ldots, N-1$.
\qen
\end{remark}

\begin{lemma}
For each $p \in \NN$, the functor
\begin{equation}
	Q^{(p)} \colon L_\infty\Agd_A \longrightarrow L_\infty\Agd_A^\dg
\end{equation}
preserves weak equivalences between $A$-cofibrant $L_\infty$-algebroids.
\end{lemma}

\begin{proof}
Given an $A$-cofibrant $L_\infty$-algebroid $\frg$, the sequence of functors $\{ Q^{(p)} \}_{p \in \NN}$ produces a functor
\begin{equation}
	Q^{(\bullet)} \frg \colon \NN \longrightarrow L_\infty\Agd_A^\dg\ ,
	\quad
	p \longmapsto Q^{(p)} \frg\ .
\end{equation}
It follows from Lemma~\ref{st:pushout square for Q^(l)} that this functor is a Reedy cofibrant object of $\Fun(\NN, L_\infty\Agd_A^\dg)$ (see the proof of~\cite[Proposition~5.17]{Nuiten:HoAlg_for_Lie_Algds} for details).
The associated graded of this sequence is~\cite[Remark~5.20]{Nuiten:HoAlg_for_Lie_Algds}
\begin{equation}
\label{eq:ass grd of Q^(*)frg}
	\mathrm{gr} (Q^{(p)} \frg)
	= \Big( F \big( \Sym^p_A (\frg[1]) [-1] \big), \dd_\frg \Big)\ ,
\end{equation}
where \smash{$F \big( \Sym^p_A (\frg[1]) [-1] \big)$} is the free Lie algebroid on the dg $A$-module of degree~$p$ polynomials in $\frg[1]$, endowed with the extension of the differential on $\frg$ by the Leibniz rule.

For each morphism $f \colon \frg \longrightarrow \frh$ in $L_\infty \Agd_A^{(\infty)}$, we obtain a natural transformation
\begin{equation}
	Q^{(\bullet)} f \colon Q^{(\bullet)} \frg \longrightarrow Q^{(\bullet)} \frh
\end{equation}
in $\Fun(\NN, L_\infty\Agd_A^\dg)$.
Suppose that $f$ is a weak equivalence between $A$-cofibrant $L_\infty$-algebroids.
By Lemma~\ref{st:detecting weqs on ass grds} and Remark~\ref{rmk:detecting weqs on ass grds} it suffices to check that $Q^{(\bullet)} f$ induces weak equivalences of dg $\RN$-modules on the associated gradeds in~\eqref{eq:ass grd of Q^(*)frg}.
This follows by the same argument as in the proof of~\cite[Proposition~5.22(a)]{Nuiten:HoAlg_for_Lie_Algds}.
\end{proof}

\begin{remark}
Let $\iota \colon L_\infty\Agd_A^\dg \longhookrightarrow  L_\infty\Agd_A$ be the canonical inclusion functor.
If $\frg$ is an $A$-cofibrant $L_\infty$-algebroid over $A$, each morphism in the chain $Q^{(p)} \frg \longrightarrow Q^{(p+1)} \frg$ is a cofibration in \smash{$L_\infty\Agd_A^\dg$}.
Since $Q^{(1)} \frg = ( F(\frg), \dd_\frg)$ is itself a cofibrant object in \smash{$L_\infty\Agd_A^\dg$}, it follows that $Q^{(p)} \frg$ is cofibrant for each $p \in \NN$.
Thus $Q^{(p)} \circ \iota$ maps $A$-cofibrant $L_\infty$-algebroids over $A$ to cofibrant $L_\infty$-algebroids; in particular, it restricts to the full subcategory \smash{$(L_\infty\Agd_A^\dg)_\cof \subset L_\infty\Agd_A^\dg$} on the cofibrant objects.
\qen
\end{remark}

\begin{proposition}\label{prop:Qp}
The functors
\begin{equation}
	Q^{(p)} \circ \iota \ ,\ Q \circ \iota
	\colon L_\infty\Agd_A^\dg \longrightarrow L_\infty\Agd_A^\dg\ ,
	\quad
	p \in \NN
\end{equation}
induce $\infty$-functors
\begin{equation}
	\overline{Q^{(p)}} \ ,\ \overline{Q} \colon L_\infty\Agd_A^\infty \longrightarrow L_\infty\Agd_A^\infty\ .
\end{equation}
\end{proposition}

\begin{proof}
By the results in this subsection thus far, the functors
\begin{equation}
	Q^{(p)} \circ \iota \ ,\ Q \circ \iota
	\colon L_\infty\Agd_A^\dg \longrightarrow L_\infty\Agd_A^\dg\ ,
	\quad
	p \in \NN
\end{equation}
canonically induce $\infty$-functors
\begin{equation}
	\widetilde{Q^{(p)}} \ ,\ \widetilde{Q} \colon (L_\infty\Agd_A^\dg)_\cof [W^{-1}] \longrightarrow (L_\infty\Agd_A^\dg)_\cof [W^{-1}]\ .
\end{equation}
By~\cite[Theorem~7.5.18]{Cisinski:HCats_HoAlg} the inclusion $\jmath \colon (L_\infty\Agd_A^\dg)_\cof \longhookrightarrow  L_\infty\Agd_A^\dg$ induces a canonical equivalence of $\infty$-categories
\begin{equation}
	\overline{\jmath} \colon (L_\infty\Agd_A^\dg)_\cof [W^{-1}] \longrightarrow L_\infty\Agd_A^\infty\ .
\end{equation}
Thus we obtain canonical $\infty$-functors
\begin{equation}
	\overline{Q^{(p)}} \coloneqq \overline{\jmath} \circ \widetilde{Q^{(p)}} \circ \overline{\jmath}^{-1} \colon L_\infty\Agd_A^\infty \longrightarrow L_\infty\Agd_A^\infty
\end{equation}
as desired, and analogously for $\overline{Q} \coloneqq \overline{\jmath} \circ \widetilde{Q} \circ \overline{\jmath}^{-1}$.
\end{proof}

We can now state the functoriality of our spaces of connections.
Let $M$ be a smooth manifold and $X \in \Sh_\infty(\InfMfd)$ a dt stack.
We obtain a commutative diagram of $\infty$-categories
\begin{equation}
\begin{tikzcd}[column sep={5.5cm,between origins}, row sep={1.75cm,between origins}]
	\textstyle{\int}\, \Con_p(-) \ar[r] \ar[d]
	& \big( L_\infty\Agd_{M}^\infty \big)_{\overline{Q^{(p)}} TM/} \ar[d] \ar[r]
	& \scS_* \ar[d]
	\\
	\Sh_\infty(\InfMfd) (M, X) \ar[r, "T(M/-)"] \ar[r] \ar[rr, bend right = 15, "\Con_p(-)" description]
	& L_\infty\Agd_{M}^\infty \ar[r, "{\Hom(\overline{Q^{(p)}} TM, -)}"]
	& \scS
\end{tikzcd}
\end{equation}

Here the rightmost vertical morphism $\scS_* \longrightarrow \scS$ is the universal left fibration.
The $\infty$-functor
\begin{equation}
	\Hom(\overline{Q^{(p)}} TM, -) \colon L_\infty\Agd_M^\infty \longrightarrow \scS
\end{equation}
classifies the left fibration
\begin{equation}
	\big( L_\infty\Agd_M^\infty \big)_{\overline{Q^{(p)}} TM/}
	 \longrightarrow L_\infty\Agd_M^\infty\ ,
\end{equation}
and so the right-hand square is a pullback square (and also a homotopy pullback square in the Joyal model structure; see for instance~\cite[Corollary~5.3.6]{Cisinski:HCats_HoAlg}).
The left-hand vertical arrow is the left fibration classified by the composition of the bottom horizontal arrows, i.e.~the $\infty$-functor
\begin{align}
	\Con_p(-) \colon \Sh_\infty(\InfMfd) (M, X) & \longrightarrow \scS\ ,
	\\
	(x \colon M \longrightarrow X) &\longmapsto \Con_p (x) = L_\infty\Agd_{M}^\infty \big( \overline{Q^{(p)}} TM,\, T(M/X) \big)\ .
\end{align}
Since the codomain of this left fibration is an $\infty$-groupoid, it is even a Kan fibration.

\begin{example}\label{ex:Bun_nablaMG}
For $X = \rmB G$ the classifying object of a Lie group $G$, we obtain 
\begin{equation}
	\int\, \Con_1(-) \simeq \Bun_\nabla(M; G)\ ,
\end{equation}
the groupoid of principal $G$-bundles with (not necessarily flat) connection on $M$.

One of the main results of the remainder of this paper is an explicit proof that also for $X = \rmB^{n}\rmU(1)$ and any $p \in \NN$, there is a canonical equivalence of \textit{spaces}
\begin{equation}
	\int \Con_p(-) \simeq \Grb^{n-1}_{\nabla|p}(M)\ ,
\end{equation}
where the right-hand side is the space of $(n{-}1)$-gerbes with $p$-form connection on the manifold $M$ (if $p > n$, then $\Grb^{n-1}_{\nabla|p}(M) = \Grb^{n-1}_\flat(M)$ is the space of $(n{-}1)$-gerbes on $M$ with \textit{flat} connection).
\qen
\end{example}

As a consequence of the above and Example~\ref{eg:BG and qcoh defthy}, we can now give a definition of a space of $G_\bullet$-bundles with $p$-form connection from first principles, for \textit{any} group object $G_\bullet$ in $\Sh_\infty(\InfMfd)$ with deformation theory and any manifold $M$.

\begin{definition}
	\label{def: G-buns with l-con}
	Let $M$ be a smooth manifold and $G_\bullet \in \Fun(\bbDelta^\opp, \Sh_\infty(\InfMfd))$  a group object whose underlying object $G=G_1$ in $\Sh_\infty(\InfMfd)$ is a formally smooth dt stack.
	For each $p \in \NN$,  the \textit{$\infty$-groupoid of $G_\bullet$-principal $\infty$-bundles on $M$ with $p$-form connection} is
	\begin{equation}
		\Bun_{\nabla|p}(M; G_\bullet) \coloneqq \int\, \Con_p(-)\ .
	\end{equation}
\end{definition}

\begin{remark}
	\label{rmk: gauge action on l-cons}
	The fact that we have constructed $\Con_p(-)$ as an $\infty$-functor automatically encodes the action of bundle morphisms (in particular gauge transformations) on the space of $p$-form connections on any $G_\bullet$-bundle $P \longrightarrow M$, enabling the study of the higher gauge theory of such connections.
	\qen
\end{remark}


\subsection{Description in terms of stacks}
\label{sec: de rham description}


Let us conclude this section by giving an equivalent description of $p$-form connections in terms of a $p$-th order version of the de Rham stack. We will first treat the case of flat connections, that is, the case where $p=\infty$.

Recall that the de Rham stack $M_{\dR}$ of a smooth manifold $M$ is the stack sending each inf-manifold $U$ to the set of maps of smooth manifolds $U_{\red}\longrightarrow M$. 
This comes with a natural map $\pi\colon M\longrightarrow M_{\dR}$ that sends each $U\longrightarrow M$ to its restriction along $U_{\red}\longhookrightarrow U$.
The \v{C}ech nerve of this map determines a groupoid object in $\Sh_\infty(\InfMfd)$ of the form
\begin{equation}
\label{eq:infinitesimal pair groupoid}
		\begin{tikzcd}
			\cdots \ \arrow[r, yshift=1.5ex] \arrow[r,yshift=0.5ex] \arrow[r,yshift=-0.5ex] \arrow[r, yshift=-1.5ex] & (M\times M\times M)^{\wedge}_M \arrow[r, yshift=1ex] \arrow[r,yshift=-1ex] \arrow[r]
			& (M\times M)^{\wedge}_M \arrow[r, yshift=0.5ex] \arrow[r,yshift=-0.5ex]
			& M\arrow[r, dashed, "\pi"] & M_{\dR} \ .
		\end{tikzcd}
\end{equation}
Here the $n$-th term of the simplicial diagram is the formal completion of $M^{\times n+1}$ at the diagonal. 
For example, $(M\times M)^{\wedge}_M$ is the stack of infinitesimally close pairs of points in $M$: a $U$-point is given by two $U$-points $x, y\colon U\longrightarrow M$ such that $x\big|_{U_{\red}}=y\big|_{U_{\red}}$.

The map $\pi$ is an effective epimorphism in $\Sh_\infty(\InfMfd)$: every map $U_{\red} \longrightarrow M$ lifts to a map $U\longrightarrow M$ since $U_{\red}\longhookrightarrow U$ admits a retraction by Corollary~\ref{cor:inf manifolds retract}. Consequently, $M_{\dR}$ is the quotient of the groupoid object \eqref{eq:infinitesimal pair groupoid}.
A principal $G$-bundle on $M_{\dR}$, i.e.~a map $M_{\dR}\longrightarrow \rmB G$, therefore corresponds to a principal $G$-bundle $P\longrightarrow M$, together with a coherent family of equivalences 
\begin{equation}
x^*P\simeq y^*P
\end{equation}
for all $x, y\colon U\longrightarrow M$ that are infinitesimally close in the sense that $x\big|_{U_{\red}}=y\big|_{U_{\red}}$. 

This kind of infinitesimal ``parallel transport'' is precisely encoded by a flat connection, via

\begin{proposition}\label{prop:flat infinitesimal transport}
Let $M$ be a smooth manifold and $X$ a dt stack, and consider the map
\begin{equation}
\pi^*\colon \Sh_\infty(\InfMfd)(M_{\dR}, X)\longrightarrow \Sh_\infty(\InfMfd)(M, X) \ .
\end{equation}
For each $x\colon M\longrightarrow X$, the fibre of $x$ is naturally equivalent to the space of flat connections on $x$, that is, the space of maps of $L_\infty$-algebroids
\begin{equation}
L_\infty\Agd_M^\infty \big(TM, T(M/X) \big)\simeq L_\infty\Agd_M^\infty \big(Q TM, T(M/X) \big) \ .
\end{equation}
\end{proposition}

\begin{proof}
By Example~\ref{ex:tangent}, the tangent bundle $TM$ is the $L_\infty$-algebroid corresponding under differentiation to $\pi\colon M\longrightarrow M_{\dR}$. The first equivalence then follows from Theorem~\ref{thm:formal stack vs lie}, and the second equivalence follows from Lemma~\ref{lem:lie algebroid resolution}.
\end{proof}

For $p$-form connections, we will use a finite order variant of the groupoid object \eqref{eq:infinitesimal pair groupoid}.

\begin{definition}\label{def:p-th order de Rham stack}
	Let $M$ be a smooth manifold and consider the simplicial diagram
\begin{equation}
\label{eq:p-th order de Rham stack}
		\begin{tikzcd}
			\cdots \ \arrow[r, yshift=1.5ex] \arrow[r,yshift=0.5ex] \arrow[r,yshift=-0.5ex] \arrow[r, yshift=-1.5ex] & (M\times M\times M)^{(p)}_M \arrow[r, yshift=1ex] \arrow[r,yshift=-1ex] \arrow[r]
			& (M\times M)^{(p)}_M \arrow[r, yshift=0.5ex] \arrow[r,yshift=-0.5ex]
			& M\arrow[r, dashed, "\pi"] & M_{\dR, \leqslant p}
		\end{tikzcd}
\end{equation}
where each $(M^{\times n+1})^{(p)}_M$ is the $p$-th order infinitesimal neighbourhood of the diagonal in $M^{\times n+1}$, as in Example \ref{ex:p-th neighbourhood}. The \emph{$p$-th order de Rham stack} of $M$ is the colimit of this simplicial diagram in the $\infty$-category $\DT$ of dt stacks.
\end{definition}

\begin{remark}\label{rem:p-th order de Rham}
	The stack $M_{\dR, \leqslant p}$ is difficult to identify explicitly, because the simplicial diagram \eqref{eq:p-th order de Rham stack} is not very well behaved: it is not a groupoid object, or even an $\infty$-groupoid object in $\Sh_\infty(\InfMfd)$, and the source and target maps \smash{$(M\times M)^{(p)}_M\longrightarrow M$} are not formally smooth.
	In particular, $M_{\dR, \leqslant p}$ cannot be calculated as the colimit in the $\infty$-category $\Sh_\infty(\InfMfd)$.
	
	However,  we will only be interested in maps from $M_{\dR, \leqslant p}$ to a dt stack $X$. The space of such maps is equivalent to the space of maps from the colimit computed in $\Sh_\infty(\InfMfd)$.
	\qen
\end{remark}

\begin{example}\label{ex:p-th order for R}
	Suppose that $M=\RN$. Let $x_0, \dots, x_n$ be the coordinates of $\RN^{n+1}$ and let $t_i=x_i-x_0$ for all $i\geqslant 1$. Then the $p$-th order infinitesimal neighbourhood of the diagonal can be identified with
	\begin{equation}
	\big( \RN^{n+1} \big)^{(p)}_{\RN} \simeq \RN \times \Spec\big( \RN[t_1, \dots, t_n]\,\big/\,(t_1, \dots, t_n)^{p+1}\big) \ .
	\end{equation}
	One can therefore view $(\RN^{n+1})^{(p)}_{\RN}$ as the inf-manifold of $(n{+}1)$-tuples of points $(x_0, x_0+t_1, \dots, x_0+t_n)$ in $\RN$ that are infinitesimally close up to order $p$, in the sense that any $(p{+}1)$-fold product of their differences $t_i$ is zero.
	\qen
\end{example}

Following Example~\ref{ex:p-th order for R}, we can think of the inf-manifold $(M^{\times n+1})^{(p)}_M$ as parametrising $(n{+}1)$-tuples of points in $M$ that are infinitesimally close up to order $p$. 
Just as for the de Rham stack, a principal $G$-bundle on $M_{\dR, \leqslant p}$, i.e.\ a map $M_{\dR, \leqslant p}\longrightarrow \rmB G$, can then be described by the following data:
\begin{myitemize}
\item a principal $G$-bundle $P\longrightarrow M$;
\item for any pair of points $(x, y)\colon U\longrightarrow (M\times M)^{(p)}_M$ that are close up to $p$-th order, an equivalence $\phi_{x, y}\colon x^*P\longrightarrow y^*P$;
\item for any triple of points $(x, y, z)\colon U\longrightarrow (M\times M\times M)^{(p)}_M$ that are close up to $p$-th order, a homotopy $\phi_{y, z}\circ \phi_{x, y}\simeq \phi_{x, z}$;
\item higher coherence data.
\end{myitemize}
In other words, this endows the principal $G$-bundle $P$ with an infinitesimal form of parallel transport up to order $p$. This corresponds precisely to a $p$-form connection on $P$ via

\begin{theorem}\label{thm:p-infinitesimal transport}
Let $M$ be a smooth manifold and $X$ a dt stack, and consider the map
\begin{equation}
\pi^*\colon \Sh_\infty(\InfMfd)(M_{\dR, \leqslant p}, X)\longrightarrow \Sh_\infty(\InfMfd)(M, X) \ .
\end{equation}
For each $x\colon M\longrightarrow X$, the fibre of $x$ is naturally equivalent to the space $\Con_p(x)$ of $p$-form connections on $x$.
\end{theorem}

The proof makes use of Lemma~\ref{lem:formal completion explicit} below, applied to the case where $N=M^{\times n+1}$.

\begin{lemma}\label{lem:formal completion explicit}
Let $i\colon M\longhookrightarrow N$ be a closed embedding of smooth manifolds. Then there is an equivalence of dt stacks
\begin{equation}
\underset{p\in\NN}\colim \, N^{(p)}_M\simeq N^{\wedge}_M
\end{equation}
where the (sequential) colimit is computed in $\Sh_\infty(\InfMfd)$. Consequently, there is an equivalence of function algebras $$\Gamma(N^\wedge_M, \mathcal{O}_{N^\wedge_M})\simeq \underset{p\in\NN}\lim\, C^\infty(N)\,\big/\,I^{p+1} \ , $$ where $I$ is the kernel of $i^*\colon C^\infty(N)\longtwoheadrightarrow C^\infty(M)$.
\end{lemma}

\begin{proof}
	Recall from Example~\ref{ex:p-th neighbourhood} that there are natural maps $N^{(p)}_M\longrightarrow N$, given at the level of the underlying smooth manifolds by the inclusion $i\colon M\longhookrightarrow N$. These give rise to natural maps $N^{(p)}_M\longrightarrow N^{\wedge}_M$ and hence to a natural map 
	\begin{equation}
	f\colon \underset{p\in\NN}\colim\, N^{(p)}_M\longrightarrow N^{\wedge}_M \ .
	\end{equation}
	To see that $f$ is an equivalence of stacks, it suffices to verify that it is an equivalence locally on $N$ and $M$. We thus reduce to the case where the closed embedding is given by $i\colon M=\RN^m\times\{0\}\longhookrightarrow \RN^{m+k}=N$.
	
	In this case, let us write $x_1, \dots, x_m, y_1, \dots, y_k$ for the coordinates on $N=\RN^{m+k}$. Then there are two kinds of infinitesimal neighbourhoods of $M$ in $N$, given by
	\begin{align*}
	N^{(p)}_M &= \Spec \big( C^\infty(\RN^{m+k})\,\big/\,(y_1, \dots, y_k)^{p+1}\big) \ , \\[4pt]
	U^{(p)} &= \Spec \big( C^\infty(\RN^{m+k})\,\big/\,(y_1^{p+1}, \dots, y_k^{p+1})\big) \ .
	\end{align*}
	For each $p$, there are natural maps $N^{(p)}_M\longhookrightarrow U^{(p)}\longhookrightarrow N^{(k\,p)}_M$. Therefore there is an equivalence of colimits
	\begin{equation}
	\underset{p\in\NN}\colim\, N^{(p)}_M\simeq \underset{p\in\NN}\colim\, U^{(p)} 
	\end{equation}
	and it is enough to show that $f\colon \colim_{p\in\NN}\, U^{(p)}\longrightarrow N^\wedge_M$ is an equivalence.
	
	Using Proposition~\ref{prop:j!}, it follows that the domain and codomain of this map are both qcdt stacks whose underlying reduced stacks are equivalent to $M$. 
	By Proposition~\ref{prop:inverse_function} and Remark~\ref{rmk: detecting equivalences for qcdt}, it therefore suffices to show that the induced map
	\begin{equation}
	\underset{p\in\NN}\colim\, T_xU^{(p)}\simeq T_x\Big(\underset{p\in\NN}\colim\, U^{(p)}\Big) \longrightarrow T_{f(x)} N^\wedge_M \simeq i^*TN=i^*T\RN^{m+k}
	\end{equation}
	is an equivalence, where $x\colon M\longrightarrow U^{(p)}$ is the canonical map.
	
	The tangent complex of $U^{(p)}$ can be computed explicitly, using the cofibrant replacement of the corresponding $C^\infty$-ring given by (this would be more complicated for \smash{$N^{(p)}_M$})
	\begin{equation}
	B^{(p)} = C^\infty(\RN^{m+k})[\eta_1, \dots, \eta_k] \qquad,\qquad |\eta_i|=-1 \ , \quad \dd\eta_i= y_i^{p+1} \ .
	\end{equation}
	The map $U^{(p)}\longrightarrow U^{(p+1)}$ then corresponds to the map $B^{(p+1)}\longrightarrow B^{(p)}$ that is the identity in degree $0$ and sends $\eta_i$ to $y_i\,\eta_i$. In addition, the map $x\colon M\longrightarrow U^{(p)}$ corresponds to the map	$B^{(p)}\longrightarrow C^\infty(\RN^m)$ sending $y_i$ and $\eta_i$ to zero for $i=1,\dots,k$.
	
	As in Examples~\ref{ex:corep tangent} and~\ref{ex:mfd are qcdt}, the tangent complex $T_xU^{(p)}$ can now be computed in terms of the $C^\infty$-algebraic cotangent complex $\Omega^1_{B^{(p)}}$. The latter is generated by $\dd x_j$, $\dd y_i$ and $\dd\eta_i$ with $1\leqslant j\leqslant m$ and $1\leqslant i\leqslant k$, so that
	\begin{equation}
	T_xU^{(p)} = \Hom_{B^{(p)}}\big(\Omega^1_{B^{(p)}}, C^\infty(\RN^m) \big) = \big(i^*T\RN^{m+k} \xrightarrow{ \ 0 \ } C^\infty(\RN^m)^{\oplus k}\big) \ .
	\end{equation}
	Here the degree $1$ part is generated by elements $\partial_{\eta_i}$ dual to $\dd\eta_i$ for $i=1,\dots,k$. The map $T_xU^{(p)}\longrightarrow T_xU^{(p+1)}$ is the identity in degree $0$ and is given in degree $1$ by
	\begin{equation}
	\partial_{\eta_i} \longmapsto y_i\,\partial_{\eta_i}+\eta_i\,\partial_{y_i}=0 \ .
	\end{equation}
	Here we used that \smash{$\Omega^1_{B^{(p+1)}}\longrightarrow \Omega^1_{B^{(p)}}$} sends $\dd\eta_i$ to $\dd(y_i\,\eta_i)=y_i\,\dd\eta_i+\eta_i\, \dd y_i$, and the fact that $y_i$ and $\eta_i$ are both sent to zero in $C^\infty(\RN^m)$.
	Taking the colimit, we conclude that $\colim_{p\in\NN}\, T_xU^{(p)}\simeq i^*T\RN^{m+k}$.
\end{proof}

\begin{proof}[Proof of Theorem~\ref{thm:p-infinitesimal transport}]
	Our goal is to prove that $M_{\dR, \leqslant p}$ is the image of the $L_\infty$-algebroid $Q^{(p)}TM$ under the formal integration functor from Theorem~\ref{thm:formal stack vs lie}. It then follows from Theorem \ref{thm:formal stack vs lie} that for any map $x\colon M\longrightarrow X$ to a dt stack, there is an equivalence between maps $M_{\dR, \leqslant p}\longrightarrow X$ under $M$ and maps of $L_\infty$-algebroids $Q^{(p)}TM\longrightarrow T(M/X)$, that is, $p$-form connections on $x$.
	
	To calculate the formal integration of $Q^{(p)}TM$, let us consider the diagram of $L_\infty$-algebroids given by
	\begin{equation}
\label{eq:tangent of pair}
\begin{tikzcd}[row sep = 1cm]
\cdots \ \arrow[r, yshift=1.5ex] \arrow[r,yshift=0.5ex] \arrow[r,yshift=-0.5ex] \arrow[r, yshift=-1.5ex] & Q^{(p)}\frg_2 \arrow[r, yshift=1ex] \arrow[r,yshift=-1ex] \arrow[r]\arrow[d]
			& Q^{(p)}\frg_1\arrow[d] \arrow[r, yshift=0.5ex] \arrow[r,yshift=-0.5ex]
			& Q^{(p)}\frg_0\simeq 0\arrow[r]\arrow[d] & Q^{(p)}TM\arrow[d]\\
			\cdots \ \arrow[r, yshift=1.5ex] \arrow[r,yshift=0.5ex] \arrow[r,yshift=-0.5ex] \arrow[r, yshift=-1.5ex] & \frg_2 \arrow[r, yshift=1ex] \arrow[r,yshift=-1ex] \arrow[r]
			& \frg_1 \arrow[r, yshift=0.5ex] \arrow[r,yshift=-0.5ex]
			& \frg_0=0\arrow[r] & TM
		\end{tikzcd}
\end{equation}
Here the bottom row is the \v{C}ech nerve of the map $0\longrightarrow TM$ in the $\infty$-category of $L_\infty$-algebroids, so that $TM\simeq \colim_{\NN_0}\, \frg_\bullet$. The top row is the image of the bottom row under the $\infty$-functor $Q^{(p)}$. It follows from \cite[Proposition 5.22(b)]{Nuiten:HoAlg_for_Lie_Algds} that $Q^{(p)}TM\simeq \colim_{\NN_0}\, Q^{(p)}\frg_\bullet$. Because the formal integration functor from Theorem~\ref{thm:formal stack vs lie} is a left adjoint, it therefore suffices to prove that the formal integration of $Q^{(p)}\frg_\bullet$ is naturally equivalent to the simplicial diagram of inf-manifolds
\[
(M^{\times \bullet+1})^{(p)}_M = \Spec\big(C^\infty(M^{\times \bullet+1})\,\big/\,I_\bullet^{p+1}\big) \ ,
	\]
	where $I_n$ denotes the kernel of the restriction along the diagonal $C^\infty(M^{\times n+1})\longtwoheadrightarrow C^\infty(M)$. 
	
	As a first step to identifying the formal integration of $Q^{(p)}\frg_\bullet$, observe that the formal integration of $\frg_\bullet$ is the \v{C}ech nerve \eqref{eq:infinitesimal pair groupoid} of the map $\pi\colon M\longrightarrow M_{\dR}$: this follows from the fact that Lie differentiation, being a right adjoint, sends the \v{C}ech nerve of $\pi$ to the \v{C}ech nerve of $0\longrightarrow TM$. Each $\frg_n$ can then be identified explicitly using Example~\ref{ex: Lie diff group under M}: it is quasi-isomorphic to the shifted normal bundle of the diagonal $M\longhookrightarrow M^{\times n+1}$ with trivial bracket and anchor. 
	This normal bundle is isomorphic to a Whitney sum of tangent bundles of $M$, so 
	\[
	\frg_n\simeq TM^{\oplus n}[-1]\ .
	\]
	Let us highlight that this equivalence is \emph{not natural}: it relies on a choice of tubular neighbourhood and is hence not compatible with the simplicial structure maps.
	
	Using Proposition~\ref{prop:Qp}, it follows that each $Q^{(p)}\frg_n\simeq Q^{(p)}TM^{\oplus n}[-1]$ is (the retract of) a cofibrant $L_\infty$-algebroid with finitely many generators in degrees $\geqslant 1$. By (the proof of) Theorem \ref{thm:formal stack vs lie}, the formal integration of such $L_\infty$-algebroids are inf-manifolds: by construction, their integration is given by the spectrum of their Chevalley--Eilenberg algebras
	\[
	\Spec\big(\ChEil^*(Q^{(p)}\frg_n)\big)\simeq \Spec\big(\ChEil^{\leqslant p}(\frg_n)\big) \ .
	\]
	We are thus left with the purely algebraic problem of giving a natural equivalence of dg $C^\infty$-rings $\ChEil^{\leqslant p}(\frg_\bullet)\simeq C^\infty(M^{\times \bullet+1})/I_\bullet^{p+1}$.
	
	To see this, first note that $\frg_n\simeq TM^{\oplus n}[-1]$ induces an equivalence of diagrams of dg $C^\infty$-rings between
	\[\begin{tikzcd}
	\ChEil^*(\frg_n)\arrow[r] & \ChEil^{\leqslant p}(\frg_n)\arrow[r] & C^\infty(M)
	\end{tikzcd}\]
	and
	\[\begin{tikzcd}
	\ChEil^*(TM^{\oplus n}[-1])\arrow[r] & \ChEil^{\leqslant p}(TM^{\oplus n}[-1])\arrow[r] & C^\infty(M) \ .
	\end{tikzcd}\]
	This shows that up to non-natural quasi-isomorphism, $\ChEil^*(\frg_n)$ is a discrete $C^\infty$-ring and $\ChEil^{\leqslant p}(\frg_n)$ is its quotient by the $(p{+}1)$-th power of the kernel of the map to $C^\infty(M)$. But there is also  a \emph{natural} model for the $C^\infty$-ring $\ChEil^*(\frg_n)$: by Theorem~\ref{thm:formal stack vs lie} and Lemma~\ref{lem:formal completion explicit}, there is a natural equivalence
	\[
	\ChEil^*(\frg_\bullet)\simeq \Gamma\big( (M^{\times \bullet+1})^\wedge_M, \mathcal{O}_{(M^{\times \bullet+1})^\wedge_M}\big) \simeq \underset{p\in\NN}\lim\, C^\infty(M^{\times \bullet+1})\,\big/\,I_\bullet^{p+1}
	\]
	of dg $C^\infty$-rings with a map to $C^\infty(M)$. The kernel of $\lim_{p\in\NN}\, C^\infty(M^{\times \bullet+1})/I_\bullet^{p+1}\longrightarrow C^\infty(M)$ is precisely the image of the ideal $I_\bullet$. Taking the quotient by the $(p{+}1)$-th power of this ideal, we therefore obtain the desired natural equivalence $\ChEil^{\leqslant p}(\frg_\bullet)\simeq C^\infty(M^{\times \bullet+1})/I_\bullet^{p+1}$.
\end{proof}


\section{$\infty$-groupoids of connections on higher $\rmU(1)$-bundles}
\label{sec: cons on higher U(1)-bundles}


The goal of this section is to prove that the notion of higher-form connections we introduced in Definition~\ref{def: Connection} reproduces the entire $\infty$-groupoids of connections on higher gerbes, i.e.~$\rmB^n \rmU(1)$-principal $\infty$-bundles.
While the connections on a classical principal bundle form a set, it has been observed in~\cite{BS:Higher_Syms_and_Deligne_Coho} that the connections on a higher bundle form a higher groupoid.
Below we start by computing two equivalent presentations of the Atiyah $L_\infty$-algebroid of any given higher gerbe (Corollary~\ref{st:At-L_oo-agd of n-gerbe from first principles} and Theorem~\ref{st:strictification of At(CG)}).
We then recall the definition of the space of $p$-form connections on an $(n{-}1)$-gerbe via \v{C}ech--Deligne cocycles and state the main result of this section in Theorem~\ref{st:finite and infinitesimal l-cons on n-gerbes}:
The \v{C}ech--Deligne picture of $p$-form connections on $(n{-}1)$-gerbes is equivalent to that following from Definition~\ref{def: Connection}.
This is remarkable because the former makes no reference to derived geometry, and so is entirely independent of the geometric methods and approach to connections we developed in the first part of this paper.


\subsection{Atiyah $L_\infty$-algebroids of higher $\rmU(1)$-bundles}
\label{subsec:Atiyahalgebroids}


Our goal will now be to calculate the $L_\infty$-algebroid associated to a higher $\rmU(1)$-bundle $$x \colon M \longrightarrow \rmB^n \rmU(1)$$ on a smooth manifold $M$.
Such a bundle is also called an \textit{$(n{-}1)$-gerbe} on $M$.
We will carry out the computation in a few steps, using the fact that the bottom-right square in the  diagram
\begin{equation}
	\label{diag:pullback U(1)-Z-R}
	\begin{tikzcd}[row sep = 1cm , column sep = 1cm]
		M \arrow[rd, "x"] \arrow[rrd, bend left=30] \arrow[rdd, bend right=30, "x_{\mathrm{top}}"{description}]
		& &
		\\
		& \rmB^n \rmU(1) \arrow[r] \arrow[d]
		& \ast\arrow[d, "0"]
		\\
		& \rmB^{n+1} \ZN \arrow[r]
		& \rmB^{n+1} \RN
	\end{tikzcd}
\end{equation}
 is a homotopy pullback square of smooth stacks. Here $$[x_{\mathrm{top}}] \ \in \ \pi_0\, \Sh_\infty(\Mfd) (M, \rmB^{n+1}\ZN)\cong \rmH^{n+1}(M; \ZN)$$ is the integral cohomology class classifying the $(n{-}1)$-gerbe. 

Theorem~\ref{thm:formal stack vs lie} implies that taking the associated $L_\infty$-algebroids preserves homotopy pullbacks.
Our strategy will therefore be to compute the $L_\infty$-algebroids associated to the three corners of the square \eqref{diag:pullback U(1)-Z-R}, as well as the morphisms between them, and then compute the homotopy pullback in the $\infty$-category of $L_\infty$-algebroids.
Let us point out that the homotopy pullback diagram we will construct is only going to live in the $\infty$-category $L_\infty\Agd_{M}^\dg[W^{-1}] = L_\infty\Agd_{M}^\infty$; its morphisms will be presented by $\infty$-morphisms of $L_\infty$-algebroids.

We start by identifying the $L_\infty$-algebroid associated to the bottom-right object in \eqref{diag:pullback U(1)-Z-R}.

\begin{proposition}\label{prop:differentiation of B^nR}
	Let $0\colon M\longrightarrow \rmB^{n+1}\RN$ be the zero map. The associated $L_\infty$-algebroid is given by
	\begin{equation}
		T(M/\rmB^{n+1}\RN) \simeq TM \oplus C^ \infty(M)[n]\ ,
	\end{equation}
	with the only non-trivial brackets given by the usual commutator bracket on $TM$ and by $[X, f] = \pounds_X(f)$ for $X\in TM$ and $f\in C^\infty(M)$.
\end{proposition}

\begin{proof}
	Let $\frg$ be a complete $L_\infty$-algebroid over $M$ and let $M\longrightarrow Y$ be the associated map to a qcdt stack (via Theorem~\ref{thm:formal stack vs lie}). Then there are natural equivalences of spaces
	\begin{align*}
		L_\infty\Agd_M^{\cpl, \infty} \big( \frg, T(M/\rmB^{n+1} \RN) \big)
		&\simeq \Sh_\infty(\InfMfd)_{M/}\big(Y, \rmB^{n+1} \RN \big)
		\\[4pt]
		& \simeq \DK \Big(\hofib\big(\Gamma(Y, \calO_Y)\longrightarrow C^\infty(M) \big) [n+1] \Big)
		\\[4pt]
		&\simeq \DK \Big(\hofib\big(\ChEil^*(\frg)\longrightarrow C^\infty(M) \big) [n+1]\Big)\ .
	\end{align*}
	Here the first and last equivalences follow from Theorem~\ref{thm:formal stack vs lie}. The second equivalence follows from Examples~\ref{ex:function algebras} and~\ref{ex:function algebras and maps to B^n RN}.
	
	The last homotopy fibre can be identified with the \emph{reduced} Chevalley--Eilenberg complex $\ol{\ChEil}{}^*(\frg)$, and so the last term corresponds under the Dold--Kan correspondence to the  cochain complex
		$$
				\begin{tikzcd}
				Z^{n+1} \ol{\ChEil}{}^*(\frg)
				& \ol{\ChEil}{}^{n}(\frg) \arrow[l]
				& \ol{\ChEil}{}^{n-1}(\frg) \arrow[l]
				& \cdots \ . \arrow[l]
			\end{tikzcd}
		$$
		Unraveling the definition of $\infty$-morphisms of $L_\infty$-algebroids, we see that the degree~$(n{+}1)$ cocycles in the reduced Chevalley--Eilenberg complex of an $L_\infty$-algebroid $\frg$ correspond precisely to the $\infty$-morphisms $\frg\longleadsto TM\oplus C^\infty(M)[n]$. Combining this with the simplicial resolution from Theorem~\ref{st:fib res of fibrant dg abelian extensions}, one finds that 
		\begin{equation}
			\DK \big( \ol{\ChEil}{}^*(\frg)[n+1] \big)
			\simeq L_\infty\Agd_M^{\cpl, \infty} \big( \frg, TM\oplus C^\infty(M)[n]\big)
		\end{equation}
		is a model for the $\infty$-categorical mapping space we are interested in. It then follows from the Yoneda lemma that $T(M/\rmB^{n+1}\RN)\simeq TM\oplus C^\infty(M)[n]$.
\end{proof}

\begin{lemma}
	\label{st:T(M/B^(n+1)ZN) and T(*/M)}
	The $L_\infty$-algebroids associated to $M\longrightarrow \ast$ and $M\longrightarrow \rmB^{n+1}\ZN$ each coincide with the tangent bundle $TM$.
\end{lemma}

\begin{proof}
	We have already proven the first case in Example~\ref{ex:tangent}. For the second case, notice that $\rmB^{n+1}\ZN$ is a discrete, or constant, stack (we will expand on this further  in the paragraphs following this proof). In particular, $\rmB^{n+1}\ZN(U)=\rmB^{n+1}\ZN(U_{\red})$ for any inf-manifold $U$. It readily follows that the formal moduli problem associated to $M\longrightarrow \rmB^{n+1}\ZN$ is the terminal one, so that the corresponding $L_\infty$-algebroid is given by $TM$.
\end{proof}

Consequently, the homotopy pullback square \eqref{diag:pullback U(1)-Z-R} gives rise to a homotopy pullback square of $L_\infty$-algebroids
$$
	\begin{tikzcd}[row sep = 1cm , column sep = 1cm]
		T \big( M/\rmB^n\rmU(1) \big) \arrow[r, rightsquigarrow] \arrow[d, rightsquigarrow]
		& TM \arrow[d, rightsquigarrow]
		\\
		TM\arrow[r, rightsquigarrow]
		& TM\oplus C^\infty(M)[n]
	\end{tikzcd}
	$$
Note that this diagram a priori only lives in the $\infty$-category $L_\infty\Agd_M$. By our computation of the mapping spaces of this $\infty$-category in terms of $\infty$-morphisms, we can present the bottom and right map by $\infty$-morphisms of $L_\infty$-algebroids. It remains to identify these two $\infty$-morphisms more explicitly and then calculate the homotopy pullback.

The proof of Proposition~\ref{prop:differentiation of B^nR} shows that there is a bijection
$$
	\Ho(L_\infty \Agd_M^\dg) \big( TM,\, TM \oplus C^\infty(M)[n] \big)
	\cong \overline{\rmH}{}^{n+1}_\dR(M)
$$
between the set of $\infty$-morphisms $TM\longleadsto TM\oplus C^\infty(M)[n]$ and the set of degree~$(n{+}1)$ reduced de~Rham cohomology classes (i.e.~the cohomology of $\Omega^{\geqslant 1}(M)$). Therefore it will now suffice to identify the two de Rham cohomology classes associated to the bottom and vertical map in \eqref{diag:pullback U(1)-Z-R}.

To achieve this we will need to take a slight formal detour. First, let $\Gamma\colon\Sh_\infty(\InfMfd)\longrightarrow \scS$ be the global sections functor, evaluating a sheaf at the final object $*$ (represented by the one-point smooth manifold). This functor admits a fully faithful left adjoint $\mathrm{cst}$, which takes \emph{constant sheaves}, and this in turn has a further left adjoint $\Pi_\infty$ taking \textit{shapes}, or \textit{homotopy types}.
Altogether we obtain a double adjunction 
\begin{equation}
	\begin{tikzcd}[column sep = 2.5cm]
		\Sh_\infty(\InfMfd) \ar[r, bend left = 30, "\Pi_\infty" {description, pos = 0.62}, "\perp"' {pos = 0.67, yshift = -0.075cm}] \ar[r, bend right = 30, "\Gamma" {description, pos = 0.61}]
		& \scS \ar[l, "\mathrm{cst}" description, "\perp" yshift = -0.1cm]
	\end{tikzcd}
\end{equation}

Using the adjoint pairs \eqref{diag:infinitesimal cohesion}, these functors can also be described as follows: for any homotopy type $K$, the $\infty$-sheaf $\mathrm{cst}(K)$ is given by $\mathrm{cst}(K)(U)=\mathrm{cst}(K)(U_{\red})$, where on the right hand side we take the constant sheaf on $K$ on the site of ordinary smooth manifolds. Likewise, $\Pi_\infty(X)=\Pi_\infty(X_{\red})$ is the shape of the reduction of $X$, i.e.~the shape of its restriction to the site of smooth manifolds. In other words, both functors $\mathrm{cst}$ and $\Pi_\infty$ are insensitive to any derived or infinitesimal structure.

By a slight abuse of notation, we will denote the composition $\mathrm{cst} \circ \Pi_\infty$ again by $\Pi_\infty$; it will be clear from the context whether the output produced by the functor is a single space or a sheaf of spaces on $\InfMfd$.
Any object $X \in \Sh_\infty(\InfMfd)$ comes naturally endowed with a map $X \longrightarrow \Pi_\infty(X)$ with the following universal property: for any map $X\longrightarrow \mathrm{cst}(K)$ to a constant sheaf, there exists a unique map
$$
	\begin{tikzcd}[column sep = 0.5cm , row sep = 1cm]
		& X \arrow[ld] \arrow[rd]
		\\
		\Pi_\infty(X) \arrow[rr, dashed] & & \mathrm{cst}(K)
	\end{tikzcd}
$$
as indicated.

\begin{lemma}
	Let $X\in\Sh_\infty(\InfMfd)$, and consider the natural map $X\longrightarrow \Pi_\infty(X)$. The natural map $\Pi_\infty(X)^\wedge_X\longrightarrow \ast^\wedge_X=X_\dR$ of formal completions (Definition~\ref{def: de Rham stack}) is an equivalence of sheaves.
\end{lemma}

\begin{proof}
	Since $\Pi_\infty(X)$ is a constant sheaf, it satisfies  $\Pi_\infty(X)(U)\simeq \Pi_\infty(X)(U_{\red})$ for all inf-manifolds $U$. We then compute
	$$
	\Pi_\infty(X)^\wedge_X(U)=\Pi_\infty(X)(U)\times_{\Pi_\infty(X)(U_\red)} X(U_{\red}) \simeq X(U_{\red})\simeq X_\dR(U)\ ,
	$$
	as required.
\end{proof}

In particular, for any sheaf $X\in\Sh_\infty(\InfMfd)$ there is a natural sequence of maps
$$
	X \longrightarrow X_\dR \longrightarrow \Pi_\infty(X)\ .
$$

\begin{lemma}
	\label{lem:poincare}
	Let $M$ be a smooth manifold. The restriction along the canonical map $M_\dR\longrightarrow \Pi_\infty(M)$ induces a natural quasi-isomorphism of function algebras
	$$
		\Gamma(\Pi_\infty(M), \calO_{\Pi_\infty(M)}) \eq \Gamma(M_\dR, \calO_{M_\dR})\ .
	$$
	Consequently, the restriction map
	\begin{equation}
		\Sh_\infty(\InfMfd)(\Pi_\infty(M), \rmB^{n+1}\RN) 
		 \longrightarrow \Sh_\infty(\InfMfd)(M_\dR, \rmB^{n+1}\RN)
	\end{equation}
	is an equivalence for all $n \geqslant 0$.
\end{lemma}

This is a version of the de Rham Theorem relating singular and de Rham cohomology.

\begin{proof}
Both $X\longmapsto X_\dR$ and $X\longmapsto \Pi_\infty(X)$ preserve homotopy colimits of $\infty$-sheaves, and $\Gamma(-, \calO)$ sends such homotopy colimits to homotopy limits. Since every smooth manifold can be written as a homotopy colimit of an open cover by contractibles, this reduces the assertion to the case where $M$ is contractible. But in that case $\Pi_\infty(M)\simeq \ast$, so that $\Gamma(\Pi_\infty(M), \calO_{\Pi_\infty(M)})\simeq \RN$ is the initial cdga, while $\Gamma(M_\dR, \calO_{M_\dR})\simeq \Omega^*(M)\simeq \RN$ by the usual de Rham Theorem and Example~\ref{ex:tangent}. In particular, any map between these cdgas is a quasi-isomorphism.
	The second part of the lemma follows from the first part and Example~\ref{ex:function algebras and maps to B^n RN}.
\end{proof}

For a smooth manifold $M$ and $p \in \NN$, we write $\Omega^p_\cl(M)$ for the abelian group of closed $p$-forms on $M$.

\begin{proposition}
	\label{prop:pullback diagram explicit}
	Let $M$ be a smooth manifold and consider the diagram
	$$
	\begin{tikzcd}[column sep = 1cm, row sep = 1cm]
		M \arrow[r] \arrow[d, "x_{\mathrm{top}}"']
		& \ast \arrow[d]
		\\
		\rmB^{n+1}\ZN \arrow[r]
		& \rmB^{n+1}\RN
	\end{tikzcd}
	$$
 of $\infty$-sheaves.	The associated diagram in $L_\infty\Agd_M$ is given (up to homotopy) by
	$$\begin{tikzcd}[row sep = 1cm , column sep = 1cm]
		& TM \arrow[d, rightsquigarrow, "0"]\\
		TM \arrow[r, rightsquigarrow, "H"] & TM \oplus C^\infty(M)[n]
	\end{tikzcd}$$
	where the right vertical $\infty$-morphism corresponds to $0 \in \Omega_\cl^{n+1}(M)$ and the bottom $\infty$-morphism corresponds to any closed $(n{+}1)$-form $H \in \Omega_\cl^{n+1}(M)$ representing the integral cohomology class $[x_\mathrm{top}]$.
\end{proposition}

\begin{proof}
	The left-bottom part of the diagram can be further decomposed as the sequence
	$$
		\begin{tikzcd}
			M \arrow[r]
			& M_\dR \arrow[r]
			& \Pi_\infty(M) \arrow[r]
			& \rmB^{n+1} \ZN \arrow[r]
			& \rmB^{n+1} \RN\ .
		\end{tikzcd}
	$$
	Since by assumption the map $\Pi_\infty(M)\longrightarrow \rmB^{n+1}\ZN$ classifies the cohomology class $[x_{\mathrm{top}}]$, Lemma~\ref{lem:poincare} implies that the homotopy class of the composite map $M_\dR\longrightarrow \rmB^{n+1}\RN$ is classified by the de Rham cohomology class associated to $x_{\mathrm{top}}$. 
	
	Now note that there are isomorphisms
	\begin{align}
		&\pi_0\, L_\infty\Agd_M^\infty \big( T(M/\rmB^{n+1}\ZN), \ T(M/\rmB^{n+1}\RN)\big)
		\\[4pt]
		&\hspace{2cm} \cong \ol{\rmH}{}^{n+1}_{\dR}(M)
		 \cong \pi_0 \Big( \Sh_\infty(\InfMfd) \big( M_\dR, \rmB^{n+1} \RN \big) \times^h_{\Sh_\infty(\InfMfd)(M, \rmB^{n+1}\RN)} \{0\} \Big)\ ,
	\end{align}
	where the middle term is the reduced de Rham cohomology (i.e.\ the cohomology of $\Omega^{\geqslant 1}(M)$---in these degrees there is no difference). The first isomorphism is the isomorphism appearing in the proof of Proposition~\ref{prop:differentiation of B^nR} (using Lemma~\ref{st:T(M/B^(n+1)ZN) and T(*/M)}), and the second isomorphism follows from Example~\ref{ex:tangent}. 	The induced morphism 
	$$
		\begin{tikzcd}
			TM = T(M/\rmB^{n+1}\ZN)\arrow[r, rightsquigarrow]
			& T(M/\rmB^{n+1}\RN) = TM \oplus C^\infty(M)[n]
		\end{tikzcd}
	$$
	can therefore be described by any representative $H$ of the image of $[x_{\mathrm{top}}]$ in de Rham cohomology.
	
	For the right-hand vertical map one argues similarly, but in this case one simply needs a representative for the \emph{zero} cohomology class.
\end{proof}

This now leads to the first main result of this section, namely an identification of the Atiyah $L_\infty$-algebroid $T(M/\rmB^n\rmU(1))$ through

\begin{corollary}
	\label{st:At-L_oo-agd of n-gerbe from first principles}
	Let $M$ be a smooth manifold and $x\colon M\longrightarrow \rmB^n\rmU(1)$ a morphism of stacks (i.e.~the classifying map of an $(n{-}1)$-gerbe on $M$). The associated $L_\infty$-algebroid is given by
	\begin{equation}
		T \big( M / \rmB^n\rmU(1) \big)
		\simeq \big( TM \oplus C^{\infty}(M)[n-1],\, H \big) \ ,
	\end{equation}
	whose underlying complex is $TM \oplus C^{\infty}(M)[n-1]$ with zero differential, whose anchor map is the canonical projection onto $TM$, whose binary bracket is given by the Lie bracket on vector fields and by the Lie derivative $[X, f] =\pounds_X(f)$ for all $X\in TM$ and $f\in C^{\infty}(M)[n-1]$, and whose $(n{+}1)$-ary bracket is given on $X_1,\dots,X_{n+1}\in TM$ by
	\begin{equation}
		[X_1, \dots, X_{n+1}]
		= H(X_1, \dots, X_{n+1})\ ,
	\end{equation}
	for any closed $(n{+}1)$-form $H \in \Omega^{n+1}_\cl(M)$ on $M$ which represents the cohomology class in $\rmH^{n+1}(M, \ZN)$ defined by $x$.
	All other brackets are zero.
\end{corollary}

\begin{proof}
	One has to compute the homotopy pullback of the diagram of $L_\infty$-algebroids and $\infty$-morphisms appearing in Proposition~\ref{prop:pullback diagram explicit}. This can be done explicitly, cf.~\cite{Nuiten:HoAlg_for_Lie_Algds}.
\end{proof}

The Atiyah $L_\infty$-algebroid $T(M/\rmB^n\rmU(1))$ in Corollary~\ref{st:At-L_oo-agd of n-gerbe from first principles} depends on the choice of a closed $(n{+}1)$-form $H$, and we abbreviate
\begin{align}
\At_{n-1}(H):=\big( TM \oplus C^{\infty}(M)[n-1],\, H \big) \ .
\end{align}
The form $H$ can be thought of as the curvature form of a full $n$-form connection on the $(n{-}1)$-gerbe classified by the morphism $x \colon M \longrightarrow \rmB^n \rmU(1)$ (see Example~\ref{ex:fullnform} below).
We conclude this section by providing an equivalent model of $T(M/\rmB^n\rmU(1))$ as a dg Lie algebroid.
It will not depend on a choice of an $(n{+}1)$-form, but only on a presentation of the abstract morphism $x \colon M \longrightarrow \rmB^n\rmU(1)$ via a choice of a good open covering $\CU$ of $M$ and a $\rmU(1)$-valued \v{C}ech $n$-cocycle $g$ subordinate to $\CU$.

\begin{remark}
Let $\Cart$ be the category of cartesian spaces $\RN^k$,
with $k\in\NN_0$, and smooth maps.	The morphism $x \colon M \longrightarrow \rmB^n\rmU(1)$ can be presented by \v{C}ech data as follows:
	there are canonical equivalences of $\infty$-categories
	\begin{equation}
		\Sh_\infty(\Mfd) \simeq \Sh_\infty(\Cart) \simeq L_W \Fun(\Cart^\opp, \sSet)\ ,
	\end{equation}
	where the first equivalence is induced by the canonical inclusion $\Cart \longhookrightarrow  \Mfd$ (see for instance~\cite[Theorem~B.8]{Bunk_R-loc_HoThy} for a presentation, or~\cite[Corollary~A.5.6]{ADH:Differential_Cohomology} for an $\infty$-categorical proof) and the last $\infty$-category is the $\infty$-localisation of the local projective model structure on $\Fun(\Cart^\opp, \sSet)$.
	
	If we let $\CU = \{U_a\}_{a \in \varLambda}$ be a good open covering of the smooth manifold $M$ and $\cC\, \CU$ its simplicial \v{C}ech nerve, then it follows that the morphism $x \colon M \longrightarrow \rmB^n \rmU(1)$ can be presented by a \v{C}ech cocycle $g \colon \cC\,\CU \longrightarrow \rmB^n\rmU(1)$.
	\qen
\end{remark}

In the following we fix an $(n{-}1)$-gerbe $x \colon M \longrightarrow \rmB^n\rmU(1)$ on $M$, together with a presentation by a \v{C}ech cocycle $$g \colon \cC\,\CU \longrightarrow \rmB^n\rmU(1) \ . $$ We denote the \v{C}ech differential by $\check{\delta}:\Omega^q(\cC_k\,\CU)\longrightarrow \Omega^q(\cC_{k+1}\,\CU)$, for $q,k\in\NN_0$.

\begin{definition}
	\label{def:cCAt(g)}
	The \textit{Atiyah $L_\infty$-algebroid $\CAt(g)$} of $g$ is given as follows:
	the non-vanishing part of its underlying cochain complex reads as
	\begin{equation}
		\begin{tikzcd}[column sep=1.8cm]
			C^\infty(\cC_0\, \CU) \ar[r, "(-1)^n\, \check{\delta}"]
			& C^\infty(\cC_1\, \CU) \ar[r, "(-1)^{n-1}\, \check{\delta}"]
			& \cdots \ar[r, "(-1)^3\, \check{\delta}"]
			& C^\infty(\cC_{n-2}\, \CU) \ar[r, "{(0, (-1)^2\, \check{\delta})}"]
			& E_{n-1}(g)
		\end{tikzcd}
	\end{equation}
	where
	\begin{equation}
		E_{n-1}(g)
		\coloneqq \big\{ (X, f) \in TM \oplus C^\infty(\cC_{n-1}\, \CU)\ \big|\ \dd\log(g)(X) - \check{\delta} f = 0 \big\}
	\end{equation}
	is situated in degree zero.
	Its anchor map is the canonical projection onto $TM$, while its only non-trivial bracket is the binary bracket which is given by
	\begin{alignat}{3}
		\big[ (X, f)\,,\, (Y, h) \big] &= \big( [X,Y]\,,\, \pounds_Y f - \pounds_X h \big)\ ,
		&&\quad
		&& \forall\, (X, f)\,,\, (Y, h) \in \CAt(g)_0 \ ,
		\\[4pt]
		\big[ (X, f)\,,\, h \big] &= \pounds_X h \ ,
		&&\quad
		&& \forall\, (X, f) \in \CAt(g)_0\,,\, h \in \CAt(g)_k\,,\, -(n-1)\leqslant k\leqslant -1 \ ,
		\\[-10pt]
		[f, h] &= 0 \ ,
		&&\quad
		&& \forall\, f \in \CAt(g)_k\,,\, h \in \CAt(g)_l\,,\, -(n-1)\leqslant l, k \leqslant -1 \ .
	\end{alignat}
\end{definition}

In particular, $\CAt(g)$ is even a fibrant dg Lie algebroid on $M$.

\begin{remark}
	The kernel of the anchor map $\CAt(g) \longrightarrow TM$ is the abelian dg Lie algebroid with zero anchor map and underlying complex
	\begin{equation}
		\begin{tikzcd}[column sep=2cm]
			C^\infty(\cC_0\, \CU) \ar[r, "(-1)^n\, \check{\delta}"]
			& C^\infty(\cC_1\, \CU) \ar[r, "(-1)^{n-1}\,\check{\delta}"]
			& \cdots \ar[r, "(-1)^3\, \check{\delta}"]
			& C^\infty(\cC_{n-2}\, \CU) \ar[r, "(-1)^2\, \check{\delta}"]
			& \ker (\check{\delta})\ .
		\end{tikzcd}
	\end{equation}
	This complex is quasi-isomorphic to $C^\infty(M)[n-1]$.
	\qen
\end{remark}

There is a notion of $p$-form connections on a \v{C}ech cocycle $g \colon \cC\,\CU \longrightarrow \rmB^n\rmU(1)$ which predates and is independent of the formalism for $p$-form connections developed in the present paper.
Giving an $n$-form connection $\CA$ of this type on $g$ is equivalent to enhancing $g$ to a full \v{C}ech--Deligne differential cocycle (we review this in more detail in Section~\ref{sec: Con_l(g) equivalence} below).
The second main result of this section is the strictification of $T(M/\rmB^n\rmU(1))$ to a dg Lie algebroid on $M$ through

\begin{theorem}
\label{st:strictification of At(CG)}
	Let $H \in \Omega^{n+1}_\cl(M)$ be a closed $(n{+}1)$-form whose class $[H] \in \rmH^{n+1}_\dR(M)$ agrees with the image of $[g] \in \rmH^{n+1}(M; \ZN)$ in de Rham cohomology.
	Any enhancement of $g$ to a full \v{C}ech--Deligne differential cocycle $(g,\CA)$ with curvature $H$ canonically induces an equivalence of $L_\infty$-algebroids on $M$:
	\begin{equation}
		\begin{tikzcd}
			\varPhi_{g, \CA} \colon \At_{n-1}(H)
			\ar[r, rightsquigarrow, "\simeq"]
			& \CAt(g)\ .
		\end{tikzcd}
	\end{equation}
\end{theorem}

We prove Theorem~\ref{st:strictification of At(CG)} in Appendix~\ref{app:Pf of strictification}.

\begin{remark}
Theorem~\ref{st:strictification of At(CG)} for $n=1$ also appears in~\cite[Proposition~5.1.3]{FRS} and for $n=2$ in~\cite[Proposition~5.2.6]{FRS}. The novelty in our approach is a universal explicit computation and strictification of the Atiyah $L_\infty$-algebroid of an $(n{-}1)$-gerbe for all $n\geqslant1$.
\qen
\end{remark}



\subsection{$\infty$-groupoids of $p$-form connections: \v{C}ech--Deligne vs.~derived geometric}
\label{sec: Con_l(g) equivalence}


Given a Lie group $G$ and a principal $G$-bundle $P \longrightarrow M$, the connections on $P$ form a \textit{set} $\Con(P)$.
This set carries an affine action of the space of 1-forms on $M$ valued in the adjoint bundle of $P$, but it does not exhibit any higher structure.
Higher $\rmU(1)$-principal bundles $\CG$ over a manifold $M$ and their connections have a description in terms of \v{C}ech--Deligne cocycles.
Using this description, it has already been observed in~\cite{BS:Higher_Syms_and_Deligne_Coho} that the $p$-form connections on $\CG$ form an $\infty$-groupoid which, in general, does have non-trivial morphisms and higher morphisms; that is, we should have the following expectation, which is a feature of higher gauge theory:
\begin{itemize}
\item[]	\textsf{\large In general, the collection of $p$-form connections on a given higher principal bundle is not a set, but a higher groupoid with its own internal morphisms.}
\end{itemize}

We emphasise that the description of $p$-form connections on $(n{-}1)$-gerbes in terms of \v{C}ech--Deligne differential cocycles is entirely unrelated to the formalism we have developed in this paper.
The main upshot of our new derived deformation theoretic approach is that, in contrast to the \v{C}ech--Deligne picture, it straightforwardly generalises to generic types of principal $\infty$-bundles on manifolds in $\Sh_\infty(\InfMfd)$.
In this section we prove that for higher gerbes $x \colon M \longrightarrow \rmB^n\rmU(1)$ on $M$ the derived geometric approach reproduces the full homotopy type of $p$-form connections on $x$ as defined independently via \v{C}ech--Deligne differential cocycles.

We first recall the definition of $p$-form connections used in~\cite{BS:Higher_Syms_and_Deligne_Coho}.
It extends the geometric picture on differential cohomology promoted in~\cite{Gajer:Geo_of_Deligne_Coho, FSS:Cech_diff_char_classes_via_L_infty, Schreiber:DCCT}.
For $0 \leqslant p \leqslant n \in \NN_0$, consider the complex of sheaves of abelian groups
\begin{equation}
	\rmb^{n-p}\,\rmb_\nabla^p \rmU(1) \ \coloneqq
	\begin{tikzcd}[column sep = 1.2cm]
		\big( 0 \ar[r]
		& \rmU(1) \ar[r, "\dd \log"]
		& \Omega^1 \ar[r, "\dd"]
		& \cdots \ar[r, "\dd"]
		& \Omega^{p-1} \ar[r, "\dd"]
		& \Omega^p \ar[r, "\dd"]
		& 0 \big)\ ,
	\end{tikzcd}
\end{equation}
where the sheaf $\rmU(1)$ of smooth $\rmU(1)$-valued functions is situated in degree $-n$.
We define a simplicial presheaf on $\Cart \subset \Mfd$ as
\begin{equation}
	\rmB^{n-p}\, \rmB_\nabla^p \rmU(1)
	\coloneqq \DK \circ \rmb^{n-p}\,\rmb_\nabla^p \rmU(1)\ ,
\end{equation}
where $\DK$ is the Dold--Kan functor.

Let $M$ be a smooth manifold and $\CU = \{U_a\}_{a \in \varLambda}$ a good open cover of $M$.
The Kan complex of \textit{$(n{-}1)$-gerbes with $p$-form connection on $M$ relative to $\CU$} is defined as the simplicially enriched hom
\begin{equation}
	\Grb^{n-1}_{\nabla|p}(\CU)
	= \ul{\Fun(\Cart^\opp, \sSet)} \big( \cC\, \CU, \rmB^{n-p}\, \rmB_\nabla^p \rmU(1) \big)\ ,
\end{equation}
where $\cC\, \CU$ denotes the \v{C}ech nerve of the good open cover $\CU$ of $M$, and where the underline indicates that we are taking the enriched hom object.
Changing the good open covering does not change the homotopy type of this space, since we can equivalently describe \smash{$\Grb^{n-1}_{\nabla|p}(\CU)$} as computing the mapping space $\Map(M, \rmB^{n-p}\, \rmB_\nabla^p \rmU(1))$ in the local projective model structure on $\Fun(\Cart^\opp, \sSet)$.
Taking the colimit over all good open covers $\CU$ of $M$, one can promote \smash{$\Grb^{n-1}_{\nabla|p}$} to a simplicial homotopy sheaf of Kan complexes on $\Mfd$ (see~\cite[Appendix~B]{Bunk_R-loc_HoThy}).
Its value on $M$ is the Kan complex of \textit{$(n{-}1)$-gerbes with $p$-form connection on $M$}.

For each $0 \leqslant m \leqslant p \leqslant n \in \NN_0$ there are canonical projection morphisms
\begin{equation}
	\rmb^{n-p}\, \rmb_\nabla^p \rmU(1)
	 \longrightarrow \rmb^{n-m}\, \rmb_\nabla^m \rmU(1)\ ,
\end{equation}
which forget the forms of highest degrees.
Since the functors $\DK$ and \smash{$\ul{\Fun(\Cart^\opp, \sSet)} (\cC\, \CU, -)$} are right Quillen functors,%
\footnote{The latter functor is right Quillen from the projective model structure on simplicial presheaves on $\Cart$, localised at \v{C}ech nerves of good open coverings, to the Kan--Quillen model structure on simplicial sets.}
they induce Kan fibrations
\begin{equation}
	\Grb^{n-1}_{\nabla|p}(\CU)
	 \longrightarrow \Grb^{n-1}_{\nabla|m}(\CU)\ .
\end{equation}
In the extreme case $p=n$, we also write $$\Grb^{n-1}_{\nabla}(\CU):=\Grb^{n-1}_{\nabla|n}(\CU)$$ and refer to this as the Kan complex of \emph{$(n{-}1)$-gerbes on $M$ with (full) connection relative to $\CU$}. At the opposite extreme $p=0$, we also write $$\Grb^{n-1}(\CU):=\Grb^{n-1}_{\nabla|0}(\CU)$$ and refer to this as the Kan complex of \emph{$(n{-}1)$-gerbes on $M$ relative to $\CU$} (without connection).

\begin{definition}
[\textbf{\cite[Definition~10.22]{BS:Higher_Syms_and_Deligne_Coho}}]
Let $\CG$ be an $(n{-}1)$-gerbe on $M$, presented by a morphism $g \colon \cC\,\CU \longrightarrow \rmB^n\rmU(1)$.
For $0 \leqslant p \leqslant n \in \NN_0$, the Kan complex of \textit{$p$-form connections on $\CG$} is the fibre
\begin{equation}
\label{eq:Con_l on n-gerbe via Deligne complex}
\begin{tikzcd}[row sep = 1cm,column sep=1cm]
	\Con_p(g) \ar[r] \ar[d]
	& \Grb^{n-1}_{\nabla|p}(\CU) \ar[d]
	\\
	\Delta^0 \ar[r, "\{g\}"']
	& \Grb^{n-1}(\CU)
\end{tikzcd}
\end{equation}
\end{definition}

Since all vertices in the cospan underlying the pullback square~\eqref{eq:Con_l on n-gerbe via Deligne complex} are Kan complexes and the right-hand vertical map is a Kan fibration, the square~\eqref{eq:Con_l on n-gerbe via Deligne complex} is even homotopy cartesian in $\sSet$.
This makes the homotopy type of $\Con_p(g)$ independent of the choice of the cover $\CU$ and cocycle $g$.

We can also compute $\Grb^{n-1}_{\nabla|p}(\CU)$ as follows:
define the complexes
\begin{equation}
	\Grb^{n-1, \ch}_{\nabla|p}(\CU)
	:= \tau_{\leqslant 0} \circ \Tot^\times \circ N \big( \rmb^{n-p}\, \rmb_\nabla^p \rmU(1) (\cC\, \CU) \big)\ ,
\end{equation}
where $N$ takes the double complex associated to a cosimplicial complex and $\Tot^\times$ is the total complex functor (based on products, rather than direct sums).
The canonical projections $	\rmb^{n-p}\, \rmb_\nabla^p \rmU(1) \longrightarrow \rmb^{n-m}\, \rmb_\nabla^m \rmU(1)$ then induce morphisms
\begin{equation}
	\pi^{n, \ch}_{p,m} \colon \Grb^{n-1, \ch}_{\nabla|p}(\CU) \longrightarrow \Grb^{n-1, \ch}_{\nabla|m}(\CU)\ ,
\end{equation}
for each $0 \leqslant m \leqslant p \leqslant n \in \NN_0$.
There are canonical isomorphisms
\begin{equation}
	\Grb^{n-1}_{\nabla|p}(\CU)
	\cong \Mod_\ZN^\dg \big( C_*(\Delta), \Grb^{n-1, \ch}_{\nabla|p}(\CU) \big)\ ,
\end{equation}
where we have used that the Dold--Kan functor $\DK$ is canonically isomorphic to taking maps in $\Mod_\ZN^\dg$ out of the simplicial chains on the standard simplices~\cite{Kan:Functors_involving_CSS_complexes}.

In this setting, we can describe the $m$-simplices of $\Grb^{n-1}_{\nabla|p}(\CU)$ explicitly as follows:
for each injective morphism $$\varphi \colon [l] \longhookrightarrow  [m] $$ in $\bbDelta$, they consist of a tuple
\begin{equation}
	\label{eq: data for p-spls in Con_l(g)}
	\CA_\varphi = \big(g_\varphi, A_\varphi^{(1)}, \ldots, A_\varphi^{(\ell(\varphi))}\big)\ ,
\end{equation}
where
\begin{align}
	g_\varphi & \colon \cC_{n-l}\, \CU \longrightarrow \rmU(1)\ ,
	\\[4pt]
	A^{(k)}_\varphi & \in \Omega^k(\cC_{n-(l+k)}\,\CU)\ ,
	\quad
	k = 1, \ldots, \ell(\varphi):=\min(n-l, p) \ .
\end{align}
Here we set
\begin{alignat}{4}
	g_\varphi &= 1\  \quad &\text{if} &&\quad 0 &> n - l  \ ,
	\\[4pt]
	A^{(k)}_\varphi &= 0\  \quad &\text{if} &&\quad k &> n - l\ .
\end{alignat}

These data further satisfy 
\begin{align}
	\label{eq:conditions for p-spls in Con_l(g)}
	\begin{split}
	\check{\delta}\, g_\varphi
	&= \prod_{r = 0}^l\, g_{\varphi \circ \partial_r}^{(-1)^r}\ ,
	\\[4pt]
	\dd \log g_\varphi + (-1)^{l+1} \check{\delta} A^{(1)}_\varphi
	&= \sum_{r = 0}^l\, (-1)^r\, A^{(1)}_{\varphi \circ \partial_r} \ ,
	\\[4pt]
	\dd A^{(k-1)}_\varphi + (-1)^{l+k}\, \check{\delta} A^{(k)}_\varphi
	&= \sum_{r = 0}^l\, (-1)^r\, A^{(k)}_{\varphi \circ \partial_r }\ ,
	\quad  k = 2, \ldots, \ell(\varphi) \ ,
	\\[4pt]
	\dd A^{(\ell(\varphi))}_\varphi
	&=  \sum\limits_{r = 0}^l\, (-1)^r\, A^{(m)}_{\varphi \circ \partial_r} \ ,
		\quad l \geqslant 1 \ ,\ p > n-l\ .
	\end{split}
\end{align}
In the cases where $p \leqslant n-l$ there is no condition on $\dd A^{(\ell(\varphi))}_\varphi$.

We describe the simplicial structure using the generating morphisms of $\bbDelta^\opp$:
let $\CA$ denote an $m$-simplex in \smash{$\Grb^{n-1}_{\nabla|p}(\CU)$}.
Then $d_i \CA$ is the $(m{-}1)$-simplex which associates to each injective morphism $\varphi^- \colon [l] \longhookrightarrow  [m-1]$ the data
\begin{subequations}
\begin{equation}
	\label{eq: spl structure on Con_l(g), d_i}
	(d_i \CA)_{\varphi^-} = \CA_{\partial^i \circ \varphi^-}\ .
\end{equation}
The $(m{+}1)$-simplex $s_i \CA$ associates to an injective morphism $\varphi^+ \colon [q] \longhookrightarrow  [m+1]$ the data
\begin{equation}
	\label{eq: spl structure on Con_l(g), s_i}
	(s_i \CA)_{\varphi^+} =
	\begin{cases}
		\ \CA_{\sigma^i \circ \varphi^+}\ , & \text{if} \ \sigma^i \circ \varphi^+ \ \text{is injective}\ ,
		\\[4pt]
		\ 0\ , & \text{otherwise} \ .
	\end{cases}
\end{equation}
\end{subequations}

Consequently, the $m$-simplices of $\Con_p(g)$ are in bijection with those morphisms of chain complexes
\begin{equation}
	\CA \colon C_*(\Delta^m) \longrightarrow \Grb^{n-1, \ch}_{\nabla|p}(\CU)
\end{equation}
such that the composition
\begin{equation}
	\pi^{n, \ch}_{p,0} \circ \CA \colon C_*(\Delta^m) \longrightarrow \Grb^{n-1, \ch}(\CU)
\end{equation}
satisfies the conditions:
\begin{myitemize}
\item $\pi^{n, \ch}_{p,0} \circ \CA$ maps each degree~zero generator of the graded $\ZN$-module underlying $C_*(\Delta^m)$ to the element $g \in \Grb^{n-1, \ch}(\CU)_0$, and

\item $\pi^{n, \ch}_{p,0} \circ \CA$ maps each generator in degree $-l$ of the graded $\ZN$-module underlying  $C_*(\Delta^m)$ to the element $1 \in \Grb^{n-1, \ch}(\CU)_{-l}$.
\end{myitemize}
In other words, an $m$-simplex of $\Con_p(g)$ is an assignment of a tuple $\CA_\varphi = \big(g_\varphi, A^{(1)}_\varphi, \ldots, A^{(p)}_\varphi\big)$ to each injective morphism $\varphi \colon [l] \longhookrightarrow  [m]$ exactly as in~\eqref{eq: data for p-spls in Con_l(g)}, but with the additional requirement
\begin{equation}
	g_\varphi = \begin{cases}
		\ g\ , & l = 0\ ,
		\\[4pt]
		\ 1\ , & l \geqslant 1\ .
	\end{cases}
\end{equation}

\begin{remark}
We also extend the definition of $\Con_p(g)$ to $p > n$:
in this case, the data of $m$-simplices is as for $p \leqslant n$, but with the additional condition 
\begin{equation}
	\dd A^{(n)}_\varphi = 0
\end{equation}
whenever $\varphi \colon [0] \longrightarrow [m]$ presents a vertex of $\Delta^m$.
In particular, it follows that
\begin{equation}
	\Con_p(g) = \Con_q(g)
	\qquad
	\forall\, p, q > n\ .
\end{equation}

It thus makes sense to define
\begin{equation}
	\Con_\flat(g) \coloneqq \Con_p(g)\ ,
\end{equation}
the space of \textit{flat connections} on the $(n{-}1)$-gerbe $\CG$, where $p$ is any integer with $p > n$.
\qen
\end{remark}

\begin{remark}
This description  of $m$-simplices in $\Con_p(g)$ based on the functors \smash{$\Mod_\ZN^\dg(C_* (\Delta), -)$} is a much larger model than that produced by the smaller model of the functor $\DK$ (see for instance~\cite[Section~II.4]{GJ:Simplicial_HoThy}).
However, we emphasise again that the two models produce isomorphic simplicial sets:
much of the information in the families $\{\CA_\varphi\, | \, \varphi \colon [l] \longhookrightarrow  [m] \}$ is redundant.
To exemplify this, observe that in the larger model a 1-simplex is specified by data assigned to each of its vertices, and its edge.
In the smaller model determined by $\DK$, we instead specify data for only one vertex and the edge.
Let us denote these schematically by $c_0$ and $c_{01}$, respectively.
The data $c_1$ assigned to the second vertex is then already determined to be $c_1 = c_0 \pm \dd c_{01}$.
\qen
\end{remark}

The main result of this section is Theorem~\ref{st:finite and infinitesimal l-cons on n-gerbes} below, which states that the entire space of $p$-form connections on any higher gerbe produced by our derived geometric formalism is weakly homotopy equivalent the space of $p$-form connections produced by the \v{C}ech--Deligne description of higher gerbes with connections.
As the \v{C}ech--Deligne formalism makes no reference to derived techniques or deformation theory, this relates two genuinely independent definitions of connections on $(n{-}1)$-gerbes.

\begin{theorem}
	\label{st:finite and infinitesimal l-cons on n-gerbes}
	Suppose that $M$ is a smooth manifold, $\CU = \{U_a\}_{a \in \varLambda}$ is a good open cover of $M$, and $g \colon \cC\, \CU \longrightarrow \rmB^n \rmU(1)$ is a \v{C}ech $n$-cocycle.
	Let $x \colon M \longrightarrow \rmB^n\rmU(1)$ be the object in $\QCDT_{M/}$ presented by $g$.
	Then for each $p \in \NN_0$ there is a canonical equivalence
	\begin{equation}
		\Con_p(g)
		\simeq \Con_p \big( \CAt(g) \big)
	\end{equation}
	between the space of $p$-form connections on $g$ defined via \v{C}ech--Deligne cocycles and the space of $p$-form connections on $x$ in the derived geometric sense of Definition~\ref{def: Connection}.
\end{theorem}

Recall that by Corollary~\ref{st:At-L_oo-agd of n-gerbe from first principles} and Theorem~\ref{st:strictification of At(CG)}, the $L_\infty$-algebroid $\CAt(g)$ on the right-hand side is indeed a model for the relative tangent complex $T\big(M/\rmB^n\rmU(1)\big)$.
The proof of Theorem~\ref{st:finite and infinitesimal l-cons on n-gerbes} is computationally involved and is presented in two appendices:
in Appendix~\ref{app:simplicial resolution} we first develop techniques for computing mapping spaces of $L_\infty$-algebroids.
In particular, we construct computationally convenient fibrant simplicial resolutions of a class of $L_\infty$-algebroids which are close to being abelian (they arise from Lie actions of $TM$ on chain complexes).
We believe that generalisations of these constructions to more general $L_\infty$-algebroids should exist.
The construction of the weak homotopy equivalence in Theorem~\ref{st:finite and infinitesimal l-cons on n-gerbes} is contained in Appendix~\ref{app:Pf of equivalence of connection spaces}.

\begin{example}\label{ex:fullnform}
	We compute the curvature of an $p$-form connection on an $(n{-}1)$-gerbe  by combining Equation~\eqref{eq:curvature of l-conn algebraically} and Diagram~\eqref{eq: diagram for Con_l(g) via L_oo}.
	Concretely, consider Diagram~\eqref{eq: diagram for Con_l(g) via L_oo} with $r = 0$ and in the left-hand column add the order $(p{+}1)$ term.
	By computations analogous to those in Appendix~\ref{app:Pf of equivalence of connection spaces}, we obtain the curvature $F^p(\CA)$ of a $p$-form connection 
	\begin{equation}
		\CA=(A^{(1)}, \ldots, A^{(p)})
	\end{equation}
	on an $(n{-}1)$-gerbe, for $1 \leqslant p \leqslant n$.
	Only terms involving an entry of degree zero contribute to the brackets in~\eqref{eq:curvature of l-conn algebraically}, and all higher-order brackets are zero since the target $L_\infty$-algebroid is a dg Lie algebroid.
	
	With a hat denoting an omitted argument, we obtain 
	\begin{align}
		F^p(\CA)(X_0, \ldots, X_p)
		&= \sum_{0 \leqslant i < j \leqslant p}\, (-1)^{i+j}\, A^{(p)} \big( [X_i, X_j], X_0, \ldots, \widehat{X_i}, \ldots, \widehat{X_j}, \ldots, X_p \big)
		\\ &\qquad + \sum_{i = 0}^p\, (-1)^i\, \big[ \big( X_i, A^{(1)}(X_i) \big)\,,\, A^{(p)}(X_0, \ldots, \widehat{X_i}, \ldots, X_p) \big]
		\\[4pt]
		&= \sum_{0 \leqslant i < j \leqslant p}\, (-1)^{i+j}\, A^{(p)} \big( [X_i, X_j], X_0, \ldots, \widehat{X_i}, \ldots, \widehat{X_j}, \ldots, X_p)
		\\ &\qquad + \sum_{i = 0}^p\, (-1)^i\, \pounds_{X_i} \big( A^{(p)}(X_0, \ldots, \widehat{X_i}, \ldots, X_p) \big)
		\\[4pt]
		&= \dd A^{(p)} (X_0, \ldots, X_p)\ .
	\end{align}
	In particular, in the case $p = n$ we recover the usual curvature of a full connection on an $(n{-}1)$-gerbe:
	\begin{align}
	F^n(\CA) = \dd A^{(n)} = H \ .
	\end{align}
	This is also in agreement with the notion of curvature from the \v{C}ech--Deligne theory of differential cocycles.
	We again emphasise that the notion of curvature of an $p$-form connection we introduced in Section~\ref{sec: curvature of higher connections} is defined in an entirely independent way of this classical curvature, and instead recovers it from fundamental principles.
	\qen
\end{example}


\section{Higher Courant algebroids and symmetries of higher $\rmU(1)$-bundles}
\label{sec: Higher Courant algebroids}


In this final section we discuss the relation between higher $\rmU(1)$-bundles, Atiyah $L_\infty$-algebras and higher Courant algebroids.
A higher Courant algebroid has a canonical extension to an $L_\infty$-algebra which we argue can be thought of as encoding the infinitesimal symmetries of a $\rmU(1)$-$n$-bundle with $(n{-}1)$-form connection, i.e. symmetries which preserve all of the connection data except the highest form degree part (Conjecture~\ref{conj:Coursym}); we prove this rigorously for $n=2$ in Example~\ref{eg: deforming gerbe with connective structure}. We compare these $L_\infty$-algebroids to the Atiyah $L_\infty$-algebras which encode the infinitesimal symmetries preserving only the underlying $(n{-}1)$-gerbe (Theorem~\ref{prop:Linftyproj} and Conjecture~\ref{conj:CourAtsym}). The relevant objects do not have deformation theory, and a detailed study of their global moduli problem naturally leads to a weakened notion of $p$-form connections in terms of certain splittings of maps of $L_\infty$-algebras. In particular, we show that a weak $1$-form connection on a $\rmU(1)$-bundle with connection is simply the original connection with the restriction that it be flat (Example~\ref{eg: deforming U(1)-bundle with connection}), while a weak $1$-form connection on a $1$-gerbe with $1$-form connection is the same as an enhanced curving in the terminology of~\cite{TD:Chern_corr_for_higher_PrBuns} (Example~\ref{eg: deforming gerbe with connective structure}). 
	

\subsection{Higher Courant brackets and their $L_\infty$-structures}
\label{sub:highercourant}


We consider higher (and lower) analogues of (split) exact Courant algebroids~\cite{Hitchin:2003cxu,Ekstrand:2009qz,Zambon:2010ka,Bi:2010ssd}. Let $M$ be a smooth manifold and $n\in\NN$. Consider the vector bundle
\begin{align}
	\CE_{n-1}\coloneqq TM\oplus\midwedge^{n-1}\,T^*M
\end{align}
endowed with the symmetric fibrewise pairing
\begin{align}
	\ip{\,\cdot\ ,\,\cdot\,}:\CE_{n-1}\times \CE_{n-1} \longrightarrow \midwedge^{n-2}\,T^*M
\end{align}
given by a symmetrization of contractions between vectors and $(n{-}1)$-forms as
\begin{align}
	\ip{(X,\alpha)\ ,\,(Y,\beta)} \coloneqq \iota_X\beta + \iota_Y\alpha \ .
\end{align}
We further endow the $C^\infty(M)$-module of sections of $\CE_{n-1}$ with the \emph{Dorfman bracket}
\begin{align}
	\dbracket{(X,\alpha)\ ,\,(Y,\beta)} \coloneqq \big([X,Y] \ ,\, \pounds_X\beta-\iota_Y\,\mathrm{d}\alpha \big) \ .
\end{align}

With a suitable smooth bracket-preserving bundle map $\rho:\CE_{n-1} \longrightarrow TM$, the Dorfman bracket obeys Jacobi and Leibniz identities
\begin{align}
	\dbracket{e_1,\dbracket{e_2,e_3}} &= \dbracket{\dbracket{e_1,e_2},e_3} + \dbracket{e_2,\dbracket{e_1,e_3}} \ , \\[4pt]
	\dbracket{e_1,f\,e_2} &= f\,\dbracket{e_1,e_2} + \big(\pounds_{\rho(e_1)}\,f\big)\,e_2 \ , \\[4pt]
	\dbracket{f\,e_1,e_2} &= f\,\dbracket{e_1,e_2} - \big(\pounds_{\rho(e_2)}\,f\big)\,e_1 + \mathrm{d}f\wedge\ip{e_1,e_2} \ ,
\end{align}
as well as a compatibility condition between the pairing and the bracket
\begin{align}
	\pounds_{\rho(e_3)}\ip{e_1,e_2} = \ip{\dbracket{e_3,e_1},e_2} + \ip{e_1,\dbracket{e_3,e_2}} \ ,
\end{align}
for all $e_1,e_2,e_3\in \Gamma(M,\CE_{n-1})=TM\oplus\Omega^{n-1}(M)$ and $f\in C^\infty(M)$. 

When the anchor map $\rho$ is the projection to $TM$, we call the vector bundle $\CE_{n-1}$ endowed with these structure maps the \emph{standard higher Courant algebroid} on $M$. In this case it is also possible to twist the Dorfman bracket by a closed $(n{+}1)$-form $H\in\Omega_{\rm cl}^{n+1}(M)$ to the \emph{$H$-twisted Dorman bracket}
\begin{align}
	\dbracket{(X,\alpha)\ ,\,(Y,\beta)}_{H} \coloneqq \dbracket{(X,\alpha)\ ,\,(Y,\beta)} + (0,\iota_X\,\iota_Y H) \ .
\end{align}

The decomposition of the Dorfman bracket into its respective antisymmetric and symmetric parts
\begin{align}
	\dbracket{e_1,e_2} = \dbracket{e_1,e_2}_{\rm C} + \big(0\,,\,\tfrac12\,\mathrm{d}\ip{e_1,e_2}\big)
\end{align}
defines the \emph{Courant bracket}
\begin{align}
	\dbracket{(X,\alpha)\ ,\,(Y,\beta)}_{\rm C} \coloneqq \big([X,Y] \ ,\, \pounds_X\beta - \pounds_Y\alpha - \tfrac12\,\mathrm{d}(\iota_X\beta-\iota_Y\alpha)\big) \ .
\end{align}
Unlike the Dorfman bracket, the Courant bracket is antisymmetric but it violates the Jacobi identity by a coboundary-like term. In other words, neither the Dorfman nor Courant brackets define Lie structures on $\CE_{n-1}$. Below we will discuss how higher Courant algebroids can instead be viewed as $L_\infty$-algebras.

The geometric subgroup preserving the Dorfman bracket, the pairing and the anchor is the semi-direct product 
\begin{align}
\cSym(\CE_{n-1}) := \cDiff(M) \ltimes \Omega^n_{\rm cl}(M)
\end{align}
of diffeomorphisms of $M$ with \emph{gauge transformations} $\mathrm{e}^{\,B}:\Gamma(M,\CE_{n-1})\longrightarrow \Gamma(M,\CE_{n-1})$ parametrised by closed $n$-forms $B\in\Omega_{\rm cl}^{n}(M)$ acting as
\begin{align}
	\mathrm{e}^{\,B}(X,\alpha) = (X,\alpha+\iota_XB) \ .
\end{align}

\begin{example} \label{ex:gentangentbundle}
	Higher Courant algebroids appear in a number of supergravity contexts:
\begin{myenumerate}
\item The standard case $n=2$ is the \emph{generalised tangent bundle}
	\begin{align}
		\CE_1 = TM\oplus T^*M \ .
	\end{align}
	It represents an isotropic splitting of the exact sequence of vector bundles
	\begin{align}
		0 \longrightarrow T^*M \xrightarrow{ \ \rho^* \ } \CE_1 \xrightarrow{ \ \rho \ } TM \longrightarrow 0
	\end{align}
	and gives rise to an exact Courant algebroid with \v{S}evera class $[H]\in\mathrm{H}^3_\dR(M)$ represented by the closed $3$-form $H$~\cite{Severa:2017oew}, called the NSNS flux in string theory. Gauge transformations are called $B$-field transformations by the Kalb--Ramond 2-form $B $. Type II supergravity can be entirely formulated in terms of the generalised geometry of exact Courant algebroids~\cite{Coimbra:2011nw}.\label{ex:gentangentbundle1}

\item The case $n=3$ is the \emph{exceptional tangent bundle}
	\begin{align}
		\CE_2 = TM\oplus\midwedge^2\,T^*M
	\end{align}
	of four-dimensional compactifications of M-theory~\cite{Hull:2007zu,Chatzistavrakidis:2019seu}. It
	represents an isotropic splitting of the exact sequence of vector bundles
	\begin{align}
		T^*M\otimes T^*M \longrightarrow \CE_2 \xrightarrow{ \ \rho \ } TM \longrightarrow 0 \ ,
	\end{align}
	giving the extension of string charges to membrane charges, with closed 4-form $G$-flux $H$ and $C$-field gauge transformations by the 3-form $B$.
	
\item For $n\geqslant 3$, $\CE_{n-1}$ is the higher (twisted) Courant algebroid of $\mathrm{GL}^+(n{+}2,\RN)$ generalised geometry that appears in $(n{+}1)$-dimensional compactifications of type~IIB supergravity~\cite{Hull:2007zu}.
\qen
\end{myenumerate}
\end{example}

\begin{remark} \label{rem:higherCourantU1}
Higher Courant algebroids arise from higher $\rmU(1)$-bundles in the following way. Let $(g,A^{(1)},\dots,A^{{(n-1)}})$ be a \v{C}ech--Deligne differential cocycle of an $(n{-}1)$-gerbe $x:M\longrightarrow \rmB^{n}\rmU(1)$ with $(n{-}1)$-form connection on a smooth manifold $M$ with a good open covering $\CU=\{U_a\}_{a\in\varLambda}$, for $n\geqslant 2$ (see Section~\ref{sec: Con_l(g) equivalence}). Over each open set $U_a$ we put the standard higher Courant algebroid \smash{$TU_a\oplus\midwedge^{n-1}\,T^*U_a$} and glue them together on overlaps $U_a\cap U_b$, for all $a,b\in\varLambda$, using the automorphism $$(X,\alpha)\longmapsto \big(X,\alpha+\iota_X\,\dd A_{ab}^{{(n-1)}}\big) \ , $$ which  is a gauge transformation because $\dd A^{{(n-1)}}\in\Omega_{\rm cl}^{n}(\check{C}_1\,\CU)$.
Applying the de~Rham differential $\dd$ to the relation
\begin{align}
\dd A^{(n-2)} - (-1)^n\,\check\delta A^{(n-1)} = 0
\end{align}
implies the cocycle condition
\begin{align}
\dd\,\check\delta A^{(n-1)} = 0 
\end{align}
in $\Omega_{\rm cl}^n(\check{C}_2\,\CU)$, which in turn gives rise to  an extension
\begin{align} \label{eq:highercourantseq}
0\longrightarrow \midwedge^{n-1}\,T^*M \longrightarrow \CE_{n-1} \xrightarrow{ \ \rho \ } TM \longrightarrow 0
\end{align}
of $TM$ in the category of vector bundles. 

Now enlarge the data $(g,A^{(1)},\dots,A^{{(n-1)}})$ by $A^{{(n)}}\in\Omega^{n}(\check{C}_0\,\CU)$ to get an $n$-form connection whose curvature $H$ represents the class $[g]\in\mathrm{H}^{n+1}(M;\ZN)$ classifying the $\rmU(1)$-$n$-bundle $x:M\longrightarrow \rmB^{n}\rmU(1)$. Then from the equation $$\check\delta A^{(n)} = -(-1)^n\,\dd A^{(n-1)}$$ we get a Lagrangian splitting $$TM \longrightarrow \CE_{n-1} \ , \qquad X\longmapsto \big(X,\iota_X A^{(n)}\big)$$ of the exact sequence \eqref{eq:highercourantseq}, making the vector bundle \smash{$\CE_{n-1}\cong TM\oplus\midwedge^{n-1}\,T^*M$} into a standard higher Courant algebroid with the $H$-twisted Dorfman bracket. 

In the case $n=2$, i.e.~the twisted generalised tangent bundle from Example~\ref{ex:gentangentbundle}\ref{ex:gentangentbundle1}, this is the well-known Hitchin construction of an exact Courant algebroid from a gerbe with connective structure in the case that the \v{S}evera class $[H]$ is integral~\cite{Hitchin:2003cxu}. 
\qen
\end{remark}

We proceed to discuss how to canonically associate to any higher Courant algebroid $\CE_{n-1}$ an $n$-term $L_\infty$-algebra, or equivalently a Lie $n$-algebra, extending the $H$-twisted Courant bracket. This $L_\infty$-structure was worked out by Zambon in~\cite[Propositions~8.1 and~8.4]{Zambon:2010ka}, generalising the $n=2$ result of Roytenberg and Weinstein~\cite{Roytenberg:1998vn,Roytenberg:2007zz,Getzler2010} who show that the sections of an $H$-twisted Courant algebroid form a 2-term $L_\infty$-algebra. To this end, recall that the \emph{Bernouilli numbers} $B_k\in\mathbbm{Q}$ for $k\geqslant0$ are defined by the coefficients in the power series expansion
\begin{align}
\sum_{k=0}^\infty \, B_k \, \frac{t^k}{k!} := \frac{t}{\mathrm{e}^{\,t}-1} = 1-\frac t2 + \frac{t^2}{6\cdot2!}-\frac{t^4}{30\cdot4!}+\frac{t^6}{42\cdot 6!}-\frac{t^8}{30\cdot8!} + \cdots \ .
\end{align} 
We summarise Zambon's construction as 

\begin{proposition}\label{prop:Zambon}
Let $(\CE_{n-1},\ip{\,\cdot\ ,\,\cdot\,},\dbracket{\,\cdot\ ,\,\cdot\,}_{H})$ be an $H$-twisted higher Courant algebroid on a smooth manifold $M$ for $n\geqslant 2$. There is an $n$-term $L_\infty$-algebra ${\rm Cour}_{n-1}(H)$ over $TM$ with underlying complex
\begin{align}
C^\infty(M)[n-1]\xrightarrow{ \ \dd \ }\Omega^1(M)[n-2]\xrightarrow{ \ \dd \ }\cdots\xrightarrow{ \ \dd \ }\Omega^{n-2}(M)[1]\xrightarrow{ \ \dd \ }\Gamma(M,\CE_{n-1}) =  TM\oplus\Omega^{n-1}(M) \ ,
\end{align}
where $\dd$ is the de~Rham differential. For $v_i=(e_i,\xi_i)$ with $e_i=(X_i,\alpha_i)\in\Gamma(M,\CE_{n-1})$ and $\xi_i\in\Omega^{\leqslant n-2}(M):=\bigoplus_{r=0}^{n-2}\,\Omega^r(M)$, the non-zero higher brackets are given as follows:
\begin{myitemize}
\item Binary brackets:
\begin{align}
[e_1,e_2] = \dbracket{e_1,e_2}_{{\rm C},H} \quad , \quad [e_1,\xi_2] = \tfrac12\,\pounds_{X_1}\,\xi_2 \ .
\end{align}
\item Ternary brackets:
\begin{align}
[e_1,e_2,e_3] &= -T_{H}(e_1,e_2,e_3) := -\tfrac1{3!}\,\big(\ip{\dbracket{e_1,e_2}_{{\rm C},H}\ ,\,e_3} + \ip{\dbracket{e_3,e_1}_{{\rm C},H}\ ,\,e_2} + \ip{\dbracket{e_2,e_3}_{{\rm C},H}\ ,\,e_1}\big) \ , \\[4pt]
[\xi_1,e_2,e_3] &= -\tfrac1{3!}\,\big(\tfrac12\,(\iota_{X_2}\,\pounds_{X_3}\,\xi_1-\iota_{X_3}\,\pounds_{X_2}\,\xi_1) + \iota_{[X_2,X_3]}\,\xi_1\big) \ .
\end{align}
\item $k$-ary brackets for all odd $5\leqslant k\leqslant n+1$:
\begin{align}
[v_1,\dots,v_k] &= \sum_{i=1}^k \, (-1)^i \, \mu_k\big((\xi_i,\alpha_i),X_1,\dots,\widehat{X_i},\dots,X_k\big) + (-1)^{\frac{k+1}2} \, k\,B_{k-1} \ \iota_{X_k}\cdots\iota_{X_1}\,H \ ,
\end{align}
where
\begin{align}
& \mu_k\big((\xi_1,\alpha_1),X_2,\dots,X_k\big) \\[4pt]
& \hspace{1cm} = (-1)^{\frac{k+1}2} \, \frac{12\,B_{k-1}}{(k-1)\,(k-2)} \\
& \hspace{3cm} \, \times \sum_{2\leqslant i<j\leqslant k} \, (-1)^{i+j+1} \ \iota_{X_k}\cdots\widehat{\iota_{X_j}}\cdots\widehat{\iota_{X_i}}\cdots\iota_{X_1} \big([\xi_1,X_i,X_j] - T_0(\alpha_1,X_i,X_j)\big) \ ,
\end{align}
and the hat indicates omission of the corresponding term.
\end{myitemize}
\end{proposition}

Gauge transformations $\mathrm{e}^{\,B}:\Gamma(M,\CE_{n-1})\longrightarrow\Gamma(M,\CE_{n-1})$ induce strict $L_\infty$-isomorphisms $\id+\iota_{\rho(-)}B:{\rm Cour}_{n-1}(H)\longrightarrow{\rm Cour}_{n-1}(H)$~\cite[Proposition~8.5]{Zambon:2010ka}.

\subsection{Higher derived brackets}
\label{sub:derivedbracket}

For later use below, let us take a brief detour to review a general construction of $L_\infty$-algebras due to Fiorenza--Manetti~\cite{FM07} and Getzler~\cite{Getzler2010}, which plays a central role in the construction of Zambon. Consider a pair consisting of a dgla and a sub-dgla
$$
\frh\subseteq \frg \ .
$$
In this situation, the homotopy fibre of the inclusion $\frh\longhookrightarrow \frg$ of the underlying complexes is quasi-isomorphic to $\frg/\frh[-1]$. Our goal is to endow the complex $\frg/\frh[-1]$ with an $L_\infty$-structure that models the homotopy fibre in the $\infty$-category of $L_\infty$-algebras $L_\infty\Alg^\infty_\RN$.

To this end, let us fix a compatible decomposition of graded vector spaces
$$
\frg \cong \frh\,\oplus \, \frg/\frh \ .
$$
The differential on $\frg$ then assumes the form
$$
\dd_\frg = \begin{pmatrix}\dd_\frh & \sfD\\ 0 & \dd_{\frg/\frh}\end{pmatrix}
$$
for some intertwining map $\sfD\colon \frg/\frh\longrightarrow \frh$ satisfying $\dd_\frg\circ\sfD=-\sfD\circ\dd_{\frg/\frh}$. Finally, let us write
$$
\pi_\frh\colon \frg\longrightarrow \frh \qquad,\qquad \pi_{\frg/\frh}\colon \frg\longrightarrow \frg/\frh
$$
for the induced projections; here $\pi_{\frg/\frh}$ preserves the differential and does not depend on the choice of splitting.

Let us now recall that in the model category of dglas, a path object for a dgla $\frg$ is given by the dgla $\frg[t, \dd t]$ of polynomial differential forms on the $1$-simplex with values in $\frg$. Using this, we obtain the  rather big model for the homotopy fibre of $\frh\longhookrightarrow \frg$ in dglas given by
$$
\hofib_\mathrm{PL}(\frh\subseteq \frg) = \big\{x(t)+y(t)\,\dd t\in \frg[t, \dd t] \ \big| \ x(0)=0\ ,\, x(1)\in \frh\big\} \ .
$$
Here the subscript ${}_\mathrm{PL}$ indicates that we are using polynomial differential forms. 

This model for the homotopy fibre comes equipped with a natural deformation retract
$$\begin{tikzcd}
\frg/\frh[-1]\arrow[r, "i", yshift=1ex] & \hofib_\mathrm{PL}(\frh\subseteq \frg) \arrow[l, "r", yshift=-1ex]\arrow[loop right, "h"]
\end{tikzcd}$$
where
\begin{align*}
i(x) &= \sfD(x)\, t + x\, \dd t \ , \\[4pt]
r\big(x(t)+y(t)\,\dd t\big) &= \pi_{\frg/\frh}\Big(\int_0^1\, y(s)\,\dd s\Big) \ , \\[4pt]
h\big(x(t)+y(t)\,\dd t\big) &= \int_0^t\, y(s)\,\dd s - t\, \pi_{\frg/\frh}\Big(\int_0^1\, y(s)\,\dd s\Big) \ .
\end{align*}
Note that the inclusion $i$ is well defined because $\sfD(x)\, t$ vanishes at $t=0$ and is contained in $\frh$ when $t=1$. Likewise, the contracting homotopy $h$ is well-defined because $h\big(x(t)+y(t)\,\dd t\big)$ vanishes at $t=0$ and its value at $t=1$ is
$$
\int_0^1\, y(s)\,\dd s - \pi_{\frg/\frh}\Big(\int_0^1\, y(s)\,\dd s\Big) = \pi_{\frh}\Big(\int_0^1\, y(s)\,\dd s\Big) \ \in \ \frh \ .
$$
Using  $\dd_\frg\big(\int_0^t\, y(s)\,\dd s\big)=-\int_0^t\, \big(\dd_{\frg}y(s)\big)\,\dd s$ by the Koszul sign rule, one readily verifies $r\,i=\id$, $\id-i\,r = \dd \, h+h\,\dd$ and $h^2=0$.

We can thus apply the Homotopy Transfer Theorem to obtain an $L_\infty$-structure on $\frg/\frh[-1]$ and an enhancement of the maps $i$ and $r$ to $\infty$-quasi-isomorphisms $$i_\infty\colon \frg/\frh[-1]\longrightarrow \hofib_\mathrm{PL}(\frh\subseteq \frg) \qquad \text{and} \qquad r_\infty\colon \hofib_\mathrm{PL}(\frh\subseteq \frg)\longrightarrow \frg/\frh[-1] \ . $$ The brackets on $\frg/\frh[-1]$ are defined combinatorially by summing over rooted trees~\cite{LodayVallette:Algebraic_Operads}, where
\begin{myitemize}
\item each $k$-ary vertex is labelled by the $k$-ary bracket of $\hofib_\mathrm{PL}(\frh\subseteq \frg)$ (for $k\geqslant 2$),
\item each leaf is labelled by the inclusion $i$,
\item each internal edge is labelled by the homotopy $h$, and
\item the root is labelled by the restriction $r$.
\end{myitemize}
The resulting compositions then define an operation on $\frg/\frh[-1]$. 

As discussed in \cite[Section 5]{FM07}, this greatly simplifies in the present situation because $\hofib_\mathrm{PL}(\frh\subseteq \frg)$  has only a binary bracket and
$$
[\mathrm{im}(h), \mathrm{im}(h)]\subseteq \ker(h)\cap \ker(r) \ , 
$$
since $\mathrm{im}(h)$ consists of $0$-forms while $h$ and $r$ vanish on $0$-forms. Up to isomorphism the only rooted tree that remains is therefore the ``left comb'' (with automorphism group the degree~$2$ symmetric group $\Sigma_2$), so that the $k$-ary bracket on $\frg/\frh[-1]$ is given by
\begin{align} \label{eq:karyformula}
[x_1, \dots, x_k]_{\frg/\frh[-1]} = \frac{1}{2}\,\sum_{\sigma\in \Sigma_{k}}\, \pm\, r\big[h\big[h\big[\cdots h\big[i(x_{\sigma(1)}), i(x_{\sigma(2)})\big], i(x_{\sigma(3)})\big], \dots\big], i(x_{\sigma(k)})\big] \ .
\end{align}
Here the brackets on the right-hand side are the Lie brackets in $\hofib_\mathrm{PL}(\frh\subseteq\frg)$, and $\pm$ is the Koszul sign of the degree $k$ permutation $\sigma$ acting on $x_1\otimes\cdots \otimes x_k\in (\frg/\frh)^{\otimes k}$ (without the shift!). In other words, if the degrees of the elements in $\frg/\frh[-1]$ are shifted by one, then the bracket becomes graded-symmetric, as it should be.

Recall that $i(x) = \sfD(x)\, t+x\, \dd t$ and that $\mathrm{im}(h)$ consists of $0$-forms in $t$. Since $h$ and $r$ vanish on $0$-forms, the  formula \eqref{eq:karyformula} simplifies to
$$
[x_1, \dots, x_k]_{\frg/\frh[-1]} = \sum_{\sigma\in \Sigma_{k}}\, \pm\, r\big[h\big[h\big[\cdots h\big[\sfD(x_{\sigma(1)})\, t, x_{\sigma(2)}\,\dd t\big], x_{\sigma(3)}\,\dd t\big], \dots\big], x_{\sigma(n)}\,\dd t\big] \ .
$$
One can make this expression more explicit in terms of the Lie bracket $[-,-]_\frg$ on $\frg$. To this end, let us introduce the notation 
$$
\{x_1, \dots, x_k\} = \pi_{\frg/\frh}\big[\big[\cdots \big[[x_{1}, x_2]_\frg, x_3\big]_\frg,\dots\big]_\frg, x_k\big]_\frg
$$
for elements $x_1, \dots, x_k\in \frg$. 

\begin{proposition}\label{prop:derived brackets}
Let $\frh\subseteq \frg$ be the inclusion of a sub-dgla, together with a splitting of graded vector spaces $\frg\cong \frh\oplus \frg/\frh$ as above. Then the transferred $L_\infty$-structure on $\frg/\frh[-1]$ is given by
\begin{align*}
[x_1, \dots, x_k]_{\frg/\frh[-1]} &= \sum_{\sigma\in \Sigma_k}\,\pm \  \sum_{i=1}^{k-1} \ \sum_{1<k_1<k_2<\dots<k_i=k} \, \frac{(-1)^{i-1}}{k_1!\,(k_2-k_1+1)!\cdots (k_i-k_{i-1}+1)!}\\
& \qquad \, \times \big\{\big\{\cdots\big\{\{\sfD(x_{\sigma(1)}), x_{\sigma(2)}, \dots, x_{\sigma(k_1)}\}, x_{\sigma(k_1+1)},\dots, x_{\sigma(k_2)}\big\}, \dots,  x_{\sigma(k_i)}\big\} \ ,
\end{align*}
where $\pm$ is the Koszul sign of the permutation $\sigma$ acting on $x_1\otimes\cdots \otimes x_k\in (\frg/\frh)^{\otimes k}$.
\end{proposition}

\begin{proof}
This is a standard result, see e.g.~\cite{Voronov05}. The homotopy $h$ consists of two terms:
\begin{myenumerate}
\item[(a)] A term that takes the primitive.
\item[(b)] A term that integrates over $[0, 1]$ and then applies  $-t\, \pi_{\frg/\frh}$.
\end{myenumerate}
We therefore obtain a sum of $2^{k-1}$ terms, corresponding to the sum over $1<k_1<\dots<k_{i-1}<k$ where we take part (b) of the homotopy $h$. Each of these times we apply $-\pi_{\frg/\frh}$, yielding the iteration of brackets and the overall sign $(-1)^{i-1}$. On the other hand, each time we apply (a) in between, we pass from a coefficient $t^l/l!$ to coefficient $t^{l+1}/(l+1)!$, and this gets reset to a linear term in $t$ when we apply (b). This leads to the formula for the denominator.
\end{proof}

Let us mention two particular cases of relevance to us in which the formula of Proposition~\ref{prop:derived brackets} simplifies.

\begin{corollary}\label{cor:derived brackets}
Suppose that $\frg\cong \frh\oplus \frg/\frh$ is as in Proposition \ref{prop:derived brackets}.
\begin{myenumerate}
\item If the bracket of $\frg$ vanishes on $\frg/\frh$, then
$$
[x_1, \dots, x_k]_{\frg/\frh[-1]} =\frac{1}{k!}\,\sum_{\sigma\in \Sigma_k}\, \pm\,  \pi_{\frg/\frh}\big[\big[\cdots\big[\big[\sfD(x_{\sigma(1)}), x_{\sigma(2)}\big]_\frg, x_{\sigma(3)}\big]_\frg, \dots\big]_\frg, x_{\sigma(k)}\big]_\frg \ .
$$
\item\label{cor:derived brackets2} If $\big[\frg/\frh, \ \frg/\frh+\sfD(\frg/\frh)\big]_{\frg}\subseteq \frg/\frh$, then
$$
[x_1, \dots, x_k]_{\frg/\frh[-1]} =-\frac{B_{k-1}}{(k-1)!}\,\sum_{\sigma\in \Sigma_k}\, \pm\, \big[\big[\cdots\big[\big[\sfD(x_{\sigma(1)}), x_{\sigma(2)}\big]_\frg, x_{\sigma(3)}\big]_\frg, \dots\big]_\frg, x_{\sigma(k)}\big]_\frg \ ,
$$
\end{myenumerate}
where $B_{k-1}$ are the Bernouilli numbers.
\end{corollary}

\begin{proof}
In the first case, any iteration of brackets $\{\{x, \dots\}, y, \dots\}$ vanishes by assumption. The only remaining term is therefore $\{x_1, \dots, x_k\}$, that is, the term corresponding to $i=1$. 

In the second case, the assumptions imply that we can replace $\pi_{\frh/\frg}[-, -]_\frg$ by $[-, -]_\frg$ in all formulas. In particular, all iterated brackets are equal to the iterated bracket 
$$
\big[\big[\cdots\big[\big[\sfD(x_{\sigma(1)}), x_{\sigma(2)}\big]_\frg, x_{\sigma(3)}\big]_\frg, \dots\big]_\frg, x_{\sigma(k)}\big]_\frg \ .
$$
It therefore suffices to identify the coefficients $C_k$ appearing in the first line of the equation in Proposition~\ref{prop:derived brackets}. Separating out the summand with $i=1$, these coefficients are given by
$$
C_k = \frac{1}{k!} + \sum_{i=2}^{k-1} \ \sum_{1<k_1<k_2<\dots<k_i=k} \, \frac{(-1)^{i-1}}{k_1!\,(k_2-k_1+1)!\cdots (k_i-k_{i-1}+1)!} \quad,\quad k\geqslant 2 \ .
$$
Let us furthermore define $C_1=-1$ (which is incompatible with this formula). 

Using this, we obtain a recurrence relation for the $C_k$ as follows. Rewriting the summation over  $1<k_1<k_2<\dots<k_i=k$ as a sum over $1<r+1<r+l_1<\dots<r+l_j=k$ with $j=i-1$, we obtain
\begin{align*}
C_k &= \frac{1}{k!}-\sum_{r=1}^{k-2}\,\frac{1}{(r+1)!} \ \sum_{j=1}^{k-r-1} \ \sum_{1<l_1<l_2<\dots<l_j=k-r} \, \frac{(-1)^{j-1}}{l_1!\,(l_2-l_1+1)!\cdots(l_j-l_{j-1}+1)!}\\[4pt]
&= -\frac{C_1}{k!} - \sum_{r=1}^{k-2}\, \frac{1}{(r+1)!}\, C_{k-r} = -\sum_{r=1}^{k-1}\,\frac{1}{(r+1)!}\, C_{k-r} \ .
\end{align*}
The proof is completed by recalling that the {Bernouilli numbers} may be completely determined from their standard recursion formula
\begin{align}
\sum_{k=0}^m \, {m+1 \choose k} \, B_k = \delta_{m,0} \ ,
\end{align}
for $m\geqslant0$, whence $C_k=-\frac{B_{k-1}}{(k-1)!}$. (Alternatively, one can use the same argument as in the proof of~\cite[Lemma~5.4]{FM07}.)
\end{proof}

\begin{example}\label{ex:higher derived brackets}
Suppose that $\frg$ is a dgla, and let $\frg^{\geqslant 0}$ be the sub-dgla of elements of non-negative degree. Then $\frg/\frg^{\geqslant 0}=\frg^{<0}$ and $\sfD(\frg^{<0})\subseteq \frg^0$. We are therefore in the setting of Corollary \ref{cor:derived brackets}\ref{cor:derived brackets2}; the resulting $L_\infty$-structure on $\frg^{<0}[-1]$ coincides with the higher derived bracket construction of Getzler \cite{Getzler2010}, up to multiplying the $k$-ary bracket by $(-1)^{k-1}$.
This $L_\infty$-structure  depends only on $\frg^{\leqslant 0}$: if $\frg\longrightarrow\mathfrak{k}$ induces an isomorphism in non-positive degrees, then there is a natural isomorphism $\frg^{<0}[-1]\cong \mathfrak{k}^{<0}[-1]$.
\qen
\end{example}

\subsection{Projection to the Atiyah $L_\infty$-algebra}
\label{sub:Courant-Atiyah}

Our goal now is to apply the construction of Section~\ref{sub:derivedbracket} to relate Zambon's $L_\infty$-algebra ${\rm Cour}_{n-1}(H)$ to the Atiyah $L_\infty$-algebras $\mathrm{At}_{n-1}(H)=\big(TM\oplus C^\infty(M)[n-1],H\big)$ over $TM$ which underlie the $L_\infty$-algebroids associated with $\rmU(1)$-$n$-bundles that we introduced in Section~\ref{subsec:Atiyahalgebroids}. For this, it is instructive to  first recall the standard representative case at the lowest level, discussed already in Examples~\ref{Example: Lie group} and~\ref{ex:conn-principal}.

\begin{example}
\label{ex:AtiyahLie}
	For $n=1$, the central extension
	\begin{align}
		\CE_0 = TM\oplus{\RN} 
	\end{align}
	of the tangent bundle $TM$ represents a splitting of the short exact sequence of vector bundles
	\begin{align} \label{eq:n=0exact}
		0\longrightarrow \RN\xrightarrow{ \ j \ }\CE_0\xrightarrow{ \ \rho \ } TM \longrightarrow 0 \ ,
	\end{align}
	where $\RN$ denotes the trivial line bundle over $M$ (regarded as a bundle of abelian Lie algebras).
	Fixing a closed 2-form $H$ (a 2-cocycle of the Lie algebroid $TM$), the $H$-twisted Dorfman bracket
	\begin{align} \label{eq:dbracketn=0}
		\dbracket{(X,f)\ ,\,(Y,g)}_H = \big([X,Y] \ ,\, \pounds_X g - \pounds_Yf - H(X,Y)\big)
	\end{align}
	makes the space of sections $\Gamma(M,\CE_0)=TM\oplus C^\infty(M)$ into a Lie algebra. Hence $\CE_0$ is a Lie algebroid with anchor $\rho$ given by the projection to $TM$. The gauge symmetries of the bracket, induced by closed 1-forms $B$,  are automorphisms of the exact sequence \eqref{eq:n=0exact} which preserve the Lie algebra structure. These can be thought of as taking tensor products of $\CE_0$ with a trivial $\rmU(1)$-bundle equipped with a flat connection~$B$. 
	
In the case where the twisting 2-form $H$ defines an integer cohomology class $[H]\in \mathrm{H}^2(M;\ZN)$, Hitchin discusses in~\cite{Hitchin:2003cxu} how this can be interpreted in terms of the Atiyah Lie algebroid $\At(P)$ associated to a principal $\rmU(1)$-bundle $\pi:P \longrightarrow M$ through the Atiyah exact sequence
	\begin{align}
		0\longrightarrow \RN\xrightarrow{ \ j \ }TP/\rmU(1)\xrightarrow{ \ \pi_* \ } TM \longrightarrow 0 \ .
	\end{align}
	A choice of splitting of this sequence gives an isomorphism $TP/\rmU(1)\cong TM\oplus\RN$, which may be used to induce a Lie algebroid structure on $TM\oplus\RN$ from that of $TP/\rmU(1)$. The splitting is equivalent to a choice of $1$-form connection $\tau:TM\longrightarrow \At(P)= TP/\rmU(1)$ on $P$ (cf. Example~\ref{ex:conn-principal} and Remark~\ref{rem:higherCourantU1}), whose curvature $H=F^1(\tau)$ measures the failure of $\tau$ from being a Lie algebra homomorphism. The Lie bracket induced by this choice of splitting coincides with \eqref{eq:dbracketn=0}, and in this case there is an isomorphism
	\begin{align}
	\At(P) \ \cong \ \At_0(H)=\mathrm{Cour}_0(H) \ .
	\end{align}
The gauge transformations with $[B]\in \mathrm{H}^1(M;\ZN)$ arise from gauge transformations of the underlying $\rmU( 1)$-bundle $P$. 
\qen
\end{example}

For $n\geqslant 2$, Example~\ref{ex:AtiyahLie} together with Remark~\ref{rem:higherCourantU1} suggest an interpretation of an $H$-twisted higher Courant algebroid  as a `higher Atiyah Lie algebroid' associated to an $(n{-}1)$-gerbe equipped with an $(n{-}1)$-form connection. In the following we make this  statement precise by generalising and extending the $n=2$ constructions of~\cite[Section~5]{FRS}. 

Let $M$ be a smooth manifold and let $\Omega^*(M)$ be its de Rham algebra. Let us write $\mathrm{Der}(\Omega^*(M))$ for the dgla of derivations of $\Omega^*(M)$. This is given in degree $k$ by the set of degree-preserving $\mathbbm{R}$-linear maps $Q\colon \Omega^*(M)\longrightarrow \Omega^*(M)[k]$ such that $$ Q(\alpha\wedge \beta) = Q(\alpha)\wedge\beta + (-1)^{k\, |\alpha|}\, \alpha\wedge Q(\beta) \ . $$ The Lie bracket is the commutator bracket and the differential $\delta$ is given by the commutator with the de Rham differential:
$$
[Q_1, Q_2] = Q_1\,Q_2 - (-1)^{|Q_1|\,|Q_2|}\, Q_2\,Q_1\qquad,\qquad \delta(Q) = [\dd,Q] = \dd\circ Q - (-1)^{|Q|}\, Q\circ \dd \ .
$$
Observe that $\mathrm{Der}(\Omega^*(M))$ is concentrated in degrees $\geqslant -1$: derivations of degree $\leqslant -2$ vanish on $0$-forms and $1$-forms, and hence on all of $\Omega^*(M)$. The derivations of degree $-1$ are all of the form $\iota_X$, the contraction along a vector $X\in TM$  on $M$. The differential sends each such $\iota_X$ to the derivation $\pounds_X$, the Lie derivative along $X$.

By definition, $\mathrm{Der}(\Omega^*(M))$ acts on $\Omega^*(M)$. We can therefore form the semi-direct product dgla
$$
\mathrm{Der}(\Omega^*(M))\ltimes \Omega^*(M)[-n]
$$
with bracket and differential
$$
\big[(Q_1,\alpha_1)\ ,\, (Q_2,\alpha_2)\big] = \big([Q_1, Q_2] \ ,\, Q_1(\alpha_2) - (-1)^{|Q_2|\,|\alpha_1|}\, Q_2(\alpha_1)\big)\qquad,\qquad \delta(Q,\alpha) = \big([\dd, Q]\ ,\,\dd\alpha\big) \ .
$$
Here $|\alpha|$ is the cohomological degree of $\alpha$, which differs (because of the shift) from its form degree.

The closed $(n{+}1)$-form $H$ defines an element of degree $1$ in this dgla. It is a Maurer--Cartan element, since $\delta(H)=\dd H=0$ and $[H,H]=0$. We can therefore twist the differential by $H$ to get the differential
$$
\delta_H(Q,\alpha) = \big([\dd, Q] \ ,\, \dd\alpha - (-1)^{|Q|}\,Q(H)\big) \ .
$$
Let us denote the resulting dgla by
$$
\frg_H = \big(\mathrm{Der}(\Omega^*(M))\ltimes \Omega^*(M)[-n]\ ,\, \delta_H\ ,\, [-, -]\big) \ .
$$

For each $0\leqslant p\leqslant n-1$, let us consider the graded subspace
$$
F^{p+1}\frg_H = \mathrm{Der}(\Omega^*(M))^{\geqslant 0} \oplus \Omega^{\geqslant p+1}(M)[-n]
$$
consisting of all derivations of degree $\geqslant 0$ and all forms of form degree $\geqslant p+1$. This is closed under the differential and the bracket, so $F^{p+1}\frg_H\subseteq \frg_H$ is a sub-dgla. Using the isomorphism $\mathrm{Der}(\Omega^*(M))^{-1}\cong TM$ which identifies a vector field $X$ with the contraction $\iota_X$, the shifted quotient $\frg_H/F^{p+1}\frg_H[-1]$ is isomorphic to the non-positively graded complex
\begin{subequations}\label{eq:cour-at-complex}
\begin{equation}\label{eq:cour-at-complex0}
\begin{tikzcd}[row sep=0pc,column sep=1.5pc]
 C^\infty(M)[n-1]\arrow[r, "0"] &  0\arrow[r,"0"] & \cdots \arrow[r,"0"] & 0 \arrow[r,"0"] & TM
\end{tikzcd} 
\end{equation}
for $p=0$, to
\begin{equation}\label{eq:cour-at-complex0n}
\begin{tikzcd}[row sep=0pc,column sep=1.5pc]
 C^\infty(M)[n-1]\arrow[r, "\dd"] & \Omega^1(M)[n-2] \arrow[r,"\dd"] & \cdots\arrow[r, "\dd"] & \Omega^{p}(M)[n-p-1]\arrow[r,"0"] & 0\arrow[r,"0"] & \cdots \arrow[r,"0"] & 0 \arrow[r,"0"] & TM 
\end{tikzcd} 
\end{equation}
for $0<p<n-1$, and to
\begin{equation}\label{eq:cour-at-complexn}
\small
\begin{tikzcd}[row sep=0pc,column sep=1.5pc]
 C^\infty(M)[n-1]\arrow[r, "\dd"] & \Omega^1(M)[n-2] \arrow[r,"\dd"] & \cdots\arrow[r, "\dd"] & \Omega^{n-2}(M)[1]\arrow[r,"\dd"] &  TM\oplus \Omega^{n-1}(M) \cong \Gamma(M, \CE_{n-1})
\end{tikzcd} 
\normalsize
\end{equation}
\end{subequations}
for $p=n-1$. Note that the complexes \eqref{eq:cour-at-complex0} and \eqref{eq:cour-at-complexn} are precisely the complexes underlying the $L_\infty$-algebras $\At_{n-1}(H)$ and ${\rm Cour}_{n-1}(H)$, respectively.

As graded vector spaces, there is a canonical splitting $\frg_H\cong F^{p+1}\frg_H\oplus \frg_H/F^{p+1}\frg_H$. The map $\sfD\colon \frg_H/F^{p+1}\frg_H\longrightarrow F^{p+1}\frg_H$ is given by 
$$
\sfD(X)= \pounds_X-(-1)^{n+1}\,\iota_XH \qquad,\qquad \sfD(\alpha)= \dd\alpha
$$
for $X\in TM$ and $\alpha\in \Omega^p(M)$. The map $\sfD$ vanishes on forms of degree $<p$. We can therefore apply Proposition \ref{prop:derived brackets} to obtain an $L_\infty$-structure on $(\frg_H/F^{p+1}\frg_H)[-1]$.

\begin{proposition}\label{prop:cour-intermediate}
Let $n\geqslant 2$ and $0\leqslant p\leqslant n-1$.
The complex \eqref{eq:cour-at-complex} admits an $L_\infty$-structure with non-trivial higher brackets given on $X_i\in TM$ and $\alpha\in \Omega^{\leqslant p}(M):=\bigoplus_{r=0}^p\,\Omega^r(M)$ as follows:
\begin{myitemize}
\item {Binary brackets:} 
\begin{align}
[X_1, X_2] &= \left\{\begin{array}{ll} [X_1, X_2]_{TM}+ \iota_{X_1}\,\iota_{X_2}H & \text{if } \ p=n-1 \ , \\[4pt]
{[X_1, X_2]_{TM}} & \text{if } \ p<n-1 \ , \end{array}\right.
\\[4pt]
 [X_1, \alpha] &= \left\{\begin{array}{ll}\pounds_{X_1}\alpha & \text{if } \ \alpha\in \Omega^{\leqslant p-1}(M) \ ,\\[4pt]
\pounds_{X_1}\alpha + \iota_{X_1}\,\dd\alpha & \text{if } \ \alpha\in \Omega^{p}(M) \ .
\end{array}\right.
\end{align}

\item $k$-ary brackets for all odd $3\leqslant k\leqslant p+2$:
\begin{align*}
[X_1, \dots, X_{k-1}, \alpha] &= (-1)^{k}\, B_{k-1}\,\iota_{X_1}\cdots\iota_{X_{k-1}}\,\dd\alpha- \frac{B_{k-1}}{k-1}\,\sum_{i=1}^{k-1}\, \iota_{X_1}\cdots \widehat{\iota_{X_i}} \dots \iota_{X_{k-1}}\, \pounds_{X_i}\alpha\\*
& \quad \, +\frac{2\,B_{k-1}}{(k-1)\,(k-2)}\,\sum_{1\leqslant i<j\leqslant k-1}\, \iota_{X_1}\cdots \widehat{\iota_{X_i}}\cdots \widehat{\iota_{X_j}} \cdots \iota_{X_{k-1}}\,\iota_{[X_i, X_j]}\,\alpha \ .
\end{align*}

\item $k$-ary brackets for  $n-p+1\leqslant k\leqslant n+1$:
\begin{align}
[X_1, \dots, X_k] = k!\ C_{k, p} \ \iota_{X_1}\cdots\iota_{X_k}H \ ,
\end{align}
where
$$
C_{k, p}=\sum_{i=1}^{k+p-n} \ \sum_{n-p<k_1<k_2<\dots<k_i=k}  \, \frac{(-1)^{i-1}}{k_1!\,(k_2-k_1+1)!\cdots (k_i-k_{i-1}+1)!} \ .
$$
\end{myitemize}
\end{proposition}

\begin{proof}
The Lie bracket of $\frg_H$ has the  properties
\begin{align*}
[\mathrm{Der}(\Omega^*(M)), \mathrm{Der}(\Omega^*(M))]&\subseteq \mathrm{Der}(\Omega^*(M)) \ , \\[4pt]
[\mathrm{Der}(\Omega^*(M)), \Omega^*(M)]&\subseteq \Omega^*(M) \ , \\[4pt]
[\Omega^*(M), \Omega^*(M)]&=0 \ .
\end{align*}
In particular, any Lie word involving two differential forms is zero. Likewise, any Lie word in derivations $\pounds_{X_i}$ and $\iota_{X_i}$ of degree $\leqslant 0$ will vanish as soon as it contains two contractions $\iota_{X_i}$ and $\iota_{X_j}$ (for degree reasons). The only non-trivial brackets that can appear as part of the formulas from Proposition \ref{prop:derived brackets} are hence of the  form
\begin{align*}
\{\pounds_{X_1}, \iota_{X_2}\} &= \iota_{[X_1, X_2]} \ , \\[4pt]
\{\alpha, \iota_{X_1}, \dots, \iota_{X_k}\} &= (-1)^{k\, |\alpha|}\,\iota_{X_1} \cdots\iota_{X_k}\alpha && \text{if } \ \alpha\in \Omega^{\leqslant p+k}(M)[-n+1] \ , \\[4pt]
\{\pounds_{X_1}, \alpha, \iota_{X_2}, \dots, \iota_{X_k}\} &= (-1)^{(k-1)\,|\alpha|}\, \iota_{X_2} \cdots\iota_{X_{k}}\,\pounds_{X_1}\alpha && \text{if } \ \alpha\in \Omega^{\leqslant p+k-1}(M)[-n+1] \ , \\[4pt]
\{\pounds_{X_1}, \iota_{X_2}, \alpha, \iota_{X_3},\dots, \iota_{X_k}\} &= (-1)^{(k-1)\,|\alpha|}\, \iota_{X_3} \cdots\iota_{X_{k}}\,\iota_{[X_1, X_2]}\,\alpha & & \text{if } \ \alpha\in \Omega^{\leqslant p+k-1}(M)[-n+1] \ .
\end{align*}

Consequently, only a few iterated brackets can appear in the formula for the $k$-ary bracket. These are of the form
\begin{align*}
\big\{\cdots\big\{\big\{\dd\alpha, \iota_{X_1}, \dots\big\}, \dots, \iota_{X_k}\big\}, \dots, \iota_{X_k}\big\} &= (-1)^{k\,(n-p)}\,\iota_{X_1}\cdots\iota_{X_{k}}\,\dd\alpha && \text{if }\alpha\in \Omega^{p}(M)[-n+1] \ , \\[4pt]
\big\{\cdots\big\{\big\{\iota_{X_1}H, \iota_{X_2}, \dots, \iota_{X_{k_1}}\big\}, \dots\big\}, \dots, \iota_{X_k}\big\} &= \iota_{X_1}\cdots\iota_{X_{k}}H && \text{if }k_1\geqslant n+1-p \ , \\[4pt]
\big\{\cdots\big\{\big\{\pounds_{X_1}, \alpha, \iota_{X_2}\dots\big\}, \dots\big\}, \dots, \iota_{X_k}\big\} &= (-1)^{(k-1)\,|\alpha|}\,\iota_{X_2}\cdots\iota_{X_k}\,\pounds_{X_1}\alpha \ ,  & \\[4pt]
\big\{\cdots\big\{\big\{\pounds_{X_1}, \iota_{X_2}, \alpha,  \dots\big\}, \dots\big\}, \dots, \iota_{X_k}\big\} &= (-1)^{(k-1)\,|\alpha|}\,\iota_{X_3}\cdots\iota_{[X_{1}, X_2]}\,\alpha \ , \\[4pt]
\big\{\cdots\big\{\big\{\pounds_{X_1} , \iota_{X_2}\big\}, \alpha, \dots\big\}, \dots, \iota_{X_k}\big\} &= (-1)^{(k-1)\,|\alpha|}\,\iota_{X_3}\cdots\iota_{[X_{1}, X_2]}\,\alpha \ , \\[4pt]
\{\pounds_{X_1}, \iota_{X_2}\} &= \iota_{[X_1, X_2]} \ ,
\end{align*}
where all remaining elements indicated by ``$\dots$'' are of the form $\iota_{X_j}$ for $X_j\in TM$ and $\alpha\in \Omega^{\leqslant p}(M)$. Plugging these into Proposition \ref{prop:derived brackets} yields the result.
\end{proof}

\begin{definition}\label{def:CourpLinfty}
For $n\geqslant 2$ and $0\leqslant p\leqslant n-1$, we write $\mathrm{Cour}_{n-1}^{\leqslant p}(H)$ for the $L_\infty$-algebra of Proposition~\ref{prop:cour-intermediate}.
\end{definition}

\begin{lemma}\label{lem:CourpAt}
$\mathrm{Cour}_{n-1}^{\leqslant n-1}(H)\cong\mathrm{Cour}_{n-1}(H)$ \ and \ $\mathrm{Cour}_{n-1}^{\leqslant 0}(H)\cong \mathrm{At}_{n-1}(H)$.
\end{lemma}

\begin{proof}
Zambon defines $\mathrm{Cour}_{n-1}(H)$ using Getzler's higher derived bracket construction (Example~\ref{ex:higher derived brackets}) applied to the dgla of shifted poly-derivations
$$
\Big(\Sym_{\Omega^*(M)}\big(\mathrm{Der}(\Omega^*(M))[n]\big)\Big)[-n] \ ,
$$
twisted by the Maurer--Cartan element $H\in \Omega^*(M)$. In degrees $\leqslant 0$, this dgla is isomorphic to the sub-dgla $\mathrm{Der}(\Omega^*(M))\oplus \Omega^*(M)[-n]$ of poly-derivations of weight $\leqslant 1$ in the symmetric algebra. By Example \ref{ex:higher derived brackets}, this implies  $\mathrm{Cour}_{n-1}^{\leqslant n-1}(H)\cong\mathrm{Cour}_{n-1}(H)$. 

On the other hand, when $p=0$ the $L_\infty$-structure from Proposition \ref{prop:cour-intermediate} only has one non-trivial higher bracket (for degree reasons), which is $[X_1, \dots, X_{n+1}]=H(X_1,\dots,X_{n+1})$.
\end{proof}

By construction, the $L_\infty$-algebras of Definition~\ref{def:CourpLinfty} fit into a diagram of the form
\begin{equation}
{\begin{tikzcd}[column sep = 0.75cm]
\mathrm{Cour}_{n-1}^{\leqslant n-1}(H)\arrow[r, dashed]\arrow[dd, "i_\infty"{swap}, "\sim", xshift=-6pt] & \cdots \arrow[r,dashed] & \mathrm{Cour}_{n-1}^{\leqslant p}(H)\arrow[r, dashed]\arrow[dd, "i_\infty"{swap}, "\sim", xshift=-6pt] & \cdots \arrow[r,dashed] & \mathrm{Cour}_{n-1}^{\leqslant 0}(H)\arrow[dd, "i_\infty"{swap}, "\sim", xshift=-6pt]\\ & & & & \\
\hofib_\textrm{PL}(F^{n}\frg_H\subseteq \frg_H)\arrow[r]\arrow[uu, "r_\infty"{swap}, xshift=6pt] & \cdots \arrow[r] & \hofib_\textrm{PL}(F^{p+1}\frg_H\subseteq \frg_H)\arrow[r]\arrow[uu, "r_\infty"{swap}, xshift=6pt] & \cdots \arrow[r] & \hofib_\textrm{PL}(F^{1}\frg_H\subseteq \frg_H)\arrow[uu, "r_\infty"{swap}, xshift=6pt]
\end{tikzcd} } 
\end{equation}
where the bottom horizontal maps are induced by the inclusions $F^{p+1}\frg_H\subseteq F^{p}\frg_{H}$. In particular, there exists a sequence of $\infty$-morphisms of $L_\infty$-algebras, indicated by the dashed horizontal arrows in the top row, by composing the bottom horizontal maps with $r_\infty$ and $i_\infty$. At the linear level, the map $\mathrm{Cour}_{n-1}^{\leqslant p}(H)\longrightarrow \mathrm{Cour}_{n-1}^{\leqslant p-1}(H)$ is the projection sending the $p$-forms to zero. 

Notably, when combined with Lemma~\ref{lem:CourpAt} this construction establishes

\begin{theorem} \label{prop:Linftyproj}
For $n\geqslant 2$, there is a natural fibration of $L_\infty$-algebras over $TM$
\begin{equation}
\begin{tikzcd}
\mathrm{Cour}_{n-1}(H) \ar[r,rightsquigarrow] & \mathrm{At}_{n-1}(H) \ ,
\end{tikzcd}
\end{equation}
which at the linear level is the natural projection of  complexes
\begin{equation} \label{eq:complexprojection}
\begin{tikzcd}[row sep = 1cm]
C^\infty(M)[n-1] \arrow[r, "\dd"] \arrow[d,"\id"] & \Omega^1(M)[n-2] \arrow[r,"\dd"] \arrow[d,"0"] & \cdots 	\arrow[r,"\dd"] & \Omega^{n-2}(M)[1]\arrow[r,"\dd"] \arrow[d,"0"] & TM\oplus\Omega^{n-1}(M) \arrow[d,"\rho"]	\\
C^\infty(M)[n-1]  \arrow[r, "0"]  & 0 \arrow[r,"0"] & \cdots 	\arrow[r,"0"] & 0 \arrow[r,"0"] & TM
			\end{tikzcd}
\end{equation}
\end{theorem}

\begin{remark} \label{rem:FRSn=1}
Theorem~\ref{prop:Linftyproj} generalises the $n=2$ result of~\cite[Proposition 5.2.3]{FRS}, which connects the Courant algebroid with the Atiyah Lie $2$-algebroid.
\qen
\end{remark}


\subsection{Global moduli problems}


We would now like to expound the view that the $L_\infty$-algebra associated to a higher Courant algebroid, detailed in Proposition~\ref{prop:Zambon}, encodes the infinitesimal symmetries of an $(n{-}1)$-gerbe with $(n{-}1)$-form connection, as suggested by Remark~\ref{rem:higherCourantU1} and through its relation to the Atiyah $L_\infty$-algebra in Theorem~\ref{prop:Linftyproj}. For this, we seek an equivalence akin to Theorem~\ref{st:At-L_oo-agd of n-gerbe from first principles}. Note, however, that $\mathrm{Cour}^{\leqslant p}_{n-1}(H)$ does \emph{not} define an $L_\infty$-algebroid for $p>0$, and so this problem should be tackled in a global setting of underlying $L_\infty$-algebras. In order to achieve this within the framework of the present paper, in this section we take a slight detour to analyse some examples of moduli problems from the global perspective discussed in Section~\ref{sec: infinitesimal symmetries}.

For a smooth manifold $M$ and stack $X \in \Sh_\infty(\InfMfd)$, we will also refer to a morphism $x \colon M \longrightarrow X$ as an $X$-structure on $M$; this is convenient when thinking of $X$ as a classifying stack for some geometric data on $M$.
For some stacks $X \in \Sh_\infty(\InfMfd)$ there are different notions of moduli stack of $X$-structures on $M$, depending on how we define an $X$-structure on products $U \times M$ of $M$ with a parameter space $U$.
For example, to define smooth families of differential forms or connections on $M$ parametrised by $U$, we can either use all differential forms on $U \times M$, or we can restrict to the vertical differential forms, i.e.~those forms which vanish as soon as they are evaluated on a tangent vector to the parameter space $U$ (this is closely related to concretification~\cite{Schreiber:DCCT, BSS:YM_fields_on_Lorentz_mfds}).
Importantly, it may happen that only one of these moduli stacks has deformation theory. 

\begin{example}
	\label{eg:moduli and tangents of p-forms}
	Let $$\Omega^p \colon \Mfd^\opp \longrightarrow \Set$$ be the sheaf assigning to each smooth manifold its set of $p$-forms.
	For any fixed manifold $M$, consider the sheaf $\Hom_v(M, \Omega^p) / \cDiff(M)$, whose value on $U$ is the action groupoid of $\cDiff(M)(U)$ on the set of \textit{vertical} $p$-forms on $U \times M$.
	
	We know from Example~\ref{eg: forms on M have no DefThy} that $\Omega^p$ and the stack of differential forms on $M$ do not have deformation theory when we define a $U$-family of differential forms on $M$ to be a form on $U \times M$. 
	However, Example~\ref{ex:vertical forms def thy} shows that $\Hom_v(M, \Omega^p)$ \textit{does} have deformation theory.
	We can thus encode a $p$-form on $M$ as a morphism $$\omega \colon * \longrightarrow \Hom_v(M, \Omega^p)\,\big/\,\cDiff(M)$$ and ask for the relative tangent complex of this morphism.
	Note that this deformation problem does not correspond to a formal moduli problem over $M$ in the sense of Theorem~\ref{thm:formal stack vs lie}, for instance.
	\qen
\end{example}

In order to compute the tangent stacks of moduli problems of this form, we use the following results:
let $A \in \Sh_\infty(\InfMfd)$ be an abelian group object, and suppose that it is endowed with a $\cDiff(M)$-action which is compatible with the group structure of $A$, thus defining a semi-direct product group $\cDiff(M) \ltimes A$.
Consider a morphism $$a \colon * \longrightarrow A \, \big/ \, \cDiff(M) \ . $$
We define a group object $\cSym(a)$ as the homotopy pullback
\begin{equation}
	\label{eq: def Sym(a) for a: * ---> A}
	\begin{tikzcd}[column sep = 2cm, row sep = 1.25cm]
		\cSym(a) \ar[r] \ar[d]
		& \cDiff(M) \ar[d, "{\phi \longmapsto (\phi,\, a^{-1}\, (\phi^{-1})^*a)}"]
		\\
		\cDiff(M) \ar[r, "{\phi \longmapsto (\phi, 1)}"']
		& \cDiff(M) \ltimes A
	\end{tikzcd}
\end{equation}

Via the canonical equivalence between maps $b \longrightarrow c$ and $1 \longrightarrow b^{-1}\, c$ in group objects in spaces, where $1$ denotes a unit element, we readily see

\begin{lemma}
	If $A = X^M$ for some group object $X$, and $a = x^\dashv$ is the morphism corresponding to a morphism $x \colon M \longrightarrow X$ with deformation theory, then $\cSym(a)$ agrees with the smooth $\infty$-group $\cSym(x) = \cDiff_{/X}(M, M)$ of symmetries of $x$ introduced in Remark~\ref{rem:loop space relative diffeo}.
\end{lemma}

\begin{proposition}
	\label{st: computing sym(a) in A//Diff(M) case}
	In the situation of~\eqref{eq: def Sym(a) for a: * ---> A}, suppose that $A$ has deformation theory and let $E$ be a cochain complex which describes the $L_\infty$-algebra of $A$; since $A$ is abelian, all of its brackets vanish.
	Then the differential graded Lie algebra over $TM$ associated to a morphism $a \colon * \longrightarrow A / \cDiff(M)$ is\footnote{The subscript ${}_\pb$ refers to the pullback square \eqref{eq: def Sym(a) for a: * ---> A}.}
	\begin{equation}
		\sym(a) = \big( TM \oplus E[-1], \dd_\pb, [-,-] \big) \ ,
	\end{equation}
	whose differential and bracket are defined as follows:
	let $$a^{-1}\, \dd a \colon TM \longrightarrow TM \ltimes E$$ denote the map induced on dglas by the right vertical morphism of group objects in~\eqref{eq: def Sym(a) for a: * ---> A}.
	
	With $|\rho| \in \ZN$ denoting the degree of an element $\rho \in E[-1]$ in the shifted complex, the differential reads as
	\begin{alignat}{5}
		\dd_\pb(\rho) &= \dd_{E[-1]} \rho
		&& \quad
		&& \in \ (E[-1])^{|\rho| + 1}\ ,
		&& \qquad
		&& |\rho| \neq -1, 0\ ,
		\\[4pt]
		\dd_\pb(\rho) &= (0, \dd_{E[-1]}\rho)
		&& \quad
		&& \in \ TM \oplus E^{-1}\ ,
		&&\qquad
		&& |\rho| = -1\ ,
		\\[4pt]
		\dd_\pb(X, \rho) &= \dd_{E[-1]} \rho + a^{-1}\, \dd a(X)
		&& \quad
		&& \in \ E^0\ ,
		&&\qquad
		&& |\rho| = 0\ .
	\end{alignat}
	The brackets read as
	\begin{alignat}{3}
		\big[ (X, \rho), (Y, \mu) \big]
		&= \big( [X, Y], X \mu - Y \rho \big)\ ,
		&&\qquad
		&& \forall\ (X, \rho), (Y, \mu) \in TM \oplus (E[-1])^0\ ,
		\\[4pt]
		\big[ (X, \rho), \mu \big]
		&= X \mu\ ,
		&&\qquad
		&& \forall\ (X, \rho) \in TM \oplus (E[-1])^0\ ,\  \mu \in (E[-1])^n \ ,\, n \neq 0\ .
	\end{alignat}
\end{proposition}

\begin{proof}
	The cartesian square~\eqref{eq: def Sym(a) for a: * ---> A} in the $\infty$-category of group objects in $\Sh_\infty(\InfMfd)$ induces a cartesian square
	\begin{equation}
		\begin{tikzcd}[column sep = 2cm, row sep = 1.25cm]
			\rmB \cSym(a) \ar[r] \ar[d]
			& \rmB \cDiff(M) \ar[d, "{\phi \longmapsto (\phi,\, a^{-1}\, (\phi^{-1})^*a)}"]
			\\
			\rmB \cDiff(M) \ar[r, "{\phi \longmapsto (\phi, 1)}"']
			& \rmB \cDiff(M) \ltimes A
		\end{tikzcd}
	\end{equation}
	in $\Sh_\infty(\InfMfd)^{*/}_{\geqslant 0}$, the $\infty$-topos of pointed connected objects in $\Sh_\infty(\InfMfd)$.
	Forming the relative tangent complex with respect to the basepoint map of each vertex, we obtain a homotopy cartesian square
	\begin{equation}
		\begin{tikzcd}[column sep = 2cm, row sep = 1.25cm]
			\sym(a) \ar[r] \ar[d]
			& TM \ar[d, "{X \longmapsto (X,\, a^{-1}\, \dd a(X))}"]
			\\
			TM \ar[r, "{X \longmapsto (X, 0)}"']
			& TM \ltimes E
		\end{tikzcd}
	\end{equation}
	of dglas.
	Here $\sym(a)$ is the dgla associated to $\cSym(a)$, which we aim to compute.
	
	We can now make $\sym(a)$ explicit by computing this homotopy pullback of dglas.
	To that end, we replace the bottom horizontal morphism by a fibration:
	consider the factorisation of the morphism $TM \longrightarrow TM \ltimes E$ as
	\begin{equation}
		\begin{tikzcd}[column sep = 1.5cm]
			TM \ar[r, "{(\id_{TM} , 0)}"]
			& TM \oplus \Cone(\id_E) \ar[r, "\pr_{TM \ltimes E}"]
			& TM \ltimes E
		\end{tikzcd}
	\end{equation}
	where $\pr$ is the canonical projection.
	The brackets on $TM \oplus \Cone(\id_E)$ are given by the $TM$-action on $E$ and on $E[-1]$.
	Then we can compute $\sym(a)$ as a strict pullback of dglas, i.e.~the diagram of dglas
	\begin{equation}
		\begin{tikzcd}[column sep = 2cm, row sep = 1.25cm]
			\sym(a) \ar[r] \ar[d]
			& TM \ar[d, "{X \longmapsto (X,\, a^{-1}\, \dd a(X))}"]
			\\
			TM \oplus \Cone(\id_E) \ar[r, "\pr_{TM \ltimes E}"']
			& TM \ltimes E
		\end{tikzcd}
	\end{equation}
 is cartesian.	One now checks by a direct computation that $\sym(a)$ is the dgla described in the statement of Proposition~\ref{st: computing sym(a) in A//Diff(M) case}.
\end{proof}

\begin{example}
	We return to the example of deforming a map (cf. Example~\ref{eg:moduli and tangents of p-forms}) $$\omega \colon * \longrightarrow \Omega^p_v(M) \, \big/ \, \cDiff(M) \ . $$
	In this case $A = \Hom_v(M, \Omega^p)$, and so likewise $E = \Hom_v(M, \Omega^p)$ since $A$ is even a vector space.
	By Proposition~\ref{st: computing sym(a) in A//Diff(M) case} we obtain 
	\begin{equation}
		\sym(\omega) =
		\big( TM \xrightarrow{ \ \pounds_{(-)}\omega \ } \Omega^p(M)\ , [-,-] \big) \ ,
	\end{equation}
	with $TM$ lying in degree zero, and the brackets given by the action of $TM$ on itself and on $\Omega^p(M)$ via the Lie derivative.
	The 2-term dgla $\sym(\omega)$ can be viewed as the derived stabiliser of $\omega$ under the Lie algebra action of $TM$ on $\Omega^p(M)$.
	In particular
	\begin{equation}
		\rmH^0\big(\sym(\omega)\big)
		= \big\{ X \in TM\, \big| \, \pounds_X \omega = 0\big\}
	\end{equation}
	is the Lie algebra of the stabiliser subgroup of $\omega$ in $\cDiff(M)$.
	\qen
\end{example}

In the remainder of this section, we will look at other examples which suggest that there is an interesting notion of 
\emph{generalised higher connections} worthy of further investigation. These are order~$p$ splittings of the map $\sym(x) \longrightarrow TM$ for
\begin{equation}
	x \colon * \longrightarrow \mathrm{Moduli}(M,X)
\end{equation}
a point in a moduli stack of $X$-structures on $M$ which is not of the form $\Hom(M,X)/\cDiff(M)$ of an internal hom where $X$ is a qcdt stack.

\subsection{Infinitesimal symmetries and weak connections}

Building on the concept of vertical families of $p$-forms on $M$, one can construct moduli stacks of $(n{-}1)$-gerbes with $p$-form connection on $M$.
We denote these stacks by \smash{$\ul{\Grb}^{n-1}_{\nabla|p}(M)_v$}.
These moduli stacks further come endowed with a canonical action of $\cDiff(M)$ via pullback; see~\cite{BS:Higher_Syms_and_Deligne_Coho} for the details of such a construction.
The fact that the classifying stacks $\rmB^{n-p}\, \rmB^p_\nabla \rmU(1)$ for $(n{-}1)$-gerbes with $p$-form connections are abelian group objects is passed on to the moduli stacks \smash{$\ul{\Grb}^{n-1}_{\nabla|p}(M)_v$}.
Furthermore, the $\cDiff(M)$-action is compatible with this abelian group structure.
Thus we can form the semi-direct product $ \cDiff(M) \ltimes \ul{\Grb}^{n-1}_{\nabla|p}(M)_v $ and apply Proposition~\ref{st: computing sym(a) in A//Diff(M) case}.
We will now describe the resulting Lie algebra.

Consider an $(n{-}1)$-gerbe with $p$-form connection on $M$ given as an element \smash{$[g, \CA] \in \ul{\Grb}^{n-1}_{\nabla|p}(M)_v(*)$}. This defines a morphism
\begin{equation}\label{eq:gApoint}
	[g, \CA] \colon * \longrightarrow \ul{\Grb}^{n-1}_{\nabla|p}(M)_v \,\big/\, \cDiff(M)\ .
\end{equation}
We fix a good open cover $\CU$ of $M$ and a representative $(g, \CA)$ of the class $[g, \CA]$.
That is, $(g, \CA)$ is a \v{C}ech--Deligne cocycle on $M$ with respect to $\CU$, as spelled out in Section~\ref{sec: Con_l(g) equivalence}.

In the notation of Proposition~\ref{st: computing sym(a) in A//Diff(M) case} we then have 
\begin{equation}
	\label{eq: E for deforming [g, CA]}
	E^l = \begin{cases}
		 \ 0\ ,
		& l < -n\ ,
		\\[4pt]
		 \ \displaystyle\bigoplus\limits_{k = 0}^l\, \Omega^k(\cC_{n+l-k}\, \CU)\ ,
		& -n \leqslant l \leqslant -n+p\ ,
		\\[12pt]
		 \ \displaystyle\bigoplus\limits_{k = 0}^p\, \Omega^k(\cC_{n+l-k}\, \CU)\ ,
		& -n+p < l < 0\ ,
		\\[12pt]
		 \ \ker(\DD)\ ,
		& l = 0\ ,
		\\[4pt]
		 \ 0 & l \geqslant 0\ ,
	\end{cases}
\end{equation}
with differential 
$$
	\dd_E = \DD = (-1)^{l+k}\, \check\delta+\dd
$$
given as the differential on the totalisation of the \v{C}ech--de~Rham double complex associated to the cover $\CU$ and the complex of sheaves of abelian groups
\begin{equation}
	\begin{tikzcd}
		\Omega^0 \ar[r, "\dd"]
		& \Omega^1 \ar[r,"\dd"]
		& \cdots \ar[r, "\dd"]
		& \Omega^p \ar[r]
		&0\ .
	\end{tikzcd}
\end{equation}
The morphism denoted $a^{-1}\, \dd a$ in Proposition~\ref{st: computing sym(a) in A//Diff(M) case} reads as
\begin{equation}
	[g, \CA]^{-1}\, \dd [g, \CA] \colon TM \longrightarrow TM \ltimes E\ ,
	\quad
	X \longmapsto \big(X,\, \dd \log g(X),\, \pounds_X A^{{(1)}}, \ldots, \pounds_X A^{{(p)}}\big)\ .
\end{equation}

Proposition~\ref{st: computing sym(a) in A//Diff(M) case} thus computes the dgla which captures the deformations of the point \eqref{eq:gApoint}.
In particular, this dgla is concentrated in degrees $-n+1, \ldots,-1,0, 1$, has underlying complex $TM \oplus E[-1]$ with $E$ described in~\eqref{eq: E for deforming [g, CA]}, and brackets given by the Lie derivative.
Its differential is the differential of $E[-1]$, with the additional term $[g, \CA]^{-1}\, \dd [g, \CA]$ in degree zero.

\begin{remark}
	In the case where $p = 0$, i.e.~where we are considering an $(n{-}1)$-gerbe $[g]$ without connection on $M$, the complex $E$ from~\eqref{eq: E for deforming [g, CA]} reads as
	\begin{equation}
		E = \left(
		\begin{tikzcd}
			\Omega^0(\cC_0\, \CU) \ar[r, "\check\delta"]
			& \Omega^0(\cC_1\, \CU) \ar[r, "\check\delta"]
			& \cdots \ar[r, "\check\delta"]
			& \Omega^0(\cC_{n-1}\, \CU) \ar[r, "\check\delta"]
			& \ker (\check\delta)
		\end{tikzcd}
		\right)\ ,
	\end{equation}
	with $\ker(\check\delta)$ sitting in degree zero.
	Thus we obtain 
	\begin{equation}
		\sym(g) = \left(
		\begin{tikzcd}[column sep = 0.75cm]
			\Omega^0(\cC_0\, \CU) \ar[r, "\check\delta"]
			& \cdots \ar[r, "\check\delta"]
			& \Omega^0(\cC_{n-2}\, \CU) \ar[r, "\check\delta"]
			& TM \oplus \Omega^0(\cC_{n-1}\, \CU) \ar[rrr, " g^{-1}\, \dd g(-) + \check\delta"]
			& & & \ker (\check\delta)
		\end{tikzcd}
		\right)\ ,
	\end{equation}
	with the term containing $TM$ now in degree zero.
	
	Since the last differential here is surjective, the canonical morphism
	\begin{equation}
		\begin{tikzcd}[column sep = 2cm, row sep = 1cm]
			\Omega^0(\cC_0\, \CU) \ar[r, "(-1)^{n-1} \id"] \ar[r] \ar[d, "(-1)^{n-1}\,\check \delta",swap]
			& \Omega^0(\cC_0\, \CU) \ar[d, "\check\delta"]
			\\
			\vdots \ar[d, "(-1)^3\,\check \delta",swap]
			& \vdots \ar[d, "\check\delta"]
			\\
			\Omega^0(\cC_{n-2}\, \CU) \ar[d, "(-1)^2\,\check \delta",swap] \ar[r, "\id"]
			& \Omega^0(\cC_{n-2}\, \CU) \ar[d, "\check\delta"]
			\\
			\big\{ (X, f) \in TM \oplus \Omega^0(\cC_{n-1}\, \CU) \, \big| \, \dd \log g(X) - \check\delta f = 0 \big\} \ar[d] \ar[r, hookrightarrow]
			& TM \oplus \Omega^0(\cC_{n-1}\, \CU) \ar[d, "g^{-1}\, \dd g(-) + \check\delta"]
			\\
			0 \ar[r]
			& \ker  (\check\delta)
		\end{tikzcd}
	\end{equation}
	where the signs of the horizontal morphisms $\id$ alternate, is a weak equivalence from the dgla $\cC \At(g)$; this morphism respects the brackets because only brackets with elements in degree zero are non-trivial.
	That is, there is a canonical equivalence
	\begin{equation}
		\sym(g) \simeq \cC\At(g) 
	\end{equation}
	of dglas.
	\qen
\label{rem:symgCAt}\end{remark}

Forgetting the respective highest form degrees of the connection data provides a chain of moduli stacks
\begin{equation}
	\ul{\Grb}^{n-1}_\nabla(M)_v \longrightarrow \ul{\Grb}^{n-1}_{\nabla|n-1}(M)_v
	 \longrightarrow \ul{\Grb}^{n-1}_{\nabla|n-2}(M)_v
	 \longrightarrow \cdots
	 \longrightarrow \ul{\Grb}^{n-1}_{\nabla|1}(M)_v
	 \longrightarrow \ul{\Grb}^{n-1}(M)_v\ .
\end{equation}
These forgetful maps are compatible with the abelian group structure and $\cDiff(M)$ actions, and so we obtain an associated chain of dglas
\begin{align}\label{eq:symchain}
	\sym (g, \CA)
	 \longrightarrow \sym^{\leqslant n-1} (g, \CA)
	& \longrightarrow \sym^{\leqslant n-2} (g, \CA) \longrightarrow  \cdots
	 \longrightarrow \sym^{\leqslant 1} (g,  \CA)
	 \longrightarrow \sym (g) \ ,
\end{align}
where $[g, \CA]$ is an $(n{-}1)$-gerbe with full connection, $\sym^{\leqslant p}(g,\CA):=\sym(g,\fgt_{p+1,\dots,n}\,\CA)$ and $\fgt_{p+1, \ldots, n}$ is the morphism which forgets form degrees $n, n-1, \ldots, p+1$ of the connection data, for $0\leqslant p\leqslant n-1$.

By combining Theorem~\ref{st:strictification of At(CG)} together with Remark~\ref{rem:higherCourantU1}, Theorem~\ref{prop:Linftyproj} and Remark~\ref{rem:symgCAt}, we are then led to formulate

\begin{conjecture} \label{conj:Coursym}
	Let $H \in \Omega_{\rm cl}^{n+1}(M)$ be the field strength of the $(n{-}1)$-gerbe connection $\CA$.
	There is a canonical equivalence
	\begin{equation}
		\mathrm{Cour}_{n-1}(H) \weq \sym^{\leqslant n-1} (g, \CA)
	\end{equation}
	of $L_\infty$-algebras over $TM$.
\end{conjecture}

Conjecture~\ref{conj:Coursym} is supported by Example~\ref{eg: deforming gerbe with connective structure} below, which shows that the conjecture is true for $n = 2$. It also agrees with the general expectations sketched in \cite[Example~3.2.8]{FRS2}.

\begin{remark}
From the discussion preceding Theorem~\ref{prop:Linftyproj}, it is tempting to further conjecture an analogous equivalence between the $L_\infty$-algebras \smash{$\mathrm{Cour}_{n-1}^{\leqslant p}(H)$} of Proposition~\ref{prop:cour-intermediate} and $\sym^{\leqslant p}(g,\CA)$ for all $0\leqslant p\leqslant n-1$. For $p=0$ the conjecture is true by Lemma~\ref{lem:CourpAt}. This would identify $\mathrm{Cour}_{n-1}^{\leqslant p}(H)$ as the $L_\infty$-algebra controlling the infinitesimal symmetries of an $(n{-}1)$-gerbe with $p$-form connection, for each $0\leqslant p\leqslant n-1$. However, in general we lack a geometric interpretation of $\mathrm{Cour}_{n-1}^{\leqslant p}(H)$ as an extension of some algebroid structure arising from the data of an $(n{-}1)$-gerbe with $p$-form connection, so there is currently no basis for such an expectation. 
\qen
\end{remark}

\begin{conjecture} \label{conj:CourAtsym}
There is a commutative diagram of $L_\infty$-algebras over $TM$
\begin{equation}
\begin{tikzcd}[row sep = 1cm]
\mathrm{Cour}_{n-1}(H) \ar[d] \ar[r] & \sym^{\leqslant n-1}(g,\CA) \ar[d] \\
\At_{n-1}(H) \ar[r] & \cC\At(g)
\end{tikzcd}
\end{equation}
where the top horizontal arrow is the weak equivalence of Conjecture~\ref{conj:Coursym}, the bottom horizontal arrow is induced by the $\infty$-equivalence from Theorem~\ref{st:strictification of At(CG)}, the left vertical arrow is the canonical $\infty$-projection of Theorem~\ref{prop:Linftyproj}, and the right vertical arrow is induced from the composition of the weak equivalence from Remark~\ref{rem:symgCAt} with the chain of morphisms in \eqref{eq:symchain}.
\end{conjecture}

\begin{remark}
For $n=2$, the statement of Conjecture~\ref{conj:CourAtsym} is similar to that of~\cite[Proposition~5.2.6]{FRS}. Conjectures~\ref{conj:Coursym} and~\ref{conj:CourAtsym} are evidently true for $n=1$ by Example~\ref{ex:AtiyahLie}.
\qen
\end{remark}

\begin{example}
	\label{eg: deforming U(1)-bundle with connection}
	Let $[g, A]$ be a $\rmU(1)$-bundle with connection on $M$.
	Its associated dgla reads as
	\begin{equation}
		\sym(g, A) =
		\left(
		\begin{tikzcd}[column sep = 0.75cm]
			TM \oplus \Omega^0(\cC_0\, \CU)
			\ar[r, "\dd_\pb"]
			& \big\{ (\omega^0_1, \omega^1_0) \in \Omega^0(\cC_1\, \CU) \oplus \Omega^1(\cC_0\, \CU)\, \big| \,
			\check\delta \omega^0_1 = 0\ ,\, \dd \omega^0_1 - \check\delta \omega^1_0 = 0 \big\}
		\end{tikzcd}
		\right)
	\end{equation}
	where the first term sits in degree zero.
	Explicitly, the differential reads as
	\begin{equation}
		\dd_\pb (X, \omega^0_0) = (-\check\delta \omega^0_0 + \dd \log g(X),\, -\dd \omega^0_0 + \pounds_X A)\ ,
	\end{equation}
	where the minus signs stem from the shift of $E$. 
	
	Note that $\sym(g,A)$ is not an $L_\infty$-algebroid over $M$, but it is an $L_\infty$-algebra over the Lie algebra $\Gamma(M, TM)$ (it lives purely on the right-hand side of the square in Theorem~\ref{thm:global}).
	Consequently, we cannot ask for $p$-form connections on $[g, A]$ in the sense of Definition~\ref{def: Connection}, but---inspired by that definition---we can still ask for $C^\infty(M)$-linear 1-splittings of the canonical projection $\sym(g, A) \longrightarrow TM$.
	These splittings could still be interpreted as 1-form connections on the bundle with 1-form connection $[g,A]$ in some weak sense.
	
	Such a 1-splitting consists of a $C^\infty(M)$-linear map $TM \longrightarrow TM \oplus \Omega^0(\cC_0\, \CU)$ which is the identity on $TM$; that is, it is given by a locally defined 1-form $a \in \Omega^1(\cC_0\, \CU)$.
	This map has to satisfy a Maurer--Cartan condition, which in this case unravels to the condition that
	\begin{equation}
		\dd_\pb \circ (\id_{TM}, a) = 0\ .
	\end{equation}
	Explicitly, for $X \in TM$ we compute
	\begin{equation}
		\dd_\pb \circ (\id_{TM}, a)(X)
		= \big( - \check\delta a(X) + \dd \log g(X),\,
		- \dd\, \iota_X a + \pounds_X A \big)\ .
	\end{equation}
	
	Vanishing of the first component is equivalent to $a \in \Omega^1(\cC_0\, \CU)$ defining a connection on the bundle $[g]$.
	For the second component, we can thus write $a = A + \check\delta \alpha$, where $\alpha \in \Omega^1(M)$ is a globally defined 1-form on $M$.
	It then remains to solve the equation
	\begin{equation}
		0 = - \dd\, \iota_X (A + \alpha) + \pounds_X A
		\qquad \Longleftrightarrow \qquad
		\iota_X\, \dd A - \dd\, \iota_X \alpha = 0\ ,
	\end{equation}
	for all vector fields $X$.
	Since the first term is $C^\infty(M)$-linear in $X$, whereas the second term is not, the only way to solve the equation is with $\dd A = 0$ and $\alpha = 0$.
	It follows that $\sym(g, A)$ admits a $C^\infty(M)$-linear 1-splitting only if its connection $A$ is flat, and in that case the only such splitting is given by $A$ itself.
	
	Heuristically, this can be understood as stating that a bundle with a 1-form connection, seen as a geometric object, admits a 1-form connection
	compatible with its existing 1-form connection only if the latter is flat.
	Moreover, even then the only admissible additional 1-form connection coincides with the existing one (which is compatible with itself precisely because it is flat).
	\qen
\end{example}

\begin{example}
	\label{eg: deforming gerbe with connective structure}
	We now consider a 1-gerbe with 1-form connection $[g, A^{{(1)}}] \in \Grb^1_{\nabla|1}(M)$; a 1-form connection on a 1-gerbe is more commonly called a \textit{connective structure}.
	By Proposition~\ref{st: computing sym(a) in A//Diff(M) case} and \eqref{eq: E for deforming [g, CA]} we find that
	\begin{equation}
		\sym(g, A^{{(1)}}) = \left(
		\begin{tikzcd}
			\Omega^0(\cC_0\, \CU) \ar[r,"\dd_\pb"]
			& TM \oplus \Omega^0(\cC_1\, \CU) \oplus \Omega^1(\cC_0\, \CU) \ar[r,"\dd_\pb"]
			& \ker(\DD)
		\end{tikzcd}
		\right)\ ,
	\end{equation}
	with the first term in degree $-1$.
	Its bracket is given by the Lie action of vector fields on functions and differential forms.
	
	We claim that this dgla is equivalent to its sub-dgla $\sym'(g, A^{{(1)}})$ with underlying cochain complex
	\begin{equation}
		\begin{tikzcd}
			\Omega^0(\cC_0\, \CU) \ar[r,"\dd_\pb"]
			& \ker(\dd_\pb)\ ,
		\end{tikzcd}
	\end{equation}
	where explicitly
	\begin{align}
		\ker(\dd_\pb)
		= \big\{ &(X, \omega^0_1, \omega^1_0) \in TM \oplus \Omega^0(\cC_1\, \CU) \oplus \Omega^1(\cC_0\, \CU)\, \big|
		\\
		&\hspace{4cm} - \check\delta \omega^0_1 + \dd \log g(X) = 0\ ,\, - \dd \omega^0_1 + \check\delta \omega^1_0 + \pounds_X A^{(1)} = 0 \big\}\ .
	\end{align}
	This amounts to showing that the map $TM \oplus \Omega^0(\cC_1\, \CU) \oplus \Omega^1(\cC_0\, \CU) \longrightarrow \ker(\DD)$ is surjective.
	
	To that end, recall that an element in $\ker(\DD)$ is a pair
	\begin{equation}
		(\omega^0_2, \omega^1_1) \ \in \ \Omega^0(\cC_2\, \CU) \oplus \Omega^1(\cC_1\, \CU)
		\end{equation}
		such that
		\begin{equation}
		\check\delta \omega^0_2 = 0
		\qquad , \qquad
		- \dd \omega^0_2 + \check\delta \omega^1_1 = 0\ .
	\end{equation}
	Since $\check\delta \omega^0_2 = 0$, there exists  $\mu^0_1 \in \Omega^0(\cC_1\, \CU)$ such that $\omega^0_2 = \check\delta \mu^0_1$.
	In particular, each pair $(\omega^0_2, \omega^1_1)$ is cohomologous to a pair of the form $(0, \omega^1_1)$ (for some new $\omega^1_1$) via the image of the element $(0, \mu^0_1, 0) \in TM \oplus \Omega^0(\cC_1\, \CU) \oplus \Omega^1(\cC_0\, \CU)$ under $\dd_\pb$.
	An element $(0, \omega^1_1) \in \ker(\DD)$ now satisfies  $\check\delta \omega^1_1 = 0$, so that we find  $\mu^1_0 \in \Omega^1(\cC_0\, \CU)$ with $\check\delta \mu^1_0 = \omega^1_1$.
	The residual equivalence relation remaining after restricting to the subspace of $\ker(\DD)$ consisting of pairs $(0, \omega^1_1)$ is generated by the elements of $TM \oplus \Omega^0(\cC_1\, \CU) \oplus \Omega^1(\cC_0\, \CU)$ of the form $(X, \omega^0_1, \omega^1_0)$ with $-\check\delta \omega^0_1 + \dd \log g(X) = 0$.
	In particular, the element $(0, 0, \mu^1_0)$ establishes the identity $[0, \omega^1_1] = 0$ in cohomology.
	
	This concludes the argument, showing that the canonical morphism
	\begin{equation}
		\begin{tikzcd}[column sep = 1cm, row sep = 1cm]
			\Omega^0(\cC_0\, \CU) \ar[r,"\id"] \ar[d, "\dd_\pb",swap]
			& \Omega^0(\cC_0\, \CU) \ar[d, "\dd_\pb"]
			\\
			\ker(\dd_\pb) \ar[r,hookrightarrow] \ar[d]
			& TM \oplus \Omega^0(\cC_1\, \CU) \oplus \Omega^1(\cC_0\, \CU) \ar[d, "\dd_\pb"]
			\\
			0  \ar[r]
			& \ker(\DD)
		\end{tikzcd}
	\end{equation}
	is a quasi-isomorphism of dglas.
	
	The dgla $\sym'(g, A^{{(1)}})$ coincides with the dgla constructed by Collier from a gerbe with connective structure in~\cite[Definition~10.38 and Proposition~10.48]{Collier:Inf_symmetries_of_gerbes}.
	Furthermore, by~\cite[Theorem~12.50]{Collier:Inf_symmetries_of_gerbes} this 2-term dgla is canonically equivalent to the 2-term $L_\infty$-algebra $\mathrm{Cour}_1(H)$ of sections of the Courant algebroid associated to $(g, A^{{(1)}})$ via the generalised tangent bundle construction of Remark~\ref{rem:higherCourantU1}.
	
	Gerbes with 1-form connection on $M$ are classified by maps $[g,A^{{(1)}}] \colon M \longrightarrow \rmB\,\rmB_\nabla\rmU(1)$.
	Since $\rmB_\nabla\rmU(1)$ is an abelian group object, these classifying maps are of the form $x \colon M \longrightarrow \rmB G$ for some group object $G$, and we can thus understand a 1-gerbe with connective structure on $M$ as a $\rmB_\nabla\rmU(1)$-principal $\infty$-bundle on $M$.
	Note, however, that $\rmB_\nabla\rmU(1)$ does not have deformation theory (see Example~\ref{eg: forms on M have no DefThy}), and so there are no connections on $[g,A^{{(1)}}]$ in the sense of Definition~\ref{def: Connection}.
	There may, however, be connections in a weaker sense:
	as in Example~\ref{eg: deforming U(1)-bundle with connection}, we conclude this example by investigating the $C^\infty(M)$-linear 1-splittings of the projection $\sym(g,A^{(1)}) \longrightarrow TM$ of $L_\infty$-algebras (rather than of $L_\infty$-algebroids).
	
	Such a 1-splitting is a $C^\infty(M)$-linear map
	\begin{equation}
		(\id_{TM}, a, b) \colon TM \longrightarrow TM \oplus \Omega^0(\cC_1\, \CU) \oplus \Omega^1(\cC_0\, \CU)
	\end{equation}
	over $TM$ satisfying the Maurer--Cartan condition, which in this case reads as
	\begin{equation}
		\dd_\pb \circ (\id_{TM}, a, b) = 0\ .
	\end{equation}
	Here $a \in \Omega^1(\cC_1\, \CU)$ and $b \in \Omega^1(\cC_0 \,\CU, T\,\CU)$, or equivalently $b$ is a family of locally defined rank~2 tensors.
	Note that, in general, $b$ has no symmetry properties (it is neither symmetric nor antisymmetric as a map $T^\vee \cC_0\, \CU \otimes T^\vee \cC_0\, \CU \longrightarrow C^\infty(\cC_0\, \CU)$).
	Explicitly, given a vector field $X \in TM$, the condition $\dd_\pb \circ (\id_{TM}, a, b) = 0$ amounts to
	\begin{equation}
		- \check\delta a(X) + \dd \log(g)(X) = 0
		\qquad , \qquad
		\check\delta \big( b(X) \big) - \dd \big( a(X) \big) + \pounds_X A^{(1)} = 0\ .
	\end{equation}
	
	The first equation states that $a$ is a connective structure on the gerbe defined by the $\rmU(1)$-cocycle $g$.
	We can thus write $a = A^{(1)} + \check\delta \alpha$, for some $\alpha \in \Omega^1(\cC_0\, \CU)$.
	The second equation then becomes
	\begin{equation}
		\check\delta \big( b(X) \big) + \dd A^{{(1)}}(X) - \dd\, \iota_X \alpha = 0\ ,
	\end{equation}
	for all vector fields $X$.
	Since the first two terms are $C^\infty(M)$-linear in $X$, but the third is not, the only way to satisfy this identity for all $X \in TM$ simultaneously is by demanding that
	\begin{equation}
		\check\delta \big( b(X) \big) + \dd A^{{(1)}}(X) = 0
		\qquad , \qquad
		 \dd\, \iota_X \alpha = 0\ .
	\end{equation}
	Here the second identity states that $\iota_X \alpha = \alpha(X)$ is locally constant, for each $X \in TM$.
	Since $\alpha$ is $C^\infty(M)$-linear, the only solution is $\alpha = 0$.
	That is, the connective structure $a$ necessarily coincides with the connective structure $A^{{(1)}}$.
	Finally, the first identity is equivalent to saying that $b$ is an \emph{enhanced curving} on $(g, A^{{(1)}})$ in the sense of~\cite[Definition~3.3]{TD:Chern_corr_for_higher_PrBuns}.
	
	This provides a new and generalised perspective on enhanced curvings.
	A particularly interesting point is that, if we split the 2-tensor $b = A^{{(2)}} + \gamma$ into its antisymmetric and symmetric parts, then $A^{{(2)}}$ is a classical curving on $(g, A^{{(1)}})$, whilst the symmetric part $\gamma$ satisfies $\check\delta \gamma = 0$ and thus defines a global symmetric 2-tensor field on $M$.
	If we impose the additional condition on $b$ that its symmetric part $\gamma$ be positive definite, then the 1-splitting---or equivalently enhanced curving---$b$ encodes simultaneously the $B$-field and Riemannian metric of, for instance, NSNS supergravity.
	This was one of the motivations for introducing enhanced curvings in~\cite{TD:Chern_corr_for_higher_PrBuns}.
	\qen
\end{example}

Recalling the discussion of Remark~\ref{rem:String2group}, Example~\ref{eg: deforming gerbe with connective structure} as well as results of~\cite{TD:Chern_corr_for_higher_PrBuns} motivate the following conjecture: 
it is known~\cite{TD:Chern_corr_for_higher_PrBuns, Waldorf:Multiplicative_Gerbes, Waldorf:String_and_Chern-Simons} that the extension $$\rmB^2 \rmU(1) \longrightarrow \String(K) \longrightarrow K$$ has an enhancement to an extension $$\rmB\, \rmB_\nabla \rmU(1) \longrightarrow \String_{\rm CS}(K) \longrightarrow K \ , $$ where $\rmB_\nabla \rmU(1)$ is the classifying stack for $\rmU(1)$-bundles with connection.
A classifying map $$P_K:M \longrightarrow \rmB K$$ thus induces a gerbe \textit{with connective structure} on the principal $K$-bundle $P_K$, which is the obstruction to a geometric string structure.

\begin{conjecture}\label{conj:hetSuGra}
	Let $\Hom_v(M, \rmB^2\, \rmB_\nabla \rmU(1))$ and $\Hom_v(M, \rmB \String_{\rm CS}(K))$ denote the full substacks of the respective internal homs where all the 1-forms induced via pullback by the connective structure of the multiplicative gerbe on $K$ are \textit{vertical}, i.e.~they only have form legs along $M$ (compare Example~\ref{ex:vertical forms def thy}).
	Then:
	\begin{myenumerate}
		\item The stacks $\Hom_v(M, \rmB^2\, \rmB_\nabla \rmU(1))/\cDiff(M)$ and $\Hom_v(M, \rmB \String_{\rm CS}(K))/\cDiff(M)$ have (possibly quasi-coherent) deformation theory.
		
		\item A lift of $P_K:M \longrightarrow \rmB K$ to $\rmB \String_{\rm CS}(K)$ exists if and only if a lift to $\rmB \String(K)$ exists.
		
		\item The $L_\infty$-algebra $T\big(* /\, (\Hom_v(M, \rmB \String_{\rm CS}(K))/\cDiff(M))\big)$ has a canonical map to the Lie algebra $TM$, and a $C^\infty(M)$-linear splitting of this map (as dg $C^\infty(M)$-modules) is the same as a geometric string structure on $P_K$ together with a symmetric 2-tensor $\gamma \in \Gamma(M, \midvee^2\, T^*M)$.
	\end{myenumerate}
\end{conjecture}

Compared to Example~\ref{eg: deforming gerbe with connective structure}, this allows one to add the data of a principal bundle connection to the $B$-field and Riemannian metric, encoding in this way the bosonic fields of heterotic supergravity; indeed, the modified Bianchi identity for the $H$-flux in this case is captured by~\cite[Lemma~3.2.4]{Waldorf:String_and_Chern-Simons}. It is also of use in 11-dimensional supergravity, wherein string connections appear as the ${\rm E}_8$-model for the 3-form $C$-field of M-theory~\cite{Diaconescu:2003bm,Freed:2004yc,Fiorenza:2012mr}.


\begin{appendix}


\numberwithin{equation}{section}


\section{Proof of Theorem~\ref{thm:formal stack vs lie}}
\label{sec:lie diff proof}


This appendix is devoted to a proof of Theorem~\ref{thm:formal stack vs lie}.
Let us start by introducing an intermediate $\infty$-category of $C^\infty$-geometric formal moduli problems, which sits between $\Sh_\infty(\InfMfd)$ and $L_\infty\Agd_M$.

\begin{definition}
	\label{def:FMP}
	Let $M$ be a smooth manifold. Let $\Art_{/C^\infty(M)}$ be the smallest full $\infty$-subcategory of $\CRing^\infty_{/C^\infty(M)}$ satisfying the following two conditions:
	\begin{myenumerate}
		\item It contains the terminal object $\id\colon C^\infty(M)\longrightarrow C^\infty(M)$.
		
		\item It is closed under square-zero extensions by shifted copies of $C^\infty(M)$ (in non-positive degree). More precisely, for any $n\geqslant 0$ and any pullback square in $\CRing^\infty_{/C^\infty(M)}$ of the form
		\begin{equation}
			\label{diag:pullback art}
			\begin{tikzcd}[row sep = 1cm , column sep = 1cm]
				R_\eta \arrow[d] \arrow[r]
				& C^\infty(M) \arrow[d, "{(\id, 0)}"]
				\\
				R \arrow[r, "\eta"]
				& {C^\infty(M) \oplus C^\infty(M)[n+1]} \arrow[r, "\pi_1"]
				& C^\infty(M)
			\end{tikzcd}
		\end{equation}
		such that $R \longrightarrow C^\infty(M)$ is in $\Art_{/C^\infty(M)}$, the composition $R_\eta\longrightarrow C^\infty(M)$ is also in $\Art_{/C^\infty(M)}$.
	\end{myenumerate}
	
	The $\infty$-category of \emph{formal moduli problems} under $M$ is the full $\infty$-subcategory $$\FMP_{M/}\subseteq\Fun(\Art_{/C^\infty(M)}, \scS)$$  on the functors $F\colon \Art_{/C^\infty(M)}\longrightarrow \scS$ with the following properties:
	\begin{myenumerate}
		\item[(a)] $F\big(C^\infty(M)\big)\simeq \ast$.
		
		\item[(b)] $F$ preserves the pullback squares \eqref{diag:pullback art}.
	\end{myenumerate}
\end{definition}

Informally, the definition of a formal moduli problem under a smooth manifold $M$ is very similar to that of a prestack with deformation theory (Definition~\ref{def:qcoh def thy}), except that a formal moduli problem is only defined on a category of infinitesimal thickenings of $M$ and not on all inf-manifolds. In particular, various results and constructions from Section~\ref{sec:infinitesimal} have analogues for formal moduli problems.

\begin{example}
	\label{ex:prorep FMP}
	Every $A\in \Art_{/C^\infty(M)}$ determines a \emph{representable} formal moduli problem, which by a slight abuse of notation we write as
	\begin{equation}
		\Spf(A) \colon \Art_{/C^\infty(M)} \longrightarrow \scS\ ,
		\qquad B \longmapsto \Art_{/C^\infty(M)}\big(A, B\big)\ .
	\end{equation}
	We think of $\Spf(A)$ as the `formal spectrum' of $A$. For instance, applying this to $C^\infty(M)$ itself yields the initial formal moduli problem under $M$, which we will simply denote by $M$. Every other formal moduli problem $F$ comes with a canonical map $x\colon M\longrightarrow F$, corresponding to the unique point of $F(C^\infty(M))$.
	
	Likewise, any filtered colimit of such representable formal moduli problems is again a formal moduli problem; following~\cite[Definition 1.5.3]{Lurie:DAG-X}, we will call these \emph{pro-representable} formal moduli problems under $M$.
	\qen
\end{example}

\begin{example}
	\label{ex:function algebras FMP}
	Similarly to Example~\ref{ex:function algebras}, there is an adjoint pair of functors
	$$
	\begin{tikzcd}
		\coSpf : \cdga_{\RN/C^\infty(M)} \arrow[r, yshift=1ex]
		& \big(\FMP_{M/}\big)^\opp : \Gamma(-, \calO)\ . \arrow[l, yshift=-1ex, "\perp"']
	\end{tikzcd}
	$$
	The left adjoint sends each $R\longrightarrow C^\infty(M)$ to the formal moduli problem
	\begin{equation}
		\coSpf(R) \colon \Art_{/C^\infty(M)} \longrightarrow \scS\ ,
		\qquad
		R' \longmapsto \cdga_{\RN/C^\infty(M)}\big(R, R'\big)\ ,
	\end{equation}
	where on the right-hand side we forget from the $C^\infty$-ring structure on $R$ to its underlying cdga.
	On the other hand,  it follows that $$\Gamma(F, \calO_F)\simeq \underset{\alpha}\holim\ R_\alpha$$ whenever$$F = \underset{\alpha}\hocolim\ \Spf(R_\alpha) $$  is a pro-representable formal moduli problem.
	\qen
\end{example}

Let $F$ be a formal moduli problem, and let $x\colon M\longrightarrow F$ be the map from the initial formal moduli problem. For every $C^\infty(M)$-module $I$ given as a finite complex of finite free $C^\infty(M)$-modules in non-positive degrees, the square-zero extension $C^\infty(M)\oplus I$ defines an object in $\Art_{/C^\infty(M)}$. By the discussion preceding Definition~\ref{def:tangent} there exists a unique object $T_xF$ in the $\infty$-category $\Mod^\infty_{C^\infty(M)}$ equipped with a natural equivalence
$$
\DK\big(T_x F \otimes^h_{C^\infty(M)} I \big)
\simeq F(C^\infty(M)\oplus I)
$$
for any such $C^\infty(M)$-module $I$. When $F = M$ is the initial formal moduli problem, this is simply the tangent bundle $TM$. Consequently, every formal moduli problem $F$ defines a map of dg $C^\infty(M)$-modules $TM \longrightarrow T_x F$, and we define
$$
T(M/F) = \hofib (TM \longrightarrow T_x F)
$$
to be its homotopy fibre. By definition this comes endowed with a natural anchor map $$T(M/F) \longrightarrow TM . $$

\begin{theorem}[{\cite[Theorem 4.2.1 and Remark 4.2.22]{Nuiten:Thesis}}]
	\label{thm:fmp vs lie}
	There is an equivalence of $\infty$-categories
	$$
	\FMP_{M/} \eq L_\infty\Agd_M^\infty\ ,
	\qquad
	F \longmapsto T(M/F)\ ,
	$$
	sending each formal moduli problem under $M$ to $T(M/F)$, endowed with a certain $L_\infty$-algebroid structure on it. 
	Furthermore, if $F$ is a formal moduli problem under $M$, and $\frg$ is the associated $L_\infty$-algebroid, then there is a natural equivalence
	$$
	\Gamma(F, \calO_F) \simeq \ChEil^*(\frg)\ .
	$$
\end{theorem}

In an algebraic setting, this also appears as~\cite[Theorem 5.1]{Nuiten:Koszul_duality_for_Lie_algebroids}; the proof from there carries over verbatum to the $C^\infty$-ring setting used here.

We will deduce Theorem~\ref{thm:formal stack vs lie} from Theorem~\ref{thm:fmp vs lie} by first relating the $\infty$-category of formal moduli problems under $M$ to the $\infty$-category of morphisms $x\colon M\longrightarrow X$ to a dt prestack $X$, and then passing to associated stacks. For the first step, we consider the canonical fully faithful functor
$$
\begin{tikzcd}
	j \colon \big(\Art_{/C^\infty(M)}\big)^{\opp} \arrow[hookrightarrow, r]
	& \PSh_\infty(\InfMfd)_{M/}\ ,
	&\big( R\longrightarrow C^\infty(M)\big) \arrow[r, mapsto]
	& \big( M\longrightarrow \Spec(R) \big)
\end{tikzcd}
$$
induced by the inclusion $\Art_{/C^\infty(M)} \longhookrightarrow  \CRing_{/C^\infty(M)}$.
It induces an adjoint pair
\begin{equation}
	\label{diag:j*}
	\begin{tikzcd}
		j_! : \PSh_\infty\big(\Art_{/C^\infty(M)}^{\opp}\big) \arrow[r, yshift=1ex, "\perp"']
		& \PSh_\infty(\InfMfd)_{M/} \arrow[l, yshift=-1ex] : j^*
	\end{tikzcd}
\end{equation}
where $j_!$ is given by left Kan extension and $j^*$ sends every object $x \colon M\longrightarrow X$ to the functor
$$
j^*X(R) = X(\Spec(R))\times^h_{X(M)} \{x\}\ .
$$

\begin{lemma}
	\label{lem:j*}
	Let $M$ be a smooth manifold.
	\begin{myenumerate}
		\item If $x\colon M\longrightarrow X$ is a map to a dt prestack, then $j^*X$ is a formal moduli problem. Furthermore, the tangent complex of $X$ at $M$ coincides with the tangent complex of the associated formal moduli problem, i.e.\ there is a natural equivalence of $C^\infty(M)$-modules
		$$
		T_xX\simeq T_x(j^*X) \ .
		$$
		
		\item Let $x\colon M\longrightarrow X$ and $f\colon X\longrightarrow Y$ be maps of qcdt stacks inducing equivalences $M\simeq X_\red\simeq Y_\red$. Then $f$ is an equivalence if and only if $j^*(f)\colon j^*X\longrightarrow j^*Y$ is an equivalence.
		In other words, the functors $(-)_\red$ and $j^*$ jointly detect equivalences in \smash{$\QCDT_{M/}$}.
	\end{myenumerate}
\end{lemma}

\begin{proof}
	Assertion~(1) follows directly from the formula for $j^*$. 
	Assertion~(2) follows from Proposition~\ref{prop:inverse_function}: indeed, every map $x'\colon N\longrightarrow X$ from a smooth manifold factors over $x\colon M\longrightarrow X$. The conditions on $X$ and $Y$ then guarantee that the map $T_{x'}X\longrightarrow T_{f(x')}Y$ arises from $T_xX\longrightarrow T_{f(x)}Y$ by tensoring up and completing (Proposition~\ref{prop:qcdt stack 1}). Since $T_xX\longrightarrow T_{f(x)}Y$ is an equivalence, the result follows.
\end{proof}

We will now study the left adjoint $j_!$ in more detail.

\begin{proposition}\label{prop:j!}
	Let $F\colon \Art_{/C^\infty(M)}\longrightarrow \scS$ be a formal moduli problem under $M$.
	\begin{myenumerate}
		\item The map $M\longrightarrow j_!F$ induces an equivalence $M\simeq (j_!F)_\red$.
		
		\item The prestack $j_!F$ is a \emph{strict} qcdt prestack (see Definition~\ref{def:qcoh def thy}).
		
		\item The unit map $F\longrightarrow j^*\,j_!F$ is an equivalence.
	\end{myenumerate}
\end{proposition}

\begin{proof}
	We will prove all assertions in increasing level of generality. Let us start by treating the case where $F$ is representable, i.e.~$F = \Spf(A)$ for some $A\in \Art_{/C^\infty(M)}$. In this case, $j_!F$ is the representable sheaf $\Spec(A)$  (recall from Proposition~\ref{prop:subcanonical} than $\InfMfd$ is subcanonical), so that $(j_!F)_\red\simeq M$ and $j^*\,j_!F\simeq F$ since $j$ is fully faithful. The representable sheaf $j_!F = \Spec(A)$ satisfies the Schlessinger condition, as it sends any pushout of inf-manifolds to a pullback. 
	
	Using Lemma~\ref{lem:qcdt}, the strong tangent representability for $j_!F$ is a consequence of the following observation: for any map $y\colon N\longrightarrow \Spec(A)$ from a smooth manifold and a bounded complete non-positively graded $C^\infty(N)$-module $I$, we need to check that the natural map
\begin{equation} \label{eq:TySpec}
	\begin{tikzcd}
		\DK \big( T_y\Spec(A) \otimes^h_{C^\infty(N)} I\big) \arrow[r]
		& \Der_{\Spec(A)}(y; I)
	\end{tikzcd}
\end{equation}
	is an equivalence. The spaces in \eqref{eq:TySpec} can be identified explicitly in terms of the $C^\infty$-algebraic cotangent complex $L_A$ of the dg $C^\infty$-ring $A$, as in Examples~\ref{ex:corep tangent} and~\ref{ex:mfd are qcdt}. Since the canonical point $x\colon M\longrightarrow \Spec(A)$ induces an isomorphism on reductions, the point $y$ factors through $x$ as $$N \longrightarrow M \longrightarrow \Spec(A) \ , $$ for some smooth map of manifolds $f\colon N\longrightarrow M$. Using this decomposition, the map \eqref{eq:TySpec} can be identified in a purely algebraic way as the Dold--Kan image of
	\begin{equation}
		\label{eq:comparison cotangent}\begin{tikzcd}
			\Hom_{C^\infty(M)} \big( L_A\otimes_A C^\infty(M), C^\infty(N) \big) \otimes^h_{C^\infty(N)} I \arrow[r]
			& \Hom_{C^\infty(M)} \big( L_A\otimes_A C^\infty(M), f_*I \big)\ .
	\end{tikzcd}\end{equation}
		
	For the \emph{canonical} basepoint $x:M\longrightarrow \Spec(A)$, there is an explicit model for $L_A\otimes_A C^\infty(M)$. Indeed, by~\cite[Corollary 4.2.18]{Nuiten:Thesis} there exists a dg Lie algebroid $\frg$ whose underlying complex is a projectively cofibrant, negatively graded complex of finite free $C^\infty(M)$-modules  such that $A\simeq \ChEil^*(\frg)$. Using this and~\cite[Lemmas 4.1.33 and 4.1.35]{Nuiten:Thesis}, one finds 
	\begin{align*}
		L_A\otimes_A C^\infty(M)
		&\cong \Omega^1(M) \oplus \frg^\vee \ ,
		\\[4pt]
		\Hom_{C^\infty(M)} \big( L_A\otimes_A C^\infty(M), C^\infty(N) \big)
		&\cong \big(TM \oplus \frg[1] \big) \otimes_{C^\infty(M)} C^\infty(N) \ ,
	\end{align*}
	with a certain differential, where $\frg^\vee$ denotes the $C^\infty(M)$-linear dual complex of $\frg$.
	This is a non-positively graded complex of finite free $C^\infty(M)$-modules and hence cofibrant. It follows that both of these complexes are projectively cofibrant complexes of finite-type projective $C^\infty(N)$-modules. If we use these models, then there is no need to derive the tensor product in \eqref{eq:comparison cotangent}, and one sees directly that the map in~\eqref{eq:comparison cotangent} becomes an isomorphism because $I$ is bounded. It follows that $\Spf(A)$ has strictly quasi-coherent deformation theory.
	
	Next suppose that $F$ is a pro-representable formal moduli problem (Example~\ref{ex:prorep FMP}). Since both $j_!$ and $j^*$ preserve filtered colimits, as does reduction, it follows that $(j_!F)_\mathrm{red}\simeq M$ and the unit map $F\longrightarrow j^*\,j_!F$ is an equivalence. Proposition~\ref{prop:def thy} implies that $j_!F$ has quasi-coherent deformation theory; its proof shows that $\kappa_!\,\kappa^*\Der_{j_!F}\longrightarrow \Der_{j_!F}$ is a filtered colimit of equivalences, so that $j_!F$ satisfies strong tangent representability.
	
	Finally, let $F$ be a general formal moduli problem. By~\cite[Proposition 1.5.8]{Lurie:DAG-X}, there exists a simplicial diagram of formal moduli problems $F_\bullet$ such that $|F_\bullet| \simeq F$ and each $F_n$ is pro-representable. Moreover, for each horn inclusion the induced map of tangent complexes
	$$\begin{tikzcd}
		T_x(F_n)\arrow[r] & T_x\big(F(\Lambda^n_i)\big)=\underset{\Delta^k\longrightarrow \Lambda^n_i}\holim\ T_x(F_k)
	\end{tikzcd}$$
	has a fibre whose cohomology groups vanish in positive degrees, because the matching maps are formally smooth and this is detected on tangent complexes (by an argument for formal moduli problems analogous to that of Proposition~\ref{prop:subm theorem}).
	This gives rise to a simplicial diagram $j_!(F_\bullet)$ of qcdt prestacks (with strong tangent representability) such that $|j_!(F_\bullet)| \simeq j_!(F)$ since $j_!$ preserves colimits.
	
	We claim that $j_!(F_\bullet)$ is a formally smooth $\infty$-groupoid, that is, for each horn inclusion the induced map
	$$
	f \colon j_!(F_n) \longrightarrow X \coloneqq \underset{\Delta^k \longrightarrow \Lambda^n_i}\holim\ j_!(F_k)
	$$
	is formally smooth. By Proposition~\ref{prop:subm theorem}, it suffices to verify that, for every point $y \colon N \longrightarrow j_!(F_n)$ from a smooth manifold, the fibre of the map $T_y\big(j_!(F_n)\big)\longrightarrow T_{f(y)}X$ has vanishing cohomology groups in positive degrees. Since the reduction of $j_!(F_n)$ is equivalent to $M$ (by the second step of this proof), the map $y$ again factors over the canonical point $M\longrightarrow j_!(F_n)$. Since both $j_!(F_n)$ and $X$ are qcdt prestacks with strong tangent representability (by Proposition~\ref{prop:def thy} and its proof), their tangent complexes at the point $y$ are obtained by pulling back their tangent complexes at the canonical point. 
	
	It thus suffices to verify that, for the canonical point $x \colon M\longrightarrow j_!(F_n)$, the map
	$$
	T_x \big( j_!(F_n) \big) \simeq T_x \big( j^*\,j_!(F_n) \big)
	 \longrightarrow \underset{\Delta^k \longrightarrow \Lambda^n_i}\holim\ T_x \big( j^*\,j_!(F_k) \big)
	\simeq T_x X 
	$$
	has a fibre with vanishing cohomology in positive degrees. Here the last equivalence uses that $j^*$ preserves homotopy limits. Since $F_\bullet\simeq j^*\,j_!F_\bullet$, this claim follows from the construction of $F_\bullet$.
	We conclude that $j_!F_\bullet$ is a formally smooth $\infty$-groupoid, so that its colimit $j_!F$ is a qcdt prestack (with strong tangent representability) by Proposition~\ref{prop:form smooth oo-groupoid qcdt}. Since reduction preserves colimits, it also follows that $M\simeq (j_!F)_\red$. To see that the unit map $F\longrightarrow j^*\,j_!F$ is an equivalence of formal moduli problems, it suffices to verify that the induced map on tangent complexes is an equivalence. This follows again from Proposition~\ref{prop:form smooth oo-groupoid qcdt}.
\end{proof}

\begin{proof}[Proof of Theorem~\ref{thm:formal stack vs lie}]
	Let us write $$\DT_{M/}\subseteq \mathrm{PreSt}^\mathrm{dt}_{M/}\subseteq \PSh_\infty(\InfMfd)_{M/}$$ for the full $\infty$-subcategories of maps $M\longrightarrow X$ where $X$ is a stack (resp. prestack) with deformation theory. By Theorem~\ref{thm:fmp vs lie} and Lemma~\ref{lem:j*}, there are functors
	\begin{equation}
		\label{eq:sheafified j*}
		\begin{tikzcd}[column sep = 1.25cm]
			T(M/-) \colon \DT_{M/} \arrow[r, hookrightarrow]
			& \mathrm{PreSt}^\mathrm{dt}_{M/} \arrow[r, "j^*"]
			& \FMP_{M/} \arrow[r, "\sim"]
			& L_\infty\Agd_{M}^\infty\ ,
		\end{tikzcd}
	\end{equation}
	where $j^*$ is the restriction along $j$ from Diagram \eqref{diag:j*}. The composition sends each stack $X$ with deformation theory  to $T(M/X)$ equipped with a certain $L_\infty$-algebroid structure on it. Proposition~\ref{prop:j!} implies that each of the  three functors in \eqref{eq:sheafified j*} admits a left adjoint, giving
	\begin{equation}
		\label{eq:sheafified j!}
		\begin{tikzcd}[column sep = 1.4cm]
			L_\infty\Agd_{M}^\infty \arrow[r, "\sim"]
			& \FMP_{M/} \arrow[r, "j_!"]
			& \mathrm{PreSt}^\mathrm{dt}_{M/} \arrow[r, "X\longmapsto \widetilde{X}"]
			& \DT_{M/}\ .
		\end{tikzcd}
	\end{equation}
	Here $j_!$ is the functor from \eqref{diag:j*}, i.e.~taking left Kan extension along $j$, and the last functor takes associated sheaves. By Proposition~\ref{prop:j!}, $j_!(F)$ is a qcdt prestack, so taking associated sheaves yields a map $$x \colon M \longrightarrow \widetilde{j_!(F)}$$ to a qcdt stack by Proposition~\ref{prop:def thy}. 
	
	The proof of Proposition~\ref{prop:qcdt stack 2} shows that $T_x\widetilde{j_!F}$ arises as the global sections of the sheaf associated to the functor which sends an open subspace $u \colon U\longhookrightarrow  M$ to
	\begin{equation}
		T_{u^*x }j_! F \simeq T_x j_! F\otimes_{C^\infty(M)} C^\infty(U)\ .
	\end{equation}
	Here we used that $j_!F$ has strictly quasi-coherent deformation theory (see Definition~\ref{def:qcoh def thy} and Proposition~\ref{prop:j!}). Consequently
	$$
	T_x \widetilde{j_!(F)} \simeq \big( T_x j_! F \big)^\cpl
	\simeq \big(T_xF\big)^\cpl
	$$ 
	is the completion of the tangent complex of $F$.
	Using this, one sees that the unit map $F\longrightarrow j^*\widetilde{j_!F}$ is an equivalence if and only if $T_xF$ is a complete complex of $C^\infty(M)$-modules. Equivalently, this means that the corresponding $L_\infty$-algebroid is complete. 
	
	We conclude that the restriction of the functor~\eqref{eq:sheafified j!} to the complete $L_\infty$-algebroids is fully faithful, with essential image contained in the full $\infty$-subcategory of qcdt stacks $X$ under $M$ inducing an equivalence  $M\simeq X_{\red}$. To show that this is an equivalence, it now remains to verify that the right adjoint~\eqref{eq:sheafified j*} detects equivalences between qcdt stacks under $M$ whose reduction is $M$.
	This can be seen as follows:
	let
	\begin{equation}
		\begin{tikzcd}[row sep = 1cm , column sep = 1cm]
			& M \ar[dl, "x"'] \ar[dr, "p(x)"] &
			\\
			Y \ar[rr, "p"'] & & X
		\end{tikzcd}
	\end{equation}
	be a morphism of qcdt stacks under $M$.
	Suppose also that the maps $x \colon M \longrightarrow Y$ and $p(x) \colon M \longrightarrow X$ induce equivalences on reductions, i.e.~$X_\red \simeq M \simeq Y_\red$, and that the induced morphism $$T(M/X) \longrightarrow T(M/Y)$$ is an equivalence of $L_\infty$-algebroids on $M$.
	We will deduce the claim from Proposition~\ref{prop:inverse_function}, by which it would suffice to show that the induced map $T_y Y \longrightarrow T_{p(y)} X$ is an equivalence in $\Mod_{C^\infty(M)}^\infty$, for every manifold point $y \colon N \longrightarrow Y$.
	
	First we observe that since the canonical morphism $M \longrightarrow Y_\red$ is an equivalence, any point $y \colon N \longrightarrow Y$ factors canonically through the map $x:M \longrightarrow Y$.
	By Remark~\ref{rmk: detecting equivalences for qcdt}, and since $Y$ and $X$ are qcdt stacks, it thus suffices to show that the morphism $T_x Y \longrightarrow T_{p(x)} X$ is an equivalence.
	We now use that the map $T(M/X) \longrightarrow T(M/Y)$ is an equivalence of $L_\infty$-algebroids.
	This induces, in particular, an equivalence of the underlying dg $C^\infty(M)$-modules (which we denote by the same symbols).
	Since $T(M/Y)$ is the homotopy fibre $\hofib(TM \longrightarrow T_x Y)$, and analogously $T(M/X) = \hofib(TM \longrightarrow T_{p(x)}X)$, we recover $T_x Y$ as the homotopy fibre of the anchor map up to a shift,
	\begin{equation}
		T_x Y \simeq \hofib(T_x Y \longrightarrow TM)[1]\ ,
	\end{equation}
	and analogously for $T_{p(x)}X$.
	Thus $T_x Y \longrightarrow T_{p(x)}X$ is an equivalence, and so the claim follows.
	
	Finally, it remains to identify the function algebra on such a qcdt stack $X$ under $M$ with the Chevalley--Eilenberg complex $\ChEil^*(T(M/X))$. To that end, let $F=j^*X$ be the corresponding formal moduli problem. By~\cite[Proposition 1.5.8]{Lurie:DAG-X}, one can write $$F=\underset{\alpha}\colim\ \Spf(A_\alpha)$$ as a sifted colimit of a diagram of corepresentable formal moduli problems. Consequently, we can write $X$ as a sifted colimit
	$$
	X \simeq \widetilde{j_!F} \simeq \underset{\alpha}\colim\ \Spec(A_\alpha)
	$$
	of representable sheaves in $\Sh_\infty(\InfMfd)$. Since the function algebras from Examples~\ref{ex:function algebras} and~\ref{ex:function algebras FMP} both send colimits to limits of cdgas, we obtain the desired equivalence
	$$
	\ChEil^*\big(T(M/X)\big)\simeq \Gamma(F, \calO_F)\simeq \underset{\alpha}\lim\ A_\alpha \simeq \Gamma(j_!F, \calO_{j_!F})\simeq \Gamma(X, \calO_X)
	$$
	where the first equivalence is provided by Theorem~\ref{thm:fmp vs lie}.
\end{proof}


\section{Proof of Theorem~\ref{thm:global}}
\label{sec:global sections proof}


We will now turn to the proof of Theorem~\ref{thm:global}. Let us start by giving a more precise definition of $\cDiff(M)$ as an $\infty$-group in $\Sh_\infty(\InfMfd)$ in terms of its classifying stack.

\begin{definition}
	\label{def:bdiff}
	Let $M$ be a smooth manifold. The stack
	$$
	\rmB \cDiff(M)\colon \InfMfd^{\opp} \longrightarrow \scS
	$$
	is the functor sending each inf-manifold $U$ to the maximal $\infty$-subgroupoid of $\InfMfd_{/U}$ on those objects $p \colon E \longrightarrow U$ such that the homotopy pullback
	$$
	\begin{tikzcd}[row sep = 1cm , column sep = 1cm]
		E \big|U_\red \arrow[r] \arrow[d]
		& E \arrow[d]
		\\
		U_{\red} \arrow[r]
		& U
	\end{tikzcd}
	$$
	defines a fibre bundle $E\big|U_\red\longrightarrow U_\red$ of (ordinary) smooth manifolds with fibres diffeomorphic to $M$. There is a canonical basepoint $\ast\longrightarrow \rmB \cDiff(M)$, classifying the trivial fibre bundle $M\longrightarrow \ast$.
\end{definition}

\begin{lemma}
	\label{lem:bdiff def thy}
	\begin{myenumerate}
		\item The basepoint $\ast\longrightarrow \rmB \cDiff(M)$ is an effective epimorphism in $\Sh_\infty(\InfMfd)$.
		
		\item There is an equivalence $\sfOmega(\rmB\cDiff(M)) \simeq \cDiff(M)$ of group objects in $\Sh_\infty(\InfMfd)$.
	\end{myenumerate}
\end{lemma}

In particular, this implies that $\rmB \cDiff(M)$ is indeed the classifying stack of the diffeomorphism group from Definition~\ref{def:diff}, and that it has deformation theory by Lemma~\ref{lem:diff def thy}.

\begin{proof}
	A morphism $f \colon Y \longrightarrow X$ in $\Sh_\infty(\InfMfd)$ is an effective epimorphism if and only if the morphism $\pi_0 f \colon \pi_0 Y \longrightarrow \pi_0 X$ it induces is an epimorphism~\cite[Proposition~7.2.1.14]{Lurie:HTT}.
	For any $Y \in \Sh_\infty(\InfMfd)$, the sheaf $\pi_0 Y$ arises as the sheafification of the $\infty$-presheaf $U \longmapsto \pi_0(Y(U))$ (which is zero-truncated, i.e.~a sheaf of sets).
	Thus to show that $f$  is an effective epimorphism in $\Sh_\infty(\InfMfd)$, it suffices to show that, for each $U \in \Sh_\infty(\InfMfd)$, the induced map $\pi_0(f_{|U}) \colon \pi_0(Y(U)) \longrightarrow \pi_0(X(U))$ is surjective.
	
	To show this for the map $* \longrightarrow \rmB\cDiff(M)$, suppose that $U$ is an inf-manifold and let $p \colon E \longrightarrow U$ be classified by a map $U \longrightarrow \rmB \cDiff(M)$.
	By the above discussion, it suffices show that $p$ is locally trivial.
	To that end, observe that since $E \times_U U_\red \longrightarrow U_\red$ is a fibre bundle, we can find a cover $\{U_a\}_{a\in\varLambda}$ of $U$ such that $E\times_U (U_a)_{\red} \longrightarrow (U_a)_{\red}$ is the trivial bundle for each $a\in\varLambda$. It then suffices to show that $E\times_U U_a\longrightarrow U_a$ is trivial as well, for each index $a$. 
	
	To see this, consider the  diagram
	$$
	\begin{tikzcd}[row sep = 1cm , column sep = 1cm]
		M\times (U_a)_{\red} \arrow[r] \arrow[d]
		& E \times_U U_a \arrow[r, dashed] \arrow[d]
		& M \times (U_a)_{\red} \arrow[d]
		\\
		(U_{a})_{\red} \arrow[r]
		& U_a\arrow[r, dashed]
		& (U_{a})_{\red}
	\end{tikzcd}
	$$
	in which the left square is a pullback. By Corollary~\ref{cor:inf manifolds retract}, there exists a dashed map $U_a\longrightarrow (U_a)_{\red}$ such that the bottom row composes to the identity. There then exists a dashed map $E\times_U U_a\longrightarrow M\times (U_a)_{\red}$ making the diagram commute, such that the top row composes to the identity: indeed, this follows from the fact that $M\times (U_a)_{\red}\longrightarrow E\times_U U_a$ is of the form $\Spec(A)\longrightarrow \Spec(A')$ where $\rmH^0(A')\longrightarrow \rmH^0(A)$ is surjective with nilpotent kernel, while the right vertical map is formally smooth since it is a submersion (see Example~\ref{ex:submersions are formally smooth}). The right-hand square then induces a map $f\colon E\times_U U_a\longrightarrow M\times U_a$ over $U_a$, whose base change along $(U_a)_{\red}\longrightarrow U_a$ is an equivalence. It follows from Lemma~\ref{lem:nil ext}(2) that $f$ is an equivalence, so that $E\times_U U_a\longrightarrow U_a$ is indeed trivial.
	This shows that $\ast\longrightarrow \rmB \cDiff(M)$ is an effective epimorphism.
	
	Unravelling the definitions, the space of maps $U\longrightarrow \sfOmega(\rmB\cDiff(M))$ is now the union of path components in $\InfMfd(U\times M, M)$ consisting of $\phi\colon U\times M\longrightarrow M$ such that $(\id, \phi)\colon U\times M\longrightarrow U\times M$ is an equivalence of inf-manifolds. This is precisely the space of $U$-points of the stack $\cDiff(M)$.
\end{proof}

Let us now turn to the construction of the right adjoint $$\Hom^\wedge\colon \DT_{M/} \longrightarrow \DT_{\ast \dslash \rmB \cDiff(M)}$$ from Theorem~\ref{thm:global}. The stack $\rmB \cDiff(M)$ comes with a universal bundle $p\colon M/\cDiff(M)\longrightarrow \rmB \cDiff(M)$ which fits into a pullback square (see for instance~\cite{NSS:oo-bundles_I, NSS_oo-Bundles_II, Bunk:Pr_oo-Bundles, Bunk:oo-bundles})
$$
\begin{tikzcd}[row sep = 1cm , column sep = 1cm]
	M \arrow[d, "{q'}"{swap}] \arrow[r, "{p'}"]
	& \ast\arrow[d, "q"]
	\\
	M/\cDiff(M) \arrow[r, "p"{swap}]
	& \rmB \cDiff(M)
\end{tikzcd}
$$
One can also describe $p$ as the associated bundle for the canonical action of $\cDiff(M)$ on $M$. 

Since $\Sh_\infty(\InfMfd)$ is locally cartesian closed, we then obtain a composite adjoint pair
\begin{equation} \label{eq:ptadjpair}
\begin{tikzcd}[column sep = 1.25cm]
	\Sh_\infty(\InfMfd)_{/\rmB \cDiff(M)} \arrow[r, yshift=1ex, "p^*", "\perp"']
	& \Sh_\infty(\InfMfd)_{/(M/\cDiff(M))}\arrow[l, yshift=-1ex, "p_*"] \arrow[r, yshift=1ex, "t_!", "\perp"']
	& \Sh_\infty(\InfMfd)\arrow[l, yshift=-1ex, "t^*"]\ .
\end{tikzcd}
\end{equation}
Here the functor $p^*$ takes the pullback along $p$, and $t_!$ post-composes with the canonical morphism $$t \colon M/\cDiff(M) \longrightarrow \ast \ . $$ The composite right adjoint $p_*\,t^*$ first sends $X$ to $X\times M/\cDiff(M)$ and then takes the internal mapping object relative to $\rmB \cDiff(M)$. 
We can equivalently describe $p_*\,t^*(X)\longrightarrow \rmB \cDiff(M)$ by the object $q^*\,p_*\,t^*(X)$ endowed with a $\cDiff(M)$-action, obtained by pulling back along the \v{C}ech nerve of the map $q$. Since $q^*\,p_*\simeq p'_*\,q'^*$ (which holds in any locally cartesian closed $\infty$-category), one sees that $q^*\,p_*\,t^*(X)\simeq\Hom(M, X)$, with $\cDiff(M)$-action given by pre-composition. In light of this, let us introduce the following notation: we set
\begin{equation}
	\Hom(M, X)\,\big/\,\cDiff(M) \coloneqq p_*\,t^*(X)
	\ \in \ \Sh_\infty(\InfMfd)_{/\rmB\cDiff(M)}\ .
\end{equation}

The functor $t_!\,p^*$ sends the canonical basepoint $q\colon \ast\longrightarrow \rmB \cDiff(M)$ to $M$. Consequently, the adjoint pair \eqref{eq:ptadjpair} induces an adjoint pair on under-categories
$$
\begin{tikzcd}
	t_!\, p^* : \Sh_\infty(\InfMfd)_{* \dslash \rmB \cDiff(M)} \arrow[r, yshift=1ex, "\perp"']
	& \Sh_\infty(\InfMfd)_{M/} : p_*\,t^* \arrow[l, yshift=-1ex]
\end{tikzcd}
$$
where $t_!\, p^*$ sends $$\ast\longrightarrow Y\longrightarrow \rmB \cDiff(M)$$ to the pullback $$M = t_!\, p^*(\ast)\longrightarrow t_!\,p^*(Y)=Y\times_{\rmB\cDiff(M)} M/\cDiff(M) \ . $$
Its right adjoint acts as
$$
p_*\,t^* \big( x \colon M \longrightarrow X \big) = \Hom(M, X)\,\big/\,\cDiff(M)
$$
with basepoint $$x^{\dashv}\colon \ast\longrightarrow \Hom(M, X)\,\big/\,\cDiff(M)$$ induced under the internal hom adjunction by the map $x$.

\begin{lemma}
	\label{lem:dt hom}
	The functors $t_!\,p^*$ and $p_*\,t^*$ both preserve dt stacks.
\end{lemma}

\begin{proof}
	Since $Y$, $\rmB\cDiff(M)$ and $M/\cDiff(M)$ are all dt stacks, $t_!\,p^*(Y)$ is a dt stack as well. If \smash{$X\in \DT_{M/}$}, then $\Hom(M, X)$ is a dt stack by Proposition~\ref{prop:def thy}.
	It follows that $\Hom(M, X)/\cDiff(M)$ is the quotient of a formally smooth $\infty$-groupoid (Definition \ref{def:form smooth oo-groupoid}) and hence a dt stack as well.
\end{proof}

\begin{corollary}
	\label{cor:completed hom adj}
	There is an adjoint pair
	$$
	\begin{tikzcd}
		t_!\, p^* : \DT_{\ast \dslash \rmB\cDiff(M)} \arrow[r, yshift=1ex, "\perp"']
		& \DT_{M/} : p_*\,t^*\ . \arrow[l, yshift=-1ex]
	\end{tikzcd}
	$$
\end{corollary}

It now remains to show that the functor $p_*\,t^*$ is compatible with Lie differentiation. Instead of constructing the diagram of right adjoints from Theorem \ref{thm:global}, our strategy will be to construct the corresponding diagram of \emph{left} adjoints. The advantage of working with left adjoints is that it will suffice to describe them on generators. 

\begin{remark}
	\label{rem:dense}
	Let $\mathscr{V}$ be a presentable $\infty$-category. A small full $\infty$-subcategory $\mathscr{C}\subseteq \mathscr{V}$ is said to be \emph{dense} if the induced colimit-preserving functor $\PSh_\infty(\mathscr{C})\longrightarrow \mathscr{V}$ is a localisation (equivalently its right adjoint is fully faithful). This implies that for any other presentable $\infty$-category $\mathscr{W}$, restriction to $\mathscr{C}$ defines a fully faithful functor $\Fun^{\rmL}(\mathscr{V}, \mathscr{W})\longhookrightarrow  \Fun(\mathscr{C}, \mathscr{W})$.
	
	This is compatible with taking over-categories in the following sense: given an adjoint pair of functors 
$$
	\begin{tikzcd}
		F : \mathscr{V} \arrow[r, yshift=-1ex]
		& \mathscr{W} : G \arrow[l, yshift=1ex, "\perp"',swap]
	\end{tikzcd}
	$$	
and an object $w\in \mathscr{W}$, consider the induced adjunction
$$
	\begin{tikzcd}
		F : \mathscr{V}_{/G(w)} \arrow[r, yshift=-1ex]
		& \mathscr{W}_{/w} : G \ . \arrow[l, yshift=1ex, "\perp"',swap]
	\end{tikzcd}
	$$
If $\mathscr{C}\subseteq \mathscr{V}$ is dense, then $\mathscr{C}_{/G(w)}\subseteq \mathscr{V}_{/G(w)}$ is dense. The functor $F\colon \mathscr{V}_{/G(w)}\longrightarrow \mathscr{W}_{/w}$ is then the unique left adjoint extending
	$$
	\mathscr{C}_{/G(w)} \simeq \mathscr{C}\times_{\mathscr{W}} \mathscr{W}_{/w}
	\xrightarrow{ \ \pi_2 \ } \mathscr{W}_{/w}\ ,
	$$
	where $\pi_2$ denotes the canonical projection functor.
	\qen
\end{remark}

\begin{proof}[Proof of Theorem~\ref{thm:global}]
	Throughout this proof, let us abbreviate $R=C^\infty(M)$.
	We are going to construct a diagram of left adjoints of the form
	\begin{equation}
		\label{eq:left adj}
		\begin{tikzcd}[column sep=2cm , row sep = 1cm]
			\DT_{M/} 
			& \DT_{\ast \dslash \rmB \cDiff(M)} \arrow[l, "t_!\,p^*"{swap}] \arrow[r, "\fgt"]
			& \DT_{\ast/}
			\\
			L_\infty\Agd_M^{\cpl, \infty} \arrow[u, "\Phi_M"]
			& \big(L_\infty\Alg^\infty_\RN \big)_{/\Gamma(M, TM)} \arrow[l, "(R\otimes(-))^{\cpl}"{swap}] \arrow[u, "\Phi_\ast"{swap}] \arrow[r, "\fgt"]
			& L_\infty\Alg^\infty_\RN\arrow[u, "\Phi_\ast"{swap}]
		\end{tikzcd}
	\end{equation}
	such that the desired square from Theorem~\ref{thm:global} is obtained from the left square by passing to right adjoints.
	
	We start by discussing the horizontal functors. The two functors on the right forget the maps to $\rmB\cDiff(M)$ and $TM$, respectively. The functor $t_!\,p^*$ is the left adjoint to the functor $p_*\,t^*$ from Corollary~\ref{cor:completed hom adj}. The bottom left horizontal functor is the left adjoint to the forgetful functor, and sends a map of $\RN$-linear $L_\infty$-algebras $\frg\longrightarrow TM$ to the completion of the $L_\infty$-algebroid $R\otimes_{\RN} \frg$ whose brackets are obtained from those of $\frg$ via the Leibniz rule. 
	
	As for the vertical functors, we define $\Phi_M$ to be the (fully faithful) left adjoint to the functor \smash{$T(M/-) \colon \DT_{M/}\longrightarrow L_\infty\Agd_M^{\cpl, \infty}$} from Theorem \ref{thm:formal stack vs lie}.
	Replacing $M$ by the point induces an analogous adjunction
	\begin{equation}\label{eq:adjunction-point}
		\begin{tikzcd}
			\Phi_* : L_\infty\Alg^\infty_\RN \arrow[r, yshift=1ex, "\perp"']
			& \DT_{*/} : T(*/-)\ . \arrow[l, yshift=-1ex]
		\end{tikzcd}
	\end{equation}
	It now remains to construct a lift of the adjoint pair $(\Phi_*, T(\ast/-))$ to an adjunction on over-categories (which we denote by the same symbols)
	\begin{equation}\label{eq:adjunction-point-lift}
		\begin{tikzcd}
			\Phi_* : (L_\infty\Alg^\infty_\RN)_{/\Gamma(M, TM)} \arrow[r, yshift=1ex, "\perp"']
			& \DT_{* \dslash \rmB\cDiff(M)} : T(*/-)\ , \arrow[l, yshift=-1ex]
		\end{tikzcd}
	\end{equation}
	that makes both squares commute.
	
	To achieve this, we use that the $\infty$-categories of $L_\infty$-algebroids and $L_\infty$-algebras over $TM$ have dense subcategories given by {good} objects \cite[Proposition 5.4]{Nuiten:Koszul_duality_for_Lie_algebroids}. Here an $L_\infty$-algebroid (resp. $L_\infty$-algebra over $TM$) is said to be \emph{good} if it is weakly equivalent to a cofibrant quasi-free $L_\infty$-algebroid (resp. $L_\infty$-algebra) on a finite number of generators, each in strictly positive degree; such an $L_\infty$-algebroid is automatically complete, because its underlying complex is quasi-isomorphic to a complex of finitely generated free $R$-modules (see Remark \ref{rem:complete}). By~\cite[Corollary~4.2.18]{Nuiten:Thesis}, taking the Chevalley--Eilenberg complex defines equivalences
	$$
	\ChEil^*\colon L_\infty\Agd_M^{\cpl, \infty, \mathrm{good}}
	\eq \Art^{\opp}_{/R}
	\qqandqq
	\ChEil^*\colon L_\infty\Alg_\RN^{\cpl, \infty, \mathrm{good}}
	\eq \Art^{\opp}_{/\RN}\ .
	$$
	
	Using this and unravelling the proof of Theorem~\ref{thm:formal stack vs lie}, one sees that the restriction
	$$
	\Phi_M\colon L_\infty\Agd_M^{\cpl, \infty, \mathrm{good}} \subseteq L_\infty\Agd_M^{\cpl, \infty}
	 \longrightarrow \DT_{M/}\ ,
	\qquad \frg \longmapsto \big( M\longrightarrow \Spec(\ChEil^*(\frg)) \big)
	$$
	sends a good $L_\infty$-algebroid $\frg$ to the spectrum of the Chevalley--Eilenberg complex of $\frg$. Similarly, for $M=\ast$ we find that $\Phi_*\colon L_\infty\Alg^{\infty, \mathrm{good}}_\RN \longrightarrow \DT_{\ast/}$ sends $\frg$ to $\Spec(\ChEil^*(\frg))$.
	
	In light of Remark \ref{rem:dense}, it now suffices to provide a dashed functor
	\begin{equation}\label{eq:restricted left}
		\begin{tikzcd}[column sep = 2.5cm, row sep = 1.25cm]
			\DT_{M/}
			& \DT_{\ast\dslash \rmB\cDiff(M)}\arrow[l, "t_!\,p^*"{swap}]\arrow[r, "\fgt"]
			& \DT_{\ast/}
			\\
			L_\infty\Agd_M^{\infty, \mathrm{good}} \arrow[u, "\Spec(\ChEil^*)"]
			& \big(L_\infty\mathrm{Alg}_\RN^{\infty,\mathrm{good}}\big)_{/\Gamma(M,TM)}
			\arrow[u, "\chi"{swap}, dashed]\arrow[l, "{R\otimes(-)}"{swap}]\arrow[r, "\fgt"]
			& 	L_\infty\mathrm{Alg}_\RN^{\infty,\mathrm{good}}\arrow[u, "\Spec(\ChEil^*)"{swap}]
		\end{tikzcd}
	\end{equation}
		such that (a) both squares commute, and (b) the right square is a pullback. Indeed, combining property~(b) and the second half of Remark \ref{rem:dense} implies that the functor $\chi$ determines an adjunction \eqref{eq:adjunction-point-lift} which is nothing but the adjunction \eqref{eq:adjunction-point} lifted to over-categories. Property (a) then implies that the resulting diagram of left adjoints~\eqref{eq:left adj} commutes.
		
		We now define $\chi$ to be the functor which sends each good $L_\infty$-algebra $\frg \longrightarrow TM$ to the morphism $\Spec(\ChEil^*(\frg))\longrightarrow \rmB \cDiff(M)$ classifying the fibre bundle
		$$
		\Spec \big( \ChEil^*(R \otimes \frg) \big)
		 \longrightarrow \Spec \big( \ChEil^*(\frg) \big)\ .
		$$
		Note that this is indeed classified by a map to $\rmB \cDiff(M)$, because the  homotopy pullback along $$\ast=\Spec\big(\ChEil^*(\frg)\big)_{\red}\longrightarrow \Spec\big(\ChEil^*(\frg)\big)$$ is equivalent to $M$: it can be computed as the spectrum of \smash{$\RN\otimes_{\ChEil^*(\frg)}^h \ChEil^*(R\otimes\frg)$}, which is quasi-isomorphic to $R=C^\infty(M)$ by \cite[Corollary~6.5]{Nuiten:Koszul_duality_for_Lie_algebroids}.
		
		With this definition of $\chi$, Diagram \eqref{eq:restricted left} commutes because $t_!\,p^*$ sends a map $Y\longrightarrow \rmB\cDiff(M)$ to the total space of the corresponding bundle, which here coincides with $\Spec(\ChEil^*(R\otimes\frg))$. The right square is a pullback square by~\cite[Theorem~4.4.1 and Proposition~4.4.12]{Nuiten:Thesis}: this asserts that the $\infty$-category of good $L_\infty$-algebras over $TM$ is equivalent to the $\infty$-category of good $L_\infty$-algebras equipped with an action on $R=C^\infty(M)$ by derivations. In turn, this is equivalent to the $\infty$-category of good $L_\infty$-algebras $\frg$ together with a map $\ChEil^*(\frg)\longrightarrow A$ of dg $C^\infty$-rings equipped with an equivalence \smash{$\RN \otimes_{\ChEil^*(\frg)}^h A\simeq R$}. Taking spectra and using Definition \ref{def:bdiff}, the latter $\infty$-category is precisely the pullback of the right square in~\eqref{eq:restricted left}.
	\end{proof}


\section{Proof of Theorem~\ref{st:strictification of At(CG)}}
\label{app:Pf of strictification}


\begin{proof}[Proof of Theorem~\ref{st:strictification of At(CG)}]
	To begin with, recall from the discussion in Section~\ref{sec: Con_l(g) equivalence} (see in particular~\eqref{eq: data for p-spls in Con_l(g)} and~\eqref{eq:conditions for p-spls in Con_l(g)}) that a full connection $\CA$ on the $(n{-}1)$-gerbe presented by the $\rmU(1)$-cocycle $g$ is an $n$-tuple
	\begin{equation}
		\CA = (A^{(1)}, \ldots, A^{(n)})\ ,
		\quad
		A^{(p)} \in \Omega^p(\cC_{n-p}\,\CU) \ ,
	\end{equation}
	satisfying the identities
	\begin{align}
	\begin{split}
		\label{eq:n-gerbe conn identities}
		\dd \log (g) - \check{\delta} A^{(1)} &= 0\ ,
		\\[4pt]
		\dd A^{(p-1)} + (-1)^p\, \check{\delta} A^{(p)} &= 0\ ,
		\quad  p = 2, \ldots, n\ .
		\end{split}
	\end{align}

By Remark~\ref{rmk:constructing cdgc-maps into CE cdgc} an $\infty$-morphism $\varPhi_{g, \CA}$ as desired is the same as a morphism
	\begin{equation}
		\phi_{g,\CA} \colon \ol{\Sym}_{C^\infty(M)} \big( \At_{n-1}(H)[1] \big)
		 \longrightarrow \CAt(g)[1]
	\end{equation}
	of graded $C^\infty(M)$-modules such that 
	\begin{myenumerate}
		\item its underlying $\RN$-linear map is a Maurer--Cartan element in
		\begin{equation}
			\Hom_\RN \Big( \ol{\Sym}_{\RN} \big( \At_{n-1}(H)[1] \big),\, \CAt(g)[1] \Big)\ ,
\end{equation}
and
		
		\item the composition
		\begin{equation}
			\begin{tikzcd}[column sep=1.8cm]
				\ol{\Sym}_{C^\infty(M)} \big( \At_{n-1}(H)[1] \big) \ar[r, "\phi_{g,\CA}"]
				& \CAt(g)[1] \ar[r, "{\rho_{\CAt(g)}[1]}"]
				& TM[1]
			\end{tikzcd}
		\end{equation}
		vanishes on all higher tensor powers of $\At_{n-1}(H)[1]$, and its restriction to $\At_{n-1}(H)[1]$ coincides with the shifted anchor map $\rho_{\At_{n-1}(H)}[1]$.
	\end{myenumerate}
	
	We begin by specifying the morphism $\phi_{g,\CA}$.
	Note that the shifted complex $\At_{n-1}(H)[1]$ exhibits vector fields $X \in TM$ in degree $-1$ and functions $f \in C^\infty(M)$ in degree $-n$.
	We can depict the data of the cochain complexes \smash{$\ol{\Sym}_\RN (\At_{n-1}(H)[1])$} and $\CAt(g)[1]$, as well as the components of the morphism $\phi_{g,\CA}$, as 
	\begin{equation}
		\label{eq:diag for At(g) proof}
		\begin{tikzcd}[column sep=2cm, row sep=1cm]
			0
			& 0
			\\
			TM \ar[r, "{(\id_{TM}, A^{(1)})}"] \ar[u]
			& E_{n-1}(g) \ar[u]
			\\
			\midwedge_\RN^2 TM \ar[r, "A^{(2)}"] \ar[u, "\delta_\ChEil"]
			& C^\infty(\cC_{n-2}\,\CU) \ar[u, "{(0, (-1)^2 \,\check{\delta})}"']
			\\
			\midwedge_\RN^3 TM \ar[r, "A^{(3)}"] \ar[u, "\delta_\ChEil"]
			& C^\infty(\cC_{n-3}\,\CU) \ar[u, "(-1)^3 \,\check{\delta}"']
			\\
			\vdots \ar[u, "\delta_\ChEil"]
			& \vdots \ar[u, "(-1)^4 \,\check{\delta}"']
			\\
			\midwedge_\RN^{n-1} TM \ar[r, "A^{(n-1)}"] \ar[u, "\delta_\ChEil"]
			& C^\infty(\cC_1\,\CU) \ar[u, "(-1)^{n-1}\,\check{\delta}"']
			\\
			\midwedge_\RN^n TM \oplus C^\infty(M) \ar[r, "A^{(n)} - \check{\delta}"] \ar[u, "\delta_\ChEil"]
			& C^\infty(\cC_0\,\CU) \ar[u, "(-1)^n \,\check{\delta}"']
			\\
			\midwedge_\RN^{n+1} TM \oplus \big( TM \otimes_\RN C^\infty(M) \big) \ar[r] \ar[u, "\delta_\ChEil"]
			& 0 \ar[u]
			\\
			\vdots \ar[u]
			& \vdots \ar[u]
		\end{tikzcd}
	\end{equation}
	
Here $\delta_{\rm CE}$ is the Chevalley--Eilenberg differential on $\ol{\Sym}_\RN(\At_{n-1}(H)[1])$ given by
\begin{align}
\delta_{\rm CE}(v_1\cdots v_k) = \sum_{j=1}^{k} \ \sum_{\sigma\in{\rm Sh}(j,k-j)} \, \epsilon(\sigma) \, [v_{\sigma(1)},\dots, v_{\sigma(j)}]_{\At_{n-1}(H),j}\, v_{\sigma(j+1)}\cdots v_{\sigma(k)}
\end{align}
for $v_j=(f_j,X_j)$ with $f_j\in C^\infty(M)$ and $X_j\in TM$, where $\mathrm{Sh}(j,k-j)$ is the set of $(j,k{-}j)$-shuffle permutations. By Corollary~\ref{st:At-L_oo-agd of n-gerbe from first principles}, the only non-zero summands are with $j=2$ and $j=n+1$.
	
	The squares in Diagram~\eqref{eq:diag for At(g) proof} do \textit{not} commute; $\phi_{g,\CA}$ is merely a map of graded $C^\infty(M)$-modules, not a morphism of chain complexes.
	We readily observe that it is compatible with the anchor maps as desired and vanishes on all higher tensor powers of $\At_{n-1}(H)[1]$.
	The top non-zero row of Diagram~\eqref{eq:diag for At(g) proof} is situated in degree $-1$.
	The top horizontal map is well-defined by~\eqref{eq:n-gerbe conn identities}.
	
We are thus left to show that $$\phi_{g, \CA}:\ol{\Sym}_\RN\big(\At_{n-1}(H)[1]\big)\longrightarrow \CAt[1]$$  satisfies the Maurer--Cartan equation
	\begin{equation}
		\label{eq:MC eqn for phi_(g,CA)}
		\dd_\Hom \phi_{g, \CA} + \sum_{l \geqslant 2}\, \frac{1}{l!}\, [\phi_{g, \CA}, \dots, \phi_{g, \CA}]_{\Hom,l} = 0\ ,
	\end{equation}
where  we use the notation
	\begin{equation}
		[\phi_{g, \CA}, \dots, \phi_{g, \CA}]_{\Hom,l}
		= [-,\dots,-]_{\CAt(g),l} \circ (\phi_{g, \CA} \otimes_\RN \cdots \otimes_\RN \phi_{g, \CA}) \circ \sfDelta^{l-1}\ .
	\end{equation}
Since $\CAt(g)$ is even a dg Lie algebroid, this bracket is zero unless $l = 2$, while $\sfDelta$ is the coproduct on $\ol{\Sym}_\RN(\At_{n-1}(H)[1])$ with
\begin{align}
\sfDelta(v_1\cdots v_k) = \sum_{j=1}^{k-1} \ \sum_{\sigma\in{\rm Sh}(j,k-j)} \,\epsilon(\sigma) \, (v_{\sigma(1)}\cdots v_{\sigma(j)})\otimes (v_{\sigma(j+1)}\cdots v_{\sigma(k)})
\end{align}
and $\sfDelta^k(v_1\cdots v_k)=0$.	
We now check \eqref{eq:MC eqn for phi_(g,CA)} in each degree $-k$ for $k\geqslant 2$, grouped into cases as follows.

\ul{$k=2$}. \ 
		Let $X_1 \wedge X_2 \in \midwedge_\RN^2 TM$.
		Then
		\begin{align}
			(\dd_\Hom \phi_{g, \CA}) (X_1 \wedge X_2)
			&= (0, (-1)^2 \,\check{\delta}) \circ A^{(2)}(X_1, X_2)
			- (\id_{TM}, A^{(1)}) \circ \delta_\ChEil (X_1 \wedge X_2)
			\\[4pt]
			&= \big( - [X_1, X_2], (-1)^2 \,\check{\delta} A^{(2)} (X_1, X_2) - A^{(1)}([X_1, X_2]) \big)\ ,
		\end{align}
		and
		\begin{align}
			&[\phi_{g, \CA}, \phi_{g, \CA}]_{\Hom, 2} (X_1 \wedge X_2)
			\\[4pt]
& \hspace{2cm}	= [-]_{\CAt(g),2} \circ (\phi_{g, \CA} \otimes \phi_{g, \CA}) (X_1 \otimes X_2 - X_2 \otimes X_1)
			\\[4pt]
			& \hspace{2cm} = [-]_{\CAt(g),2} \big( (X_1, A^{(1)}(X_1)) \otimes (X_2, A^{(1)}(X_2)) - (X_2, A^{(1)}(X_2)) \otimes (X_1, A^{(1)}(X_1))\big)
			\\[4pt]
			& \hspace{2cm} = 2\, \big( [X_1, X_2], \pounds_{X_1} (A^{(1)}(X_2)) - \pounds_{X_2} (A^{(1)}(X_1)) \big)\ .
		\end{align}
		Therefore in total we obtain
		\begin{align}
			&\Big( \dd_\Hom \phi_{g, \CA} + \sum_{l \geqslant 2} \,\frac{1}{l!}\, [\phi_{g, \CA}, \dots, \phi_{g, \CA}]_{\Hom,l} \Big) (X_1 \wedge X_2)
			\\[4pt]
			&\hspace{4cm} = \big( - [X_1, X_2], (-1)^2 \,\check{\delta}A^{(2)} (X_1, X_2) - A^{(1)}([X_1, X_2]) \big)
			\\
			&\hspace{4cm} \qquad \, + \big( [X_1, X_2], \pounds_{X_1} (A^{(1)}(X_2)) - \pounds_{X_2} (A^{(1)}(X_1)) \big)
			\\[4pt]
			&\hspace{4cm} = \big( 0, ((-1)^2 \,\check{\delta} A^{(2)} + \dd A^{(1)}) (X_1, X_2) \big)
			\\
			&\hspace{4cm} = 0\ .
		\end{align}
		
\ul{$3 \leqslant k \leqslant n-1$}. \ 
		Let $X_1 \wedge \cdots \wedge X_k \in \midwedge_\RN^k TM$.
		Note that
		\begin{equation}
			\delta_{\ChEil} (X_1 \wedge \cdots \wedge X_k)
			= \sum_{1 \leqslant i < j \leqslant k}\, (-1)^{i+j-1}\, [X_i, X_j] \wedge X_1 \wedge \cdots \wedge \widehat{X_i} \wedge \cdots \wedge \widehat{X_j} \wedge \cdots \wedge X_k
		\end{equation}
		 (with the $-1$ in the exponent guaranteeing that $\delta_{\ChEil} (X_1 \wedge X_2) = [X_1, X_2]$). The hat denotes omission of a tensor factor. Then
		\begin{align}
			(\dd_\Hom \phi_{g, \CA}) (X_1 \wedge \cdots \wedge X_k)
			&= (-1)^k \,\check{\delta} \big( A^{(k)}(X_1, \ldots, X_k) \big)
			- A^{(k-1)} \circ \delta_\ChEil (X_1 \wedge \cdots \wedge X_k)
			\\[4pt]
			&= (-1)^k\, (\check{\delta} A^{(k)} ) (X_1, \ldots, X_k)
			\\
			&\qquad - \sum_{1 \leqslant i < j \leqslant k}\, (-1)^{i+j-1}\, A^{(k-1)} \big( [X_i, X_j], X_1, \ldots, \widehat{X_i}, \ldots, \widehat{X_j}, \ldots, X_k \big)\ .
		\end{align}
		
		Moreover, we compute
		\begin{align}
			&[\phi_{g, \CA}, \phi_{g, \CA}]_{\Hom, 2} (X_1 \wedge \cdots \wedge X_k)
			\\[4pt]
			&\quad = [-]_{\CAt(g),2} \circ (\phi_{g, \CA} \otimes \phi_{g, \CA})
			\bigg( \sum_{l = 1}^k \ \sum_{\sigma \in \mathrm{Sh}(l, k-l)}\, \epsilon(\sigma)\, (X_{\sigma(1)} \wedge \cdots \wedge X_{\sigma(l)}) \otimes (X_{\sigma(l+1)} \wedge \cdots \wedge X_{\sigma(k)}) \bigg)
			\\[4pt]
			& \quad = [-]_{\CAt(g),2} \circ (\phi_{g, \CA} \otimes \phi_{g, \CA})
			\bigg( \sum_{\sigma \in \mathrm{Sh}(1, k-1)}\, \epsilon(\sigma)\, X_{\sigma(1)} \otimes (X_{\sigma(2)} \wedge \cdots \wedge X_{\sigma(k)})
			\\
			& \hspace{6cm} + \sum_{\sigma \in \mathrm{Sh}(k-1, 1)}\, \epsilon(\sigma)\, (X_{\sigma(1)} \wedge \cdots \wedge X_{\sigma(k-1)}) \otimes X_{\sigma(k)} \bigg)
			\\[4pt]
			& \quad = [-]_{\CAt(g),2} \bigg( \sum_{\sigma \in \mathrm{Sh}(1, k-1)}\, \epsilon(\sigma)\, \big( X_{\sigma(1)}, A^{(1)}(X_{\sigma(1)}) \big) \otimes A^{(k-1)}(X_{\sigma(2)} \wedge \cdots \wedge X_{\sigma(k)})
			\\
			&\hspace{4cm} + \sum_{\sigma \in \mathrm{Sh}(1, k-1)}\, \epsilon(\sigma)\, A^{(k-1)}(X_{\sigma(1)} \wedge \cdots \wedge X_{\sigma(k-1)}) \otimes \big( X_{\sigma(k)}, A^{(1)}(X_{\sigma(k)}) \big) \bigg)
			\\[4pt]
			& \quad = 2\, \sum_{\sigma \in \mathrm{Sh}(1, k-1)}\, \epsilon(\sigma)\, \pounds_{X_{\sigma(1)}} \big( A^{(k-1)}(X_{\sigma(2)} \wedge \cdots \wedge X_{\sigma(k)}) \big)
			\\[4pt]
			& \quad = 2\, \sum_{i=1}^k\, (-1)^{i-1}\, \pounds_{X_i} \big( A^{(k-1)}(X_1 \wedge \cdots \wedge \widehat{X_i} \wedge \cdots \wedge X_k) \big)\ .
		\end{align}
		Here we have used that $[-]_{\CAt(g),2}$ gives zero unless at least one of its arguments is in degree one (this allows us to discard all other splittings of the exterior product $X_1 \wedge \cdots \wedge X_k$).
		Subsequently, we have used the graded antisymmetry of the binary bracket to combine the two sums into one with a factor of $2$.
		
		Combining these results, we obtain 
		\begin{align}
			&\Big( \dd_\Hom \phi_{g, \CA} + \sum_{l \geqslant 2}\, \frac{1}{l!}\, [\phi_{g, \CA}, \dots, \phi_{g, \CA}]_{\Hom,l} \Big) (X_1 \wedge \cdots \wedge X_k)
			\\[4pt]
			& \hspace{4cm} = (-1)^k\, (\check{\delta} A^{(k)}) (X_1, \ldots, X_k)
			\\*
			&\hspace{4cm} \qquad - \sum_{1 \leqslant i < j \leqslant k}\, (-1)^{i+j-1}\, A^{(k-1)} \big( [X_i, X_j], X_1, \ldots, \widehat{X_i}, \ldots, \widehat{X_j}, \ldots, X_k \big)
			\\*
			&\hspace{4cm} \qquad + \sum_{i=1}^k\, (-1)^{i-1}\, \pounds_{X_i} \big( A^{(k-1)}(X_1 \wedge \cdots \wedge \widehat{X_i} \wedge \cdots \wedge X_k) \big)
			\\[4pt]
			& \hspace{4cm} = \big( (-1)^k \,\check{\delta} A^{(k)} + \dd A^{(k-1)} \big) (X_1, \ldots, X_k)
			\\
			& \hspace{4cm} = 0\ ,
		\end{align}
		where we have used the global (or coordinate-free) formula for the de Rham differential (see for instance~\cite[Proposition~14.32]{Lee:Smooth_Mfds_v2}).
		
\ul{$k=n$}. \ 
		We have to evaluate the $\RN$-linear map
		\begin{equation}
			\dd_\Hom \phi_{g, \CA} + \sum_{l \geqslant 2}\, \frac{1}{l!}\, [\phi_{g, \CA}, \dots, \phi_{g, \CA}]_{\Hom,l}
			\
			\in	\ \Hom_\RN \Big( \ol{\Sym}_{\RN} \big( \At_{n-1}(H)[1]), \CAt(g) \Big)
		\end{equation}
		on elements of the form $X_1 \wedge \cdots \wedge X_n + f$, where $X_1, \ldots, X_n \in TM$ and $f \in C^\infty(M)$.
		By linearity, we can split this into two computations.
		For arguments of the form $X_1 \wedge \cdots \wedge X_n \in \midwedge_\RN^n TM$, the computation is exactly as in the case $3\leqslant k\leqslant n-1$ above; that is, we obtain 
		\begin{align}
			\Big( \dd_\Hom \phi_{g, \CA} + \sum_{l \geqslant 2}\, \frac{1}{l!}\, [\phi_{g, \CA}, \dots, \phi_{g, \CA}]_{\Hom,l} \Big) (X_1 \wedge \cdots \wedge X_n)
			&= \big( (-1)^n \,\check{\delta} A^{(n)} + \dd A^{(n-1)} \big) (X_1, \ldots, X_n) \\[4pt]
			&= 0\ .
		\end{align}
		It remains to consider arguments of the form $f \in C^\infty(M)$.
		The Chevalley--Eilenberg differential on $\ChEil_*(\At_{n-1}(H))$ satisfies $\delta_\ChEil(f) = 0$ (see Corollary~\ref{st:At-L_oo-agd of n-gerbe from first principles}).
		Since $\check{\delta}^2 = 0$ and also $\sfDelta(f) = 0$, it follows that
		\begin{equation}
			\Big( \dd_\Hom \phi_{g, \CA} + \sum_{l \geqslant 2}\, \frac{1}{l!}\, [\phi_{g, \CA}, \dots, \phi_{g, \CA}]_{\Hom,l} \Big) (f)
			= 0\ .
		\end{equation}
		
\ul{$k=n+1$}. \ 
		We again use linearity to split the computation into two cases.
		First, we evaluate
		\begin{align}
			\Big( \dd_\Hom \phi_{g, \CA} + \sum_{l \geqslant 2}\, \frac{1}{l!}\, [\phi_{g, \CA}, \dots, \phi_{g, \CA}]_{\Hom,l} \Big) (X_1 \wedge \cdots \wedge X_{n+1})\ .
		\end{align}
		We obtain 
		\begin{align}
			&(\dd_\Hom \phi_{g, \CA}) (X_1 \wedge \cdots \wedge X_{n+1})
			\\[4pt]
			& \qquad = - \sum_{1 \leqslant i < j \leqslant n+1}\, (-1)^{i+j-1}\, A^{(n)} \big( [X_i, X_j], X_1, \ldots, \widehat{X_i}, \ldots, \widehat{X_j}, \ldots, X_{n+1} \big)
			+ (\check{\delta} H) (X_1, \ldots, X_{n+1})\ ,
		\end{align}
		and, as in the case $3\leqslant k\leqslant n-1$ above,
		\begin{align}
			[\phi_{g, \CA}, \phi_{g, \CA}]_{\Hom, 2} (X_1 \wedge \cdots \wedge X_{n+1})
			= 2\, \sum_{i=1}^{n+1}\, (-1)^{i-1}\, \pounds_{X_i} \big( A^{(n)} (X_1 \wedge \cdots \wedge \widehat{X_i} \wedge \cdots \wedge X_{n+1}) \big)\ .
		\end{align}
		Together this yields
		\begin{align}
			&\Big( \dd_\Hom \phi_{g, \CA} + \sum_{l \geqslant 2}\, \frac{1}{l!}\, [\phi_{g, \CA}, \dots, \phi_{g, \CA}]_{\Hom,l} \Big) (X_1 \wedge \cdots \wedge X_{n+1})
			= (\check{\delta}H - \dd A^{(n)}) (X_1, \ldots, X_n)
			= 0\ .
		\end{align}
		This uses precisely that $H$ is the curvature of the $(n{-}1)$-gerbe connection $\CA$.
		
		Finally, we compute
		\begin{align}
			\Big( \dd_\Hom \phi_{g, \CA} + \sum_{l \geqslant 2}\, \frac{1}{l!}\, [\phi_{g, \CA}, \dots, \phi_{g, \CA}]_{\Hom,l} \Big) (X_1 \vee f)\ .
		\end{align}
		We find 
		\begin{align}
			(\dd_\Hom \phi_{g, \CA}) (X_1 \vee f)
			= \check{\delta} (\pounds_{X_1} f)
			= \dd \,\check{\delta}(f) (X_1)\ ,
		\end{align}
		and
		\begin{align}
			[\phi_{g, \CA}, \phi_{g, \CA}]_{\Hom, 2} (X_1 \vee f)
			&= [-]_{\CAt(g),2} \circ (\phi_{g, \CA} \otimes \phi_{g, \CA}) \circ \sfDelta (X_1 \vee f)
			\\[4pt]
			&= [-]_{\CAt(g),2} \circ (\phi_{g, \CA} \otimes \phi_{g, \CA}) (X_1 \otimes f + (-1)^n\, f \otimes X_1)
			\\[4pt]
			&= [-]_{\CAt(g),2} \big( (X_1, A^{(1)}(X_1)) \otimes (- \check{\delta} f) + (-1)^n\, (- \check{\delta} f) \otimes (X_1, A^{(1)}(X_1)) \big)
			\\[4pt]
			&= - 2\, \pounds_{X_1} \,\check{\delta} f
			\\
			&= -2\,\dd\,\check\delta(f)(X_1) \ .
		\end{align}
		These two contributions cancel each other, and we obtain 
		\begin{align}
			\Big( \dd_\Hom \phi_{g, \CA} + \sum_{l \geqslant 2}\, \frac{1}{l!}\, [\phi_{g, \CA}, \dots, \phi_{g, \CA}]_{\Hom,l} \Big) (X_1 \vee f)
			= 0\ .
		\end{align}

All higher degree terms vanish since $\CAt(g)[1]$ is concentrated in degrees ranging from $-n$ to $-1$.
	This completes the proof that $\phi_{g, \CA}$ satisfies the Maurer--Cartan equation~\eqref{eq:MC eqn for phi_(g,CA)}.
	
	Finally, we have to check that the $\infty$-morphism $\phi_{g, \CA}$ is indeed a weak equivalence of $L_\infty$-algebroids over $C^\infty(M)$.
	By definition (see for instance~\cite[Definition~3.10]{KS:Intro_to_L_oo-algebras}) this is the case if and only if the linear part of the morphism
	\begin{equation}
		\label{eq:linear part of phi_(g,CA)}
		\begin{tikzcd}[column sep=2cm, row sep=1cm]
			0
			& 0
			\\
			TM \ar[r, "{(\id_{TM}, A^{(1)})}"] \ar[u]
			& E_{n-1}(g) \ar[u]
			\\
			0 \ar[u] \ar[r]
			& C^\infty(\cC_{n-2}\,\CU) \ar[u, "{(0, (-1)^2 \,\check{\delta})}"']
			\\
			\vdots \ar[u]
			& \vdots \ar[u, "(-1)^3 \,\check{\delta}"']
			\\
			0 \ar[r] \ar[u]
			& C^\infty(\cC_1\,\CU) \ar[u, "(-1)^{n-1} \,\check{\delta}"']
			\\
			C^\infty(M) \ar[r, "\check{\delta}"] \ar[u]
			& C^\infty(\cC_0\,\CU) \ar[u, "(-1)^n \,\check{\delta}"']
			\\
			0 \ar[u]
			& 0 \ar[u]
		\end{tikzcd}
	\end{equation}
	is a quasi-isomorphism of cochain complexes.
	
	First, we claim that the morphism $(\id_{TM}, A^{(1)})$ induces an isomorphism on $\rmH^{-1}$ (recall that the complexes in \eqref{eq:linear part of phi_(g,CA)} are shifted by one degree).
	Let $(X, f)$ and $(X', f')$ be elements of $E_{n-1}(g)$.
	These are cohomologous if and only if $X = X'$ and there exists some $h \in C^\infty(\cC_{n-2} \,\CU)$ such that $\check{\delta} h = f' - f$.
	
	Consider the maps
	\begin{alignat}{3}
		[\id_{TM}, A^{(1)}] & \colon TM  \longrightarrow \rmH^{-1} \big( E_{n-1}(g) \big)\ ,
		& \quad
		X &\longmapsto \big[ X, A^{(1)}(X) \big]\ ,
		\\[4pt]
		\pr_{TM} & \colon \rmH^{-1} \big( E_{n-1}(g) \big) \longrightarrow TM\ ,
		&\quad
		[Y,f] &\longmapsto Y\ .
	\end{alignat}
	These are isomorphisms:
	we directly see that $\pr_{TM} \circ [\id_{TM}, A^{(1)}]$ is the identity on $TM$, and 
	\begin{equation}
		[\id_{TM}, A^{(1)}] \circ \pr_{TM} [Y,f]
		= \big[ Y, A^{(1)}(Y) \big]
		= [Y, f] + \big[ 0, A^{(1)}(Y) - f \big]\ .
	\end{equation}
	Recall that the \v{C}ech cohomology of the sheaf of smooth real-valued functions $C_M^\infty$ satisfies 
	\begin{equation}
		\label{eq:Cech coho of real fcts}
		\check{\rmH}^{p-1} \big( M; C_M^\infty \big) = 0
	\end{equation}
for all $ p \geqslant 1$, and that it may be computed using any good open cover $\CU$ of the manifold $M$.
	Since moreover
	\begin{equation}
		\check{\delta} \big( A^{(1)}(Y) - f \big)
		= (\check{\delta} A^{(1)})(Y) - \check{\delta} f
		= \dd \log(g) (Y) - \dd \log(g) (Y)
		= 0\ ,
	\end{equation}
	it follows that $[\id_{TM}, A^{(1)}] \circ \pr_{TM}$ is the identity on cohomology.
	
	In degrees $-2, \ldots, -(n-1)$ the morphism of cochain complexes~\eqref{eq:linear part of phi_(g,CA)} is an isomorphism on cohomology groups by~\eqref{eq:Cech coho of real fcts}.
	Finally, it also induces an isomorphism of cohomology groups in degree $-n$; this follows from the fact that the \v{C}ech complex $C_*(\CU, C^\infty_M)$ is an acyclic resolution of the sheaf $C^\infty_M$ on $M$.
\end{proof}


\section{A simplicial resolution for dg abelian extensions of $L_\infty$-algebroids}
\label{app:simplicial resolution}

Let $A$ be a cdga over $\RN$ with differential $\dd_A$, and $\frg$ a  dg Lie algebroid over $A$ with differential $\dd_\frg$, bracket $[-,-]_\frg$ and anchor map $\rho \colon \frg \longrightarrow T_A$.
In this appendix we construct an explicit fibrant simplicial resolution of $\frg$ in \smash{$L_\infty\Agd^\dg_A$}, under certain abelianness assumptions on $\frg$ (see Theorem~\ref{st:fib res of fibrant dg abelian extensions} below).
Essentially, the Atiyah $L_\infty$-algebroid of any higher bundle whose structure group is an $\bbE_\infty$-group satisfies this assumption.
The fibrant simplicial resolution we construct may be of independent interest, but in this paper it will serve as a key ingredient in the proof of Theorem~\ref{st:finite and infinitesimal l-cons on n-gerbes} which we give in Appendix~\ref{app:Pf of equivalence of connection spaces}.

Recall the notion of a dg Lie algebroid representation on a dg $A$-module (see for instance~\cite[p.~29]{Nuiten:HoAlg_for_Lie_Algds}).

\begin{lemma}
	For each $l \in \ZN$, setting
	\begin{equation}
		[-,-] \colon \frg \otimes_\RN \ker(\rho)[l]  \longrightarrow \ker(\rho)[l]\ ,
		\quad
		(x, v) \longmapsto [x,v]_\frg
	\end{equation}
	defines a dg Lie algebroid representation of $\frg$ on $\ker(\rho)[l]$, satisfying 
	\begin{equation}
		[a\,x,v] = a\, [x,v]
		\qqandqq
		[x,a\,v] = (-1)^{|a|\,|x|}\, a\,[x, v] + \pounds_{\rho(x)}(a) \, v
		\ ,
	\end{equation}
	for all $a\in A$, $x\in\frg$ and $v\in\ker(\rho)[l]$.
\end{lemma}

\begin{proof}
	Since $\ker(\rho) \subset \frg$ is a Lie ideal, the action map $[-,-] \colon \frg \otimes_\RN  \ker(\rho)[l] \longrightarrow \ker(\rho)[l]$ is well-defined.
	Compatibility of the action with the differential and the $A$-module structure follow from the fact that $\frg$ is a dg Lie algebroid over $A$.
	The degree shift on $\ker(\rho)$ does not enter in these conditions, since the differential and the $A$-action only have to be commuted with elements $x \in \frg$, i.e.~the entries of the first argument of the bracket $[-,-]$.
\end{proof}

Let $\ch_\RN \colon L_\infty \Agd^\dg_A \longrightarrow \Mod^\dg_\RN$ denote the forgetful functor.\footnote{The notation reminds us that we keep only the underlying cochain complex.}
We define a simplicial object in $\Mod_A^\dg$ as the pullback
\begin{equation}
	\label{eq:ch_k(hat(frg)) as pullback}
	\begin{tikzcd}[column sep=1.75cm, row sep=1cm]
		\Hom_{T_A} \big( C_* (\Delta; \RN), \ch_\RN(\frg) \big) \ar[r] \ar[d,"\widehat\rho",swap]
		& \Hom_\RN \big( C_* (\Delta; \RN), \ch_\RN(\frg) \big) \ar[d, "\rho_*"]
		\\
		\sfc\, \ch_\RN(T_A) \cong \Hom_\RN \big( \sfc\, C_* (\Delta^0; \RN), \ch_\RN(T_A) \big) \ar[r]
		& \Hom_\RN \big( C_* (\Delta; \RN), \ch_\RN(T_A) \big)
	\end{tikzcd}
\end{equation}
of simplicial objects in $\Mod^\dg_A$, where $\sfc$ denotes the constant diagram functor.

Explicitly, we find the following:
first, at the level of underlying graded modules 
\begin{equation}
	\grmod_\RN \big( C_* (\Delta^n; \RN) \big)
	= \bigoplus_{l \in \NN_0} \ \bigoplus_{[l] \longhookrightarrow  [n]}\, \RN[l]\ ,
\end{equation}
that is,
\begin{equation}
	C_{-l} (\Delta^n; \RN) = \bigoplus_{[l] \longhookrightarrow  [n]}\, \RN\ ,
\end{equation}
where the sum runs over all injective morphisms $[l] \longhookrightarrow  [n]$ in $\bbDelta$.
We denote a pure element in the summand indexed by $\psi \colon [l] \longhookrightarrow  [n]$ as a pair $(\lambda, \psi)$, where $\lambda \in \RN$.
Note that by our convention to treat all differentials as having degree $+1$, the complexes of simplicial chains are concentrated in non-positive degrees.

At the level of underlying graded $\RN$-modules, we thus obtain 
\begin{align}
	\grmod_\RN\, \Hom_\RN \big( C_* (\Delta^n; \RN), \ch_\RN(\frg) \big)
	&= \bigoplus_{l \in \NN_0} \ \bigoplus_{[l] \longhookrightarrow  [n]}\, \grmod_\RN\, \ch_\RN(\frg)[-l]
	\\[4pt]
	&= \bigoplus_{l \in \NN_0} \ \bigoplus_{[l] \longhookrightarrow  [n]}\, \grmod_\RN\, \ch_\RN \big( \frg[-l] \big)\ .
\end{align}
We denote the pure elements of these graded modules as pairs $(x, \varphi)$, where $\varphi \colon [l] \longhookrightarrow  [n]$ and $x \in \frg$.
The degree of a pure element $(x, \varphi)$ is
\begin{equation}
	\label{eq:degree in Hom_k(CD, frg)}
	\big| (x, \varphi) \big| = |x| + l\ .
\end{equation}
We also refer to $l \in \NN_0$ as the \textit{weight} of $(x, \varphi)$.
The element $(x, \varphi)$ represents the morphism of $\RN$-modules
\begin{equation}
	C_{-l}(\Delta^n; \RN) \longrightarrow \ch_\RN (\frg)\ ,
	\qquad
	(\lambda, \psi) \longmapsto
	\begin{cases}
		0\ , & \psi \neq \varphi\ ,
		\\[4pt]
		(\lambda\, x, \varphi)\ , & \psi = \varphi\ .
	\end{cases}
\end{equation}
The pairs $(x, \varphi)$ generate the graded $\RN$-module $\grmod_\RN\, \Hom_\RN (C_* (\Delta^n; \RN), \ch_\RN(\frg))$.

The differential on $\Hom_\RN (C_* (\Delta^n; \RN), \ch_\RN(\frg))$ is the usual differential on the  hom  complex, but we will need to describe it explicitly.
It will be useful to denote a morphism $\varphi \colon [l] \longrightarrow [n]$ in $\bbDelta$ by its ordered sequence of values; that is, we write
\begin{equation}
	\varphi =\big(  0 \leqslant \varphi(0) < \cdots < \varphi(l) \leqslant n \big)\ .
\end{equation}
The differential acts on a pair $(x, \varphi) = (x, 0 \leqslant \varphi(0) < \cdots < \varphi(l) \leqslant n)$ as
\begin{align}
	\label{eq:dd on Hom(CDelta, frg)}
	\begin{split}
	\dd (x, \varphi)
	&=  (\dd_\frg x, \varphi) - (-1)^{|x| + l}\, \sum_{r = 0}^{l+1}\, (-1)^r \  \sum_{\substack{\psi \colon [l+1] \longhookrightarrow  [n]\\ \psi \circ \partial_r = \varphi}}\, (x, \psi)
	\\[4pt]
	&= (\dd_\frg x, \varphi) \\
	& \qquad \, - \sum_{r = 0}^{l+1}\, (-1)^{r+|x|+l} \ \sum_{j  = \varphi(r-1) +1}^{\varphi(r)  - 1} \big( x, \varphi(0) < \cdots < \varphi(r-1) < j < \varphi(r) < \cdots < \varphi(l+1) \big)\ ,
	\end{split}
\end{align}
where we adopt the convention $\varphi(-1) = 0$ and $\varphi(l+1) = n$.

The simplicial structure takes the following form:
a morphism $\alpha\colon [m] \longrightarrow [n]$ acts as
\begin{equation}
	\label{eq:simplical mps on ch(frg)}
	\Hom_\RN \big( \alpha, \ch_\RN(\frg) \big) (x, \varphi)
	= \sum_{\substack{\psi \colon [l] \longhookrightarrow  [m] \\ \alpha \circ \psi = \varphi}}\, (x, \psi)\ .
\end{equation}
Indeed, one can check that for each $(\lambda, \psi) \in C_{-l}(\Delta^m;\RN)$,
\begin{equation}
	\label{eq:spl structure of Hom(CDelta, frg)}
	\Big( \Hom_\RN \big( \alpha, \ch_\RN(\frg) \big) (x, \varphi) \Big) (\lambda, \psi)
	= (x, \varphi) (\lambda, \alpha \circ \psi)
	=
	\begin{cases}
		0\ , & \alpha \circ \psi \neq \varphi\ ,
		\\[4pt]
		(\lambda\, x,\varphi)\ , & \alpha \circ \psi = \varphi\ .
	\end{cases}
\end{equation}

The simplicial graded vector space underlying $\Hom_\RN (C_* (\Delta; \RN), \ch_\RN(\frg))$ further carries an $A$-module structure.
It is given as
\begin{align}
	A \otimes_\RN \grmod_\RN\, \Hom_\RN \big( C_* (\Delta; \RN), \ch_\RN(\frg) \big)
	& \longrightarrow \grmod_\RN\, \Hom_\RN \big( C_* (\Delta; \RN), \ch_\RN(\frg) \big)\ ,
	\\
	a \otimes_\RN (x, \varphi)
	&\longmapsto a\, (x, \varphi) = (a\,x, \varphi)\ .
\end{align}
One checks the compatibility with the differential,
\begin{align}
	\dd \big( a \, (x, \varphi) \big)
	&= \dd (a\,x, \varphi)
	\\[4pt]
	&= \big( \dd_\frg (a\,x), \varphi \big) - \sum_{r = 0}^{l+1}\, (-1)^{r+|a|+|x|+l} \ \sum_{\substack{\psi \colon [l+1] \longhookrightarrow  [n]\\ \psi \circ \partial_r = \varphi}}\, (a\,x, \psi)
	\\[4pt]
	&= \big( \dd_A a \, x + (-1)^{|a|}\, a \, \dd_\frg x, \varphi \big) - \sum_{r = 0}^{l+1}\, (-1)^{r+|a|+|x|+l} \ \sum_{\substack{\psi \colon [l+1] \longhookrightarrow  [n]\\ \psi \circ \partial_r = \varphi}}\, (a\,x, \psi)
	\\[4pt]
	&= \dd_A a \, (x, \varphi) + (-1)^{|a|}\, a \, \dd (x, \varphi)\ .
\end{align}
That is, $\Hom_\RN (C_* (\Delta; \RN), \ch_\RN(\frg))$ is even a simplicial dg $A$-module.

We can now describe the simplicial dg $A$-module in the top-left corner in~\eqref{eq:ch_k(hat(frg)) as pullback} as follows:
its underlying graded $\RN$-module in simplicial degree $n \in \NN_0$  is the subcomplex
\begin{equation}
	\Hom_{T_A} \big( C_* (\Delta^n; \RN), \ch_\RN(\frg) \big)
	\subset \Hom_\RN \big( C_* (\Delta^n; \RN), \ch_\RN(\frg) \big)
\end{equation}
with weight decomposition
\begin{align}
\begin{split}
	\label{eq:frg: ul k-module}
	\grmod_\RN\, \Hom_{T_A} \big( C_* (\Delta^n; \RN), \ch_\RN(\frg) \big)
	& = \Big\{ (x_0, \ldots, x_n) \in {\bigoplus\limits_{i = 0}^n}\, \frg\, \Big| \, \rho(x_0) = \cdots = \rho(x_n) \in T_A \Big\}
	\\
	& \qquad \, \oplus \ \bigoplus_{l \in \NN} \ \bigoplus_{[l] \longhookrightarrow  [n]}\, \grmod_\RN\, \ch_\RN \big( \ker(\rho) \big)[-l]
	\\[4pt]
	&\eqqcolon \widehat{\frg}_{n,(0)} \ \oplus \ \bigoplus_{l \in \NN}\, \widehat{\frg}_{n, (l)} =: \widehat{g}_n \ .
	\end{split}
\end{align}
Here we have introduced notation for the \textit{weight~$l$ graded submodule} $\widehat{\frg}_{n,(l)}$, for $l \in \NN_0$.

We also write
\begin{equation}
	\widehat{\frg}_{n, (+)} \coloneqq \bigoplus_{l \in \NN}\, \widehat{\frg}_{n, (l)}
	\qqandqq
	\widehat{\frg}_{n, (l)} = \bigoplus_{\varphi \colon [l] \longhookrightarrow  [n]}\, \widehat{\frg}_{n, (\varphi)}\ .
\end{equation}
We emphasise that, in general, each of the graded $\RN$-modules $\widehat{\frg}_{n, (l)}$ may still be non-zero in each degree.
The anchor map of $\widehat{\frg}_{n,(l)}$, i.e.~the left vertical arrow in~\eqref{eq:ch_k(hat(frg)) as pullback}, vanishes on the positive weight component and acts as
\begin{equation}
	\widehat{\rho}(x_0, \ldots, x_n) = \rho(x_0) = \cdots = \rho(x_n)
\end{equation}
on the weight zero part.

We check that the differential restricts to this graded $\RN$-submodule:
since the anchor map $\rho$ of $\frg$ is a morphism of dg $A$-modules, the differential automatically restricts to the positive weight component.
It remains to consider its action on elements
\begin{equation}
	\xi = \sum_{i = 0}^n\, \big( x_i, i \colon [0] \longhookrightarrow  [n] \big)
	\ \in \ \widehat{\frg}_{n, (0)}\ .
\end{equation}
Note that since the anchor map $\rho$ preserves the $\frg$-degree, the condition $\rho(x_i) = \rho(x_j)$ for all $i, j = 0, \ldots, n$ implies that $|x_i| = |x_j|$ for all $i, j = 0, \ldots, n$.
In particular, $\xi$ is a pure element, i.e.~it has a well-defined degree $|\xi| = |x_i|$, for any $i = 0, \ldots, n$.

From~\eqref{eq:dd on Hom(CDelta, frg)} we obtain 
\begin{align}
\begin{split}
	\label{eq:dd_(hat(frg)) on weight zero}
	\dd_{\widehat{\frg}_n} \xi
	&= \sum_{i = 0}^n\, \dd \big( x_i, i \colon [0] \longhookrightarrow  [n] \big)
	\\[4pt]
	&= \sum_{i = 0}^n\, \bigg( \big( \dd_\frg x_i, i \colon [0] \longhookrightarrow  [n] \big) \\
	& \hspace{2cm}
	- \sum_{j = 0}^{i-1}\, (-1)^{|x_i|}\, \big( x_i, 0 \leqslant j < i \leqslant n \big)
	+ \sum_{j = i+1}^{n}\, (-1)^{|x_i|}\, \big( x_i, 0 \leqslant i < j \leqslant n \big) \bigg)
	\\[4pt]
	&= \underbrace{\sum_{i = 0}^n\, \big( \dd_\frg x_i, i \colon [0] \longhookrightarrow  [n] \big)}_{\in\, \widehat{\frg}_{n, (0)}}
	\ - \ \underbrace{\bigoplus_{\varphi \colon [1] \longhookrightarrow  [n]}\, \big( \underbrace{ (-1)^{|x_{\varphi(1)}|}\, x_{\varphi(1)} - (-1)^{|x_{\varphi(0)}|}\, x_{\varphi(0)}}_{\in\, \ker(\rho)}, \varphi \big)}_{\in\, \widehat{\frg}_{n, (+)}}\ .
	\end{split}
\end{align}
Since $\widehat{\rho}$ vanishes on the positive weight part, it also follows that $\widehat{\rho} \colon \widehat{\frg}_n \longrightarrow T_A$ is a morphism of dg $\RN$-modules.

\begin{remark}\label{rem:dgraise}
	We observe that the differential restricts to maps
	\begin{equation}
		\dd_{\widehat{\frg}_n} \colon \widehat{\frg}_{n, (l)} \longrightarrow \widehat{\frg}_{n, (l)} \oplus \widehat{\frg}_{n, (l+1)}\ .
	\end{equation}
	It has a weight preserving component and a component which increases the weight by one.
	\qen
\end{remark}

\begin{definition}\label{def:ab ext}
	A fibrant $L_\infty$-algebroid $\frg$ over $A$ is a \textit{dg abelian extension of $T_A$} if $\ker(\rho)$ is an abelian $L_\infty$-algebroid, i.e.~it has only trivial brackets.
\end{definition}

\begin{lemma}
	\label{st:dgLAgd action on ker(rho) for ab ext}
	Let $\frg$ be a dg abelian extension of $T_A$.
	Then the action of $\frg$ on $\ker(\rho)$ factors through $T_A$.
	That is, $[x_0, v] = [x_1, v]$ for each $x_0, x_1 \in \frg$ with $\rho(x_0) = \rho(x_1)$ and each $v \in \ker(\rho)$.
\end{lemma}

\begin{proof}
	This follows because the bracket on $\ker(\rho)$ is trivial and $x_1 - x_0 \in \ker(\rho)$.
\end{proof}

For a dg abelian extension $\frg$ of $T_A$ we thus obtain an induced action of $T_A$ on $\ker(\rho)$, which we denote by
\begin{equation}
	\label{eq:Lie derivative on ker(rho)}
	\pounds_{(-)} (-) \colon T_A \otimes_\RN \ker(\rho) \longrightarrow \ker(\rho)\ ,
	\qquad
	(X,v)\longmapsto \pounds_X v \coloneqq [\hat{X}, v]\ ,
\end{equation}
where $\hat{X} \in \frg$ is any element with $\rho(\hat{X}) = X$.

From now on we assume that $\frg$ is a dg abelian extension of $T_A$.
Under this assumption, we now endow the simplicial dg $A$-module $\Hom_{T_A} (C_* (\Delta; \RN), \ch_\RN(\frg))$ defined in~\eqref{eq:ch_k(hat(frg)) as pullback} with the structure of a simplicial dg Lie algebroid over $A$.
Note that the splitting
\begin{equation}
	\Hom_{T_A} \big( C_* (\Delta^n; \RN), \ch_\RN(\frg) \big)
	= \bigoplus_{l \in \NN_0} \ \bigoplus_{\varphi \colon [l] \longhookrightarrow  [n]}\, \widehat{\frg}_{n, (\varphi)}
\end{equation}
persists at the level of dg $A$-modules (in particular, the $A$-action changes degrees, but preserves weights and hence these dg $\RN$-submodules).

For each $n \in \NN_0$, we introduce a bracket
\begin{equation}
	[-,-]_{\widehat{\frg}_n} \colon \Hom_{T_A} \big( C_* (\Delta^n; \RN), \ch_\RN(\frg) \big) \otimes_\RN \Hom_{T_A} \big( C_* (\Delta^n; \RN), \ch_\RN(\frg) \big) \longrightarrow \Hom_{T_A} \big( C_* (\Delta^n; \RN), \ch_\RN(\frg) \big)
\end{equation}
as follows:
\begin{subequations}\label{eq:hat(frg) backet}
	\begin{myenumerate}
		\item The bracket returns zero whenever both its arguments have non-zero \textit{weight}:
		\begin{equation}
			\label{eq:hat(frg) backet w=(+,+)}
			[\xi, \zeta]_{\widehat{\frg}_n} = 0 \ ,
			\qquad
			\forall\, \xi, \zeta \in \widehat{\frg}_{n, (+)}\ .
		\end{equation}
		
		\item Given two elements
		\begin{equation}
			\xi = \sum_{i = 0}^n\, \big( x_i, i \colon [0] \longhookrightarrow  [n] \big)
			\ \in \ \widehat{\frg}_{n, (0)}
			\qqandqq
			\zeta = \sum_{i = 0}^n\, \big( y_i, i \colon [0] \longhookrightarrow  [n] \big)
			\ \in \ \widehat{\frg}_{n, (0)}
		\end{equation}
		of weight zero, we set
		\begin{equation}
			\label{eq:hat(frg) backet w=(0,0)}
			[\xi, \zeta]_{\widehat{\frg}_n}
			= \sum_{i = 0}^n\, \big( [x_i, y_i]_\frg,\, i \colon [0] \longhookrightarrow  [n] \big)
			\ \in \ \widehat{\frg}_{n, (0)}\ .
		\end{equation}
		
		\item Finally, consider a pair of elements
		\begin{equation}
			\xi = \sum_{i = 0}^n\, \big( x_i, i \colon [0] \longhookrightarrow  [n] \big)
			\ \in \ \widehat{\frg}_{n, (0)}
			\qqandqq
			\big( y, \psi \colon [l] \longhookrightarrow  [n] \big) \ \in \ \widehat{\frg}_{n, (l)}\ ,
		\end{equation}
		where $l > 0$.
		It follows from~\eqref{eq:frg: ul k-module} that there exists $X \in T_A$ such that $\rho(x_i) = X$, for all $i = 0, \ldots, n$.
		We set
		\begin{equation}
			\label{eq:hat(frg) backet w=(0,+)}
			\big[ \xi, (y, \psi) \big]_{\widehat{\frg}_n}
			= \frac{1}{l+1}\, \sum_{i = 0}^l\, \big( [x_{\psi(i)}, y]_\frg,\, \psi \big)
			= (\pounds_X y, \psi) \ \in \ \widehat{\frg}_{n,(l)} \ ,
		\end{equation}
		where we have used Lemma~\ref{st:dgLAgd action on ker(rho) for ab ext} and the notation introduced in~\eqref{eq:Lie derivative on ker(rho)}. We also write
		\begin{equation} \label{eq:0+bracketrho}
		\big[ \xi, (y, \psi) \big]_{\widehat{\frg}_n}
			= (\pounds_{\widehat{\rho}(\xi)} y, \psi) \ .
		\end{equation}
	\end{myenumerate}
\end{subequations}

We extend this bracket antisymmetrically and $\RN$-bilinearly to a morphism
\begin{equation}
	[-,-]_{\widehat{\frg}_n} \colon \Hom_{T_A} \big( C_* (\Delta^n; \RN), \ch_\RN(\frg) \big) \otimes_\RN \Hom_{T_A} \big( C_* (\Delta^n; \RN), \ch_\RN(\frg) \big) \longrightarrow \Hom_{T_A} \big( C_* (\Delta^n; \RN), \ch_\RN(\frg) \big)
\end{equation}
of graded $\RN$-modules.
For later use we record
\begin{equation}
	\big[ (y, \psi), \xi \big]_{\widehat{\frg}_n}
	= (-1)^{(|y| + l)\, |\xi| +1}\, \big[ \xi, (y, \psi) \big]_{\widehat{\frg}_n}
	= \frac{(-1)^{l\, |\xi|}}{l+1}\, \sum_{i = 0}^l\, \big( [y, x_{\psi(i)}]_\frg,\, \psi \big)
	= (-1)^{(|y|+l)\, |\xi|+1}\, (\pounds_{\widehat\rho(\xi)} y, \psi)\ ,
\end{equation}
in the situation of~\eqref{eq:hat(frg) backet w=(0,+)}.

\begin{remark}
	\label{rmk:frg_(n,(+)) is dg Lie-ideal}
	We see from~\eqref{eq:hat(frg) backet w=(0,+)} that the non-zero weight dg $\RN$-submodule is an ideal for the bracket $[-,-]_{\widehat{\frg}_n}$, for each $n \in \NN_0$.
	\qen
\end{remark}

\begin{lemma}
	\label{st:hat(frg)_n is dgla over k}
	For each $n \in \NN_0$, the bracket~\eqref{eq:hat(frg) backet} endows the dg $\RN$-module $\Hom_{T_A} (C_* (\Delta^n; \RN), \ch_\RN(\frg))$ with the structure of a dg Lie algebra over $\RN$.
\end{lemma}

\begin{proof}
	The bracket is (graded) antisymmetric and bilinear by construction.
	
	\underline{Leibniz rule.} \ 
	We check the interaction of the bracket with the differential $\dd_{\widehat{\frg}_n}$ by going through the three cases in~\eqref{eq:hat(frg) backet}.
	As $\dd_{\widehat{\frg}_n}$ does not lower the weight of an element  by Remark~\ref{rem:dgraise}, we see that
	\begin{align}
		\dd_{\widehat{\frg}_n} [\xi, \zeta]_{\widehat{\frg}_n} &= 0\ ,
		\\[4pt]
		[\dd_{\widehat{\frg}_n} \xi, \zeta]_{\widehat{\frg}_n} + (-1)^{|\xi|} \, [\xi, \dd_{\widehat{\frg}_n} \zeta]_{\widehat{\frg}_n} &= 0\ ,
	\end{align}
	whenever $\xi$ and $\zeta$ are elements of $\widehat{\frg}_{n, (+)}$.
	
	Next we consider elements
	\begin{equation}
		\xi = \sum_{i = 0}^n\, \big( x_i, i \colon [0] \longhookrightarrow  [n] \big)
		\ \in \ \widehat{\frg}_{n, (0)}
		\qqandqq
		\zeta = \sum_{i = 0}^n\, \big( y_i, i \colon [0] \longhookrightarrow  [n] \big)
		\ \in \ \widehat{\frg}_{n, (0)}\ ,
	\end{equation}
	and using~\eqref{eq:dd on Hom(CDelta, frg)} compute
	\begin{align}
		\dd_{\widehat{\frg}_n} [\xi, \zeta]_{\widehat{\frg}_n}
		&= \dd_{\widehat{\frg}_n}\, \bigg( \sum_{i = 0}^n \big( [x_i, y_i]_\frg,\, i \colon [0] \longhookrightarrow  [n] \big) \bigg)
		\\[4pt]
		&= \sum_{i = 0}^n\, \bigg( \big( [\dd_\frg x_i, y_i]_\frg + (-1)^{|\xi|}\, [x_i, \dd_\frg y_i]_\frg,\, i \colon [0] \longhookrightarrow  [n] \big)
		\\
		& \hspace{2cm} - (-1)^{|\xi| + |\zeta|}\, \sum_{j = 0}^{i-1}\, \big( [x_i, y_i]_\frg,\, 0 \leqslant j < i \leqslant n \big)
		\\
		& \hspace{2cm} + (-1)^{|\xi| + |\zeta|}\, \sum_{j = i+1}^{n}\, \big( [x_i, y_i]_\frg,\, 0 \leqslant i < j \leqslant n \big) \bigg)
		\\[4pt]
		&= \sum_{i = 0}^n\, \big( [\dd_\frg x_i, y_i]_\frg + (-1)^{|\xi|}\, [x_i, \dd_\frg y_i]_\frg,\, i \colon [0] \longhookrightarrow  [n] \big)
		\\
		& \qquad\, + \sum_{0 \leqslant i < j \leqslant n}\, (-1)^{|\xi| + |\zeta|}\, \big( [x_i, y_i]_\frg - [x_j, y_j]_\frg,\, 0 \leqslant i < j \leqslant n \big)\ ,
	\end{align}
	where we have used~\eqref{eq:dd_(hat(frg)) on weight zero} and that $\frg$ is a dg Lie algebra over $\RN$.
	
	We also compute one half of the right-hand side of the Leibniz identity:
	\begin{align}
		[\dd_{\widehat{\frg}_n} \xi, \zeta]_{\widehat{\frg}_n}
		&= \sum_{i,r = 0}^n\, \big[ \dd_{\widehat{\frg}_n} \big( x_i, i \colon [0] \longhookrightarrow  [n] \big)\,,\,
		\big( y_r, r \colon [0] \longhookrightarrow  [n] \big) \big]_{\widehat{\frg}_n}
		\\[4pt]
		&= \sum_{i,r = 0}^n\, \bigg( \Big[ \big( \dd_\frg x_i, i \colon [0] \longhookrightarrow  [n] \big)\,,\,
		\big( y_r, r \colon [0] \longhookrightarrow  [n] \big) \Big]_{\widehat{\frg}_n}
		\\
		& \hspace{2cm} - \sum_{j = 0}^{i-1}\, (-1)^{|\xi|}\, \Big[ \big( x_i, 0 \leqslant j < i \leqslant n \big)\,,\, \big( y_r, r \colon [0] \longhookrightarrow  [n] \big) \Big]_{\widehat{\frg}_n}
		\\
		& \hspace{2cm} + \sum_{j = i+1}^{n}\, (-1)^{|\xi|}\, \Big[ \big(x_i, 0 \leqslant i < j \leqslant n \big)\,,\, \big( y_r, r \colon [0] \longhookrightarrow  [n] \big) \Big]_{\widehat{\frg}_n} \bigg)
		\\[4pt]
		&= \sum_{i = 0}^n\, \big( [\dd_\frg x_i, y_i]_\frg , i \colon [0] \longhookrightarrow  [n] \big)
		\\
		& \qquad\, - \sum_{i = 0}^n \ \sum_{j = 0}^{i-1}\, (-1)^{|\xi| + |\zeta|} \, \frac{1}{2}\, \big( [x_i, y_j]_\frg + [x_i, y_i]_\frg, 0 \leqslant j < i \leqslant n \big)
		\\
		&\qquad\, + \sum_{i = 0}^n \ \sum_{j = i+1}^{n}\, (-1)^{|\xi| + |\zeta|} \, \frac{1}{2} \, \big( [x_i, y_i]_\frg + [x_i, y_j]_\frg, 0 \leqslant i < j \leqslant n \big)
		\\[4pt]
		&= \sum_{i = 0}^n\, \big( [\dd_\frg x_i, y_i]_\frg , i \colon [0] \longhookrightarrow  [n] \big)
		\\
		&\qquad\, + \sum_{0 \leqslant i < j \leqslant n}\, \frac{(-1)^{|\xi| + |\zeta|}}{2} \, \big( [x_i, y_i]_\frg + [x_i, y_j]_\frg - [x_j, y_i]_\frg - [x_j, y_j]_\frg, 0 \leqslant i < j \leqslant n \big)\ .
	\end{align}
		
	The degree $|\zeta|$ which appears in the sign factor in the third step here arises since we compute a bracket of the type~\eqref{eq:hat(frg) backet w=(0,+)}, but with the order of arguments reversed:
	the first argument is of weight one and the second of weight zero.
	To compute the bracket in this order, we recall the formula~\eqref{eq:degree in Hom_k(CD, frg)} for the degree in $\Hom_{T_A} (C_* (\Delta^n; \RN), \ch_\RN(\frg))$ and compute
	\begin{align}
		&\sum_{r = 0}^n\, \Big[ \big(x_i, 0 \leqslant i < j \leqslant n \big)\,,\, \big( y_r, r \colon [0] \longhookrightarrow  [n] \big) \Big]_{\widehat{\frg}_n}
		\\[4pt]
		& \hspace{4cm} = \sum_{r = 0}^n\, (-1)^{(|\xi| + 1)\, |\zeta| + 1}\, \Big[ \big( y_r, r \colon [0] \longhookrightarrow  [n] \big)\,,\, \big(x_i, 0 \leqslant i < j \leqslant n \big) \Big]_{\widehat{\frg}_n}
		\\[4pt]
		& \hspace{4cm} = (-1)^{(|\xi| + 1)\, |\zeta| + 1}\,\tfrac12\, \big( [y_i, x_i]_\frg + [y_j, x_i]_\frg, 0 \leqslant i < j \leqslant n \big)
		\\[4pt]
		& \hspace{4cm} = (-1)^{|\zeta|} \, \tfrac12\, \big( [x_i, y_i]_\frg + [x_i, y_j]_\frg, 0 \leqslant i < j \leqslant n \big)\ .
	\end{align}
	
	Noting that the degree of $\zeta = \sum_{i = 0}^n\, (y_i, i \colon [0] \longhookrightarrow  [n]) \in \widehat{\frg}_{n,(0)}$ is $|\zeta| = |y_i|$, for each $i = 0, \ldots, n$, we also compute the other half of the right-hand side of the Leibniz identity:
	\begin{align}
		(-1)^{|\xi|}\, [\xi, \dd_{\widehat{\frg}_n} \zeta]_{\widehat{\frg}_n}
		&= \sum_{i = 0}^n\, (-1)^{|\xi|}\, \big( [x_i, \dd_\frg y_i]_\frg , i \colon [0] \longhookrightarrow  [n] \big)
		\\[4pt]
		&\qquad \, - \sum_{i = 0}^n \ \sum_{j = 0}^{i}\, \frac{(-1)^{|\xi| + |\zeta|}}{2}\, \big( [x_j, y_i]_\frg + [x_i, y_i]_\frg, 0 \leqslant j < i \leqslant n \big)
		\\
		&\qquad \, + \sum_{i = 0}^n \ \sum_{j = i+1}^{n}\, \frac{(-1)^{|\xi| + |\zeta|}}{2}\, \big( [x_i, y_i]_\frg + [x_j, y_i]_\frg, 0 \leqslant i < j \leqslant n \big)
		\\[4pt]
		&= \sum_{i = 0}^n\, (-1)^{|\xi|}\, \big( [x_i, \dd_\frg y_i]_\frg , i \colon [0] \longhookrightarrow  [n] \big)
		\\
		&\qquad \, + \sum_{0 \leqslant i < j \leqslant n}\, \frac{(-1)^{|\xi| + |\zeta|}}{2}\, \big( [x_i, y_i]_\frg + [x_j, y_i]_\frg - [x_i, y_j]_\frg - [x_j, y_j]_\frg, 0 \leqslant i < j \leqslant n \big)\ .
	\end{align}
	
	Combining the above results, we indeed find that the Leibniz rule
	\begin{equation}
		\dd_{\widehat{\frg}_n} [\xi, \zeta]_{\widehat{\frg}_n}
		= [\dd_{\widehat{\frg}_n} \xi, \zeta]_{\widehat{\frg}_n} + (-1)^{|\xi|}\, [\xi, \dd_{\widehat{\frg}_n} \zeta]_{\widehat{\frg}_n}
	\end{equation}
	is satisfied, for all
	\begin{equation}
		\xi = \sum_{i = 0}^n\, \big( x_i, i \colon [0] \longhookrightarrow  [n] \big)
		\ \in \ \widehat{\frg}_{n, (0)}
		\qqandqq
		\zeta = \sum_{i = 0}^n\, \big( y_i, i \colon [0] \longhookrightarrow  [n] \big)
		\ \in \ \widehat{\frg}_{n, (0)}\ .
	\end{equation}
	
	It remains to check the Leibniz rule for pairs of the form
	\begin{equation}
		\xi = \sum_{i = 0}^n\, \big( x_i, i \colon [0] \longhookrightarrow  [n] \big)
		\ \in \ \widehat{\frg}_{n, (0)}
		\qqandqq
		\big( y, \varphi \colon [l] \longhookrightarrow  [n] \big) \ \in \ \widehat{\frg}_{n, (l)}\ ,
	\end{equation}
	with $l > 0$.
	We find 
	\begin{align}
		\dd_{\widehat{\frg}_n} \big[ \xi, (y, \varphi) \big]_{\widehat{\frg}_n}
		&= \frac{1}{l+1}\, \sum_{i = 0}^l\, \dd_{\widehat{\frg}_n} \big( [x_{\varphi(i)}, y]_\frg,\, \varphi \big)
		\\[4pt]
		&= \frac{1}{l+1}\, \sum_{i = 0}^l\, \big( \dd_\frg [x_{\varphi(i)}, y]_\frg, \varphi \big)
		\\
		&\qquad\, - \frac{(-1)^{|\xi| + |y| + l}}{l+1}\, \sum_{i = 0}^l \ \sum_{r = 0}^{l+1}\, (-1)^r \ \sum_{\substack{\psi \colon [l+1] \longhookrightarrow  [n]\\ \psi \circ \partial_r = \varphi}}\, \big( [x_{\varphi(i)}, y]_\frg,\, \psi \big)
		\\[4pt]
		&= \frac{1}{l+1}\, \sum_{i = 0}^l\, \big( [\dd_\frg x_{\varphi(i)}, y]_\frg + (-1)^{|x_i|}\, [x_{\varphi(i)}, \dd_\frg y]_\frg, \varphi \big)
		\\
		&\qquad\, - \frac{(-1)^{|\xi| + |y| + l}}{l+1}\, \sum_{i = 0}^l \ \sum_{r = 0}^{l+1}\, (-1)^r \ \sum_{\substack{\psi \colon [l+1] \longhookrightarrow  [n]\\ \psi \circ \partial_r = \varphi}}\, \big( [x_{\varphi(i)}, y]_\frg,\, \psi \big)
		\\[4pt]
		&= (\pounds_{\dd_\Hom \widehat\rho(\xi) }\, y + (-1)^{|\xi|}\, \pounds_{\widehat\rho(\xi)}\, \dd_\frg y, \varphi)
		\\
		& \qquad \, - (-1)^{|\xi| + |y| + l}\, \sum_{r = 0}^{l+1}\, (-1)^r \  \sum_{\substack{\psi \colon [l+1] \longhookrightarrow  [n]\\ \psi \circ \partial_r = \varphi}}\, ( \pounds_{\widehat\rho(\xi)}\,y,\, \psi \big)\ .
	\end{align}
	In the last step we have used the notation from~\eqref{eq:0+bracketrho}.
	
	On the other hand, using~\eqref{eq:hat(frg) backet w=(0,+)} and~\eqref{eq:hat(frg) backet w=(+,+)} we find 
	\begin{align}
		\big[ \dd_{\widehat{\frg}_n} \xi, (y, \varphi) \big]_{\widehat{\frg}_n}
		&= \sum_{i = 0}^n\, \big[ (\dd_\frg x_i, i)\,,\, (y, \varphi) \big]_{\widehat{\frg}_n}
		= \frac{1}{l+1}\, \sum_{i = 0}^l\, \big( [\dd_\frg x_{\varphi(i)}, y]_\frg,\, \varphi \big)
		= (\pounds_{\dd_\Hom \widehat{\rho}(\xi)}\, y, \varphi)\ ,
	\end{align}
	where in the first step we have used that the bracket vanishes on $\widehat{\frg}_{n,(+)} \otimes_\RN \widehat{\frg}_{n, (+)}$.
	Using~\eqref{eq:dd on Hom(CDelta, frg)}, we also compute
	\begin{align}
		(-1)^{|\xi|}\, \big[ \xi, \dd_{\widehat{\frg}_n} (y, \varphi) \big]_{\widehat{\frg}_n}
		&= (-1)^{|\xi|}\, \big[ \xi, (\dd_\frg y, \varphi) \big]_{\widehat{\frg}_n}
		- (-1)^{|\xi| + |y| + l}\, \sum_{r = 0}^{l+1}\, (-1)^r \  \sum_{\substack{\psi \colon [l+1] \longhookrightarrow  [n]\\ \psi \circ \partial_r = \varphi}} \, \big[ \xi, (y, \psi) \big]_{\widehat{\frg}_n}
		\\[4pt]
		&= (-1)^{|\xi|}\, \sum_{i = 0}^l\, \frac{1}{l+1}\, \big( [x_{\varphi(i)}, \dd_\frg y]_\frg, \varphi \big)
		\\
		&\qquad\, - (-1)^{|\xi| + |y| + l}\, \sum_{r = 0}^{l+1} \  \sum_{\substack{\psi \colon [l+1] \longhookrightarrow  [n]\\ \psi \circ \partial_r = \varphi}} \ \sum_{i = 0}^{l+1}\, \frac{(-1)^r}{l+2}\, \big( [x_{\psi(i)}, y]_\frg, \psi \big)
		\\[4pt]
		&= (-1)^{|\xi|}\, (\pounds_{\widehat\rho(\xi)}\, \dd_\frg y, \varphi)
		- (-1)^{|\xi| + |y| + l}\, \sum_{r = 0}^{l+1} \  \sum_{\substack{\psi \colon [l+1] \longhookrightarrow  [n]\\ \psi \circ \partial_r = \varphi}}\, (-1)^r\, (\pounds_{\widehat\rho(\xi)}\, y, \psi)\ .
	\end{align}
	Thus we see that the Leibniz rule is satisfied in this case as well.
	
	\underline{Jacobi identity.} \ 
	The Jacobi identity on triples $\xi_0, \xi_1, \xi_2 \in \widehat{\frg}_n$ of weight zero follows readily from~\eqref{eq:hat(frg) backet w=(0,0)} and the Jacobi identity for $[-,-]_\frg$.
	By~\eqref{eq:hat(frg) backet w=(+,+)} and~\eqref{eq:hat(frg) backet w=(0,+)}, the Jacobi identity is trivially satisfied on triples of elements in $\widehat{\frg}_n$ where at least two elements lie in the positive weight submodule $\widehat{\frg}_{n, (+)}$.
	It thus remains to check the Jacobi identity on elments $\xi, \zeta \in \widehat{\frg}_{n, (0)}$ and $(y, \varphi) \in \widehat{\frg}_{n, (+)}$.
	We compute
	\begin{align}
		\big[ \xi, [\zeta, (z, \varphi)]_{\widehat{\frg}_n} \big]_{\widehat{\frg}_n}
		&= \big[ \xi, (\pounds_{\widehat{\rho}(\zeta)} z, \varphi) \big]_{\widehat{\frg}_n}
		= \big( \pounds_{\widehat{\rho}(\xi)} (\pounds_{\widehat{\rho}(\zeta)} z), \varphi \big)\ ,
		\\[4pt]
		\big[ [\xi, \zeta]_{\widehat{\frg}_n}, (z, \varphi) \big]_{\widehat{\frg}_n}
		&= \Big[ \sum_{i = 0}^n\, \big( [x_i, y_i]_\frg, i \big), (z, \varphi) \Big]_{\widehat{\frg}_n}
		= (\pounds_{[\widehat{\rho}(\xi), \widehat{\rho}(\zeta)]} \,z, \varphi)\ ,
		\\[4pt]
		(-1)^{|\zeta| \, |\xi|} \, \big[ \zeta, [\xi, (z, \varphi)]_{\widehat{\frg}_n} \big]_{\widehat{\frg}_n}
		&= (-1)^{|\zeta|\, |\xi|}\, \big( \pounds_{\widehat{\rho}(\zeta)} (\pounds_{\widehat{\rho}(\xi)} z), \varphi \big)\ .
	\end{align}
	The Jacobi identity then follows from that of $[-,-]_\frg$; it implies that the action map $\frg \otimes_\RN \ker(\rho) \longrightarrow \ker(\rho)$ is compatible with the bracket on $\frg$, that is
	\begin{equation}
		\pounds_{[\widehat{\rho}(\xi), \widehat{\rho}(\zeta)]} \,z
		= \pounds_{\widehat{\rho}(\xi)}\, \pounds_{\widehat{\rho}(\zeta)} z
		- (-1)^{|\xi| \, |\zeta|} \, \pounds_{\widehat{\rho}(\zeta)}\, \pounds_{\widehat{\rho}(\xi)} z\ ,
	\end{equation}
as required.
\end{proof}

\begin{lemma}
	The $A$-module structure on $\widehat{\frg}_n$ inherited from that on $\frg$ via the pullback~\eqref{eq:ch_k(hat(frg)) as pullback} is compatible with the anchor map and brackets on $\widehat{\frg}_n$, for each $n \in \NN_0$.
	That is, for each $n \in \NN_0$, the dgla $\widehat{\frg}_n$ over $\RN$ carries the structure of a dg Lie algebroid over $A$.
\end{lemma}

\begin{proof}
	The $A$-action reads as
	\begin{equation}
		a\, (y, \varphi) = (a\,y, \varphi)\qquad ,
		\qquad
		a\, \xi = a \, \sum_{i = 0}^n\, \big( x_i, i \colon [0] \longhookrightarrow  [n] \big)
		= \sum_{i = 0}^n\, \big( a\, x_i, i \colon [0] \longhookrightarrow  [n] \big)\ .
	\end{equation}
	Here the $A$-action inside the brackets is the original $A$-action on $\frg$.
	Since $\widehat{\rho} = 0$ in positive weight, there is no compatibility condition to check in that case.
	In weight zero we check that
	\begin{align}
		\widehat{\rho} (a\, \xi)
		&= \widehat{\rho} \Big( \sum_{i = 0}^n\, \big( a\, x_i, i \colon [0] \longhookrightarrow  [n] \big) \Big)
		= a\, \rho(x_0)
		= a\, \widehat{\rho}(\xi)\ .
	\end{align}
	
	It remains to check that the bracket on $\widehat{\frg}_n$ satisfies the anchored Leibniz rule with respect to the $A$-action.
	In the case where both arguments of the bracket have positive weight, the bracket is zero and there is nothing to check.
	Thus let $a \in A$, $\xi \in \widehat{\frg}_{n, (0)}$ and $(y, \varphi) \in \widehat{\frg}_{n, (l)}$.
	Then
	\begin{align}
		\big[ \xi, a\, (y, \varphi) \big]_{\widehat{\frg}_n}
		&= \big( \pounds_{\widehat{\rho}(\xi)} (a\,y), \varphi \big)
		\\[4pt]
		&= ( a\, \pounds_{\widehat{\rho}(\xi)} y + \pounds_{\widehat{\rho}(\xi)}(a) \, y, \varphi)
		\\[4pt]
		&= a \, (\pounds_{\widehat{\rho}(\xi)} y, \varphi) + \pounds_{\widehat{\rho}(\xi)}(a) \, (y, \varphi)\ ,
	\end{align}
	where we have used that the bracket $[-,-]_\frg$ on $\frg$ satisfies the anchored Leibniz rule with respect to the $A$-action on $\frg$.
	We also check that
	\begin{align}
		\big[(y, \varphi), a\, \xi \big]_{\widehat{\frg}_n}
		&= (-1)^{(|a|+|\xi|)\,(|y | + l) + 1}\, (\pounds_{a\, \widehat{\rho}(\xi)} y, \varphi)
		\\[4pt]
		&= (-1)^{(|a|+|\xi|)\,(|y| + l) + 1}\, (a \, \pounds_{\widehat{\rho}(\xi)} y, \varphi)
		\\[4pt]
		&= (-1)^{(|a|+|\xi|)\,(|y| + l) + 1}\, a \, (\pounds_{\widehat{\rho}(\xi)} y, \varphi)
		\\[4pt]
		&= (-1)^{|a|\, (|y| + l)}\, a\, \big[(y, \varphi), \xi \big]_{\widehat{\frg}_n}\ ,
	\end{align}
	where we have used $\widehat{\rho}(y, \varphi) = 0$.
	
	Finally, we consider the case of two elements $\xi,\zeta\in\widehat{\frg}_{n, (0)}$.
	We compute
	\begin{align}
		[\xi, a\, \zeta]_{\widehat{\frg}_n}
		&= \sum_{i = 0}^n\, \big( [x_i, a\, y_i]_\frg,\, i \colon [0] \longhookrightarrow  [n] \big)
		\\[4pt]
		&= \sum_{i = 0}^n\, \big( a \, [x_i, y_i]_\frg + \pounds_{\rho(x_i)}(a) \, y_i,\, i \colon [0] \longhookrightarrow  [n] \big)
		\\[4pt]
		&= a \, [\xi, \zeta]_{\widehat{\frg}_n} + \pounds_{\widehat{\rho}(\xi)}(a) \, \zeta\ .
	\end{align}
	Here we have used  $\rho(x_0) = \cdots = \rho(x_n) = \widehat{\rho}(\xi)$.
\end{proof}

\begin{lemma}
	The dg Lie algebroids $\widehat{\frg}_n$, for $n \in \NN_0$, assemble into a simplicial dg Lie algebroid over $A$.
	That is, the structure of a dg Lie algebroid over $A$ on each $\widehat{\frg}_n$, for $n \in \NN_0$, is compatible with the simplicial structure maps induced by the pullback~\eqref{eq:ch_k(hat(frg)) as pullback}.
\end{lemma}

\begin{proof}
	Let $\alpha \colon [m] \longrightarrow [n]$ be a morphism in $\bbDelta$.
	We denote the morphism in \smash{$\Mod_A^\dg$} induced by $\alpha$ as $\alpha^* \colon \widehat{\frg}_n \longrightarrow \widehat{\frg}_m$ (it is described explicitly in~\eqref{eq:simplical mps on ch(frg)}).
	From the description of the $A$-action it is evident that $\alpha^*$ is a morphism of $A$-modules.
	Furthermore, $\alpha^*$ is compatible with the differential by construction of the cochain complex underlying $\widehat{\frg}_n$ as a pullback in $\Mod_\RN^\dg$ (see~\eqref{eq:ch_k(hat(frg)) as pullback}; however, this compatibility can also be seen directly, since $\alpha^*$ acts on $\varphi$ by postcomposition and the differential by precomposition).
	
Let $\xi = \sum_{i = 0}^n\, (x_i, i)$ be an element in $\widehat{\frg}_{n, (0)}$. It will be helpful to observe that
	\begin{equation}
		\alpha^*\xi
		= \alpha^* \bigg( \sum_{i = 0}^n\, (x_i, i) \bigg)
		= \sum_{i = 0}^n \ \sum_{\substack{j \colon [0] \longhookrightarrow  [m] \\ \alpha(j) = i}}\, (x_i, j)
		= \sum_{i = 0}^n \  \sum_{j \in \alpha^{-1}(i)}\, (x_i, j)
		= \sum_{j = 0}^m\, (x_{\alpha(j)}, j)\ .
	\end{equation}
	
	Using this identity, we can check the compatibility of $\alpha^*$ with brackets between elements of weight zero and of positive weight.
	We compute
	\begin{align}
		\big[ \alpha^*\xi, \alpha^*(y, \varphi) \big]_{\widehat{\frg}_m}
		&= \bigg[ \sum_{j = 0}^m\, (x_{\alpha(j)}, j),\, \sum_{\substack{\psi \colon [l] \longhookrightarrow  [m] \\ \alpha \circ \psi = \varphi}}\, (y, \psi) \bigg]_{\widehat{\frg}_m}
		\\[4pt]
		&= \sum_{\substack{\psi \colon [l] \longhookrightarrow  [m] \\ \alpha \circ \psi = \varphi}} \ \sum_{j = 0}^m\, \big[ (x_{\alpha(j)}, j), (y, \psi) \big]_{\widehat{\frg}_m}
		\\[4pt]
		&= \sum_{\substack{\psi \colon [l] \longhookrightarrow  [m] \\ \alpha \circ \psi = \varphi}}\, (\pounds_{\widehat{\rho}(\alpha^*\xi)} y, \psi)
		\\[4pt]
		&= \sum_{\substack{\psi \colon [l] \longhookrightarrow  [m] \\ \alpha \circ \psi = \varphi}}\, (\pounds_{\widehat{\rho}(\xi)} y, \psi)
		\\[4pt]
		&= \alpha^* \big[ \xi, (y, \varphi) \big]_{\widehat{\frg}_n}\ ,
	\end{align}
where we used $\widehat{\rho}(\alpha^*\xi) = \widehat{\rho}(\xi)$, which follows either directly or by construction of $\widehat{\frg}_n$ via the pullback in~\eqref{eq:ch_k(hat(frg)) as pullback}.
	
	Next we check the compatibility of $\alpha^*$ with brackets between two elements of weight zero.
	On the one hand, we compute
	\begin{align}
		\big[ \alpha^*\xi, \alpha^*\zeta \big]_{\widehat{\frg}_m}
		= \bigg[ \sum_{j = 0}^m\, (x_{\alpha(j)}, j),\,
		\sum_{s = 0}^m\, (y_{\alpha(s)}, s) \bigg]_{\widehat{\frg}_m}
		= \sum_{j = 0}^m\, \big( [x_{\alpha(j)}, y_{\alpha(j)}]_\frg,\, j \big)\ .
	\end{align}
	On the other hand, we obtain
	\begin{align}
		\alpha^* [\xi, \zeta]_{\widehat{\frg}_n}
		= \alpha^* \bigg( \sum_{i = 0}^n\, \big( [x_i, y_i]_\frg,\, i \big) \bigg)
		= \sum_{j = 0}^m\, \big( [x_{\alpha(j)}, y_{\alpha(j)}]_\frg,\, j \big)\ ,
	\end{align}
as required.
\end{proof}

To summarise, we have shown 

\begin{proposition}
	\label{st:hat(frg) is dg Lie agd}
	For each dg abelian extension $\frg$ of $T_A$, the simplicial dg $A$-module defined by the pullback~\eqref{eq:ch_k(hat(frg)) as pullback} carries the structure of a simplicial dg Lie algebroid $\widehat{\frg}$ over $A$, with brackets as defined in~\eqref{eq:hat(frg) backet}.
\end{proposition}

Next we will enhance this result to

\begin{theorem}
	\label{st:fib res of fibrant dg abelian extensions}
	For each dg abelian extension $\frg$ of $T_A$, the simplicial dg Lie algebroid $\widehat{\frg}$ is a fibrant simplicial resolution of $\frg$ in the semi-model structure on \smash{$L_\infty \Agd_A^\dg$}.
\end{theorem}

\begin{corollary}
	Let $\frh$ be an $A$-cofibrant $L_\infty$-algebroid over $A$.
	Suppose that $\frg$ is a  dg abelian extension of $T_A$.
	Then the derived space of morphisms $\frh \longrightarrow \frg$ in \smash{$L_\infty\Agd_A^\dg$} can be computed as
	\begin{equation}
		\mathbf{R}\Hom_{L_\infty\Agd_A^\dg}(\frh, \frg)
		\cong L_\infty\Agd_A^\dg(Q\frh, \widehat{\frg})
		\cong L_\infty\Agd_A(\frh, \widehat{\frg})\ ,
	\end{equation}
	where the isomorphisms are in the homotopy category $\Ho\sSet$.
\end{corollary}

\begin{proof}
	This follows by combining Remark~\ref{rmk:computing Con_l(frg)}, Proposition~\ref{st:hat(frg) is dg Lie agd} and Theorem~\ref{st:fib res of fibrant dg abelian extensions}.
\end{proof}

\begin{proof}[Proof of Theorem~\ref{st:fib res of fibrant dg abelian extensions}]
	We need to check that
	\begin{myenumerate}
		\item the functor \smash{$\widehat{\frg} \colon \bbDelta^\opp \longrightarrow L_\infty \Agd_A^{\dg}$} is Reedy fibrant, and
		
		\item the canonical morphism $\sfc\, \frg \longrightarrow \widehat{\frg}$ from the constant simplicial object on $\frg$ is an objectwise weak equivalence.
	\end{myenumerate}
	
	We start with (1), the Reedy fibrancy.
	By~\cite[Theorems~3.1 and 3.3]{Nuiten:HoAlg_for_Lie_Algds} there is a sequence of Quillen adjunctions
	\begin{equation}
		\begin{tikzcd}[column sep=1.25cm]
			\Mod_\RN^\dg \ar[r, shift left = 0.15cm, "\perp"' yshift=0.05cm]
			& (\Mod_\RN^\dg)_{/T_A} \ar[r, shift left = 0.15cm, "\perp"' yshift=0.05cm] \ar[l, shift left = 0.15cm]
			& (\Mod_A^\dg)_{/T_A} \ar[r, shift left = 0.15cm, "\perp"' yshift=0.05cm] \ar[l, shift left = 0.15cm]
			& L_\infty\Agd_A^\dg\ , \ar[l, shift left = 0.15cm]
		\end{tikzcd}
	\end{equation}
	where the right adjoints are forgetful functors.
	Here $\Mod_\RN^\dg$ is endowed with the projective model structure, and each of the other (semi-)model structures is transferred from this projective model structure along the respective right adjoints.
	However, the leftmost right adjoint does not preserve final objects.
	
	It thus suffices to check that the underlying functor $$U \widehat{\frg} \colon \bbDelta^\opp \longrightarrow (\Mod_\RN^\dg)_{/T_A}$$ is Reedy fibrant.
	This is the functor sending $[n] \in \bbDelta$ to the morphism
	\begin{equation}
		\ch_\RN(\widehat{\frg})
		= \Hom_{T_A} \big( C_* (\Delta^n; \RN), \ch_\RN(\frg) \big)
		 \longrightarrow T_A
	\end{equation}
	which appears as the left-hand vertical arrow in the pullback square~\eqref{eq:ch_k(hat(frg)) as pullback}.
	We thus have to check that the matching map
	\begin{equation}
		(U\widehat{\frg})_n \longrightarrow M_n U\widehat{\frg} \times_{M_n \sfc\,\id_{T_A}} \id_{T_A}
	\end{equation}
	is a fibration in \smash{$(\Mod_\RN^\dg)_{/T_A}$}, for each $n \in \NN_0$.
	
	The indexing category for cospans is connected, and so the functor \smash{$(\Mod_\RN^\dg)_{/T_A} \longrightarrow \Mod_\RN^\dg$} preserves and reflects pullbacks.
	Since the model structure on \smash{$(\Mod_\RN^\dg)_{/T_A}$} is transferred from that on \smash{$\Mod_\RN^\dg$}, we equivalently need to check that the matching map
	\begin{equation}
		\ch_\RN(\widehat{\frg}_n) \longrightarrow M_n \ch_\RN (\widehat{\frg}) \times_{M_n \sfc\, T_A} T_A
	\end{equation}
	is a fibration in \smash{$\Mod_\RN^\dg$}, for each $n \in \NN_0$.
	In other words, we are checking that the left-hand vertical arrow in~\eqref{eq:ch_k(hat(frg)) as pullback} is a Reedy fibration in \smash{$\Fun(\bbDelta^\opp, \Mod_\RN^\dg)$}.
	
	The stability of fibrations under pullbacks will imply this, provided we can show that the right-hand vertical morphism in~\eqref{eq:ch_k(hat(frg)) as pullback} is a Reedy fibration in \smash{$\Fun(\bbDelta^\opp, \Mod_\RN^\dg)$}.
	There are canonical isomorphisms
	\begin{equation}
		M_n \Hom_\RN \big( C_* (\Delta; \RN), \ch_\RN(\frg) \big)
		\cong \Hom_\RN \big( L_n C_* (\Delta; \RN), \ch_\RN(\frg) \big)
	\end{equation}
	where the matching object $M_n$ is computed for a $\bbDelta^\opp$-shaped diagram, and the latching object $L_n$ is computed for the $\bbDelta$-shaped diagram \smash{$C_*(\Delta;\RN) \colon \bbDelta \longrightarrow \Mod_\RN^\dg$}.
	Property~(1) will now follow if we can show that the pullback product map induced by the commutative square
	\begin{equation}
		\begin{tikzcd}[row sep = 1cm , column sep = 1cm]
			\Hom_\RN \big( C_* (\Delta^n; \RN), \ch_\RN(\frg) \big) \ar[r] \ar[d]
			& \Hom_\RN \big( L_n C_* (\Delta; \RN), \ch_\RN(\frg) \big) \ar[d]
			\\
			\Hom_\RN \big( C_* (\Delta^n; \RN), \ch_\RN(T_A) \big) \ar[r]
			& \Hom_\RN \big( L_n C_* (\Delta; \RN), \ch_\RN(T_A) \big)
		\end{tikzcd}
	\end{equation}
	is a fibration, for each $n \in \NN_0$.
		
	Since \smash{$\Mod_\RN^\dg$} is a symmetric monoidal model category, this will follow if the anchor map $\rho \colon \frg \longrightarrow T_A$ is a fibration and the canonical morphism $C_*(\Delta^n;\RN) \longrightarrow L_n C_*(\Delta;\RN)$ is a cofibration in $\Mod_\RN^\dg$.
	The first condition is satisfied by the assumption that $\frg$ is fibrant.
	For the second condition, recall that there is a sequence of Quillen adjunctions
	\begin{equation}
		\begin{tikzcd}[column sep=1.25cm]
			\sSet \ar[r, shift left = 0.15cm, "\perp"' yshift=0.05cm, "{\ZN[-]}"]
			& \Fun(\bbDelta^\opp, \Ab) \ar[r, shift left = 0.15cm, "\perp"' yshift=0.05cm, "N"] \ar[l, shift left = 0.15cm]
			& \Mod_{\ZN, \leqslant 0}^\dg \ar[l, shift left = 0.15cm, "\varGamma"] \ar[r, hookrightarrow, shift left = 0.15cm, "\perp"' yshift=0.05cm]
			& \Mod_\ZN^\dg \ar[l, shift left = 0.15cm, "\tau^{\leqslant 0}"] \ar[r, shift left = 0.15cm, "\perp"' yshift=0.05cm, "{(-) \otimes_\ZN \RN}"]
			& \Mod_\RN^\dg \ar[l, shift left = 0.15cm]
		\end{tikzcd}
	\end{equation}
	The composition of the left adjoints is the singular chains functor with coefficients in $\RN$,
	\begin{equation}
		C_*(-;\RN) \colon \sSet \longrightarrow \Mod_\RN^\dg\ ,
	\end{equation}
	where we view the singular chain complex as an unbounded chain complex.
	It can be written as a composition of left Quillen functors and so is a left Quillen functor itself.
	The cosimplicial diagram $\Delta \colon \bbDelta \longrightarrow \sSet$, $[n] \longmapsto \Delta^n$ is Reedy cofibrant and, moreover, a Reedy cofibrant resolution of the $0$-simplex $\Delta^0$ in $\sSet$~\cite[Corollaries~15.9.11 and 15.9.12]{Hirschhorn:MoCats}.
	By~\cite[Proposition~15.4.1]{Hirschhorn:MoCats} it follows that the functor \smash{$C_*(\Delta;\RN) \colon \bbDelta \longrightarrow \Mod_\RN^\dg$} is Reedy cofibrant.
	That is, the latching map \smash{$C_*(\Delta^n; \RN) \longrightarrow L_n C_*(\Delta;\RN)$} is a cofibration in $\Mod_\RN^\dg$, for each $n \in \NN_0$.
	This completes the proof of (1), i.e.~that \smash{$\widehat{\frg} \colon \bbDelta^\opp \longrightarrow L_\infty\Agd_A^\dg$} is Reedy fibrant.
	
	It remains to show (2), i.e.~that the canonical morphism $\sfc\, \frg \longrightarrow \widehat{\frg}$ is an objectwise weak equivalence of simplicial diagrams in \smash{$L_\infty\Agd_A^\dg$}.
	Since the semi-model structure on \smash{$L_\infty\Agd_A^\dg$} is transferred from the projective model structure on \smash{$\Mod_\RN^\dg$}, it suffices to check this at the level of the underlying dg $\RN$-modules.
	We again use properties of the pullback square~\eqref{eq:ch_k(hat(frg)) as pullback}.
	The bottom horizontal morphism in that diagram is an objectwise weak equivalence of functors \smash{$\bbDelta^\opp \longrightarrow \Mod_\RN^\dg$}.
	We have shown above that the right-hand vertical morphism in~\eqref{eq:ch_k(hat(frg)) as pullback} is a Reedy fibration.%
	\footnote{In fact, for this part of the proof it suffices that it is a fibration in the projective model structure on $\Fun(\bbDelta^\opp, \Mod_\RN^\dg)$, which holds since \smash{$\Mod_\RN^\dg$} is a monoidal model category and the anchor map $\rho \colon \frg \longrightarrow T_A$ is a fibration by assumption.}
	Since every object in $\Mod_\RN^\dg$ is fibrant, the square~\eqref{eq:ch_k(hat(frg)) as pullback} is homotopy cartesian in \smash{$\Mod_\RN^\dg$}.
	We infer that the top horizontal morphism in~\eqref{eq:ch_k(hat(frg)) as pullback} is an objectwise weak equivalence as well.
	
	Now consider the augmentation of the diagram~\eqref{eq:ch_k(hat(frg)) as pullback} given by
	\begin{equation}
		\begin{tikzcd}[column sep=1.75cm, row sep=1cm]
			\sfc \, \ch_\RN(\frg) \ar[ddr, bend right=30] \ar[rrd, bend left = 15, "\sim"] \ar[rd, dashed, "(\sim)"]
			& &
			\\
			& \Hom_{T_A} \big( C_* (\Delta; \RN), \ch_\RN(\frg) \big) \ar[r, "\sim"] \ar[d, twoheadrightarrow, "\widehat\rho", swap]
			& \Hom_\RN \big( C_* (\Delta; \RN), \ch_\RN(\frg) \big) \ar[d, twoheadrightarrow, "\rho_*"]
			\\
			& \sfc \, \ch_\RN(T_A)  \ar[r, "\sim"']
			& \Hom_\RN\big( C_* (\Delta; \RN), \ch_\RN(T_A) \big)
		\end{tikzcd}
	\end{equation}
	The left-hand bent arrow is the anchor map of $\frg$, and the top bent arrow arises from the collapsing map $\Delta \longrightarrow \sfc \, \Delta^0$.
	The dashed arrow is induced by the universal property of the pullback.
	It is an objectwise weak equivalence in $\Mod_\RN^\dg$ by the two-out-of-three property.
\end{proof}


\section{Proof of Theorem~\ref{st:finite and infinitesimal l-cons on n-gerbes}}
\label{app:Pf of equivalence of connection spaces}


\begin{proof}[Proof of Theorem~\ref{st:finite and infinitesimal l-cons on n-gerbes}]
	As a consequence of Remark~\ref{rmk:computing Con_l(frg)} and Example~\ref{eg:T_(C^oo(M)) is C^oo(M)-cofibrant} we can compute the space $\Con_p(\cC\At(g))$ as
	\begin{equation}
		\Con_p \big( \cC\At(g) \big)
		\simeq L_\infty\Agd_M \big( Q^{(p)} TM, \widehat{\frg} \big)\ ,
	\end{equation}
	where $\widehat{\frg}$ is any fibrant simplicial resolution of $\cC\At(g)$ in \smash{$L_\infty\Agd_{M}^\dg$}.
	In particular, since $\cC\At(g)$ is a  dg abelian extension of $TM$, we can use the simplicial resolution of $\cC\At(g)$ derived in Theorem~\ref{st:fib res of fibrant dg abelian extensions}.
	
	Thus, from now on, let $\frg = \cC\At(g)$ and let $\widehat{\frg}$ be that fibrant simplicial resolution of $\frg$.
	We will show that there is a canonical isomorphism of simplicial sets
	\begin{equation}
		\Con_p(g)
		\cong L_\infty\Agd_M \big( Q^{(p)} TM, \widehat{\frg} \big)\ .
	\end{equation}
	The kernel of the anchor map $\rho:\cC\At(g)\longrightarrow TM$ is given by
	\begin{equation}
		\ker(\rho)
		= \tau^{\leqslant 0} \big( C^\infty(\cC_*\,\CU), (-1)^n\, \check{\delta} \big)\ .
	\end{equation}
	Recall the notation $E_{n-1}(g)$ from Definition~\ref{def:cCAt(g)} and that the elements of $L_\infty\Agd_M(Q^{(p)} \frh, \frg)$ are in bijection with certain Maurer-Cartan elements in $\Hom_{C^\infty(M)}(\Sym^{\leqslant p}(\frh[1]), \frg[1])$ (see Lemma~\ref{st:describing order~l-mps of L_oo-Agds}).
	
	We first consider the case $p \leqslant n$ and compute the set of $l$-simplices, i.e.~the set
	\begin{equation}
		L_\infty\Agd_M \big( Q^{(p)} TM, \widehat{\frg}_l \big)\ ,
	\end{equation}
	using the description of $\infty$-morphisms of $L_\infty$-algebroids in Remark~\ref{rmk:constructing cdgc-maps into CE cdgc} and the paragraphs preceding it.
	To that end, we set
	\begin{equation}
		E_{n-1,l}(g)
		\coloneqq \bigg\{ \bigoplus_{i = 0}^l\, (X_i, f_i) \,
		\bigg| \, (X_i, f_i) \in E_{n-1}(g)\ ,\ X_i = X_j\ \forall\, i,j \in [l] \bigg\}\ .
	\end{equation}
	The $\RN$-linear map underlying any morphism $Q^{(p)}TM \longrightarrow \widehat{\frg}_l$ is given by a sequence of  $C^\infty(M)$-linear maps corresponding to the horizontal arrows in Diagram~\eqref{eq: diagram for Con_l(g) via L_oo} shown in Figure~\ref{fig: diagram for Con_l(g) via L_oo}.
		
	\begin{figure}[]
	\small
		\begin{equation}
			\label{eq: diagram for Con_l(g) via L_oo}
			\begin{tikzcd}[column sep=2.5cm, row sep=1cm]
				\vdots
				& \vdots
				\\
				0
				\ar[u] \ar[r]
				& \displaystyle{\bigoplus_{\varphi \colon [1] \longhookrightarrow  [l]}}\, \{f \in C^\infty(\cC_{n-1}\,\CU) \, | \, \check{\delta}f = 0\} \ 
				\oplus \ \displaystyle{\bigoplus_{r \geqslant 2}} \  \displaystyle{\bigoplus_{\varphi \colon [r] \longhookrightarrow  [l]}}\, C^\infty(\cC_{n-r}\,\CU)
				\ar[u]
				\\
				TM
				\ar[r, "{(\id_{TM}, A^{(1)}})"] \ar[u]
				&  E_{n,l}(g) \ 
				\oplus \ \displaystyle{\bigoplus_{r \geqslant 1}} \  \displaystyle{\bigoplus_{\varphi \colon [r] \longhookrightarrow  [l]}}\, C^\infty(\cC_{n-(1+r)}\,\CU)
				\ar[u,"\dd_{\widehat{\frg}_l}",swap]
				\\
				\midwedge_\RN^2 TM
				\ar[r, "A^{(2)}"] 	\ar[u, "\delta_\ChEil"]
				& \displaystyle{\bigoplus_{r \in \NN_0}} \  \displaystyle{\bigoplus_{\varphi \colon [r] \longhookrightarrow  [l]}}\, C^\infty(\cC_{n-(2+r)}\,\CU)
				\ar[u,"\dd_{\widehat{\frg}_l}",swap]
				\\
				\midwedge_\RN^3 TM
				\ar[r, "A^{(3)}"] \ar[u, "\delta_\ChEil"]
				& \displaystyle{\bigoplus_{r \in \NN_0}} \  \displaystyle{\bigoplus_{\varphi \colon [r] \longhookrightarrow  [l]}}\, C^\infty(\cC_{n-(3+r)}\,\CU)
				\ar[u,"\dd_{\widehat{\frg}_l}",swap]
				\\
				\vdots \ar[u, "\delta_\ChEil"]
				& \vdots \ar[u,"\dd_{\widehat{\frg}_l}",swap]
				\\
				\midwedge_\RN^{p-1} TM
				\ar[r, "A^{(p-1)}"] \ar[u, "\delta_\ChEil"]
				& \displaystyle{\bigoplus_{r \in \NN_0}} \  \displaystyle{\bigoplus_{\varphi \colon [r] \longhookrightarrow  [l]}}\, C^\infty(\cC_{n-(p-1+r)}\,\CU)
				\ar[u,"\dd_{\widehat{\frg}_l}",swap]
				\\
				\midwedge_\RN^p TM
				\ar[r, "A^{(p)}"] \ar[u, "\delta_\ChEil"]
				& \displaystyle{\bigoplus_{r \in \NN_0}} \  \displaystyle{\bigoplus_{\varphi \colon [r] \longhookrightarrow  [l]}}\, C^\infty(\cC_{n-(p+r)}\,\CU)
				\ar[u,"\dd_{\widehat{\frg}_l}",swap]
				\\
				0
				\ar[r] \ar[u]
				& \displaystyle{\bigoplus_{r \in \NN_0}} \  \displaystyle{\bigoplus_{\varphi \colon [r] \longhookrightarrow  [l]}}\, C^\infty(\cC_{n-(p+1+r)}\,\CU)
				\ar[u,"\dd_{\widehat{\frg}_l}",swap]
				\\
				\vdots \ar[u]
				& \vdots \ar[u]
			\end{tikzcd}
		\end{equation}
		\normalsize
		\caption{Maps from $\overline{\Sym}{}_\RN^{\leqslant p}(TM[1])$ to the cochain complex underlying the dg Lie algebroid $\widehat{\frg}_l$ of the fibrant simplicial resolution $\widehat{\frg}$ of $\cC\At(g)$.}
		\label{fig: diagram for Con_l(g) via L_oo}
	\end{figure}
	
	We will be referring to Diagram~\eqref{eq: diagram for Con_l(g) via L_oo} throughout this proof.
	The vector spaces in the right-hand column are obtained from~\eqref{eq:frg: ul k-module}.
	Note that, whilst each column of this diagram is a complex, we do not aim to construct a cochain morphism, but rather a Maurer--Cartan element in the mapping $L_\infty$-algebra; that is, Diagram~\eqref{eq: diagram for Con_l(g) via L_oo} is \textit{not} a commutative diagram.
	Each horizontal arrow is labelled by a $C^\infty(M)$-linear map, which decomposes as a sum of maps into each component of the direct sum in its codomain, that is,
	\begin{equation}
		A^{(q)} = \bigoplus_{r \in \NN_0} \ \bigoplus_{\varphi \colon [r] \longhookrightarrow  [l]}\, A^{(q)}_\varphi\ .
	\end{equation}
	Reassembling these data by $\varphi$, instead of by form degree $q$, yields tuples $\CA_\varphi = (A^{(1)}_\varphi, \ldots, A^{(p-r)}_\varphi)$ indexed by all injective maps with codomain $[l]$ in $\bbDelta$.
	Our goal is to show that these tuples are precisely the same as the tuples of differential forms making up $l$-simplices in $\Con_p(g)$ as detailed in Section~\ref{sec: Con_l(g) equivalence}.
	
	Remark~\ref{rmk:constructing cdgc-maps into CE cdgc}(2) amounts to the top horizontal map having the identity on $TM$ in its first component, as indicated.
	It remains to compute the conditions on the linear map imposed by the Maurer--Cartan condition in Remark~\ref{rmk:constructing cdgc-maps into CE cdgc}(1).
	We check this at each level in $\overline{\Sym}{}^{\leqslant p}_\RN (TM[1])$.
	Note that
	\begin{equation}
		A^{(q)}_\varphi = 0
		\qquad
		\text{whenever} \quad n - (q + r) < 0\ ,
	\end{equation}
	where $\varphi \colon [r] \longhookrightarrow  [l]$.
	
	The fact that the weight~zero component of the top non-trivial horizontal morphism $(\id_{TM}, A^{(1)})$ in~\eqref{eq: diagram for Con_l(g) via L_oo} is required to take values in $E_{n,l}(g)$ is equivalent to the condition
	\begin{equation}
		- \check{\delta} A^{(1)}_i(X) + \dd \log(g) = 0\ ,
		\qquad
		\forall\, i \in \{0, \ldots, l\}\,,\, X \in TM\ .
	\end{equation}
	
	\underline{$\midwedge_\RN^1 TM$.} \ 
	Let $X_1 \in TM$ be a vector field.
	It is sent to an element of the form
	\begin{equation}
		\bigoplus_{i = 0}^l\, \big( X_1, A^{(1)}_i(X_1) \big) \ 
		\oplus \ \bigoplus_{r \geqslant 1} \ \bigoplus_{\varphi \colon [r] \longhookrightarrow  [l]}\, A^{(1)}_\varphi(X_1)\ .
	\end{equation}
	Here the coproduct and the differential in $\overline{\Sym}{}^{\leqslant p}_\RN (TM[1])$ act as zero.
	If $\dd_{\widehat{\frg}_l}$ denotes the differential on $\widehat{\frg}_l$, we thus have to ensure that
	\begin{equation}
		\label{eq:hat(dd) CA = 0 on Lambda^1 frX(M)}
		\dd_{\widehat{\frg}_l} \bigg( \bigoplus_{i = 0}^l \big( X_1, A^{(1)}_i(X_1) \big) \ 
		\oplus \ \bigoplus_{r \geqslant 1} \ \bigoplus_{\varphi \colon [r] \longhookrightarrow  [l]}\, A^{(1)}_\varphi(X_1) \bigg)
		= 0\ .
	\end{equation}
	
	Using the explicit expression~\eqref{eq:dd on Hom(CDelta, frg)} for the differential, we evaluate the left-hand side as
	\begin{align}
		&\dd_{\widehat{\frg}_l} \bigg( \bigoplus_{i = 0}^l \big( X_1, A^{(1)}_i(X_1) \big) \ 
		\oplus \ \bigoplus_{r \geqslant 1} \ \bigoplus_{\varphi \colon [r] \longhookrightarrow  [l]} A^{(1)}_\varphi(X_1) \bigg)
		\\[4pt]
		& \hspace{1cm}= - \sum_{i = 0}^l \ \sum_{s = 0}^1\, (-1)^s \  \sum_{\substack{\psi \colon [1] \longhookrightarrow  [l]\\ \psi \circ \partial_s = i}}\, \big( A^{(1)}_i(X_1), \psi \big)
		+ \sum_{r \geqslant 1} \ \sum_{\varphi \colon [r] \longhookrightarrow  [l]}\, \big( (-1)^{r+1}\, \check{\delta} A^{(1)}_\varphi(X_1), \varphi \big)
		\\
		&\hspace{1cm} \qquad - \sum_{r \geqslant 1}\  \sum_{\varphi \colon [r] \longhookrightarrow  [l]} \ \sum_{s = 0}^{r+1}\, (-1)^s \ \sum_{\substack{\psi \colon [r+1] \longhookrightarrow  [l]\\ \psi \circ \partial_s = \varphi}}\, \big( A^{(1)}_\varphi(X_1), \psi \big)
		\\[4pt]
		& \hspace{1cm} = \sum_{r \geqslant 1} \ \sum_{\psi \colon [r] \longhookrightarrow  [l]}\, \big( (-1)^{r+1}\, \check{\delta} A^{(1)}_\psi(X_1), \psi \big)
		- \sum_{r \geqslant 0} \ \sum_{\varphi \colon [r] \longhookrightarrow  [l]} \ \sum_{s = 0}^{r+1}\, (-1)^s \  \sum_{\substack{\psi \colon [r+1] \longhookrightarrow  [l]\\ \psi \circ \partial_s = \varphi}} \big( A^{(1)}_\varphi(X_1), \psi \big)
		\\[4pt]
		& \hspace{1cm} \overset{(\ast)}{=} \sum_{r \geqslant 1} \ \sum_{\psi \colon [r] \longhookrightarrow  [l]}\, \big( (-1)^{r+1}\, \check{\delta} A^{(1)}_\psi(X_1), \psi \big)
		- \sum_{r \geqslant 1} \ \sum_{\psi \colon [r] \longhookrightarrow  [l]} \ \sum_{s = 0}^{r+1}\, (-1)^s \ \big( A^{(1)}_{\psi \circ \partial_s}(X_1), \psi \big)
		\\[4pt]
		& \hspace{1cm} = \sum_{r \geqslant 1} \  \sum_{\psi \colon [r] \longhookrightarrow  [l]}\, \bigg( (-1)^{r+1}\, \check{\delta} A^{(1)}_\psi(X_1)
		- \sum_{s = 0}^{r+1}\, (-1)^s\, A^{(1)}_{\psi \circ \partial_s}(X_1)\,,\, \psi \bigg)\ .
	\end{align}
	
	In the step labelled~$(\ast)$ we have reorganised the second term (with the quadruple sum) as follows:
	fix  $r \in \{0, \ldots, l\}$.
	Prior to step~$(\ast)$ there is a triple sum over all injective maps $\varphi \colon [r] \longhookrightarrow  [l]$, all $s \in [r{+}1]$ and all those injective maps $\psi \colon [r{+}1] \longhookrightarrow  [l]$ such that $\psi$ extends the map $\varphi$ along $\partial_s \colon [r] \longhookrightarrow  [r{+}1]$.
	We assemble all terms in this triple sum where a fixed map $\psi \colon [r{+}1] \longhookrightarrow  [l]$ appears.
	These are precisely those terms in the original sum where $\varphi = \psi \circ \partial_s$.
	In this case, the map $\varphi$ is completely determined by $\psi$ and $s$.
	Finally, observe that each injective morphism $\psi \colon [r{+}1] \longhookrightarrow  [l]$ appears in the original sum (given $\psi$, consider any term indexed by $\varphi = \psi \circ \partial_s$, for any $s \in \{0, \ldots, r+1\}$).
	In other words, the indexing sets of the original and new sum are each canonically isomorphic to the set of all commutative diagrams of injective maps
	\begin{equation}
		\begin{tikzcd}[column sep = 1cm, row sep = 1cm]
			& {[r{+}1]} \ar[dr, hookrightarrow, "\psi"] &
			\\
			{[r]} \ar[rr, hookrightarrow, "\varphi"'] \ar[ur, hookrightarrow, "\partial_s"]
			& & {[l]}
		\end{tikzcd}
	\end{equation}
	At the same time, in the new sum we have shifted the index $r$ by one, i.e.~we perform the sum over $(r \geqslant 1, \psi \colon [r] \longhookrightarrow  [p])$ instead of $(r \geqslant 0, \psi \colon [r{+}1] \longhookrightarrow  [p])$.
	
	We thus read off that the condition~\eqref{eq:hat(dd) CA = 0 on Lambda^1 frX(M)} is equivalent to the equations
	\begin{equation}
		(-1)^{r+1}\, \check{\delta} A^{(1)}_\psi(X_1)
		- \sum_{s = 0}^{r+1}\, (-1)^s\, A^{(1)}_{\psi \circ \partial_s}(X_1)
		= 0\ ,
		\qquad
		\forall\, \psi \colon [r] \longhookrightarrow  [l]\,,\, r \geqslant 1\ .
	\end{equation}
	
	\underline{$\midwedge_\RN^m TM \ , \ 2 < m < p$.} \ 
	Let $X_1, \ldots, X_m \in TM$ be vector fields on $M$.
	We compute
	\begin{align}
		&\dd_{\widehat{\frg}_l} \bigg( \bigoplus_{r \in \NN_0} \  \bigoplus_{\varphi \colon [r] \longhookrightarrow  [l]}\, A^{(m)}_\varphi(X_1, \ldots, X_m) \bigg)
		\\[4pt]
		& \hspace{1cm} = \sum_{r \in \NN_0} \ \sum_{\psi \colon [r] \longhookrightarrow  [l]}\, \big( (-1)^{m+r} \,\check{\delta} A^{(m)}_\psi(X_1, \ldots, X_m), \psi \big)
		\\*
		&\hspace{1cm} \qquad - \sum_{r \in \NN_0} \ \sum_{\varphi \colon [r] \longhookrightarrow  [l]} \ \sum_{s = 0}^{r+1}\, (-1)^s \  \sum_{\substack{\psi \colon [r+1] \longhookrightarrow  [l]\\ \psi \circ \partial_s = \varphi}}\, \big( A^{(m)}_\varphi(X_1, \ldots, X_m), \psi \big)
		\\[4pt]
		& \hspace{1cm} = \sum_{\psi \colon [0] \longhookrightarrow  [l]}\, \big( (-1)^m \,\check{\delta} A^{(m)}_\psi(X_1, \ldots, X_m), \psi \big)
		\\
		& \hspace{1cm} \qquad + \sum_{r \geqslant 1} \ \sum_{\psi \colon [r] \longhookrightarrow  [l]}\, \bigg( (-1)^{m+r} \,\check{\delta} A^{(m)}_\psi(X_1, \ldots, X_m)
		- \sum_{s = 0}^{r+1}\, (-1)^s\, A^{(m)}_{\psi \circ \partial_s}(X_1, \ldots, X_m),\, \psi \bigg)\ ,
	\end{align}
	where we have used the same reorganisation of the quadruple sum as in the computation for $\midwedge_\RN^1 TM$ above.
	
	Next we compute
	\begin{align}
		&\bigoplus_{r \in \NN_0} \ \bigoplus_{\varphi \colon [r] \longhookrightarrow  [l]}\, A^{(m-1)}_\varphi \circ \delta_\ChEil(X_1, \ldots, X_m)
		\\[4pt]
		&\hspace{1cm} = \sum_{r \in \NN_0} \ \sum_{\varphi \colon [r] \longhookrightarrow  [l]} \ \sum_{1 \leqslant i < j \leqslant m}\, (-1)^{i+j-1}\, \Big( A^{(m-1)}_\varphi \big( [X_i, X_j]_{TM}, X_1, \ldots, \widehat{X_i}, \ldots, \widehat{X_j}, \ldots, X_m \big),\, \varphi \Big)\ .
	\end{align}
	
	Finally, we compute the bracket term in the Maurer--Cartan equation
	\begin{align}
		&\big[ A^{(m)}, A^{(m)} \big]_{\Hom, 2} (X_1 \wedge \cdots \wedge X_m)
		\\[4pt]
		& \hspace{1cm} = [-]_{\widehat{\frg}_l} \circ (A^{(m)} \otimes A^{(m)})
		\bigg( \sum_{\sigma \in \mathrm{Sh}(1, m-1)}\, \epsilon(\sigma)\, X_{\sigma(1)} \otimes (X_{\sigma(2)} \wedge \cdots \wedge X_{\sigma(m)})
		\\
		&\hspace{6cm} + \sum_{\sigma \in \mathrm{Sh}(m-1, 1)}\, \epsilon(\sigma)\, (X_{\sigma(1)} \wedge \cdots \wedge X_{\sigma(m-1)}) \otimes X_{\sigma(m)} \bigg)
		\\[4pt]
		&\hspace{1cm} = \sum_{\sigma \in \mathrm{Sh}(1, m-1)}\, \epsilon(\sigma)\, \Big[ \big( X_{\sigma(1)},\, A^{(1)}(X_{\sigma(1)}) \big)\,,\, A^{(m-1)}(X_{\sigma(2)} \wedge \cdots \wedge X_{\sigma(m)}) \Big]_{\widehat{\frg}_l}
		\\
		& \hspace{1cm} \qquad + \sum_{\sigma \in \mathrm{Sh}(m-1, 1)} \  \epsilon(\sigma)\, \Big[ A^{(m-1)} (X_{\sigma(1)} \wedge \cdots \wedge X_{\sigma(m-1)})\,,\, \big( X_{\sigma(m)},\, A^{(1)} (X_{\sigma(m)}) \big) \Big]_{\widehat{\frg}_l}
		\\[4pt]
		& \hspace{1cm} = 2\, \sum_{\sigma \in \mathrm{Sh}(1, m-1)}\, \epsilon(\sigma)\, \Big[ \big( X_{\sigma(1)}\,,\, A^{(1)}(X_{\sigma(1)}) \big)\,,\, A^{(m-1)}(X_{\sigma(2)} \wedge \cdots \wedge X_{\sigma(m)}) \Big]_{\widehat{\frg}_l}
		\\[4pt]
		&\hspace{1cm} =2\, \sum_{\sigma \in \mathrm{Sh}(1, m-1)} \  \epsilon(\sigma) \\
		& \hspace{1cm} \qquad \times \sum_{i = 0}^l \ \sum_{r \in \NN_0} \  \sum_{\varphi \colon [r] \longhookrightarrow  [l]}\,
		\Big[ \big( X_{\sigma(1)},\, (A^{(1)}_i(X_{\sigma(1)}),\, i ) \big)\,,\, \big( A^{(m-1)}_\varphi(X_{\sigma(2)} \wedge \cdots \wedge X_{\sigma(m)}),\, \varphi \big) \Big]_{\widehat{\frg}_l}
		\\[4pt]
		& \hspace{1cm} = 2\, \sum_{\sigma \in \mathrm{Sh}(1, m-1)} \, \epsilon(\sigma) \ \sum_{r \in \NN_0} \ \sum_{\varphi \colon [r] \longhookrightarrow  [l]}\,
		\Big(\pounds_{X_{\sigma(1)}} \big( A^{(m-1)}_\varphi(X_{\sigma(2)} \wedge \cdots \wedge X_{\sigma(m)}) \big)\,,\,\varphi\Big)
		\\[4pt]
		& \hspace{1cm} = 2\, \sum_{r \in \NN_0} \ \sum_{\varphi \colon [r] \longhookrightarrow  [l]} \ \sum_{i = 1}^m\, (-1)^{i-1}\,
		\Big(\pounds_{X_i} \big( A^{(m-1)}_\varphi (X_1 \wedge \cdots \wedge \widehat{X_i} \wedge \cdots \wedge X_m) \big)\,,\,\varphi\Big) \ .
	\end{align}
	
	Altogether we obtain 
	\begin{align}
		&\big( \dd_\Hom A^{(m)} + \tfrac{1}{2}\, \big[ A^{(m)}, A^{(m)} \big]_{\Hom, 2} \big) (X_1 \wedge \cdots \wedge X_m)
		\\[4pt]
		& \hspace{1cm} = \sum_{i = 0}^l\, \bigg( (-1)^m \,\check{\delta} A^{(m)}_i(X_1, \ldots, X_m)
		\\
		&\hspace{3cm} - \sum_{1 \leqslant s < t \leqslant m}\, (-1)^{s+t-1}\, A^{(m-1)}_i \big( [X_s, X_t]_{TM}, X_1, \ldots, \widehat{X_s}, \ldots, \widehat{X_t}, \ldots, X_m \big)
		\\
		&\hspace{3cm} + \sum_{t = 1}^m\, (-1)^{t+1}\,
		\pounds_{X_t} \big( A^{(m-1)}_i(X_1 \wedge \cdots \wedge \widehat{X_t} \wedge \cdots \wedge X_m)\big)
		,\, i \bigg)
		\\
		&\hspace{1cm} \quad\, + \sum_{r \geqslant 1} \ \sum_{\psi \colon [r] \longhookrightarrow  [l]}\, \bigg( (-1)^{m+r} \,\check{\delta} A^{(m)}_\psi(X_1, \ldots, X_m)
		- \sum_{s = 0}^{r+1}\, (-1)^s\, A^{(m)}_{\psi \circ \partial_s}(X_1, \ldots, X_m)
		\\
		&\hspace{4cm}\qquad - \sum_{1 \leqslant i < j \leqslant m}\, (-1)^{i+j-1}\, A^{(m-1)}_\psi \big( [X_i, X_j]_{TM}, X_1, \ldots, \widehat{X_i}, \ldots, \widehat{X_j}, \ldots, X_m \big)
		\\
		&\hspace{4cm}\qquad + \sum_{i = 1}^m\, (-1)^{i+1}\,
		\pounds_{X_i} \big( A^{(m-1)}_\psi(X_1 \wedge \cdots \wedge \widehat{X_i} \wedge \cdots \wedge X_m)\big)
		,\, \psi \bigg)
		\\[4pt]
		& \hspace{1cm} = \sum_{i = 0}^l\, \bigg( \big( (-1)^m \,\check{\delta} A^{(m)}_i + \dd A^{(m-1)}_i \big) (X_1 \wedge \cdots \wedge X_m)
		,\, i \bigg)
		\\
		& \hspace{1cm} \quad + \sum_{r \geqslant 1} \ \sum_{\psi \colon [r] \longhookrightarrow  [l]}\,
		\bigg( \Big( (-1)^{m+r} \,\check{\delta} A^{(m)}_\psi
		+ \dd A^{(m-1)}_\psi
		- \sum_{s = 0}^{r+1}\, (-1)^s\, A^{(m)}_{\psi \circ \partial_s} \Big) (X_1 \wedge \cdots \wedge X_m)
		,\, \psi \bigg)\ .
	\end{align}
	
	\underline{$\midwedge_\RN^2 TM$.} \ 
	This case is analogous to the previous computation upon replacing the bracket~\eqref{eq:hat(frg) backet w=(0,+)} with the bracket~\eqref{eq:hat(frg) backet w=(0,0)}.
	
	In summary, we obtain 
	\begin{equation}
		\dd_\Hom A^{(m)} + \tfrac{1}{2}\, \big[ A^{(m)}, A^{(m)} \big]_{\Hom, 2}
		= 0\ ,
	\end{equation}
	that is, $A^{(m)}$ is an order $p$ morphism $TM \longrightarrow \widehat{g}_l$ of $L_\infty$-algebroids on $M$ if and only if
	\begin{equation}
		(-1)^{m+r} \,\check{\delta} A^{(m)}_\psi + \dd A^{(m-1)}_\psi
		=
		\begin{cases}
			\ 0\ ,
			& \psi \colon [0] \longhookrightarrow  [l]\ ,
			\\[4pt]
			\ \displaystyle\sum\limits_{s = 0}^{r+1}\, (-1)^s\, A^{(m)}_{\psi \circ \partial_s}\ ,
			& \psi \colon [r] \longhookrightarrow  [l]\ ,\ r \geqslant 1\ .
		\end{cases}
	\end{equation}
	In other words, there is a canonical bijection between $l$-simplices of $\Con_p(g)$, defined via \v{C}ech--Deligne cocycles, and $l$-simplices of $\Con_p(\cC\At(g))$, defined via our derived geometric formalism.
	Furthermore, these isomorphisms $\Con_p(\cC\At(g))_l \longrightarrow \Con_p(g)_l$ are compatible with the simplicial structure maps, as is evident from~\eqref{eq:simplical mps on ch(frg)}, \eqref{eq: spl structure on Con_l(g), d_i} and~\eqref{eq: spl structure on Con_l(g), s_i}.
	We thus obtain an isomorphism in $\sSet$ as claimed.
	
	Finally, if $p \geqslant n$ we obtain an additional non-trivial condition at $\midwedge_\RN^{n+1} TM$.
	However, we observe that \smash{$A^{(n)}_\varphi = 0$} whenever $\varphi \colon [r] \longhookrightarrow  [l]$ with $r > 0$ and \smash{$A^{(m)} = 0$} whenever $m > n$ (since the target complex is truncated).
	In this case, the second and third computations from the case $\midwedge_\RN^m TM$ above yield the additional condition 
	\begin{equation}
		\dd A^{(n)}_i = 0\ ,
		\quad
		\forall\, i \in [l]\ .
	\end{equation}
	That is, in that situation the vertices of the $l$-simplex $\CA_\varphi$ are \textit{flat} $n$-form connections.
\end{proof}

\end{appendix}

\begin{small}

\makeatletter

\interlinepenalty=10000

\makeatother

\bibliographystyle{alphaurl}
\addcontentsline{toc}{section}{References}
\bibliography{Higher_Inf_Sym_Bib}

\end{small}

\end{document}

%% file: Higher_Inf_Sym_Defs.tex
\newcommand{\PT}{{\mathrm{PT}}}
\newcommand{\Der}{{\mathrm{Der}}}
\newcommand{\QCDT}{\mathrm{St}^{\mathrm{qcdt}}}
\newcommand{\DT}{\mathrm{St}^{\mathrm{dt}}}
\newcommand{\fgt}{{\mathrm{fgt}}}
\newcommand{\pb}{{\mathrm{pb}}}
\newcommand{\cpl}{{\mathrm{cpl}}}

\newcommand{\Cone}{{\mathrm{Cone}}}
\newcommand{\DK}{{\mathrm{DK}}}

\newcommand{\cAut}{\mathsf{Aut}}
\newcommand{\Ad}{{\mathrm{Ad}}}

\newcommand{\cof}{{\mathrm{cof}}}

\newcommand{\dg}{{\mathrm{dg}}}
\newcommand{\Tot}{{\mathrm{Tot}}}
\newcommand{\rmb}{{\mathrm{b}}}
\newcommand{\Con}{{\mathrm{Con}}}

\newcommand{\ChEil}{{\mathrm{CE}}}
\newcommand{\At}{{\mathrm{At}}}
\newcommand{\CAt}{{\check{C}\mathrm{At}}}

\newcommand{\pr}{\mathrm{pr}}

\newcommand{\String}{{\mathrm{String}}}

\newcommand{\rmH}{\mathrm{H}}
\newcommand{\rmU}{{\mathrm{U}}}

\newcommand{\rmB}{{\mathrm{B}}}

\newcommand{\Sym}{{\mathrm{Sym}}}
\newcommand{\cSym}{{\mathsf{Sym}}}

\newcommand{\Map}{\mathrm{Map}}

\newcommand{\dR}{\mathrm{dR}}

\newcommand{\dd}{\mathrm{d}}
\newcommand{\DD}{\mathrm{D}}

\newcommand{\fib}{\mathrm{fib}}
\newcommand{\Hom}{{\operatorname{Hom}}}
\newcommand{\id}{{\mathrm{id}}}

\newcommand{\cDiff}{\mathsf{Diff}}

\newcommand{\opp}{\mathrm{op}}

\newcommand{\cl}{\mathrm{cl}}

\newcommand{\red}{\mathrm{red}}
\newcommand{\holim}{\operatorname{holim}}
\newcommand{\hocolim}{\operatorname{hocolim}}
\newcommand{\hofib}{{\mathrm{hofib}}}

\newcommand{\colim}{\operatorname{colim}}

\newcommand{\Ho}{\operatorname{Ho}}

\newcommand{\rmL}{{\mathrm{L}}}

\renewcommand{\lim}{\operatorname{lim}}
\newcommand{\ch}{{\operatorname{ch}}}
\newcommand{\grmod}{{\operatorname{grmod}}}
\newcommand{\FMP}{{\mathrm{FMP}}}
\newcommand{\cdga}{{\mathrm{cdga}}}

\newcommand{\scC}{\mathscr{C}}

\newcommand{\scD}{\mathscr{D}}

\newcommand{\scS}{\mathscr{S}}
\newcommand{\scM}{{\mathscr{M}}}

\newcommand{\scX}{{\mathscr{X}}}

\newcommand{\CG}{\pazocal{G}}
\newcommand{\CI}{\pazocal{I}}

\newcommand{\CU}{\pazocal{U}}

\newcommand{\CE}{\pazocal{E}}

\newcommand{\CA}{\pazocal{A}}

\newcommand{\calO}{{\mathcal{O}}}

\newcommand{\sfD}{\mathsf{D}}

\newcommand{\sfc}{\mathsf{c}}

\newcommand{\NN}{\mathbbm{N}}
\newcommand{\RN}{\mathbbm{R}}
\newcommand{\ZN}{\mathbbm{Z}}

\newcommand{\bbE}{\mathbbm{E}}

\newcommand{\frh}{\mathfrak{h}}
\newcommand{\frg}{{\mathfrak{g}}}

\newcommand{\sym}{{\mathfrak{sym}}}

\newcommand{\ul}[1]{\underline{#1}}
\newcommand{\ol}[1]{\overline{#1}}

\newcommand{\cC}{{\check{C}}}

\newcommand{\Ab}{{\mathscr{A}\mathrm{b}}}
\newcommand{\Grb}{{\mathscr{G}\mathrm{rb}}}

\newcommand{\Set}{{\mathscr{S}\mathrm{et}}}
\newcommand{\sSet}{{\mathscr{S}\mathrm{et}_{\hspace{-.03cm}\Delta}}}

\newcommand{\Fun}{{\operatorname{Fun}}}

\newcommand{\Cat}{{\mathscr{C}\mathrm{at}}}

\newcommand{\Mfd}{{\mathscr{M}\hspace{-0.02cm}\mathrm{fd}}}

\newcommand{\Cart}{{\mathscr{C}\mathrm{art}}}

\newcommand{\Sh}{{\mathscr{S}\mathrm{h}}}
\newcommand{\PSh}{{\mathscr{PS}\mathrm{h}}}

\newcommand{\Bun}{{\mathscr{B}\hspace{-0.02cm}\mathrm{un}}}

\newcommand{\Agd}{{\mathscr{A}\mathrm{gd}}}
\newcommand{\Alg}{{\mathscr{A}\mathrm{lg}}}
\newcommand{\Lie}{{\mathscr{L}\hspace{-0.03cm}\mathrm{ie}}}

\newcommand{\InfMfd}{{\mathscr{I}\mathrm{nf}\hspace{-0.01em}\mathscr{M}\mathrm{fd}}}
\newcommand{\Art}{{\mathscr{A}\mathrm{rt}}}
\newcommand{\Spec}{{\mathscr{S}\hspace{-0.03em}\mathrm{pec}}}
\newcommand{\Spf}{{\mathscr{S}\hspace{-0.03em}\mathrm{pf}}}
\newcommand{\CRing}{{\mathrm{C}^\infty \hspace{-0.06em} \mathscr{R}\mathrm{ing}}}
\newcommand{\Perf}{{\mathscr{P}\mathrm{erf}}}
\newcommand{\Mod}{{\mathscr{M}\mathrm{od}}}
\newcommand{\coSpec}{{\mathrm{co}\mathscr{S}\mathrm{pec}}}
\newcommand{\coSpf}{{\mathrm{co}\mathscr{S}\mathrm{pf}}}

\newcommand{\midwedge}{\text{\Large$\wedge$}}
\newcommand{\midvee}{\text{\Large$\vee$}}

\newcommand{\ip}[1]{\left\langle #1 \right\rangle}
\newcommand{\dbracket}[1]{[\![ #1 ]\!]}

\def\swnk1{{\textrm{\tiny $(n{-}k{+}1)$}}}
\def\swk1{{\textrm{\tiny $(k{-}1)$}}}

\newcommand{\sfOmega}{{\mathsf{\Omega}}}
\newcommand{\sfDelta}{{\mathsf{\Delta}}}

\newcommand{\weq}{\xrightarrow{ \ \sim \ }}
\newcommand{\eq}{\xrightarrow{ \ \simeq \ }}

\newcommand{\longhookrightarrow}{\,\lhook\joinrel\longrightarrow\,}  
\newcommand{\longleadsto}{\,\relbar\joinrel\leadsto\,}
\newcommand{\longtwoheadrightarrow}{\,\relbar\joinrel\twoheadrightarrow\,}

\newenvironment{myitemize}{\begin{itemize}[itemsep=-0.1cm, leftmargin=*, topsep=0cm]}{\end{itemize}}

\newenvironment{myenumerate}{\begin{enumerate}[topsep=-\parskip+0.1cm, itemsep=0.1cm, parsep=0cm, label=(\arabic*), leftmargin=*]}{\end{enumerate}}

\newenvironment{mysubenumerate}{\begin{enumerate}[topsep=-\parskip+0.1cm, itemsep=0.1cm, parsep=0cm, label=(\roman*), leftmargin=*]}{\end{enumerate}}

\newcommand{\qen}{\hfill$\triangleleft$}

\newcommand{\qandq}{\quad \text{and} \quad}
\newcommand{\qqandqq}{\qquad \text{and} \qquad}

\newcommand{\dslash}{{/\hspace{-0.1cm}/}}

\newcommand{\bbDelta}{{{\mathbbe{\Delta}}}} 